**w. b. vasantha kandasamy**

# smarandache rings

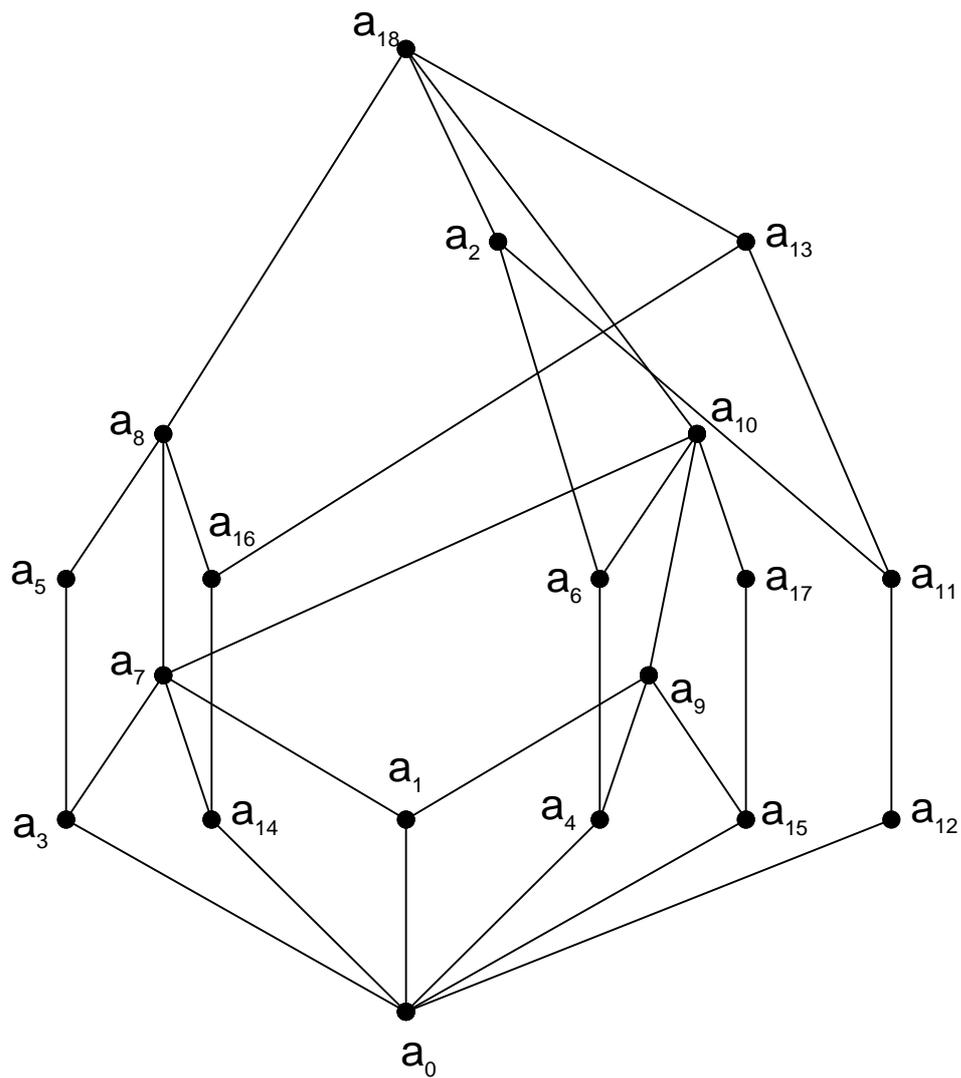

# Smarandache Rings

**W. B. Vasantha Kandasamy**

Department of Mathematics
Indian Institute of Technology, Madras
Chennai – 600036, India
e-mail: *vasantha@iitm.ac.in*
web: *http://mat.iitm.ac.in/~wbv*

2002

*The picture on the cover is the lattice representation of the S-ideals of the Smarandache mixed direct product ring $R = Z_3 \times Z_{12} \times Z_7$. This is a major difference between a ring and a Smarandache ring. For, in a ring the lattice representation of ideals is always a modular lattice but we see in case of S-rings the lattice representation of S-ideal need not in general be modular.*

# CONTENTS

# PREFACE

Over the past 25 years, I have been immersed in research in Algebra and more particularly in ring theory. I embarked on writing this book on Smarandache rings (S-rings) specially to motivate both ring theorists and Smarandache algebraists to develop and study several important and innovative properties about S-rings.

Writing this book essentially involved a good deal of reference work. As a researcher, I felt that it will be a great deal better if we thrust importance on results given in research papers on ring theory rather than detail the basic properties or classical results that the standard textbooks contain. I feel that such a venture, which has consolidated several ring theoretic concepts, has made the current book a unique one from the angle of research.

One of the major highlights of this book is by creating the Smarandache analogue of the various ring theoretic concepts we have succeeded in defining around 243 Smarandache concepts.

As it is well known, studying any complete structure is an exercise in unwieldiness. On the other hand, studying the same properties locally makes the study easier and also gives way to greater number of newer concepts. Also localization of properties automatically comes when Smarandache notions are defined. So the Smarandache notions are an excellent means to study local properties in rings.

Two levels of Smarandache rings are defined. We have elaborately dealt in case of Smarandache ring of level I, which, by default of notion, will be called as Smarandache ring. The Smarandache ring of level II could be constructed mainly by using Smarandache mixed direct product. The integral domain Z failed to be a Smarandache ring but it is one of the most natural Smarandache ring of level II.

This book is organized into five chapters. Chapter one is introductory in nature and introduces the basic algebraic structures. In chapter two some basic results and properties about rings are given. As we expect the reader to have a strong background in ring theory and algebra we have recollected for ready reference only the basic results. Chapter three is completely devoted to the introduction, description and analysis of the Smarandache rings — element-wise, substructure-wise and also by localizing the properties. The fourth chapter deals with mixed direct product of rings, which paves way for the more natural expression for Smarandache rings of level II. It is important to mention that unlike in rings where the two sided ideals form a modular lattice, we see in case of Smarandache rings the two sided ideals in general do not form a modular lattice which is described in the cover page of this book. This is a marked difference, which distinguishes a ring and a Smarandache ring. The fifth chapter contains a collection of suggested problems and it contains 200 problems in ring theory and Smarandache ring theory. It is pertinent to mention here that some problems, specially the zero divisor conjecture find several equivalent formulations. We have given many equivalent formulations, for this conjecture that has remained open for over 60 years.

I firstly wish to put forth my sincere thanks and gratitude to Dr. Minh Perez. His making my books on Smarandache notions into an algebraic structure series, provided me the necessary enthusiasm and vigour to work on this book and other future titles.

It gives me immense happiness to thank my children Meena and Kama for single-handedly helping me by spending all their time in formatting and correcting this book.

I dedicate this book to be my beloved mother-in-law Mrs. Salagramam Alamelu Ammal, whose only son, an activist-writer and crusader for social justice, is my dear husband. She was the daughter of Sakkarai Pulavar, a renowned and much-favoured Tamil poet in the palace of the King of Ramnad; and today when Meena writes poems in English, it reminds me that this literary legacy continues.

**Chapter One**

# PRELIMINARY NOTIONS

This chapter is devoted to the introduction of basic notions like, groups, semigroups, lattices and Smarandache semigroups. This is mainly done to make this book self-sufficient. As the book aims to give notions mainly on Smarandache rings, so it anticipates the reader to have a good knowledge in ring theory. We recall only those results and definitions, which are very basically needed for the study of this book.

In section one we introduce certain group theory concepts to make the reader understand the notions of Smarandache semigroups, semigroup rings and group rings. Section two is devoted to the study of semigroups used in building rings viz. semigroup rings. Section three aims to give basic concepts in lattices. The final section on Smarandache semigroups gives the definition of Smarandache semigroups and some of its properties, as this would be used in a special class of rings.

## 1.1 Groups

In this section we just define groups for we would be using it to study group rings. As the book assumes a good knowledge in algebra for the reader, we give only some definitions, notations and results with the main motivation to make the book self-contained; atleast for the basic concepts. We give examples and ask the reader to solve the problems at the end of each section, as it would help the student when she/he proceeds into the study of Smarandache rings and Smarandache notions about rings; not only for comparison of these two concepts, but to make them build more and more Smarandache structures.

**DEFINITION 1.1.1**: *A set G that is closed under a given operation '.' is called a group if the following axioms are satisfied.*

1. *The set G is non-empty.*
2. *If a, b, c ∈ G then a(bc) = (ab) c.*
3. *There are exists in G an element e such that*
   *(a) For any element a in G, ea = ae = a.*
   *(b) For any element a in G there exists an element $a^{-1}$ in G such that $a^{-1}a = aa^{-1} = e$.*

A group, which contains only a finite number of elements, is called a finite group, otherwise it is termed as an infinite group. By the order of a finite group we mean the number of elements in the group.

It may happen that a group G consists entirely elements of the from $a^n$, where a is a fixed element of G and n is an arbitrary integer. In this case G is called a cyclic group and the element a is said to generate G.

***Example 1.1.1***: Let Q be the set of rationals. Q\{0} is a group under multiplication. This is an infinite group.

***Example 1.1.2***: $Z_p$ = {0, 1, 2, … , p − 1}, p a prime be the set of integers modulo p. $Z_p$\{0} is a group under multiplication modulo p. This is a finite cyclic group of order p-1.

**DEFINITION 1.1.2**: *Let G be a group. If a . b = b . a for all a, b $\in$ G, we call G an abelian group or a commutative group.*

The groups given in examples 1.1.1 and 1.1.2 are both abelian.

**DEFINITION 1.1.3**: *Let X = {1, 2, … , n}. Let $S_n$ denote the set of all one to one mappings of the set X to itself. Define operation on $S_n$ as the composition of mappings denote it by 'o'. Now ($S_n$, o) is a group, called the permutation group of degree n. Clearly ($S_n$, o) is a non-abelian group of order n!. Throughout this text $S_n$ will denote the symmetric group of degree n.*

***Example 1.1.3***: Let X={1, 2, 3}. $S_3$ = {set of all one to one maps of the set X to itself} . The six mappings of X to itself is given below:

$$
p_0 : \begin{array}{ccc} 1 & \rightarrow & 1 \\ 2 & \rightarrow & 2 \\ 3 & \rightarrow & 3 \end{array}
$$

$$
p_1 : \begin{array}{ccc} 1 & \rightarrow & 1 \\ 2 & \rightarrow & 3 \\ 3 & \rightarrow & 2 \end{array} \qquad p_2 : \begin{array}{ccc} 1 & \rightarrow & 3 \\ 2 & \rightarrow & 2 \\ 3 & \rightarrow & 1 \end{array}
$$

$$
p_3 : \begin{array}{ccc} 1 & \rightarrow & 2 \\ 2 & \rightarrow & 1 \\ 3 & \rightarrow & 3 \end{array} \qquad p_4 : \begin{array}{ccc} 1 & \rightarrow & 2 \\ 2 & \rightarrow & 3 \\ 3 & \rightarrow & 1 \end{array}
$$

$$
\text{and} \quad p_5 : \begin{array}{ccc} 1 & \rightarrow & 3 \\ 2 & \rightarrow & 1 \\ 3 & \rightarrow & 2 \end{array}
$$

$S_3 = \{p_0, p_1, p_2, p_3, p_4, p_5\}$ is a group of order $6 = 3!$

Clearly $S_3$ is not commutative as

$$p_1 \circ p_2 = \begin{matrix} 1 & \rightarrow & 3 \\ 2 & \rightarrow & 1 \\ 3 & \rightarrow & 2 \end{matrix} = p_5$$

$$p_2 \circ p_1 = \begin{matrix} 1 & \rightarrow & 2 \\ 2 & \rightarrow & 3 \\ 3 & \rightarrow & 1 \end{matrix} = p_4.$$

Since $p_1 \circ p_2 \neq p_2 \circ p_1$, $S_3$ is a non-commutative group.

Denote $p_0, p_1, p_2, \ldots, p_5$ by

$$\begin{pmatrix} 1 & 2 & 3 \\ 1 & 2 & 3 \end{pmatrix}, \begin{pmatrix} 1 & 2 & 3 \\ 1 & 3 & 2 \end{pmatrix}, \begin{pmatrix} 1 & 2 & 3 \\ 3 & 2 & 1 \end{pmatrix}, \ldots, \begin{pmatrix} 1 & 2 & 3 \\ 3 & 1 & 2 \end{pmatrix}$$

respectively. We would be using mainly this notation.

**DEFINITION 1.1.4**: *Let (G, o) be a group. H a non-empty subset of G. We say H is a subgroup if (H, o) is a group.*

**Example 1.1.4**: Let $G = \langle g \, / \, g^8 = 1 \rangle$ be a cyclic group of order 8. $H = \{g^2, g^4, g^6, 1\}$ is subgroup of G.

**Example 1.1.5**: In the group $S_3$ given in example 1.1.3, $H = \{1, p_4, p_5\}$ is a subgroup of $S_3$.

Just we shall recall the definition of normal subgroups.

**DEFINITION 1.1.5**: *Let G be a group. A non-empty subset H of G is said to be a normal subgroup of G, if Ha = aH for every a in G or equivalently $H = \{a^{-1}ha \, / \, for \, every \, a \, in \, G \, and \, every \, h \in H\}$. If G is an abelian group or a cyclic group then all of its subgroups are normal in G.*

**Example 1.1.6**: The subgroup $H = \{1, p_4, p_5\}$ given in example 1.1.5 is a normal subgroup of $S_3$.

**Notation**: Let $S_n$ be the symmetric group of degree n. Then for $n \geq 5$, each $S_n$ has only one normal subgroup, $A_n$ which is of order $\frac{n!}{2}$ called the alternating group.

**DEFINITION 1.1.6**: *If G is a group, which has no normal subgroups then we say G is simple.*

**DEFINITION 1.1.7**: *A subnormal series of a group G is a finite sequence $H_0$, $H_1$, ..., $H_n$ of subgroups of G such that $H_i$ is a normal subgroup of $H_{i+1}$ with $H_0 = \{e\}$ and $H_n = G$.*

A normal series of G is a finite sequence $H_0$, $H_1$, ... , $H_n$ of normal subgroups of G such that $H_i \subset H_{i+1}$, $H_0 = \{e\}$ and $H_n = G$.

**Example 1.1.7**: Let $Z_{11} \setminus \{0\} = \{1, 2, \ldots, 10\}$ be the group under multiplication modulo 11. $Z_{11} \setminus \{0\}$ is a group. This has no subgroups or normal subgroups.

**Example 1.1.8**: Let G= $\langle g / g^{12} = 1 \rangle$ be the cyclic group of order 12. The series $\{1\} \subseteq \{g^6, 1\} \subseteq \{1, g^3, g^6, g^9\} \subseteq$ G. The series $\{1\} \subseteq \{1, g^6\} \subseteq \{1, g^2, g^4, g^6, g^8, g^{10}\} \subseteq$ G.

**DEFINITION 1.1.8**: *Let G be a group with identity e. We say an element $x \in$ G to be a torsion free element, if for no finite integer n, $x^n = e$. If every element in G is torsion free we say G is a torsion free group.*

**Example 1.1.9**: Let G = Q $\setminus \{0\}$; Q the field of rationals. G is a torsion free abelian group.

A torsion free group is of infinite order; by the very definition of it. The reader is requested to read more about, the composition series in groups as it would be used in studying the concept of A.C.C and D.C.C for rings in the context of Smarandache notions.

**PROBLEMS:**

1. Find all the normal subgroups in $S_n$.
2. Find all subgroups of the symmetric group $S_8$.
3. Find only cyclic subgroups of $S_9$.
4. Can $S_9$ have non-cyclic subgroups?
5. Find all abelian subgroups of $S_{12}$.
6. Find all subgroups in the dihedral group; $D_{2n} = \{a, b/a^2 = b^n = 1 \text{ and } bab = a\}$.
7. Is $D_{2.3} = \{a, b / a^2 = b^3 = 1 \text{ and } bab = a\}$ simple?

8.   Find the subnormal series of $S_n$.

9.   Find the normal series of $D_{2n}$.

10.  Find the subnormal series of $G = \{g / g^{2n} = 1\}$.

11.  Can $G = \langle g / g^p = 1$, p a prime$\rangle$ have a normal series?

12.  Find the normal series of $G = \langle g / g^{30} = 1\rangle$.

## 1.2 Semigroups

In this section we introduce the concept of semigroups mainly to study the two concepts; Smarandache semigroups and semigroup rings. Several types of semigroups are defined and their substructures like ideals and subsemigroups are also defined and illustrated with several examples. We expect the reader to have a strong background of algebra.

**DEFINITION 1.2.1**: *A semigroup is a set S together with an associative closed binary operation '.' defined on it. We shall call (S, .) a semigroup or S a semigroup.*

**Example 1.2.1**: $(Z^+ \cup \{0\}, \times)$; the set of positive integers with zero under multiplication is a semigroup.

**Example 1.2.2**: $S_{n \times m} = \{(a_{ij})/a_{ij} \in Z\}$ be the set of all $n \times m$ matrices under addition. $S_{n \times m}$ is a semigroup.

**Example 1.2.3**: $S_{n \times n} = \{(a_{ij}) / a_{ij} \in Z^+\}$ be the set of all $n \times n$ matrices under multiplication. $S_{n \times n}$ is a semigroup.

**Example 1.2.4**: Let $S(n) = \{$set of all maps from a set $X = \{x_1, x_2, \dots, x_n\}$ to itself$\}$. $S(n)$ under composition of maps is a semigroup.

**Example 1.2.5**: $Z_{15} = \{0, 1, 2, \dots, 14\}$ is the semigroup under multiplication modulo 15.

**DEFINITION 1.2.2**: *Let S be a semigroup. For a, b $\in$ S, if we have a . b = b . a, we say S is a commutative semigroup.*

**DEFINITION 1.2.3**: *Let S be a semigroup. If an element e $\in$ S such that a . e = e . a = a for all a $\in$ S, we say S is a semigroup with identity or a monoid.*

If the number of elements in a semigroup is finite we say S is a finite semigroup; otherwise S is an infinite semigroup. The semigroup given in examples 1.2.1 and 1.2.2 are commutative monoids of infinite order. The semigroup given in example 1.2.3 is an infinite semigroup which is non-commutative.

Example 1.2.4 is a non-commutative monoid of finite order. The semigroup in example 1.2.5 is a commutative monoid of finite order.

**DEFINITION 1.2.4**: *Let (S, .) be a semigroup. A non-empty subset P of S is said to be a subsemigroup if (P, .) is a semigroup.*

***Example 1.2.6***: Let $Z_{12}$ = {0, 1, 2, … , 11} be the monoid under multiplication modulo 12. P = {0, 2, 4, 8} is a subsemigroup and P is not a monoid.

Several such examples can be easily got.

**DEFINITION 1.2.5**: *Let S be a semigroup. A non-empty subset P of S is said to be a right(left) ideal of S if for all p ∈ P and s ∈ S we have ps ∈ P (sp ∈ P). If P is simultaneously both a right and a left ideal we call P an ideal of the semigroup S.*

**DEFINITION 1.2.6**: *Let S be a semigroup under multiplication. We say S has zero divisors provided 0 ∈ S and a.b = 0 for a ≠ 0, b ≠ 0 in S.*

***Example 1.2.7***: Let $Z_{16}$ = {0, 1, 2, … , 15} be the semigroup under multiplication. $Z_{16}$ has zero divisors given by

$$2.8 \equiv 0 \pmod{16}$$
$$4.4 \equiv 0 \pmod{16}$$
$$8.8 \equiv 0 \pmod{16}$$
$$4.8 \equiv 0 \pmod{16}.$$

Now we will define idempotents in semigroups.

**DEFINITION 1.2.7**: *Let S be a semigroup under multiplication. An element s ∈ S is said to be an idempotent in the semigroup if $s^2$ = s.*

***Example 1.2.8***: Let $Z_{10}$ = {0, 1, 2, … , 9} be the semigroup under multiplication modulo 10. Clearly 5 ∈ $Z_{10}$ is such that $5^2 \equiv 5 \pmod{10}$, also $6^2 \equiv 6 \pmod{10}$. Thus $Z_{10}$ has non-trivial idempotents in it.

**DEFINITION 1.2.8**: *Let S be a semigroup with unit 1 i.e., a monoid, we say an element x ∈ S is invertible if there exists a y ∈ S such that xy =1.*

***Example 1.2.9***: Let $Z_{12}$ = {0, 1, 2, … , 11} be the semigroup under multiplication modulo 12. Clearly 1 ∈ $Z_{12}$ and

$$11.11 \equiv 1 \pmod{12}$$
$$5.5 \equiv 1 \pmod{12}$$
$$7.7 \equiv 1 \pmod{12}.$$

Thus $Z_{12}$ has invertible elements.

We give some problems for the reader to solve.

**Notation**: Throughout this book S(n) will denote the set of all mapping of a set X with cardinality n to itself. i.e., X = {1, 2, ... , n}; S(n) under the composition of mappings is a semigroup. Clearly the number of elements in S(n) = $n^n$. S(n) will be addressed in this text as a symmetric semigroup.

For example the semigroup S(3) has $3^3$ i.e., 27 elements in it and S(3) is a non-commutative monoid

$$i = \begin{pmatrix} 1 & 2 & 3 \\ 1 & 2 & 3 \end{pmatrix}$$

acts as the identity. Now

$$S(2) = \left\{ \begin{pmatrix} 1 & 2 \\ 1 & 2 \end{pmatrix}, \begin{pmatrix} 1 & 2 \\ 2 & 1 \end{pmatrix}, \begin{pmatrix} 1 & 2 \\ 1 & 1 \end{pmatrix} \text{ and } \begin{pmatrix} 1 & 2 \\ 2 & 2 \end{pmatrix} \right\}$$

is a semigroup under composition of maps, in fact a monoid of order 4. We will call S(n) the symmetric semigroup of order $n^n$ by default of terminology.

## PROBLEMS:

1. Let $S = \left\{ \begin{pmatrix} a & 0 \\ 0 & 0 \end{pmatrix} / a \in Z_7 \setminus \{0\} \right\} \cup \left\{ \begin{pmatrix} 1 & 0 \\ 0 & 1 \end{pmatrix} \right\}$. Is S a semigroup under multiplication? What is the order of S?
2. Find a non-commutative semigroup of order 6.
3. Can a semigroup of order 3 be non-commutative?
4. Find the smallest non-commutative semigroup.
5. Is all semigroups of order p, p a prime, a commutative semigroup? Justify.
6. Find all subsemigroups of the symmetric semigroup S(6).
7. Find all right ideals of the symmetric semigroup S(9).
8. Find only ideals of the symmetric semigroup S(10).
9. Find a semigroup of order 26. (different from $Z_{26}$).

10. Let $S_{3 \times 3} = \{ (a_{ij}) \, / \, a_{ij} \in Z_2 \}$ i.e., set of all $3 \times 3$ matrices with entries from $Z_2 = \{0, 1\}$. Is $S_{3 \times 3}$ a semigroup? Find ideals and subsemigroups in $S_{3 \times 3}$. Does $S_{3 \times 3}$ have idempotents? Does $S_{3 \times 3}$ have zero divisors? Find units in $S_{3 \times 3}$.

11. For the semigroup $Z_{12} = \{0, 1, 2, 3, \ldots, 11\}$ under multiplication modulo 12. Find

    i. Subsemigroups which are not ideals.
    ii. Ideals.
    iii. Zero divisors.
    iv. Idempotents.
    v. Units.

12. Find in the semigroup $S(21)$ right and left ideals. Does $S(21)$ have subsemigroups which are not ideals?

## 1.3 Lattices

In this section we mainly introduce the concept of lattices as we have a well known result in ring theory which states that "the set of all two sided ideals of a ring form a modular lattice". As our main motivation for writing this book is to obtain all possible Smarandache analogous in ring we want to see how the collection of Smarandache ideals and Smarandache subrings look like. Do they form a modular lattice? We answer this question in chapter four. So we devote this section to introduce lattices and modular lattices.

**DEFINITION 1.3.1**: *Let A and B be two non-empty sets. A relation R from A to B is a subset of $A \times B$. Relations from A to A are called relation on A, for short. If (a, b) $\in$ R then we write aRb and say that a is in relation R to b. Also if a is not in relation R to b we write a$\mathbb{R}$b. A relation R on a nonempty set may have some of the following properties:*

*R is reflexive if for all a in A we have aRa.*

*R is symmetric if for all a, b in A, aRb implies bRa. R is anti symmetric if for all a,b in A, aRb and bRa imply a = b.*

*R is transitive if for all a,b,c in A aRb and bRc imply aRc. A relation R on A is an equivalence relation, if R is reflexive, symmetric and transitive.*

*In this case, [a] = {b $\in$ A / aRb} is called the equivalence class of a for any a $\in$ A.*

**DEFINITION 1.3.2**: *A relation R on a set A is called a partial order (relation) if R is reflexive, anti symmetric and transitive. In this case (A, R) is called a partially ordered set or poset.*

**DEFINITION 1.3.3**: *A partial order relation ≤ on A is called total order or lattice order if for each a, b ∈ A either a ≤ b or b ≤ a; (A, ≤) is then called a chain or a totally ordered set.*

For example {-7, 3, 2, 5, 11} is a totally ordered set under the order ≤.

Let (A, ≤) be a poset. We say a is a greatest element "if all other elements are smaller. More precisely a ∈ A is called the greatest element of A if for all x ∈ A we have x ≤ a. The element b in A is called a smallest element of A if b ≤ x for all x ∈ A. The element c ∈ A is called a maximal element of A if c ≤ x implies c = x for all x ∈ A; similarly d ∈ A is called a minimal element of A if x ≤ d implies x = d for all x ∈ A.

It can be shown that (A, ≤) has almost one greatest and one smallest element. However there may be none, one or several maximal or minimal elements. Every greatest element is maximal and every smallest element is minimal.

**DEFINITION 1.3.4**: *Let (A, ≤) be a poset and B ⊆ A.*

   a) *a ∈ A is called an upper bound of B if and only if for all b ∈ B; b ≤ a.*

   b) *a ∈ A is called a lower bound of B if and only if for all b ∈ B a ≤ b.*

   c) *The greatest amongst the lower bounds whenever it exists is called the infimum of B, and is denoted by inf B.*

   d) *The least upper bound of B, whenever it exists, is called the supremum of B and is denoted by sup B.*

**DEFINITION 1.3.5**: *A poset (L, ≤) is called lattice ordered if for every pair x, y of elements of L, the sup {x, y}and inf {x, y} exists.*

**DEFINITION 1.3.6**: *An algebraic lattice (L, ∪, ∩) is a nonempty set L with two binary operation ∩ (meet) and ∪ (join), which satisfy the following results:*

$$L_1 \quad x \cap y = y \cap x \qquad\qquad x \cup y = y \cup x$$
$$L_2 \quad x \cap (y \cap z) = (x \cap y) \cap z \qquad (x \cup y) \cup z = x \cup (y \cup z)$$
$$L_3 \quad x \cap (y \cup x) = x \qquad\qquad x \cup (x \cap y) = x$$

*Two applications of ($L_3$) namely $x \cap x = x \cap (x \cup (x \cap x)) = x$ lead to $x \cap x = x$ and $x \cup x = x$. $L_1$ is the commutative law, $L_2$ is the associative law, $L_3$ is the absorption law and $L_4$ is the idempotent law.*

**DEFINITION 1.3.7**: *A lattice L is called modular if for all x, y, z ∈ L.*

*x ≤ z imply x ∪ (y ∩ z) = (x ∪ y) ∩ z (modular equation).*

***Result 1.3.1***: The lattice given in the following figure is known as pentagon lattice:

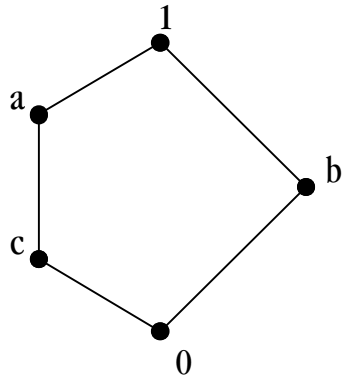

Figure 1.3.1

which is not modular.

***Result 1.3.2***: The lattice known as diamond lattice (given by figure 1.3.2) is modular.

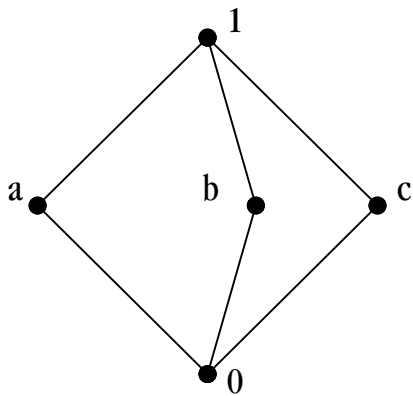

Figure 1.3.2

**DEFINITION 1.3.8**: *A lattice L is called distributive if either of the following conditions hold for all x, y, z in L.*

$$x \cup (y \cap z) = (x \cup y) \cap (x \cup z)$$
$$x \cap (y \cup z) = (x \cap y) \cup (x \cap z).$$

The lattice given in Figure 1.3.2 is the smallest modular lattice which is not distributive.

**DEFINITION 1.3.9**: *A non-empty subset S of a lattice L is called a sublattice of L if S is a lattice with respect to the restriction of $\cap$ and $\cup$ of L onto S.*

**<u>Result 1.3.3</u>**: Every distributive lattice is modular.

Proof is left for the reader as an exercise.

**<u>Result 1.3.4</u>**: A lattice is modular if and only if none of its sublattices is isomorphic to the pentagon lattice.

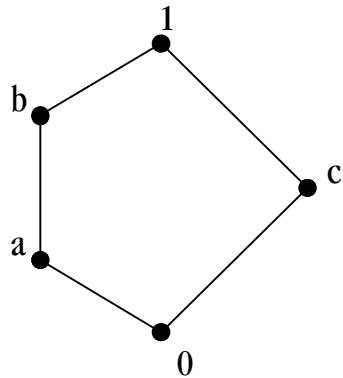

Figure 1.3.3

We leave the proof as an exercise to the reader.

Now we give some problems:

**<u>PROBLEMS:</u>**

1. Prove the lattice given in figure 1.3.4 is distributive.

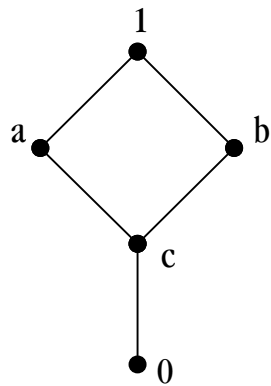

Figure 1.3.4

2. Prove the lattice given by Figure 1.3.5. is non-modular.

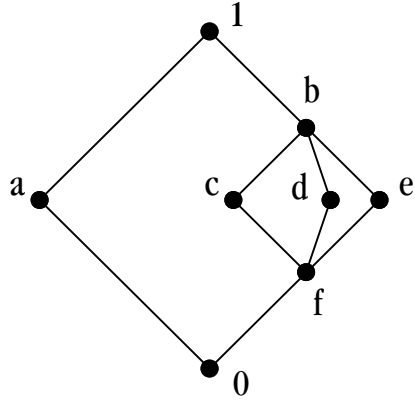

Figure 1.3.5

3. Is this lattice

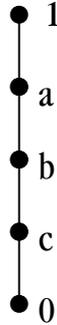

Figure 1.3.6

modular ?

4.

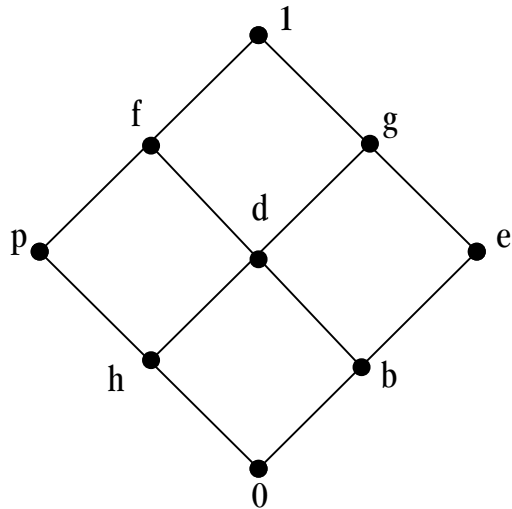

Figure 1.3.7

Is this lattice modular? distributive?

5. Give a modular lattice of order nine which is non-distributive.

## 1.4 Smarandache semigroups

In this section we introduce the notion of Smarandache semigroups (S-semigroups) and illustrate them with examples. The main aim of this is that we want to define which of the group rings and semigroup rings are Smarandache rings, while doing so we would be needing the concept of Smarandache semigroups. As the study of S-semigroups is very recent one, done by F. Smarandache, Padilla Raul and W.B. Vasantha Kandasamy [73, 60, 154, 156], we felt it is appropriate that the notion of S-semigroups is substantiated with examples.

DEFINITION [73, 60]: *A Smarandache semigroup (S-semigroup) is defined to be a semigroup A such that a proper subset A is a group (with respect to the induced operation on A).*

DEFINITION [154, 156]: *Let A be a S-semigroup. A is said to be a Smarandache commutative semigroup (S-commutative semigroup) if the proper subset of A which is a group is commutative. If A is a commutative semigroup and if A is a S-semigroup then A is obviously a S-commutative semigroup.*

*Example 1.4.1*: Let $Z_{12} = \{0, 1, 2, \ldots, 11\}$ be the semigroup under multiplication modulo 12. It is a S-semigroup as the proper subset P = {3, 9} is a group with 9 as unit; that is the multiplicative identity. That is P is a cyclic group of order 2.

*Example 1.4.2*: Let S(5) be the symmetric semigroup is a S-semigroup, as $S_5 \subset S(5)$ is the proper subset that is a symmetric group of degree 5. Further S(5) is a S-commutative semigroup as the element

$$p = \begin{pmatrix} 1 & 2 & 3 & 4 & 5 \\ 2 & 3 & 4 & 5 & 1 \end{pmatrix}$$

generates a cyclic group of order 5.

DEFINITION [154, 156]: *Let S be a S-semigroup. A proper subset X of S which is a group under the operations of S is said to be a Smarandache normal subgroup (S-normal subgroup) of the S-semigroup, if aX $\subseteq$ X and Xa $\subseteq$ X or aX = {0} and Xa = {0} for all x $\in$ S, if 0 is an element in S.*

*Example 1.4.3*: Let $Z_{10} = \{0, 1, 2, \ldots, 9\}$ be the S-semigroup of order 10 under multiplication modulo 10. The set X = {2, 4, 6, 8} is a subgroup of $Z_{10}$ which is a S-normal subgroup of $Z_{10}$.

## PROBLEMS:

1. Show $Z_{15}$ is a S-semigroup. Can $Z_{15}$ have S-normal subgroups?
2. Let $S(8)$ be the symmetric semigroup, prove $S(8)$ is a S-semigroup. Can $S(8)$ have S-normal subgroups?
3. Find all S-normal subgroups of $Z_{24} = \{0, 1, 2, \ldots , 23\}$, the semigroup of order 24 under multiplication modulo 24.
4. Give an example of a S-non-commutative semigroup.
5. Find the smallest S-semigroup which has nontrivial S-normal subgroups.
6. Is $M_{3\times3} = \{(a_{ij}) \ / \ a_{ij} \in Z_3 = \{0,1,2\}\}$ a semigroup under matrix multiplication; a S-semigroup?
7. Can $M_{3\times3}$ given in problem 6 have S-normal subgroup? Substantiate your answer.
8. Give an example of a S-semigroup of order 18 having S-normal subgroup.
9. Can a semigroup of order 19 be a S-semigroup having S-normal subgroups?
10. Give an example of a S-semigroup of order p, p a prime.

**Chapter two**

# RINGS AND THEIR PROPERTIES

In this chapter we recollect some of the basic properties of rings. This Chapter is organized into seven sections. In section one we just recall the definition of ring and give some examples. Section two is devoted to the study of special elements like zero divisors, units, idempotents nilpotents etc. Study of substructures like subrings, ideals and Jacobson radical are introduced in section three. Recollection of the concept of homomorphisms and quotient rings are carried out in section four. Special rings like polynomical rings, matrix rings, group rings etc are defined in section five. Section six introduces modules and the final section is completely devoted to the recollection of the rings which satisfy chain conditions. Every section ends with a list of problems to be solved by the reader. Finally no claim is made that we have recaptured all facts about rings we do not do it in fact the reader is expected to be well versed in ring theory.

## 2.1 Definition and Examples

In this section we recall the definition of rings and their basic properties and illustrate them with examples. Also the definition of field, integral domain and division ring are given.

**DEFINITION 2.1.1**: *A non-empty set R is said to be an associative ring if in R are defined two binary operations '+' and '.' respectively such that*

1. *(R, +) is an additive abelian group.*
2. *(R, .) is a semigroup.*
3. *a . (b + c) = a . b + a . c and*
   *(a + b) . c = a . c + b . c for*
   *all a, b, c ∈ R (the two distributive laws).*

*It may very well happen that (R, .) is a monoid, that is there is an element 1 in R such that a . 1 = 1 . a = a for every a ∈ R, in such cases we shall describe R as a ring with unit element.*

*If the multiplication in R is such that a . b = b . a for every a, b in R, then we call R a commutative ring, if a . b ≠ b . a atleast for a pair in R then R is a non-commutative ring.*

Henceforth, we simply represent a . b by ab.

***Example 2.1.1***: Let Z be the set of integers, positive, negative and 0; Z is a commutative ring with 1.

***Example 2.1.2***: Let $Z_n = \{0, 1, 2, \dots, n-1\}$ be the ring of integers modulo n. $Z_n$ is a ring under modulo addition and multiplication. $Z_n$ is a commutative ring with unit.

***Example 2.1.3***: Let $M_{n \times n} = \{(a_{ij}) / a_{ij} \in Z\}$, the set all $n \times n$ matrices with matrix addition and multiplication. $M_{n \times n}$ is a non-commutative ring with unit element.

**DEFINITION 2.1.2**: *Let (R, +, .) be a ring, if (R \ {0}, . ) is an abelian group we call R a field.*

**Notation**: Z − denotes the set of integers positive, negative and zero. Q − denotes the set of positive and negative rationals with zero R − denotes the set of reals, positive, negative with zero. $Z_n$ − set of integers modulo n. $Z_n = \{0, 1, 2, \dots, n-1\}$, $Z_p$ − set of integers modulo p, p − prime, Set of complex number of the from a+ib, a, b ∈ R or Q or Z is denoted by C.

**DEFINITION 2.1.3**: *If a ring R has a finite number of elements we say R is a finite ring, otherwise R is an infinite ring.*

**DEFINITION 2.1.4**: *Let R be a ring if mx = x + ...+ x (m-times) is zero for every x ∈ R, m a positive integer then we say characteristic of R is m. If for no m the result is true we say the characteristic of R is 0, denoted by characteristic R is 0 or characteristic R is m.*

<u>Note</u>: The rings given in examples 2.1.1 and 2.1.3 are of characteristic zero where as the ring in example 2.1.2 is of characteristic n.

***Example 2.1.4***: Let $Z_9 = \{0, 1, 2, \dots, 8\}$. This is a commutative finite ring of characteristic 9 with unit 1.

**DEFINITION 2.1.5**: *Let R be a ring, we say a ≠ 0 ∈ R is a zero divisor, if there exists b ∈ R, b ≠ 0, such that a.b = 0.*

***Example 2.1.5***: The ring $Z_{15} = \{0, 1, 2, \dots, 14\}$ is of characteristic 15. Clearly for $3 \neq 0 \in Z_{15}$ we have $5 \in Z_{15}$ such that $3.5 \equiv 0 \mod(15)$ thus $Z_{15}$ has zero divisor.

But the ring given in example 2.1.1 has no zero divisors.

**DEFINITION 2.1.6**: *Let R be a commutative ring with unit. If R has no zero divisors we say R is an integral domain. (The presence of unit is not a must).*

The ring Z given in example 2.1.1 is an integral domain.

**DEFINITION 2.1.7**: *Let R be a non-commutative ring in which the non-zero elements form a group under multiplication, then R is a division ring.*

***Example 2.1.6***: Let P be the set of symbols of the form $\alpha_0 + \alpha_1 i + \alpha_2 j + \alpha_3 k$ where all the numbers $\alpha_0$, $\alpha_1, \alpha_2$ and $\alpha_3$ are real numbers. We declare two such symbols $\alpha_0 + \alpha_1 i + \alpha_2 j + \alpha_3 k$ and $\beta_0 + \beta_1 i + \beta_2 j + \beta_3 k$ to be equal if and only if $\alpha_t = \beta_t$ for t = 0, 1, 2, 3. In other words to make P into a ring we must define a '+' and a '.' for its elements.

To this end for any $X = \alpha_0 + \alpha_1 i + \alpha_2 j + \alpha_3 k$ and $Y = \beta_0 + \beta_1 i + \beta_2 j + \beta_3 k$ define $X + Y = (\alpha_0 + \alpha_1 i + \alpha_2 j + \alpha_3 k) + (\beta_0 + \beta_1 i + \beta_2 j + \beta_3 k) = (\alpha_0 + \beta_0) + (\alpha_1 + \beta_1)i + (\alpha_2 + \beta_2)j + (\alpha_3 + \beta_3)k$ and $X \cdot Y = (\alpha_0 + \alpha_1 i + \alpha_2 j + \alpha_3 k) (\beta_0 + \beta_1 i + \beta_2 j + \beta_3 k) = (\alpha_0\beta_0 - \alpha_1\beta_1 - \alpha_2\beta_2 - \alpha_3\beta_3) + (\alpha_0\beta_1 + \alpha_1\beta_0 + \alpha_2\beta_3 - \alpha_3\beta_2)i + (\alpha_0\beta_2 + \alpha_2\beta_0 + \alpha_3\beta_1 - \alpha_1\beta_3)j + (\alpha_0\beta_3 + \alpha_3\beta_0 + \alpha_1\beta_2 - \alpha_2\beta_1)k$.

We use in the product the following relation $i^2 = j^2 = k^2 = -1 = ijk$, $ij = -ji = k$, $jk = -kj = i$, $ki = -ik = j$ Notice $\pm i$, $\pm j$, $\pm k$, $\pm 1$ form a non-abelian group of order 8 under multiplication. $0 + 0i + 0j + 0k = 0$ acts as the additive identity 0 of P. $1 = 1 + 0i + 0j + 0k$ serves as the unit. If $X = \alpha_0 + \alpha_1 i + \alpha_2 j + \alpha_3 k$ then its inverse

$$Y = \frac{\alpha_0}{\beta} - \frac{\alpha_1 i}{\beta} - \frac{\alpha_2 j}{\beta} - \frac{\alpha_3 k}{\beta}$$

where $\beta = \alpha_0^2 + \alpha_1^2 + \alpha_2^2 + \alpha_3^2$. Clearly $X \cdot Y = 1$. Thus it can be verified as $ij \neq ji$, we get a division ring.

<u>**Result 1**</u>: Every commutative division ring is a field. Left for the reader to prove.

<u>**Result 2**</u>: A finite integral domain is a field.

It is left for the reader to verify these results.

***Example 2.1.7***: Let $S_{3 \times 3} = \{(a_{ij}) / a_{ij} \in Z_{10}\}$ be the set of all $3 \times 3$ matrices; $S_{3 \times 3}$ is not a division ring but is a non-commutative ring of characteristic 10.

<u>**PROBLEMS:**</u>

1. Give an example of a commutative ring of order 12. (where by order of the ring R we mean the number of elements in R, denote by o(R) or |R|).
2. What is the order of the smallest non-commutative ring?

3. Can a ring or order 11 be non-commutative?
4. Find the zero divisors in the ring $Z_{30} = \{0, 1, 2, \ldots, 29\}$.
5. How many elements does the ring $M_{2 \times 2} = \{(a_{ij}) \; / \; a_{ij} \in Z_4 = \{0, 1, 2, 3\}\}$; (i.e., set of all $2 \times 2$ matrices with entries from $Z_4$) contain?
   a. Find zero divisors in $M_{2 \times 2}$.
   b. Find units in $M_{2 \times 2}$.
   c. Show $AB \neq BA$ atleast for a pair $A, B \in M_{2 \times 2}$.
6. Give an example of a ring of characteristic 0 which has zero divisors.
7. Find a non-commutative ring of finite order other than the matrix ring.
8. Does there exist a division ring of characteristic 0?
9. Does there exist a division ring of characteristic n, n a non-prime?
10. Find all zero divisors in the ring $Z_{25}$.
11. Find a ring of order 10 which has no unit.

## 2.2 Special Elements in Rings.

In this section we mainly introduce the concept of units, idempotents, zero divisors and regular elements, we just recall the definition of these concepts and illustrate them with examples. All properties and results related to these concepts are left for the reader to refer, books on ring theory.

**DEFINITION 2.2.1**: *Let R be a ring, an element $x \in R \setminus \{0, 1\}$ (if 1 is in R) is called an idempotent in R if $x^2 = x$ for $x \in R$.*

**Example 2.2.1**: Let $Z_{12} = \{0, 1, 2, \ldots, 11\}$ be the ring of integers modulo 12. We see $4^2 \equiv 4 \pmod{12}$ is an idempotent in it.

**Example 2.2.2**: Let $M_{n \times n} = \{(a_{ij}) \; / \; a_{ij} \in Q$ − the field of rationals$\}$; $M_{n \times n}$ is a ring under matrix addition and matrix multiplication. We have matrices $A \in M_{n \times n}$ such that $A^2 = A$.

For example take n = 3,

$$A = \begin{bmatrix} 1 & 0 & 0 \\ 0 & 0 & 0 \\ 0 & 0 & 1 \end{bmatrix}$$

is such that $A^2 = A$. Thus we have seen idempotents in case of both a commutative and a non-commutative rings.

**_Result_**: Let R be a ring with 1. If R has a nontrivial idempotent then we have nontrivial divisors of zero.

For let $x \in R \setminus \{0, 1\}$ such that $x^2 = x$ so $x^2 - x = 0$ i.e., $x(x - 1) = 0$ as $x \neq 0$ and $x \neq 1$, we have nontrivial zero divisors. We call an element nilpotent if $x^n = 0$ where $x \neq 0 \in R$ and $n \geq 2$.

**DEFINITION 2.2.2**: *Let R be a ring with 1. If for $x \in R \setminus \{0\}$ there exists a y in R with x.y = 1 we say R has units or invertible elements.*

**_Example 2.2.3_**: Let Q be the field of rationals every element in $Q \setminus \{0\}$ is a unit.

**_Example 2.2.4_**: Let $Z_{15} = \{0, 1, 2, \ldots, 14\}$ be the ring of integers modulo 15, we see $14^2 \equiv 1 \pmod{15}$, $4^2 \equiv 1 \pmod{15}$, $8.2 \equiv 1 \pmod{15}$. Thus $Z_{15}$ has nontrivial units but not all elements in $Z_{15}$ are units.

**_Example 2.2.5_**: Let $M_{5 \times 5} = \{(a_{ij}) \, / \, a_{ij} \in Q\}$ be the ring of matrices. Clearly all matrices $A \in M_{5 \times 5}$ are such that A is non-singular that is $|A| \neq (0)$ are invertible.

**DEFINITION 2.2.3**: *Let R be a ring if for $s \in R$ we have $sr \neq 0$ and $rs \neq 0$ for all $r \neq 0 \in R$; then we say s is a regular element of R.*

For instance all elements in an integral domain or a field are regular elements.

**PROBLEMS:**

1. Find all idempotents, zero divisors and units in $Z_{35}$.
2. Find the zero divisors and regular elements of the ring $M_{2 \times 2} = \{(a_{ij}) \, / \, a_{ij} \in Z_2 = \{0, 1\}\}$; where $M_{2 \times 2}$ is the matrix ring.
3. Find all the regular elements in $Z_{24}$.
4. Find only the idempotent matrices of $M_{3 \times 3} = \{(a_{ij}) \, / \, a_{ij} \in Z_3\}$.
5. How many regular elements are in $M_{2 \times 2}$ given in problem 2?
6. Does $Z_{16} = \{0, 1, 2, \ldots, 15\}$ have nilpotents of order 6?
7. Can a matrix A in $M_{3 \times 3}$ given in problem 4 have nilpotent elements of order 5? Justify your answer.
8. Give zero divisors in $Z_{12}$, which are not nilpotents. (for example $6^2 \equiv 0 \pmod{12}$).
9. Can a ring R have only nilpotent element as zero divisor? Justify your answer.
10. Find all regular elements, nilpotents, zero divisors, idempotents and units of the ring $Z_{210}$.

## 2.3 Substructures of a Ring.

In this section we introduce the concept of ideals, subrings and radicals for rings. We only recall the very basic definitions and illustrate them with examples.

**DEFINITION 2.3.1**: *Let R be a ring, a proper subset S or R is said to be a subring of R if S itself under the operations of R is a ring. Clearly {0} is a subring.*

**Example 2.3.1**: Let $Z_{15} = \{0, 1, 2, \ldots, 14\}$ be the ring of integers modulo 15. S = $\{3, 6, 9, 12, 0\}$ is a subring of $Z_{15}$.

**Example 2.3.2**: Let Z be the ring of integers, $nZ = \{0, \pm n, \pm 2n, \ldots\}$, is a subring of Z, n any positive integer.

**Example 2.3.3**: Let $M_{3 \times 3} = \{(a_{ij}) \,/\, a_{ij} \in Z_4 = \{0, 1, 2, 3\}\}$. Clearly

$$P = \left\{ \begin{pmatrix} a_{ij} & 0 & 0 \\ 0 & 0 & 0 \\ 0 & 0 & 0 \end{pmatrix} \middle/ a_{ij} \in Z_4 = \{0,1,2,3\} \right\}$$

is a subring of $M_{3 \times 3}$.

**DEFINITION 2.3.2**: *Let R be a ring. A nonempty subset I of R is said to be the right (left) ideal of R if*

      *1. I is a subgroup of R under addition.*
      *2. For all $r \in R$ and $s \in I$; $rs \in I$. ($sr \in I$).*

*I is called an ideal; if I is simultaneously both a right and a left ideal of R. If R is a commutative ring naturally the concept of right and left ideals coincide.*

**Example 2.3.4**: Let Z be the ring of integers; $pZ = \{0, \pm p, \pm 2p, \ldots\}$ is an ideal of Z. It is to be noted that in any ring R, (0) is an ideal of R; we will call (0) and R as trivial ideals of R.

**Example 2.3.5**: Let $Z_{22} = \{0, 1, 2, \ldots, 21\}$ be the ring of integers modulo 22.

Clearly $I = \{0, 11\}$ is an ideal of $Z_{22}$. Also $P = \{0, 2, 4, 6, 8, \ldots, 20\}$ is an ideal of $Z_{22}$.

**Example 2.3.6**: Let $Z_7 = \{0, 1, 2, \ldots, 6\}$, this is a ring. Clearly $Z_7$ has no ideals as $Z_7$ is a prime field of characteristic 7.

The student is expected to note that fields F have no nontrivial ideals. The only trivial ideals of F are {0} and F.

***Example 2.3.*7**: Let $M_{2\times2}$ = {($a_{ij}$) / $a_{ij} \in Z_2$ = {0, 1}} be the ring. Can $M_{2\times2}$ have ideals? It is left as an exercise to find ideals in $M_{2\times2}$.

**DEFINITION 2.3.3**: *Let R be a ring. I an ideal of R, I is said to be a principal ideal of R, if it is generated by a single element.*

***Example 2.3.8***: Let Z be the ring of integers, every element p in Z generates an ideal pZ, which is principal.

***Example 2.3.9***: Let $Z_{25}$ = {0, 1, 2, … , 24} be the ring of integers modulo 25. $\langle 5 \rangle$ = {0, 5, 10, 15, 20} is an ideal of $Z_{25}$ ('$\langle 5 \rangle$' denotes the ideal generated by 5.) which is principal.

**DEFINITION 2.3.4**: *Let R be a ring, I an ideal of R. I is said to be a maximal ideal of R; if J is an ideal of R such that I $\subset$ J $\subset$ R, then either I = J or J = R. We similarly define an ideal P of a ring R to be minimal, if S is an ideal of R such that (0) $\subset$ S $\subset$ P then either (0) = S or S = P.*

*A proper ideal P of a ring R is called prime if for xy $\in$ P we have x $\in$ P or y $\in$ P.*

***Example 2.3.10***: Let Z be the ring of integers. P = 8Z is an ideal. P is not a prime ideal as 4.2 $\in$ P but both 2 and 4 are not in P.

**DEFINITION 2.3.5**: *Let R be a ring. The intersection of all maximal ideals of a commutative ring is called the radical of the ring R denoted by rad (R). This is called the Jacobson radical of R. rad R={r $\in$ R / 1 - rx is a unit for all x$\in$R}. Thus the radical is the largest ideal K of R such that for all r $\in$ R, 1 − r is a unit.*

**DEFINITION 2.3.6**: *An ideal I of a ring R is said to be a nil ideal of R if every element of I is nilpotent. An ideal I is nilpotent if $I^n$ = 0 for some n $\geq$ 1 by $I^n$ = I . I… I, $I^2$ = I. I = {$\Sigma x_i y_i / x_i, y_i \in I$} similarly for any power of n.*

**DEFINITION 2.3.7**: *A ring R is simple if it has no two sided ideals other than (0) and R. It is interesting to note that all fields are trivially simple rings.*

**PROBLEMS:**

1. Find all ideals of $Z_{124}$.
2. Can the ring $Z_{24}$ have Jacobson radical?

3. Find all maximal ideals of $Z_{125}$.
4. Find a ring in which an ideal which is simultaneously maximal and minimal.
5. Find two right ideals of $M_{n \times n} = \{(a_{ij}) / a_{ij} \in Z_{12}\}$ which are not left ideals.
6. Let $Z_{210} = \{0, 1, 2, \ldots, 209\}$ be the ring of integers modulo 210 find
   a. Jacobson radical of $Z_{210}$.
   b. Maximal ideal.
   c. Minimal ideal.
   d. Is every ideal principal?
   e. Does $Z_{210}$ have prime ideals?
7. Find subrings which are not ideals in Q.
8. Can $Z_{210}$ given in problem 6 have subrings which are not ideals?
9. Find ideals and subrings of $Z_{25}$. Are they identical?
10. Find subrings which are not ideals in $M_{3 \times 3} = \{(a_{ij}) / a_{ij} \in Z_6 = \{0, 1, \ldots, 5\}\}$.

## 2.4 Homomorphism and Quotient Rings

In this section we recall the basic concepts of homomorphism and quotient rings and give some examples.

**DEFINITION 2.4.1**: *Let R and S be two rings. A mapping f: R → S is called a homomorphism of rings if for all a, b ∈ R we have (1) f(a+b) = f(a) + f(b) and (2) f(ab) = f(a) f(b). If f is a homomorphism, it is easy to verify f(0)=0, f(-x) = -f(x), and f($1_R$)= $1_S$; in case both rings have identity. In case f($1_R$)= $1_S$; we say the map f is unitary.*

**DEFINITION 2.4.2**: *Let f: R → S be a ring homomorphism, the kernel of the homomorphism f is defined to be the set = {x ∈ R / f(x) = 0} and is denoted by Ker f = {x ∈ R / f(x)=0}.*

*A ring homomorphism f: R → S is called*

   *1) a monomorphism if f is injective*
   *2) an epimorphism if f is surjective*
   *3) an isomorphism if f is bijective*
   *4) an endomorphism if R = S and*
   *5) an automorphism if R = S and f is an isomorphism or equivalently; we can say a homomorphism f: R → S is a monomorphism if and only if Ker (f) = {0}.*

*Clearly if f: R → S is a ring homomorphism.*

It is left as an exercise for the reader to verify that ker f is always a two sided ideal of R.

**DEFINITION 2.4.3**: *Let R be a ring, I be a two sided ideal of R, we make R / I = {a + I / a ∈ R} into a ring called the quotient ring of R by defining operations '+' and '.' as follows.*

$$(a+I) + (b+I) = (a+b) + I \text{ for all } a, b \in R$$
$$(a+I) + (-a+I) = I.$$

*So R/I is a group under addition, a+I = I for all a ∈ I so I is the additive identity of R/I.*

$$(r+I) I = I (r+I) = I \text{ for all } r \in R$$
$$(a+I) (b+I) = ab + I.$$

*So (R / I, +, .) is a ring called the quotient ring.* (Here the distributive laws are left for the reader to verify).

For any ring homomorphism f: R→S, kernel f denoted by ker f is an ideal of R and R/ker f is a ring.

Several properties about quotient rings exists the nice among them is R/I is a field if and only if I is a maximal ideal in R. If I is a maximal ideal we call R/I the residue field of R at I.

**DEFINITION 2.4.4**: *A ring R with 1 is called a local ring if the set of all non-units in R is an ideal.*

All division rings are local rings.

<u>**PROBLEMS:**</u>

1. Find a ring homomorphism $\phi$ between $Z_{20}$ and $Z_{18}$ such that the ker $\phi \neq \{0\}$.

2. Let f: $Z_{25} \rightarrow Z_{16}$ be a ring homomorphism find the quotient ring $\dfrac{Z_{25}}{\ker f}$.

3. Let $Z_{36}$ = {0, 1, 2, …, 35} be the ring of integers modulo 36. Let I = {2, …, 34, 0} and J = {3, 9, …, 33, 0} be ideals of $Z_{36}$. Find the quotient rings $Z_{36}$/I and $Z_{36}$/J.

4. Let $Z_{21}$ = {0, 1, 2, …, 20}. Find an ideal I of $Z_{21}$ such that $Z_{21}$/I is a field.

5. Prove for $Z_{12}$={0, 1, 2, …, 11}, the ring integers with I = {0, 6}, the quotient ring $Z_{12}$/ I is not a field.

6. Show in any ring R we can have several quotient rings related to different ideals. Illustrate them by an example.

7. How many quotient rings can be constructed for the ring Z?.

8. Give an example of a finite local ring.

9. Find a homomorphism from f: $Z_{27} \rightarrow Z_{18}$ such that $Z_{27}$ / ker $\phi$ is isomorphic to $Z_{18}$.

10. Let $Z_{23}$ and $Z_{19}$ be two rings. Is it possible to find a homomorphism $\phi$ from $Z_{23}$ to $Z_{19}$ such that $Z_{23}$ /ker $\phi \cong Z_{19}$. Justify your answer.

## 2.5 Special Rings

In this section we just recall the four types of rings which are specially formed and illustrate them with examples. They are polynomial rings, matrix rings, direct product of rings, ring of Gaussian integers, group rings and semigroup rings. Examples of these rings will help in the study of Smarandache ring. Throughout this section by a ring we mean only ring with unit, which is commutative.

**DEFINITION 2.5.1**: *Let R be a commutative ring with unit 1, x be an indeterminate,*

$$R[x] = \left\{ \sum_{i=0}^{n} a_i x^i \middle/ a_i \in R, \ n \ \text{can be a finite or an infinite integer} \ i \geq 0 \right\}. \quad x^0 \ \text{is}$$

*defined to be 1.*

*Let $p(x) = p_0 + p_1 x + \ldots + p_n x^n$ and $q(x) = q_0 + q_1 x + \ldots + q_m x^m$. be elements in R[x]. We say p(x) = q(x) if and only if m = n and $p_i = q_i$ for all i, $0 \leq i \leq r$. In particular $a_0 + a_1 x + \ldots + a_n x^m = 0$ if and only if each $a_i = 0$.*

*Define addition in R[x] by $p(x) + q(x) = (p_0 + p_1 x + \ldots + p_n x^n) + (q_0 + q_1 x + \ldots + q_m x^m) = (p_0 + q_0) + (p_1 + q_1)x + \ldots + q_m x^m$ if m > n*

*$p(x) q(x) = (p_0 + p_1 x + \ldots + p_n x^n) (q_0 + q_1 x + \ldots + q_m x^m) = p_0 q_0 = (p_0 q_1 + p_1 q_0)x + (p_0 q_2 + p_1 q_1 + p_2 q_0) x^2 + \ldots + p_n q_m x^{m+n}$.*

It can be easily verified that R[x] is a commutative ring with unit 1. Elements of R[x] are called polynomials and R[x] is a polynomial ring in the indeterminate x with coefficients from R.

***Example 2.5.1***: Let Z be the ring of integers. Z[x] is a polynomial ring in the variable x. Z[x] is an integral domain.

***Example 2.5.2***: Let $Z_n$ be the ring of integers modulo n. $Z_n[x]$ is a polynomial ring with coefficients from $Z_n$. $Z_n[x]$ is a commutative ring with 1 and has zero divisors if n is a non-prime.

It is interesting to note in the polynomial ring Z[x] every ideal is principal.

Polynomial rings in several variables can also be defined in a similar way. For if R[x] is a polynomial ring. Suppose y is another indeterminate then (R[x]) [y] is a polynomial ring using the commutative ring R[x] as the ring of quotients for the inderterminate y, (we assume xy = yx) denoted by R[x, y]. Suppose $x_1, \ldots, x_n$ are n variables then the polynomial ring in n variables is R[$x_1, x_2, \ldots, x_n$] where we assume $x_i x_j = x_j x_i$ for $1 \leq i, j \leq n$.

**DEFINITION 2.5.2**: *Let R be a commutative ring with 1, and $n \geq 1$ be an integer.*

*$M_{n \times n} = \{(a_{ij}) / a_{ij} \in R; \ 1 \leq i, j \leq n\}$ be the set of all $n \times n$ matrices with entries in R where*

$$(a_{ij}) = \begin{pmatrix} a_{11} & a_{12} & . & . & . & a_{1n} \\ a_{21} & a_{22} & . & . & . & a_{2n} \\ . & . & & & & . \\ . & . & & & & . \\ . & . & & & & . \\ a_{n1} & a_{n2} & . & . & . & a_{nn} \end{pmatrix}$$

*Define addition and multiplication in $M_{n \times n}$ as follows: $(a_{ij}) + (b_{ij}) = (a_{ij} + b_{ij})$ and $(a_{ij}) \cdot (b_{ij}) = (c_{ij})$ where $c_{ij} = \sum_{k=1}^{n} a_{ik} b_{kj}$ for all i, j; $1 \leq i, j \leq n$. It is easily verified that $M_{n \times n}$ is a ring called the matrix ring of order with entries from R.*

$$I_{n \times n} = \begin{pmatrix} 1 & 0 & 0 & . & . & . & 0 \\ 0 & 1 & 0 & . & . & . & 0 \\ . & & & & & & . \\ . & & & & & & . \\ . & & & & & & . \\ 0 & 0 & 0 & . & . & . & 1 \end{pmatrix}$$

*is called the $n \times n$ identity matrix of $M_{n \times n}$.*

It is interesting to note that matrix ring $M_{n \times n}$ is non-commutative and has zero divisors idempotents, nilpotents and also units.

***Example 2.5.3***: Let $M_{2 \times 2} = \{(a_{ij}) \, / \, a_{ij} \in Z_5 = \{0, 1, 2, 3, 4\}\}$ is a ring . Find units, idempotents and zero divisors in $M_{2 \times 2}$.

**DEFINITION 2.5.3**: *Let R and S be any two rings (not necessarily both R and S should be commutative rings with unit). P = R $\times$ S be the cartesian product of R and S. Define addition and multiplication on P. (x, y) + (x$_1$, y$_1$) = (x + x$_1$, y + y$_1$) and (x, y) (x$_1$, y$_1$) = (x x$_1$, y y$_1$) under these operations it is easily verified P is a ring called the direct product of the rings R and S.*

*If we take n rings say R$_1$, R$_2$, …, R$_n$ define P = R$_1$ $\times$… $\times$R$_n$ = {(r$_1$, r$_2$, …, r$_n$) / r$_i$ $\in$ R$_n$, i = 1, 2, …, r$_n$} is the direct product of the n rings R$_1$, R$_2$, …, R$_n$.*

***Example 2.5.4***: Let $P = Z_2 \times Z_5 = \{(a, b) \, / \, a \in Z_2 \text{ and } b \in Z_5\}$. P is a direct product of rings with 10 elements and has nontrivial zero divisors.

***Example 2.5.5***: Let $S = Z_2 \times Z \times Z_9$ the direct product of ring. S is an infinite ring with zero divisors, S is a commutative ring.

**DEFINITION 2.5.4**: *Consider the subset of C (the complex field) given by Z[i] = {a + ib / a, b $\in$ Z}. This is the set of integral points whose both coordinates are integers. It is easily verified Z[i] is a ring called the ring of Gaussian integers where addition and multiplication are given by (a + ib) + (c + id) = (a + c, i(b + d)) and (a + ib) (c + id) = (ac – bd, i(ad + bc)). The unity of Z[i] is 1.*

**DEFINITION 2.5.5**: *Consider the set Q[i] = {a + ib / a, b $\in$ Q}, Q the field of rationals. It is easily verified Q[i] is a ring called the ring of Gaussian numbers; in fact Q[i] is a field. It is easy to verify. Z $\subset$ Z[i] $\subset$ Q[i] $\subset$ R [i] = C.*

For more about properties of Gaussian rings refer [15].

The reader may just recall the concept of integral quaternions which was introduced in section 2.1. Now we introduce the concept of group rings and semigroup rings.

**DEFINITION 2.5.6**: *Let R be a commutative ring with unit 1 and G be a multiplicative group. The group ring, RG of the group G over the ring R consists of all finite formal sums of the form $\sum_i \alpha_i g_i$ (i-runs over a finite number) where $\alpha_i \in R$ and $g_i \in G$ satisfying the following conditions:*

i)   $$\sum_{i=1}^{n} \alpha_i g_i = \sum_{i=1}^{n} \beta_i g_i \iff \alpha_i = \beta_i \text{ for } i = 1, 2, \ldots, n, g_i \in G.$$

ii)  $$\left(\sum_{i=1}^{n} \alpha_i g_i\right) + \left(\sum_{i=1}^{n} \beta_i g_i\right) = \sum_{i=1}^{n} (\alpha_i + \beta_i) g_i \, ; g_i \in G.$$

iii) $$\left(\sum_{i} \alpha_i g_i\right)\left(\sum_{j} \beta_i g_i\right) = \sum_{k} \gamma_k m_k \text{ where } \gamma_k = \sum \alpha_i \beta_j, \ g_i b_j = m_k.$$

iv)  $$r_i m_i = m_i r_i \text{ for all } r_i \in R \text{ and } m_i \in G.$$

v)   $$r \sum_{i=1}^{n} r_i g_i = \sum_{i=1}^{n} (r r_i) g_i \text{ for } r_i \, r \in R \text{ and } \sum r_i g_i \in RG.$$

*RG is a ring with $0 \in R$ as its additive identity. Since $1 \in R$ we have $G = 1.G \subset G$ and $R.e = R \subseteq RG$ where $e$ is the identity of $G$. Clearly if we replace the group $G$ by a semigroup $S$ we say $RS$ is the semigroup ring of the semigroup $S$ over the ring $R$.*

***Example 2.5.6***: Let $Z_2 = \{0, 1\}$ be the ring and $G = \langle g \, / \, g^5 = 1 \rangle$, the group ring $Z_2 G$ is a ring which is commutative and has zero divisors. For $g^5 + 1 = (g + 1) (1 + g + g^2 + g^3 + g^4) = 0$.

It is now important to mention if R is a ring and G is any finite group or has elements of finite order than the group ring RG has nontrivial zero divisors.

If G is a torsion free abelian group and K a field of characteristic zero, the group ring KG has no zero divisors. It is pertinent to mention here till date i.e., even after 60 years the problem if K is a field of characteristic zero and G a torsion free non-abelian group; can the group ring KG have zero divisors remains open, proposed in 1940 by G Higman [33].

PROBLEMS:

1. Find ideals in $Z_7[x]$, the polynomial ring in the variable x.
2. Can $Z_8[x]$ have zero divisors? Find a maximal ideal in $Z_8[x]$.
3. Let $M_{2\times 2} = \{(a_{ij}) \, / \, a_{ij} \in Z_8\}$ be the matrix ring. Find
   a. Idempotents in $M_{2\times 2}$.
   b. Ideals (right only).

     c. Zero divisors.

     d. Units.

     e. Subrings which are not ideals.

4. Let $G = S_3$ and $Z_3 = \{0, 1, 2\}$ find in the group ring $Z_3S_3$.

     a. Zero divisors.

     b. Ideals.

     c. Units.

     d. Left ideals.

     e. Idempotents.

     f. What is the order of $Z_3S_3$?

5. Let $G = S(4)$ be the semigroup $Z_2 = \{0, 1\}$ be the field of characteristic 2. Let $Z_2G$ be the semigroup ring of the semigroup G over $Z_2$. Find

     a. Number of elements in $Z_2G$.

     b. Idempotents in $Z_2G$.

     c. Ideals in $Z_2G$.

     d. Quotient ring $\dfrac{Z_2G}{I}$ for any ideal I of $Z_2G$.

## 2.6 Modules

In this section we just recall the definition of modules and some of its basic properties and illustrate them by examples.

**DEFINITION 2.6.1**: *Let R be a ring. An R-module or a left R-module is an additive abelian group M having R as a left operator domain such that in addition to the requirement r(x + y) = rx + ry. (r ∈ R, x, y ∈ M); for all groups with operators, we also have (a+b) x = ax + by, (ab) x = a(bx), $I_R x = x$ for a, b ∈ R and x ∈ M. The elements of M are called vectors and those of the ring R are called scalars. The mapping (a, x) → ax of A × M → M is called the scalar multiplication in the R-module M. We can define a similar notion called right R-modules where R acts on the right side of M.*

**Example 2.6.1**: All the additive abelian groups over the ring of integers Z is a Z-module.

**DEFINITION 2.6.2**: *Let M be an R-module. A subgroup S of the additive group M is a submodule, if S itself is an R-module.*

**DEFINITION 2.6.3**: *Let M and N be any two R-modules. An R-module homomorphism is a mapping ϕ from M to N such that ϕ(x + y) = ϕ(x) + ϕ(y) and ϕ(αx) = αϕ(x) for all x, y ∈ M and α ∈ R.*

We illustrate this by examples and problems.

***Example 2.6.2***: Let R be a ring say $Z_8$. Now $M = Z_8 \times Z_8$ is an abelian group under addition, M is a $Z_8$ – module over $Z_8$.

**DEFINITION 2.6.4**: *Let M be a module, M is called a simple module if $M \neq (0)$ and the only submodules of M are (0) and M.*

***Example 2.6.3***: Let R be a ring $S = R \times R$ is an R-module. (show $M = R \times \{0\}$ and $N = \{0\} \times R$ are not isomorphic as S-modules).

PROBLEMS:

1. Let A and B be two submodules of a module M; prove $A \cap B$ is a submodule of M.
2. $M = Z \times Z \times Z \times Z \times Z$ is a module over Z.
    1. Find submodules of M.
    2. Find two submodules which are isomorphic in M.
3. Let $S = R \times R \times R$ be a module over R. Can S have submodules which are isomorphic?
4. Is S given in example 3, a simple module over R?
5. Give an example of a simple module.

## 2.7 Rings with chain conditions

In this section we recall the concept of chain conditions in rings; that is the concept of Artinian rings and Noetherian rings and illustrate them by examples.

**DEFINITION 2.7.1**: *Let R be a ring. R is said to be left Noetherian or Noetherian if the R-modules, $R_t$ is Noetherian. Since the submodules of $R_t$ are the same as left ideals of R; therefore this is the same as to say that the ring R satisfies the following equivalent conditions.*

   a. *Ascending chain conditions, (A.C.C) on left ideals: If every ascending chain $M_1 \subset M_2 \subset \ldots$ of left ideals of R is stationary.*

   b. *Maximum condition on left ideals: If every non-empty collection of left ideals has a maximal member.*

   c. *Finite generations of left ideals every left ideal of R is finitely generated.*

*Similarly a ring R is said to be left Artinian or Artinian if R-modules $R_t$ is Artinian; i.e., if the ring R satisfies the following equivalent conditions.*

1. *Descending chain conditions (D.C.C) on left ideals: If every descending chain $N_1 \supset N_2 \supset \ldots$ of left ideals of R is stationary.*

2. *Minimum condition on left ideals: every non-empty collection of left ideals of R has a minimal member. Finally we say that R is right Noetherian (respectively right Artinian) if the right R-modules $A_i$ is Noetherian (respectively Artinian).*

**Example 2.7.1**: Let Z be the ring of integers; Z is Noetherian but it is not Artinian because we have $(2) \supset (4) \supset (8) \supset \ldots (n) \ldots$; (n) denotes ideals of Z where n $\in$ Z.

**Example 2.7.2**: The ring of polynomials in finitely many variables over a Noetherian ring is Noetherian. (This is known as Hilbert basis theorem).

Every right Artirian ring is right Noetherian. The converse does not hold good.

<u>PROBLEMS:</u>

1. Show the finite direct product of Noetherian ring is Noetherian.
2. Is every factor ring of a right Artinian ring, Artinian? Justify your answer.
3. Show that the ring of all $2 \times 2$ matrices $\begin{pmatrix} a & b \\ 0 & c \end{pmatrix}$ where a, b, c $\in$ Q is right Noetherian but not left Noetherian.
4. Show that the ring of all $2 \times 2$ matrices $\begin{pmatrix} a & b \\ 0 & c \end{pmatrix}$; a is rational b and c are reals is right Artinian but not left Artinian.
5. Is the ring P = Z $\times$ Z $\times$ Z $\times$ Z $\times$ Z Artinian?
6. Prove R[x] is right Noetherian if R is Noetherian.
7. Show R $[x_1, \ldots, x_n]$ is right Noetherian if R is right Noetherian.

# SMARANDACHE RINGS AND ITS PROPERTIES

This is the main chapter of this book. Here we introduce several new concepts in Smarandache rings and recall the definition of Smarandache rings, Smarandache ideals and Smarandache subrings as given by Florentin Smarandache and Padilla Raul. We do not indulge in proving any of the classical results in ring theory. For in our opinion as there are many texts on ring theory any interested reader can develop all classical results and theorem to Smarandache rings.

Several new concepts from ring theory that are not found in textbooks but have appeared only as research papers are introduced in this chapter. So at this juncture the author makes it very clear that most of the ring theory concepts given in undergraduate texts are ignored as they can be treated as exercises once this book is mastered. The author felt that when several innovative concepts in ring theory – about elements and substructures in rings – which are found only in research papers are given Smarandache equivalents, certainly it would be of interest to both Smarandache algebraists and ring theorists. Hence this book incorporates both the unique concepts of ring theory and their Smarandache analogues. It contains several definitions propounded by various authors and also provides an extensive bibliography of these papers thereby making it an important piece of work on Smarandache rings.

This chapter is organized into ten sections. Section one defines Smarandache rings of level I and II, explains with examples and introduces the concept of Smarandache commutative rings. In section two, three, four the author introduces the special elements in a ring viz Smarandache units, Smarandache zero divisors and Smarandache idempotents. Several important results are given with many illustrative examples. The main substructure like S-ideals and S-subring are studied in section five leading to the definition of Smarandache simple rings, Smarandache pseudo simple rings. Smarandache modules are introduced in section six. Just the Smarandache analogue of D.C.C and A.C.C are given in section seven.

In section eight we define Smarandache group rings and Smarandache semi group rings as they serve as concrete examples in almost all illustrations. Special elements like Smarandache nilpotents, Smarandache semi idempotent, Smarandache pseudo commutative pair, S-quasi commutative elements, Smarandache semi nilpotent element etc. are introduced in the ninth section. The tenth section is the longest section and the main section of this chapter. It introduces over 70 Smarandache notions and gives around 40 theorems with 55 illustrative examples.

In several places the author leaves the proof of certain result for the reader as only by solving these at each stage can make a researcher well versed in Smarandache ring theory. This section ends with 70 problems which can be easily worked as exercise by any studious researcher. Each section starts with a brief introduction.

## 3.1 Definition of Smarandache Ring with Examples

In this section we recall the definition of the Smarandache rings and illustrate it with examples. Smarandache rings were introduced by Florentin Smarandache and Padilla Raul [73, 60] in the year 1998. Several researchers have been working on these Smarandache concepts. As we have several books on ring theory and no book on Smarandache rings here we venture to write a book on Smarandache-Rings distinctly different from usual ring theory books.

**DEFINITION [73, 60]** : *A Smarandache ring (S-ring) is defined to be a ring A, such that a proper subset of A is a field with respect to the operations induced. By proper subset we understand a set included in A different from the empty set, from the unit element if any and from A.*

**Example 3.1.1**: Let $F[x]$ be a polynomial ring over a field F. $F[x]$ is a S-ring.

**Example 3.1.2**: Let $Z_{12} = \{0, 1, 2, \ldots, 11\}$ be a ring. $Z_{12}$ is a S-ring as A = $\{0, 4, 8\}$ is a field with 4 acting as the unit element.

**Example 3.1.3**: $Z_6 = \{0, 1, 2, \ldots, 5\}$ is a S-ring; for take A = $\{0, 2, 4\}$ is a field with 4 acting as the unit of A.

It is interesting to note that we do not demand the unit of the ring to be the unit of the field; further we do not expect all rings to be S-rings.

From now onwards we will call these S-rings as S-ring I. For these rings do not help us to define Smarandache commutative ring or like concepts. So we are forced to opt for the second level of S-ring.

**DEFINITION 3.1.1**: *Let R be a ring, R is said to be a Smarandache ring of level II (S-ring II) if R contains a proper subset A (A ≠ 0) such that*

   *1. A is an additive abelian group.*
   *2. A is a semigroup under multiplication.*
   *3. For a, b, ∈ A; a.b = 0 if and only if a = 0 or b = 0.*

**THEOREM 3.1.1**: *Let R be S-ring I then R is a S-ring II.*

*Proof*: By the very definition of S-ring I and S-ring II we see every S-ring I is a S-ring II for it obviously satisfies all conditions of S-ring II.

**THEOREM 3.1.2**: *Every S-ring II need not in general be a S-ring I.*

*Proof*: Take Z[x] the polynomial ring. Z[x] is S-ring II for Z ⊂ Z[x] but Z[x] is not a S-ring I.

Thus we have the class of S-ring I to be strictly contained in the class of S-ring II.

**DEFINITION 3.1.2**: *Let R be a ring, R is said to be a Smarandache commutative ring II (S-commutative ring II) if R is a S-ring and there exists at least a proper subset A of R which is a field or an integral domain i.e. for all a, b ∈ A we have ab = ba. If R has no proper subset A (A ⊂ R) which is a field or an integral domain then we say R is a Smarandache non-commutative ring II (S-non-commutative ring II).*

Thus we can simply say R is a S-non-commutative ring II if no proper subset of R is an integral domain or a field.

**THEOREM 3.1.3**: *Let R be a ring, R is a S-commutative ring II if and only if R has atleast a proper subset, which is an integral domain.*

*Proof*: Given R is a S-commutative ring II, so R has a proper subset A, which is an integral domain.

Conversely suppose R has a proper subset A which is an integral domain by the very definition of S-ring II, R is a S-commutative ring II.

**THEOREM 3.1.4**: *Let R be a ring, R is said to be a S-non-commutative ring II if R has no proper subset A, which is an integral domain, but R has only proper subsets, which are division rings.*

*Proof*: For if R has atleast one proper subset which is an integral domain then R will be a S-commutative ring II but for R to be a S-ring II, R must have atleast a proper subset which is a division ring. Hence the claim.

From these definitions and results we see even if R is a non-commutative ring still R can be a S-commutative ring II.

***Example 3.1.4***: Let QR={$\alpha_0$ + $\alpha_1$i + $\alpha_2$j + $\alpha_3$k / $\alpha_0$, $\alpha_1$, $\alpha_2$, $\alpha_3$ ∈ Q – the field of rationals, $i^2 = j^2 = k^2 = -1 = ijk$, ij = – ji = k, jk = – kj = i, ki = –ik = j} be the ring of quaternions.

Clearly QR is a S-commutative ring II and QR is also a S-ring I. (QR is non-commutative) as it has Q ⊂ QR to be a commutative field.

**DEFINITION 3.1.3**: *Let R be a ring, R is a S-ring I (or II), we say the Smarandache characteristics (S-characteristic) of R is the characteristic of the field which is a proper subset of R (and or) the characteristic of the integral domain which is a proper subset of R or the characteristic of a division ring which is a proper subset of R.*

Thus for a ring R which is a S-ring I or S-ring II we can have several S-characteristics associated with it.

**THEOREM 3.1.5**: *Let R be a commutative finite ring. If R is a S-ring II then R is a S-ring I.*

*Proof*: By the very definition of S-ring I and S-ring II we see they are identical in a finite commutative ring as "Every finite integral domain is a field". Hence the claim.

**THEOREM 3.1.6**: *If R is a S-ring I (or S-ring II) and R[x] is a polynomial ring in the indeterminate x over R, then R[x] is a S-ring I (or S-ring II).*

*Proof*: Now R is a S-ring I (S-ring II) so $A \subset R$ (A is a field or an integral domain or a division ring) so $A[x] \subset R[x]$ is an integral domain or a division ring, hence R is S-ring II or $A \subset R[x]$, so if R is a S-ring I so is R[x].

**THEOREM 3.1.7**: *Let F be a field and G any group. Then the group ring FG is a S-ring I.*

*Proof*: The result is true as the field F is such that $F \subset FG$. Hence FG is a S-ring I.

**THEOREM 3.1.8**: *Let F be a field and S any semigroup with unit. The semigroup ring FS is a S-ring I.*

*Proof*: Left for the reader to prove.

**THEOREM 3.1.9**: *Let Z be the ring of integers and G any group, then the group ring ZG is a S-ring II and not a S-ring I.*

*Proof*: Obvious from the fact Z is only an integral domain and $Z \subset ZG$; hence ZG is a S-ring II.

**COROLLARY**: *Let Z be the ring of integers and G a non-commutative group (S a non-commutative monoid) then the group ring ZG (the semigroup ring ZS) is a S-ring II.*

*Proof*: Left for the reader to verify.

**Theorem: 3.1.10**: *Let* $M_{n \times n} = \{(a_{ij}) \, / \, a_{ij} \in Z\}$ *be the ring of matrices.* $M_{n \times n}$ *is a S-ring II.*

*Proof*: Consider the matrix, $A = \{(a_{ii}) \, / \, a_{ii} \in Z \setminus \{0\}$ and $a_{ii} = 0, i \neq j\} \cup (0)$, (where $(0)$ is the zero matrix) that is, A consists of only diagonal matrices. Then A is an integral domain, so $M_{n \times n}$ is a S-ring II and not a S-ring I.

## PROBLEMS:

1. Give an example of a S-ring II, which is not a S-ring I.
2. Can a ring with zero divisors be a S-ring I? Justify your answer with examples.
3. Give an example of a finite S-ring I of order 64.
4. What is the order of the smallest S-ring I which is non-commutative?
5. Give an example of a smallest S-ring I.
6. Find a S-ring I using the semigroup S(5). (By constructing suitable semigroup rings).
7. Let $Z_3S(4)$ be the S-ring I. Is $Z_3 S(4)$ a S-commutative ring I?
8. Let $Z_{24} S_3$ be the group ring of the group $S_3$ over the ring $Z_{24}$. How many proper subsets in $Z_{24} S_3$ are fields? Is $Z_{24} S_3$, S-commutative?
9. Is $Z_{12}G$ where $G = \langle g \, / \, g^{12} = 1 \rangle$, a S-ring II? Justify your answer.
10. $ZS(n)$ be the semigroup ring. Is $Z S(n)$ a S-ring I? Justify your claim.

## 3.2 Smarandache units in Rings

In this section we introduce the notion of Smarandache units (S-units) in rings. For introducing S-units we don't require S-ring. S-units are defined for any arbitrary ring and interesting results are obtained. We prove that units of the form $x^2 = 1$ can never be S-units. We prove every unit in the field of rationals and reals are S-units.

**Definitions 3.2.1**: *Let R be a ring with unit. We say $x \in R \setminus \{1\}$ is a Smarandache unit (S-unit) if there exists $y \in R$ with*

1. *$xy = 1$.*
2. *There exists $a, b$ in $R \setminus \{x, y, 1\}$.*
   - i) *$xa = y$ or $ax = y$ or*
   - ii) *$yb = x$ or $by = x$ and*
   - iii) *$ab = 1$.*

*(2(i) or 2(ii) is satisfied it is enough to make a S-unit).*

***Example 3.2.1***: Let $Z_9 = \{0, 1, 2, \ldots, 8\}$ be the ring under multiplication modulo 9, $2 \in Z_9$ is a S-unit for $5 \in Z_9$ is such that $2.5 \equiv 1 \pmod 9$ and $7, 4 \in Z_9$ is such that $2.7 \equiv 5 \pmod 9$ and $5.4 \equiv 2 \pmod 9$ with $7.4 \equiv 1 \pmod 9$.

***Example 3.2.2***: Let $Z_5 = \{0, 1, 2, 3, 4\}$ be the ring of integers modulo 5. Clearly $3 \in Z_5$ is a S-unit in $Z_5$ as $2.3 \equiv 1 \pmod 5$ and $4 \in Z_5$ is such that $2.4 \equiv 3 \pmod 5$ and $3.4 \equiv 2 \pmod 5$ and $4^2 \equiv 1 \pmod 5$.

**THEOREM 3.2.1**: *Every S-unit in a ring is a unit but all units in a ring need not in general be S-units.*

*Proof*: Clearly by the very definition of S-unit we see it is a unit, but every unit need not be a S-unit. We prove this by an example: Consider the ring. $Z_9 = \{0, 1, 2, \ldots, 8\}$ of modulo integers. Clearly 7 is a unit as $7.4 \equiv 1 \pmod 9$; 7 is not a S-unit in $Z_9$ as we cannot find $a, b \in Z_9 \setminus \{7, 4\}$ such that $7a \equiv 4 \pmod 9$ or $4.b \equiv 7 \pmod 5$ with $ab \equiv 1 \pmod 9$.

***Example 3.2.3***: Let $Z_{15} = \{0, 1, 2, \ldots, 14\}$ be the ring of integers modulo 15. Now $2 \in Z_{15}$ is a S-unit for $2.8 \equiv 1 \pmod{15}$, $4^2 \equiv 1 \pmod{15}$ and $2.4 \equiv 8 \pmod{15}$.

It is important to note that when we say x is a S-unit in R we do not say there exist $y \neq x$ in R, but we will prove that $x^2 = 1$ can never be a S-unit so it is not essential to say $y \neq x$ in the definition.

Similarly when we take $a, b \in R \setminus \{x, y, 1\}$ we do not demand a and b to be distinct, a = b can also occur. We illustrate this by an example.

***Example 3.2.4***: $Z_{15} = \{0, 1, 2, \ldots, 14\}$ be the ring of integers modulo 15. We have $4^2 \equiv 1 \pmod{15}$ as we cannot find $a, b \in Z_{15}$ such that $4a \equiv 4 \pmod{15}$ or $4b \equiv 4 \pmod{15}$ with $a.b = 1 \pmod{15}$.

**THEOREM 3.2.2**: *Let R be a ring with unit 1, if $x \in R \setminus \{1\}$ is a S-unit, with xy = 1 then $x \neq y$.*

*Proof*: Let $x \in R \setminus \{1\}$ is a S-unit, so by definition we have $xy = 1$ such that $a, b \in R \setminus \{x, y, 1\}$ with $xa = y$ or $ax = y$, (by = x or $yb = x$) and $ab = 1$; if $x = y$ then we have $x^2 = 1$, $xa = x$, i.e., $x^2 a = x^2$ that is $a = 1$, a contradiction to the very definition of S-unit.

***Example 3.2.5***: Let $R_{2 \times 2} = \{(a_{ij}) \, / \, a_{ij} \in Z_2 = \{0, 1\}\}$ be the collection of all $2 \times 2$ matrices with entries from $Z_2 = \{0, 1\}$. $R_{2 \times 2}$ is a ring under the matrix multiplication and matrix addition.

In R$_{2\times2}$ we have

$$\begin{pmatrix} 0 & 1 \\ 1 & 0 \end{pmatrix}, \begin{pmatrix} 1 & 0 \\ 0 & 1 \end{pmatrix}, \begin{pmatrix} 1 & 1 \\ 0 & 1 \end{pmatrix}, \begin{pmatrix} 1 & 0 \\ 1 & 1 \end{pmatrix}, \begin{pmatrix} 0 & 1 \\ 1 & 1 \end{pmatrix} \text{ and } \begin{pmatrix} 1 & 1 \\ 1 & 0 \end{pmatrix}$$

to be units with

$$\begin{pmatrix} 0 & 1 \\ 1 & 0 \end{pmatrix} \begin{pmatrix} 0 & 1 \\ 1 & 0 \end{pmatrix} = \begin{pmatrix} 1 & 0 \\ 0 & 1 \end{pmatrix},$$

$$\begin{pmatrix} 1 & 1 \\ 0 & 1 \end{pmatrix} \begin{pmatrix} 1 & 1 \\ 0 & 1 \end{pmatrix} = \begin{pmatrix} 1 & 0 \\ 0 & 1 \end{pmatrix},$$

$$\begin{pmatrix} 1 & 0 \\ 1 & 1 \end{pmatrix} \begin{pmatrix} 1 & 0 \\ 1 & 1 \end{pmatrix} = \begin{pmatrix} 1 & 0 \\ 0 & 1 \end{pmatrix}$$

and

$$\begin{pmatrix} 0 & 1 \\ 1 & 1 \end{pmatrix} \begin{pmatrix} 1 & 1 \\ 1 & 0 \end{pmatrix} = \begin{pmatrix} 1 & 0 \\ 0 & 1 \end{pmatrix}$$

where

$$\begin{pmatrix} 1 & 0 \\ 0 & 1 \end{pmatrix}$$

is the multiplicative identity. To find which of these are S-units. Clearly it can be easily verified that

$$\begin{pmatrix} 0 & 1 \\ 1 & 0 \end{pmatrix}, \begin{pmatrix} 1 & 1 \\ 0 & 1 \end{pmatrix} \text{ and } \begin{pmatrix} 1 & 0 \\ 1 & 1 \end{pmatrix}$$

are not S-units. It can be verified that none of the elements in R$_{2\times2}$ are S-units but R$_{2\times2}$ has 5 distinct units.

From this example we have made the following observations:

$$\begin{pmatrix} 0 & 1 \\ 1 & 0 \end{pmatrix}, \begin{pmatrix} 0 & 1 \\ 1 & 0 \end{pmatrix} = \begin{pmatrix} 1 & 0 \\ 0 & 1 \end{pmatrix};$$

clearly

$$x = y = \begin{pmatrix} 0 & 1 \\ 1 & 0 \end{pmatrix}.$$

Now if we do not assume a, b ∈ R \ {x, y, 1}; a = 1 or b = 1 can occur in which case

$$\begin{pmatrix} 0 & 1 \\ 1 & 0 \end{pmatrix} \begin{pmatrix} 1 & 0 \\ 0 & 1 \end{pmatrix} = \begin{pmatrix} 0 & 1 \\ 1 & 0 \end{pmatrix},$$

thus every element x in a ring R such that $x^2 = 1$ will become a S-unit. Further

$$\begin{pmatrix} 0 & 1 \\ 1 & 1 \end{pmatrix} \begin{pmatrix} 1 & 1 \\ 1 & 0 \end{pmatrix} = \begin{pmatrix} 1 & 0 \\ 0 & 1 \end{pmatrix}$$

is an example; if we do not assume a, b ∈ R \ {x, y, 1} and if we take

$$x = \begin{pmatrix} 0 & 1 \\ 1 & 1 \end{pmatrix} \text{ and } y = \begin{pmatrix} 1 & 1 \\ 1 & 0 \end{pmatrix}$$

we have by taking

$$x = a = \begin{pmatrix} 0 & 1 \\ 1 & 1 \end{pmatrix},$$

$$\begin{pmatrix} 0 & 1 \\ 1 & 1 \end{pmatrix} \begin{pmatrix} 0 & 1 \\ 1 & 1 \end{pmatrix} = \begin{pmatrix} 1 & 1 \\ 1 & 0 \end{pmatrix} = y.$$

Hence the assumptions in the definition 3.2.1 are important for S-units to be distinctly different from units.

**THEOREM 3.2.3**: *Every unit in the ring $Z_n$ = {0, 1, …, n–1} is not a S-unit.*

*Proof*: Given $Z_n$ is the ring of integers modulo n. We have n − 1 ∈ $Z_n$ is such that (n − 1)(n − 1) ≡ 1 (mod n) is a unit, which is not a S-unit by theorem 3.2.2. Thus we have in a prime field of characteristic p, p a prime every element is a unit but every element in $Z_p$ is not a S-unit contrary to prime fields of characteristic 0.

**THEOREM 3.2.4**: *Let Q be the field of rationals, every unit in Q is a S-unit.*

*Proof*: Q is the field of rationals. To prove every unit in Q is a S-unit in Q. Let m be an integer, we know $m \times \frac{1}{m} = 1$, $m \times \frac{1}{m^2} = \frac{1}{m}$ and $m^2 \times \frac{1}{m} = m$ and $\frac{1}{m^2} \times m^2 = 1$. If $m = \frac{p}{q} (q \neq 0)$ then $m_1 = \frac{q}{p}$ is such that $m_1 \ m \ = \ 1$. Now $\frac{p}{q} \times \frac{q^2}{p_2} = \frac{q}{p}$ and $\frac{q}{p} \times \frac{p^2}{q^2} = \frac{p}{q}$ and $\frac{p^2}{q^2} \times \frac{q^2}{p^2} = 1$. Hence every unit in Q is a S-unit.

In view of this we have the following theorem:

**THEOREM 3.2.5**: *If F is a prime field of characteristic 0 every unit is a S-unit.*

*Proof*: Since all prime fields of characteristic 0 are isomorphic to Q we have the result.

***Example 3.2.6***: Let Q be the field of characteristic 0 and G = {g | $g^2$ = 1}. The group ring QG = {$\alpha + \beta g$ | $\alpha, \beta \in$ Q}. Now g $\in$ QG and $g^2$ = 1 but g is not a S-unit.

**DEFINITION 3.2.2**: *Let S be a ring, if every element in S is a S-unit then we say S is a Smarandache unit domain (S-unit domain).*

*If S has no S-units, S is said to be a Smarandache unit free ring (S-unit free ring).*

***Example 3.2.7***: Q is a S-unit domain.

***Example 3.2.8***: R is a S-unit domain.

***Example 3.2.9***: $Z_p$, p a prime is a S-unit free domain.

**PROBLEMS:**

1. Find all S-units in $Z_{210}$.
2. Find all S-units of the group ring $Z_2 S_3$.
3. How many S-units does the semigroup ring $Z_4 S(3)$ have?
4. Find those units, which are not S-units in $Z_{24}$.
5. Does $M_{3 \times 3}$ = {$(a_{ij})$ / $a_{ij} \in Z_4$ = {0, 1, 2, 3}}, the ring of 3 × 3 matrices have S-units? Justify your answer.
6. Find all S-units in $QS_8$, the group ring of the group $S_8$ over the rational field Q.
7. Find the S-units in ZS(7); the semigroup ring of the semigroup S(7) over the ring of integers Z.

8. Find units in the semigroup ring $ZS(7)$ given in problem 7 which are not S-units.

9. Find the S-units of the group ring $Z_{11}G$ where G is the dihedral group given by $G = \{a, b \mid a^2 = b^9 = 1, bab = a\}$.

10. Find all units in $Z_{11}G$ in problem 9 which are not S-units.

11. Can the group ring $Z_3G$ where $G = \langle g \mid g^p = 1 \rangle$, p a prime $p > 3$ have S-units? Justify your answer.

12. Can the group ring $Z_pG$ where $G = \langle g \mid g^p = 1 \rangle$, have S-units? Does $Z_pG$ have units, which are not S-units?

## 3.3 Smarandache Zero Divisors in Rings

In this section we introduce the concept of Smarandache zero divisors (S-zero divisors) in rings and show that every S-zero divisor is a zero divisor but all zero divisors are not S-zero divisors.

**DEFINITION 3.3.1**: *Let R be a ring, we say x and y in R is said to be a Smarandache zero divisor (S-zero divisor) if xy = 0 and there exists a, b ∈ R \ {0, x, y} with*

    *1. xa = 0 or ax = 0.*
    *2. yb = 0 or by = 0.*
    *3. ab ≠ 0 or ba ≠ 0.*

**Example 3.3.1**: Let $Z_{20} = \{0, 1, 2, \ldots, 19\}$ be the ring of integers modulo 20. Clearly 10, 16 is a S-zero divisor, consider $5, 6 \in Z_{20} \setminus \{0\}$

$$5 \times 16 \equiv 0 \pmod{20}$$
$$6 \times 10 \equiv 0 \pmod{20}$$
$$6 \times 5 \not\equiv 0 \pmod{20}.$$

**Example 3.3.2**: Let $Z_{10} = \{0, 1, \ldots, 9\}$ be the ring of integers modulo 10. Clearly $2.5 \equiv 10 \equiv 0 \pmod{10}$ is a zero divisor but is not a S-zero divisor.

**THEOREM 3.3.1**: *Let R be a ring. Every S-zero divisor is a zero divisor but a zero divisor in general is not a S-zero divisor.*

*Proof*: By the very definition of S-zero divisor we see if x, y is a S-zero divisor, it is a zero divisor. But by example 3.3.2 we see $2.5 \equiv 0 \pmod{10}$ is a zero divisor in $Z_{10}$ but it is not a S-zero divisor.

***Example 3.3.3***: Let $S_{2\times 2} = \left\{ \begin{pmatrix} a & b \\ c & d \end{pmatrix} \middle/ a,b,c,d \in Z_2 = \{0,1\} \right\}$ be the set of all $2 \times 2$

matrices with entries from the ring of integers $Z_2$. Clearly $S_{2\times 2}$ is the matrix ring. Now

$$\begin{pmatrix} 1 & 0 \\ 0 & 0 \end{pmatrix}, \begin{pmatrix} 0 & 0 \\ 0 & 1 \end{pmatrix} \in S_{2\times 2}$$

is zero divisor of $S_{2\times 2}$ as

$$\begin{pmatrix} 1 & 0 \\ 0 & 0 \end{pmatrix}\begin{pmatrix} 0 & 0 \\ 0 & 1 \end{pmatrix} = \begin{pmatrix} 0 & 0 \\ 0 & 0 \end{pmatrix} \text{ and } \begin{pmatrix} 0 & 0 \\ 0 & 1 \end{pmatrix}\begin{pmatrix} 1 & 0 \\ 0 & 0 \end{pmatrix} = \begin{pmatrix} 0 & 0 \\ 0 & 0 \end{pmatrix}.$$

Now take

$$x = \begin{pmatrix} 0 & 1 \\ 0 & 0 \end{pmatrix} \text{ and } y = \begin{pmatrix} 0 & 0 \\ 1 & 0 \end{pmatrix}$$

in $S_{2\times 2}$. We have

$$\begin{pmatrix} 0 & 1 \\ 0 & 0 \end{pmatrix}\begin{pmatrix} 1 & 0 \\ 0 & 0 \end{pmatrix} = \begin{pmatrix} 0 & 0 \\ 0 & 0 \end{pmatrix}$$

but

$$\begin{pmatrix} 1 & 0 \\ 0 & 0 \end{pmatrix}\begin{pmatrix} 0 & 1 \\ 0 & 0 \end{pmatrix} = \begin{pmatrix} 0 & 1 \\ 0 & 0 \end{pmatrix} \neq \begin{pmatrix} 0 & 0 \\ 0 & 0 \end{pmatrix}$$

$$\begin{pmatrix} 0 & 0 \\ 1 & 0 \end{pmatrix}\begin{pmatrix} 0 & 0 \\ 0 & 1 \end{pmatrix} = \begin{pmatrix} 0 & 0 \\ 0 & 0 \end{pmatrix}$$

but

$$\begin{pmatrix} 0 & 0 \\ 0 & 1 \end{pmatrix}\begin{pmatrix} 0 & 0 \\ 1 & 0 \end{pmatrix} = \begin{pmatrix} 0 & 0 \\ 1 & 0 \end{pmatrix} \neq \begin{pmatrix} 0 & 0 \\ 0 & 0 \end{pmatrix}.$$

Finally

$$\begin{pmatrix} 0 & 1 \\ 0 & 0 \end{pmatrix}\begin{pmatrix} 0 & 0 \\ 1 & 0 \end{pmatrix} = \begin{pmatrix} 1 & 0 \\ 0 & 0 \end{pmatrix} \neq \begin{pmatrix} 0 & 0 \\ 0 & 0 \end{pmatrix}, \begin{pmatrix} 0 & 0 \\ 1 & 0 \end{pmatrix}\begin{pmatrix} 0 & 1 \\ 0 & 0 \end{pmatrix} = \begin{pmatrix} 0 & 0 \\ 0 & 1 \end{pmatrix} \neq \begin{pmatrix} 0 & 0 \\ 0 & 0 \end{pmatrix}.$$

Hence

$$\begin{pmatrix} 0 & 1 \\ 0 & 0 \end{pmatrix}, \begin{pmatrix} 0 & 0 \\ 1 & 0 \end{pmatrix}$$

is a S-zero divisor of the ring $S_{2\times 2}$.

***Example 3.3.4***: Let $R_{3\times 3} = \{(a_{ij}) \ / \ a_{ij} \in Z_4 = \{0, 1, 2, 3\}\}$ be the collection of all $3 \times 3$ matrices with entries from $Z_4$. Now $R_{3\times 3}$ is a ring under matrix multiplication modulo 4. We have

$$\begin{pmatrix} 1 & 0 & 0 \\ 0 & 0 & 0 \\ 0 & 0 & 2 \end{pmatrix}, \begin{pmatrix} 0 & 0 & 0 \\ 0 & 1 & 0 \\ 0 & 2 & 2 \end{pmatrix}$$

in $R_{3\times 3}$ is a zero divisor of $R_{3\times 3}$.

$$\begin{pmatrix} 1 & 0 & 0 \\ 0 & 0 & 0 \\ 0 & 0 & 2 \end{pmatrix} \begin{pmatrix} 0 & 0 & 0 \\ 0 & 1 & 0 \\ 0 & 2 & 2 \end{pmatrix} = \begin{pmatrix} 0 & 0 & 0 \\ 0 & 0 & 0 \\ 0 & 0 & 0 \end{pmatrix},$$

$$\begin{pmatrix} 0 & 0 & 0 \\ 0 & 1 & 0 \\ 0 & 2 & 2 \end{pmatrix} \begin{pmatrix} 1 & 0 & 0 \\ 0 & 0 & 0 \\ 0 & 0 & 2 \end{pmatrix} = \begin{pmatrix} 0 & 0 & 0 \\ 0 & 0 & 0 \\ 0 & 0 & 0 \end{pmatrix}.$$

Consider

$$\begin{pmatrix} 0 & 0 & 0 \\ 0 & 3 & 2 \\ 0 & 0 & 2 \end{pmatrix}, \begin{pmatrix} 0 & 0 & 0 \\ 0 & 0 & 0 \\ 0 & 2 & 2 \end{pmatrix} \in R_{3\times 3}$$

$$\begin{pmatrix} 1 & 0 & 0 \\ 0 & 0 & 0 \\ 0 & 0 & 2 \end{pmatrix} \begin{pmatrix} 0 & 0 & 0 \\ 0 & 3 & 2 \\ 0 & 0 & 2 \end{pmatrix} = \begin{pmatrix} 0 & 0 & 0 \\ 0 & 0 & 0 \\ 0 & 0 & 0 \end{pmatrix}$$

$$\begin{pmatrix} 0 & 0 & 0 \\ 0 & 3 & 2 \\ 0 & 0 & 2 \end{pmatrix} \begin{pmatrix} 1 & 0 & 0 \\ 0 & 0 & 0 \\ 0 & 0 & 2 \end{pmatrix} = \begin{pmatrix} 0 & 0 & 0 \\ 0 & 0 & 0 \\ 0 & 0 & 0 \end{pmatrix}$$

$$\begin{pmatrix} 0 & 0 & 0 \\ 0 & 1 & 0 \\ 0 & 2 & 2 \end{pmatrix} \begin{pmatrix} 0 & 0 & 0 \\ 0 & 0 & 0 \\ 0 & 2 & 2 \end{pmatrix} = \begin{pmatrix} 0 & 0 & 0 \\ 0 & 0 & 0 \\ 0 & 0 & 0 \end{pmatrix}$$

but

$$\begin{pmatrix} 0 & 0 & 0 \\ 0 & 0 & 0 \\ 0 & 2 & 2 \end{pmatrix} \begin{pmatrix} 0 & 0 & 0 \\ 0 & 1 & 0 \\ 0 & 2 & 2 \end{pmatrix} = \begin{pmatrix} 0 & 0 & 0 \\ 0 & 0 & 0 \\ 0 & 2 & 0 \end{pmatrix} \neq \begin{pmatrix} 0 & 0 & 0 \\ 0 & 0 & 0 \\ 0 & 0 & 0 \end{pmatrix}$$

$$\begin{pmatrix} 0 & 0 & 0 \\ 0 & 3 & 2 \\ 0 & 0 & 2 \end{pmatrix} \begin{pmatrix} 0 & 0 & 0 \\ 0 & 0 & 0 \\ 0 & 2 & 2 \end{pmatrix} = \begin{pmatrix} 0 & 0 & 0 \\ 0 & 0 & 0 \\ 0 & 0 & 0 \end{pmatrix}$$

$$\begin{pmatrix} 0 & 0 & 0 \\ 0 & 0 & 0 \\ 0 & 2 & 2 \end{pmatrix} \begin{pmatrix} 0 & 0 & 0 \\ 0 & 3 & 2 \\ 0 & 0 & 2 \end{pmatrix} = \begin{pmatrix} 0 & 0 & 0 \\ 0 & 0 & 0 \\ 0 & 2 & 0 \end{pmatrix} \neq \begin{pmatrix} 0 & 0 & 0 \\ 0 & 0 & 0 \\ 0 & 0 & 0 \end{pmatrix} .$$

So

$$\begin{pmatrix} 1 & 0 & 0 \\ 0 & 0 & 0 \\ 0 & 0 & 2 \end{pmatrix} \text{ and } \begin{pmatrix} 0 & 0 & 0 \\ 0 & 1 & 0 \\ 0 & 2 & 2 \end{pmatrix}$$

are S-zero divisors in $R_{3 \times 3}$.

**THEOREM 3.3.2**: *Let R be a non-commutative ring. x, y $\in$ R \ {0} be a S-zero divisor with a, b $\in$ R \ {0, x, y} satisfying the following conditions:*

1. *ax = 0 and xa $\neq$ 0.*
2. *yb = 0 and by $\neq$ 0.*
3. *ab = 0 and ba $\neq$ 0.*

*Then (xa + by)$^2$ = 0, i.e., xa + by is a nilpotent element of R.*

*Proof*: Given x, y $\in$ R \ {0} is a S-zero divisor such that xy = 0 = yx. We have a, b $\in$ R \ {0, x, y} with ax = 0 and xa $\neq$ 0 and yb = 0 and by $\neq$ 0 with ab = 0 and ba $\neq$ 0. Consider (xa + by)$^2$ = xaby + byxa + xaxa + byby; using xy = yx = 0, ab = 0, yb = 0 and ax = 0 we get xa + by to be a nilpotent element of order 2.

Now in view of this we have the following nice definition:

**DEFINITION 3.3.2**: *Let R be a commutative ring. If R has no S-zero divisors we say R is a Smarandache integral domain (S- integral domain) (Thus we may have zero divisors in R what we need is R should not have S-zero divisors).*

**THEOREM 3.3.3**: *Let R be an integral domain. Then R is a S- integral domain.*

*Proof*: Obvious by the very definition of S-integral domain.

**DEFINITION 3.3.3**: *Let R be a non-commutative ring. If R has no S-zero divisors then we say R is a Smarandache division ring (S-division ring).*

(Here also a S-division ring may have zero divisors). We will discuss and use these concepts in later chapters.

**Examples 3.3.5**: Clearly $Z_4 = \{0, 1, 2, 3\}$ is a S-integral domain but is not an integral domain as $2.2 \equiv 0 \pmod{4}$.

**THEOREM 3.3.4**: *Every S integral domain in general is not an integral domain.*

*Proof*: By example 3.3.5, $Z_4 = \{0, 1, 2, 3\}$ is not in integral domain but is a S-integral domain.

**COROLLARY**: *All division rings are S-division rings.*

*Proof*: By very definition of S-division rings.

Finally the author wishes to state that all zero divisors which occur are only from finite zeros. Finite zeros are zeros which occur in finitely constructed structure. $0 \in Z \subset Q \subset R$ is not a finite zero. For more about these please refer [159].

**PROBLEMS:**

1. Find whether $Z_{24}$ has S-zero divisors?
2. Does $Z_{14}$ have zero divisors, which are not S-zero divisors?
3. Find whether the group ring $Z_2 S_3$ has S-zero divisors?
4. Does the semigroup ring $Z_{12} S(5)$ have S-zero divisors? Can $Z_{12} S(5)$ have zero divisors which are not S-zero divisors?
5. Is $Z_{25} = \{0, 1, 2, \ldots, 24\}$ an S-integral domain? Justify your answer.
6. Give an example of S-division ring, which is not a division ring.
7. Can $Z_{12}(3)$ have nilpotent elements of order 2?

8.   Find all zero divisors in the semigroup ring $Z_2 S(4)$.
9.   Find all S-zero divisors of the group ring $Z_2 S_4$.
10.  Which ring $Z_2 S_4$ or $Z_2 S(4)$ will have more number of S-zero divisors? ($Z_2 S_4$ and $Z_2 S(4)$ given in problems in 8 and 9).

## 3.4 Smarandache idempotents in Rings

In this section we introduce the concept of Smarandache idempotents and Smarandache co-idempotents in rings and prove, if a ring has Smarandache idempotents then it has at least two divisors of zero. We prove if G is a finite group and K a field of characteristic zero then the group ring KG has nontrivial Smarandache idempotents. Finally we show group ring KG of a torsion free abelian group G over a field K of characteristic 0 has no Smarandache idempotents.

**DEFINITION 3.4.1**: *Let R be a ring. An element $0 \neq x \in R$ is a Smarandache idempotent (S-idempotent) of R if*

   *1. $x^2 = x$.*
   *2. There exists $a \in R \setminus \{x, 1, 0\}$.*
       *i.    $a^2 = x$ and*
       *ii.   $xa = a$ (ax = a) or $ax = x$ (xa = x)*

*or in (2, ii) is in the mutually exclusive sense.*

**DEFINITION 3.4.2**: *Let $x \in R \setminus \{0, 1\}$ be a Smarandache idempotent of R i.e., $x^2 = x$ and there exists $y \in R \setminus \{0, 1, x\}$ such that $y^2 = x$ and $yx = x$ or $xy = y$. We call 'y' the Smarandache co-idempotent (S-co-idempotent) and denote the pair by (x, y).*

***Example 3.4.1***: Let $Z_6 = \{0, 1, 2, \ldots, 5\}$ be the ring of integers modulo 6, then $4 \in Z_6$ is a S-idempotent of $Z_6$ for $4^2 \equiv 4 \pmod 6$ and $2 \in Z_6 \setminus \{4\}$ is such that $2^2 \equiv 4 \pmod 6$ $2.4 \equiv 2 \pmod 6$. Now $3 \in Z_6$ is such that $3^2 \equiv 3 \pmod 6$ but 3 is an idempotent of $Z_6$ but is not a S-idempotent of $Z_6$.

**THEOREM 3.4.1**: *Every S-idempotent is an idempotent but every idempotent in general is not a S-idempotent.*

*Proof*: By the very definition of S-idempotents we see every S-idempotent is an idempotent of the ring R. We see in example 3.4.1, in the ring $Z_6 = \{0, 1, 2, \ldots, 5\}$, 3 $\in Z_6$ is such that $3^2 \equiv 3 \pmod 6$ but is an idempotent which is not an S-idempotent of $Z_6$.

***Example 3.4.2***: Let $Z_{10} = \{0, 1, 2, \ldots, 9\}$ be the ring of integers modulo 10. Now the idempotents in $Z_{10}$ are 5 and 6 for $5^2 \equiv 5 \pmod{10}$ and $6^2 \equiv 6 \pmod{10}$. 5 is not a S-idempotent but 6 is a S-idempotent, $6^2 \equiv 6 \pmod{10}$ and $4 \in Z_{10}$ is such that $4^2 \equiv 6 \pmod{10}$ and $4.6 \equiv 4 \pmod{10}$.

**THEOREM 3.4.2**: *Let R be a ring. If R has a S-idempotent then R has atleast 2 nontrivial zero divisors.*

*Proof*: Let $a \in$ R be a S-idempotent, hence $a^2 = a$ and there exists $b \in$ R \ $\{a, 0, 1\}$ such that $b^2 = a$ and $ab = b$ which in turn implies $(a - 1) b = 0$. Further $a^2 = a$ implies $a (a - 1) = 0$ (as $a \neq 1$ or 0). Clearly $b \neq 0$ and $a \neq 1$. Hence the claim.

**COROLLARY**: *If R is a commutative ring and if R has S-idempotents then R has atleast 3 zero divisors.*

*Proof*: From Theorem 3.4.2 we have two zero divisors. Now $a^2 - b^2 = 0$ as $a^2 = a$ and $b^2 = a$ so $(a - b)(a + b) = 0$ is another zero divisor as $a \neq b$. The three zero divisors are distinct as $b \neq 1$ and $a \neq b$. Hence the theorem.

***Example 3.4.3***: Now consider the ring of integers modulo 12 given by $Z_{12} = \{0, 1, 2, 3, \ldots, 11\}$. Clearly this ring has two nontrivial idempotents viz 4 and 9, both of them are S-idempotents as $4^2 \equiv 4 \pmod{12}$ and $8 \in Z_{12}$ is such that $8^2 \equiv 4 \pmod{12}$ and $4.8 \equiv 8 \pmod{12}$. Now $9^2 \equiv 9 \pmod{12}$; $3 \in Z_{12}$ is such that $3^2 \equiv 9 \pmod{12}$ and $9.3 \equiv 3 \pmod{12}$. Hence the claim. Thus we have rings in which every idempotent is also a S-idempotent.

***Example 3.4.4:*** Let $Z_{15} = \{0, 1, 2, \ldots, 14\}$ be the ring of integers modulo 15. The only nontrivial idempotents of $Z_{15}$ are 6 and 10. Clearly 6 is a S-idempotent of $Z_{15}$, as $9 \in Z_{15}$ is such that $9 . 6 \equiv 9 \pmod{15}$ and $9^2 \equiv 6 \pmod{15}$ but 10 is a S-idempotent of $Z_{15}$, as $10^2 \equiv 10 \pmod{15}$ and $5^2 \equiv 10 \pmod{15}$ and $5.10 \equiv 5 \pmod{15}$. This example says even in the ring of modulo integers every idempotent is not a S-idempotent.

**THEOREM 3.4.3**: *Let $Z_n$ be the ring of integers modulo n. $Z_n$ has idempotents which are not S-idempotents when n = 2p, where p is a prime.*

*Proof*: Given $Z_n = \{0, 1, 2, \ldots, n - 1\}$ is the ring of integers modulo n and n = 2p, where p is an odd prime. Now $p^2 \equiv p \pmod{2p}$ by simple number theoretic arguments $p^2 \equiv p \pmod{2p}$ means as $p^2 + p \equiv 0 \pmod{2p}$ i.e., $p(p + 1) \equiv 0 \pmod{2p}$. Now $p \in Z_{2p}$ is an idempotent but p is not a S-idempotent for there does not exist a m $\in$ $Z_{2p}$ \ {p, 0, 1} such that $m^2 \equiv p$ and mp $\equiv$ m. But if p is a prime $m^2 = p$ is impossible. Thus in $Z_{2p}$, when p is prime, p is an idempotent which is not a S-idempotent.

***Example 3.4.5***: Let $Z_{30}$ = {0, 1, 2, …, 29} be the ring of integers modulo 30. Now this ring has 6, 10, 15, 16, 21 and 25 as non-trivial idempotents. 6 $\in$ $Z_{30}$ is a S-idempotent as $6^2 \equiv 6 \pmod{30}$ and $24^2 \equiv 6 \pmod{30}$, $6.24 \equiv 24 \pmod{30}$.

Similarly 10 is a S-idempotent as 20 serves the role of b $\in$ $Z_{30}$ \ {0, 1, 10}. 15 is not a S-idempotent. 16 is a S-idempotent with 14 acting as b for $16^2 \equiv 16 \pmod{30}$, $14^2 \equiv 16 \pmod{30}$ and $14.16 \equiv 14 \pmod{30}$ for 21, 9 serves as the element to make 21 a S-idempotent. For 25 is a S-idempotent as 5 serves the role of b.

Now we observe the following from this example.

1. All idempotents are not S-idempotents in $Z_{30}$ as 15 $\in$ $Z_{30}$ is not a S-idempotent.
2. Idempotents taken in certain pairs adds to 31. For (6, 25), (10, 21) and (15, 16). The sum of S-idempotent and S-co-idempotents adds to 30. (6, 24), (10, 20) and (16, 14).
3. For the idempotent to be S-idempotents we need a, b $\in$ $Z_{30}$ in all case with $a^2 \equiv a$, $b^2 = a$, ab = b we have 'b' such that a + b = 30.

These observations leads to a certain open problems which is given in Chapter 5.

***Example 3.4.6***: Let $Z_4$ = {0, 1, 2, 3} be the ring of integers modulo 4. $Z_4$ has no idempotents hence no S-idempotents.

***Example 3.4.7***: Let $Z_{16}$ = {0, 1, 2, …, 15} be the ring of integers modulo 16. $Z_{16}$ has no idempotents hence no S-idempotents.

***Example 3.4.8***: Let $Z_{27}$ = {0, 1, 2, …, 26}, be the ring of integers modulo 27. Clearly it can be verified $Z_{27}$ has no idempotents and so has no S-idempotents.

***Example 3.4.9***: Let $Z_3$ = {0, 1, 2} be the prime field of characteristic 3. G = $\langle g / g^2 = 1 \rangle$ be the cyclic group of order 2. $Z_3$G = {0, 1, 2, g, 2g, 1 + g, 1 + 2g, 2 + g, 2+2g}. Clearly 2 + 2g is a S-idempotent of $Z_3$G as $(2 + 2g)^2 = 2 + 2g$ and (2 + 2g) (1 + g) = 2 + 2g. Hence the claim.

***Example 3.4.10***: Let G = $\langle g / g^2 = 1 \rangle$ be the cyclic group of order 2 and Q be the field of rationals. QG be the group ring of the group G over Q, $\frac{1}{2}(1+g)$ is an idempotent of QG. $b = \dfrac{-1}{2}$ $(1 + g) \in$ QG is such that $b^2 = \dfrac{1}{2}(1 + g)$ and $ab = b$. So $\dfrac{1}{2}(1 + g)$ is a S-idempotent.

***Example 3.4.11***: Let $S_3 = \{1, p_1, p_2, p_3, p_4, p_5\}$ be the symmetric group of degree 3 where

$$1 = \begin{pmatrix} 1 & 2 & 3 \\ 1 & 2 & 3 \end{pmatrix}, \; p_1 = \begin{pmatrix} 1 & 2 & 3 \\ 1 & 3 & 2 \end{pmatrix}, \; p_2 = \begin{pmatrix} 1 & 2 & 3 \\ 3 & 2 & 1 \end{pmatrix},$$
$$p_3 = \begin{pmatrix} 1 & 2 & 3 \\ 2 & 1 & 3 \end{pmatrix}, \; p_4 = \begin{pmatrix} 1 & 2 & 3 \\ 2 & 3 & 1 \end{pmatrix} \text{ and } p_5 \begin{pmatrix} 1 & 2 & 3 \\ 3 & 1 & 2 \end{pmatrix}.$$

$Z_2 = \{0, 1\}$ be the prime field of characteristic two. Clearly the group ring $Z_2 S_3$ has idempotents which are S-idempotents. Now $(1 + p_4 + p_5)^2 = (1 + p_4 + p_5) = a$ take b $= 1 + p_1 + p_2 \in Z_2 S_3$. Now $b^2 = 1 + p_4 + p_5$, $ab = a$. Hence the claim. Now if take for the idempotent $1 + p_4 + p_5$ the element $p_1 + p_2 + p_3$ we will get $(p_1 + p_2 + p_3)^2 = 1 + p_4$ $+ p_5$ and $(1 + p_4 + p_5)(p_1 + p_2 + p_3) = p_1 + p_2 + p_3$.

This leads us to an interesting result that S-co-idempotents are not unique for a given idempotent.

**THEOREM 3.4.4**: *Let R be a ring. a $\in$ R be a S-idempotent. The S-co-idempotents of a in general is not unique.*

*Proof*: By an example. Consider example 3.4.11 the S-co-idempotent of $1 + p_4 + p_5$ is not unique.

***Example 3.4.12***: Let $Z_{105} = \{0, 1, 2, \ldots, 104\}$ be the ring of integers modulo 105. $(105 = 3 \times 5 \times 7)$. The idempotent in $Z_{105}$ are 15, 21, 36, 70, 85 and 91. It can be verified that all these idempotents are S-idempotents.

Now the S-co-idempotent for 15 is 90, for 21 is 84, 36 it is 69 for 70 the S-co-idempotent is 35, for 85 it is 20 and for 91 the S-co-idempotent is 14.

**THEOREM 3.4.5**: *Let F be a field. F has no S-idempotents.*

*Proof*: Since a field has no nontrivial idempotents so a field has no S-idempotents.

**THEOREM 3.4.6**: *Let F be a field of characteristic zero and G any group of finite order; the group ring FG has S-idempotents.*

*Proof*: Let FG be the group ring of G over F. Given G is of finite order. Two possibilities arise; order of G is prime or order of G is not a prime. Let order of G be a prime say p then $\alpha = \dfrac{+1}{p}(1 + g + g^2 + \ldots + g^{p-1})$ is such that $\alpha^2 = \alpha$ is an idempotent of FG.

For take $\alpha = \dfrac{-1}{p}(1 + g + g^2 + \ldots + g^{p-1})$. Clearly $b^2 = \alpha$ and $\alpha b = b$. Hence the claim. If the order of G is finite and not a prime then G has a subgroup say H or order m. Then by taking

$$a = \frac{1}{m}\sum_{h_i \in H} h_i$$

is such that $a^2 = a$ and take

$$b = \frac{-1}{m}\sum_{h_i \in H} h_i$$

is such that $b^2 = a$ and $ab = b$. Hence the claim.

**THEOREM 3.4.7**: *Let F be a field of characteristic 0 and G be a group having elements of finite order then the group ring FG has idempotents which are S-idempotents of FG.*

*Proof*: Let FG be the group ring of the group G over F. Given G has elements of finite order i.e., $g \in G$ is such that $g^m = 1$ ($m < \infty$). Take $a = \dfrac{1}{m}(1 + g + \ldots + g^{m-1})$ is such that $a^2 = a$ and if we take $b = \dfrac{-1}{m}(1 + g + g^2 + \ldots + g^{m-1})$ then $b^2 = a$ and $ab = b$. Hence the claim.

**THEOREM 3.4.8**: *Let G be a torsion free abelian group. F a field of characteristic zero. The group ring FG has no S-idempotents.*

*Proof*: Given G is a torsion free abelian group and F a field of characteristic zero. The group ring FG has no zero divisor, but for a ring to have S-idempotents it is guaranteed that the ring should have atleast two zero divisors. So this group ring cannot have S-idempotents as FG is a domain.

<u>**PROBLEMS:**</u>

1. Find all S-idempotents in $Z_{243} = \{0, 1, 2, \ldots, 242\}$.
2. Can $Z_{49} = \{0, 1, 2, \ldots, 48\}$ have nontrivial S-idempotents?

3. Find all S-idempotents in the group ring $Z_7G$ where $G = S_7$ the symmetric group of degree 7.

4. Can $Z_3S(4)$ the semigroup ring of the semigroup $S(4)$ over the ring $Z_3$ have S-idempotents? Justify your answer.

5. How many S-idempotents does $Z_2S(3)$ the semigroup ring of the semigroup $S(3)$ over $Z_2$ have?

6. Find all idempotents in $Z_{210} = \{0, 1, 2, \ldots, 209\}$, which are not S-idempotents.

7. Let $M_{3\times3} = \{(a_{ij}) \ / \ a_{ij} \in Z_6\}$ be the ring of all $3 \times 3$ matrices. Find all S-idempotents in $M_{3\times3}$.

8. Can $M_{5\times5} = \{(a_{ij}) \ / \ a_{ij} \in Z_{11}\}$, the ring of $5 \times 5$ matrices have idempotents? S-idempotents? Substantiate your answer.

9. Find a ring R which has idempotents but not S-idempotents.

10. Give an example of a ring R in which every idempotent is an S-idempotent.

## 3.5 Substructures in S-rings

In this section we introduce substructures in S-ring I and S-ring II. The notion of Smarandache rings, Smarandache ideals and Smarandache pseudo ideals is introduced in these S-rings I and II and they are illustrated by examples. Some interesting results about them are also obtained in this section.

**DEFINITION 3.5.1**: *Let S be a ring. A proper subset A of S is said to be a Smarandache subring (S-subring) of S if A has a proper subset B which is a field and A is a subring of S.*

**THEOREM 3.5.1**: *Let S be a ring. If S has S-subring then S is a S-ring I.*

*Proof*: Obvious from the fact that the ring has a S-subring A implies A contains a subfield which is also a subfield in S, so S is a S-ring I.

Suppose S is a S-ring I, it may not be always possible to obtain a S-subring in S or to be more precise every subring of a S-ring I need not in general be a S-subring of S.

***Example 3.5.1***: Let $Z_6 = \{0, 1, 2, 3, 4, 5\}$ be the ring of integers modulo 6. $Z_6$ is a S-ring I but $Z_6$ has no S-subring.

Clearly $Z_8 = \{0, 1, 2, 3, 4, 5, 6, 7\}$ has no subsets which are fields so $Z_8$ is not even a S-ring I.

***Example 3.5.2***: Let $Z_{12} = \{0, 1, 2, \ldots, 10, 11\}$ be the ring of integers modulo 12. $Z_{12}$ is a S-ring. In fact $Z_{12}$ has S-subring for take $P = \{0, 2, 4, 6, 8, 10\}$ and $A = \{0, 4, 8\}$ is a subfield of $Z_{12}$. So $Z_{12}$ has a S-subring.

**THEOREM 3.5.2**: *Let R be S-ring I, R may have subrings but R may not have S-subrings.*

*Proof*: Let R be a S-ring I, say $Z_6 = \{0, 1, 2, 3, 4, 5\}$ be ring of integers modulo 6. Clearly $Z_6$ is a S-ring I which has no S-subrings but has subrings $S_1 = \{0, 3\}$ and $S_2 = \{0, 4, 2\}$.

**DEFINITION [73, 60]**: *The Smarandache ideal is defined as an ideal A such that a proper subset of A is a field (with respect with the same induced operations). By proper subset we understand a set included in A, different from the empty set, from the unit element – if any and from A.*

***Example 3.5.3***: Let $Z_6$ be the S-ring i.e., $Z_6 = \{0, 1, 2, 3, 4, 5\}$. Clearly $I = \{0, 3\}$ and $J = \{0, 2, 4\}$ are ideals of $Z_6$ but none of them are S-ideals of $Z_6$.

**THEOREM 3.5.3**: *Let R be a ring if R has S-ideal then R is a S-ring. Conversely if R is a S-ring we cannot say every ideal in R is an S-ideal of R.*

*Proof*: Let R be a ring. If R has a S-ideal then we know R has a proper subset A which is a field, so R becomes a S-ring.

Now let R be a S-ring to show ideals of R need not be S-ideals of R. We prove by an example. Consider $Z_6 = \{0, 1, 2, 3, 4, 5\}$. This is a S-ring having ideals none of them are S-ideals of R.

***Example 3.5.4***: Let $Z_{10} = \{0, 1, 2, \ldots, 9\}$ be the ring of integers modulo 10. Clearly $Z_{10}$ is a S-ring having no S-ideals.

Now we proceed onto define Smarandache pseudo ideals in a S-ring.

**DEFINITION 3.5.2**: *Let (A, +, .) be a S-ring. B be a proper subset of A(B ⊂ A) which is a field. A non-empty subset S of A is said to be Smarandache pseudo right ideal (S-pseudo right ideal) of A related to B if*

> *1. (S, +) is an additive abelian group.*
> *2. For b ∈ B and s ∈ S we have s . b ∈ S.*

On similar lines we define Smarandache pseudo left ideal (S-pseudo left ideal). A non-empty subset S of A is said to be a Smarandache pseudo ideal (S-pseudo ideal), if S is both a S-pseudo right ideal and S-pseudo left ideal.

**_Remark_**: It is important to note that the phrase 'related to B' is important for if the field B is changed to $B^1$ the same S may not in general be a S-pseudo ideal related to $B^1$ also. Thus the S-pseudo ideals are different from usual ideal defined in a ring. Further we define S-pseudo ideal only when the ring itself is a S-ring I, otherwise we don't define S-pseudo ideal; for in case of S-ideals the ring by the very definition becomes a S-ring.

**_Example 3.5.5_**: Let Q[x] be the polynomial ring over the rationals. Clearly Q[x] is a S-ring. Consider S = $\langle n(x^2+1) / n \in Q \rangle$ be the set generated under addition. Now Q.S $\subset$ S and S.Q $\subset$ S, so S is a pseudo ideal of Q[x] related to Q.

**THEOREM 3.5.4**: *Let R be any S-ring. Any ideal of R is a S-pseudo ideal of R but in general, every S-pseudo ideal of R need not be an ideal of R*

*Proof*: Given R is a S-ring. So $\phi \neq$ B, B $\subset$ R, B is a field. Now I is an ideal of R, so IR $\subset$ I and RI $\subset$ I. Since B $\subset$ R we have BI $\subset$ I and IB $\subset$ I. Hence I is a S-pseudo ideal related to B.

To prove the converse, consider the S-ring given in example 3.5.5. Clearly S is a S-pseudo ideal but S is not an ideal of Q[x] as x.S not contained in S. Hence the claim.

**_Example 3.5.6_**: Let R be the field of reals. R[x] be the polynomial ring. Clearly R[x] is a S-ring. Now Q $\subset$ R[x] and R $\subset$ R[x] are fields contained in R[x]. Consider S= $\langle n(x^2+1) / n \in Q \rangle$ a group generated additively. Now S is a S-pseudo ideal relative to Q but is not a S-pseudo ideal relative to R. Thus this leads to the following result:

**THEOREM 3.5.5**: *Let R be a S-ring. Suppose A and B are two subfields of R; and S be a S-pseudo ideal related to A. S need not in general be a S-pseudo ideal related to B.*

*Proof*: By an example; in example 3.5.6 we see the set S is a S-pseudo ideal for the field Q and is not a S-pseudo ideal under the field of reals R.

**_Example 3.5.7_**: $Z_{12}$ = {0, 1, 2, …, 11} be the ring of integers modulo 12. Clearly $Z_{12}$ is a S-ring for A = {0, 4, 8} is a field in $Z_{12}$ with $4^2 \equiv 4$ (mod 12) acting as the multiplicative identity. Now S = {0, 6} is the S-pseudo ideal related to A. But S is also an ideal of $Z_{12}$. Every ideal of $Z_{12}$ is also a S-pseudo ideal of $Z_{12}$ related to A.

***Example 3.5.8***: Let $M_{2 \times 2}$ be the set of $2 \times 2$ matrices with entries from the prime field $Z_2 = \{0, 1\}$.

$$M_{2 \times 2} = \left\{ \begin{matrix} \begin{pmatrix} 0 & 0 \\ 0 & 0 \end{pmatrix}, \begin{pmatrix} 0 & 1 \\ 0 & 0 \end{pmatrix}, \begin{pmatrix} 1 & 0 \\ 0 & 0 \end{pmatrix}, \begin{pmatrix} 0 & 0 \\ 1 & 0 \end{pmatrix} \\ \begin{pmatrix} 0 & 0 \\ 0 & 1 \end{pmatrix}, \begin{pmatrix} 1 & 1 \\ 0 & 0 \end{pmatrix}, \begin{pmatrix} 0 & 0 \\ 1 & 1 \end{pmatrix}, \begin{pmatrix} 1 & 0 \\ 1 & 0 \end{pmatrix} \\ \begin{pmatrix} 0 & 1 \\ 0 & 1 \end{pmatrix}, \begin{pmatrix} 1 & 0 \\ 0 & 1 \end{pmatrix}, \begin{pmatrix} 0 & 1 \\ 1 & 0 \end{pmatrix}, \begin{pmatrix} 1 & 1 \\ 0 & 1 \end{pmatrix} \\ \begin{pmatrix} 1 & 0 \\ 1 & 1 \end{pmatrix}, \begin{pmatrix} 1 & 1 \\ 1 & 0 \end{pmatrix}, \begin{pmatrix} 0 & 1 \\ 1 & 1 \end{pmatrix} \text{and} \begin{pmatrix} 1 & 1 \\ 1 & 1 \end{pmatrix} \end{matrix} \right\}$$

be the ring of matrices under matrix addition and multiplication modulo 2.

Now $M_{2 \times 2}$ is a S-ring for

$$A = \left\{ \begin{pmatrix} 0 & 0 \\ 0 & 0 \end{pmatrix}, \begin{pmatrix} 1 & 0 \\ 0 & 0 \end{pmatrix} \right\}$$

is a field of $M_{2 \times 2}$. Let

$$S = \left\{ \begin{pmatrix} 0 & 0 \\ 0 & 0 \end{pmatrix}, \begin{pmatrix} 1 & 1 \\ 0 & 0 \end{pmatrix} \right\},$$

S is a S-pseudo left ideal related to A but S is not a S-pseudo right ideal related to A for

$$\begin{pmatrix} 1 & 1 \\ 0 & 0 \end{pmatrix} \cdot \begin{pmatrix} 1 & 0 \\ 0 & 0 \end{pmatrix} = \begin{pmatrix} 1 & 0 \\ 0 & 0 \end{pmatrix} \text{as} \begin{pmatrix} 1 & 0 \\ 0 & 0 \end{pmatrix} \notin S.$$

Now

$$B = \left\{ \begin{pmatrix} 0 & 0 \\ 0 & 0 \end{pmatrix}, \begin{pmatrix} 0 & 0 \\ 0 & 1 \end{pmatrix} \right\}$$

is also a field.

$$S = \left\{ \begin{pmatrix} 0 & 0 \\ 0 & 0 \end{pmatrix}, \begin{pmatrix} 1 & 1 \\ 0 & 0 \end{pmatrix} \right\}$$

is a left ideal related to B but not a right ideal related to B.

$$C = \left\{ \begin{pmatrix} 0 & 0 \\ 0 & 0 \end{pmatrix}, \begin{pmatrix} 0 & 0 \\ 1 & 1 \end{pmatrix} \right\}$$

is a field. Clearly S is not a S-pseudo left ideal with respect to C. But S is a S-pseudo right ideal with respect to C. Thus from the above example we obtain the following observation which is important to be noted.

***Remark***: A set S can be a S-pseudo ideal relative to more than one field. For

$$S = \left\{ \begin{pmatrix} 0 & 0 \\ 0 & 0 \end{pmatrix}, \begin{pmatrix} 1 & 1 \\ 0 & 0 \end{pmatrix} \right\}$$

is a S-pseudo left ideal related to both A and B. The same set S is not a S-pseudo left ideal with respect to the related field

$$C = \left\{ \begin{pmatrix} 0 & 0 \\ 0 & 0 \end{pmatrix}, \begin{pmatrix} 0 & 0 \\ 1 & 1 \end{pmatrix} \right\}$$

but S is a S-pseudo right ideal related to C.

Thus the same set S can be S-pseudo left ideal or right ideal depending on the related field. Clearly S is a S-pseudo ideal related to the field

$$D = \left\{ \begin{pmatrix} 0 & 0 \\ 0 & 0 \end{pmatrix}, \begin{pmatrix} 1 & 0 \\ 0 & 1 \end{pmatrix} \right\}.$$

**DEFINITION 3.5.3**: *Let R be a ring. I a S-ideal of R; we say I is a Smarandache minimal ideal (S-minimal ideal) of R if we have a $J \subset I$ where J is another S-ideal of R then J = I is the only ideal.*

**DEFINITION 3.5.4**: *Let R be a S-ring and M be a S-ideal of R, we say M is a Smarandache maximal ideal (S-maximal ideal) of R if we have another S-ideal N such that $M \subset N \subset R$ then the only possibility is M = N or N = R.*

***Example 3.5.9***: Let $Z_{15} = \{0, 1, 2, \ldots, 13, 14\}$ be the ring of integers modulo 15. Clearly I = $\{0, 3, 6, 9, 12\}$ is a S-ideal of $Z_{15}$ which is also a S-maximal ideal of $Z_{15}$.

***Example 3.5.10***: $Z_{14} = \{0, 1, 2, 3, 4, \ldots, 11, 12, 13\}$ be the ring of integers modulo 14. I = $\{0, 2, 4, 6, 8, 10, 12\}$ is a S-maximal ideal of $Z_{14}$.

***Example 3.5.11***: Let $Z_{12} = \{0, 1, 2, \ldots, 11\}$ be the ring of integers modulo 12. Now $I = \{0, 2, 4, 6, 8, 10\}$ is a S-ideal in fact S-maximal ideal. $J = \{0, 4, 8\}$ is an ideal which is a minimal ideal. Thus we have the concept of S-maximal ideal and no S-minimal ideal in the ring $Z_{12}$.

**DEFINITION 3.5.5**: *Let R be a S-ring and I be a S-pseudo ideal related to A. A $\subset$ R (A is a field). I is said to be a Smarandache minimal pseudo ideal (S-minimal pseudo ideal) of R if $I_1$, is another S-pseudo ideal related to A and (0) $\subset I_1 \subset I$ implies $I = I_1$ or $I_1 = (0)$ The minimality may vary with the different related fields.*

**DEFINITION 3.5.6**: *Let R be a S-ring. M is said to be Smarandache maximal pseudo ideal (S-maximal pseudo ideal) related to the field A, A $\subset$ R if $M_1$ is another S-pseudo ideal related to A and if $M \subset M_1$ then $M = M_1$.*

**DEFINITION 3.5.7**: *Let R be a S-ring, a S-pseudo ideal I related to a field, A, A $\subset$ R is said to be a Smarandache cyclic pseudo-ideal (S-cyclic pseudo-ideal) related to a field A, if I can be generated by a single element.*

**DEFINITION 3.5.8**: *Let R be a S-ring, a S-pseudo ideal I of R related to a field A is said to be a Smarandache prime pseudo ideal (S-prime pseudo-ideal) related to A if x. y $\in$ I implies x $\in$ I or y $\in$ I.*

***Example 3.5.12***: Let $Z_2 = \{0, 1\}$ be the prime field of characteristic 2. $Z_2[x]$ be the polynomial ring of degree less than or equal to 3, that is $Z_2[x] = \{0, 1, x, x^2 \ldots, 1 + x, 1 + x^2, \ldots, 1 + x + x^2 + x^3\}$. Clearly $Z_2[x]$ is a S-ring as it contains a field $Z_2$. $S = \{0, (1 + x), (1 + x^3), (x + x^3)\}$ is a S-pseudo ideal related to $Z_2$ and not related to $Z_2[x]$.

***Example 3.5.13***: Let $Z_2 = \{0, 1\}$ be the prime field of characteristic two. $S_3 = \{1, p_1, p_2, p_3, p_4, p_5\}$ be the symmetric group of degree 3. $Z_2 S_3$ be the group ring of the group $S_3$ over $Z_2$. $Z_2 S_3$ is a S-ring. $A = \{0, p_4 + p_5\}$ is a field. Let $S = \{0, 1 + p_1 + p_2 + p_3 + p_4 + p_5\}$ be the subset of $Z_2 S_3$. S is a S-pseudo ideal related to A and S is also a S-pseudo ideal related to $Z_2$.

**THEOREM 3.5.6**: *Let $Z_2 = \{0, 1\}$ be the prime field of characteristic 2, G any finite group of order n. Then $Z_2 G$ has S-pseudo ideals which are ideals of $Z_2 G$.*

*Proof*: Take $Z_2 = \{0, 1\}$ a field of characteristic two and the group ring $Z_2 G$ is a S-ring. Let $G = \{g_1, g_2, \ldots, g_{n-1}, 1\}$ be the set of all elements of G. $S = \{0, 1 + g_1 + \ldots + g_{n-1}\}$ is a S-pseudo ideal related to $Z_2$ and S is also an ideal of $Z_2 G$.

**DEFINITION 3.5.9**: *Let R be a S-ring II. A is a proper subset of R is a Smarandache subring II (S-subring II) of R if A is a subring and A itself is a S-ring II.*

**Example 3.5.14**: Let Z be the ring of integers; Z is a S-ring II and Z has S-subring II. Clearly Z is never a S-ring I or has a S-subring I.

**DEFINITION 3.5.10**: *Let Z[x] be the polynomial ring. Z[x] is a S-ring II. Also Z[x] has a S-subring II.*

**Example 3.5.15**: Let $pZ = \{0, \pm p, \dots \pm np, \dots\}$ be the ring (p > 3 and p a prime) $2pZ \subset pZ$ and $2pZ$ is a S-subring II.

**DEFINITION 3.5.11**: *Let R be a S-ring II, a non-empty subset I of R, is said to a Smarandache right (left) ideal II (S-right (left) ideal II) of R if*

1. *I is a S-subring II*
2. *Let $A \subset I$ be an integral domain or a division ring in I, then $ai \in I$ ($ia \in I$) for all $a \in A$ and $i \in I$. If I is simultaneously S-right ideal II and S-left ideal II then I is a Smarandache ideal II (S-ideal II) of R related to A.*

**DEFINITION 3.5.12**: *Let R be a ring if R is a S-ring I and has no S-ideals then we say R is a Smarandache simple ring I (S-simple ring I).*

**DEFINITION 3.5.13**: *Let R be a S-ring if R has no S-pseudo ideals, then we say R is a Smarandache pseudo simple ring (S-pseudo simple ring).*

**DEFINITION 3.5.14**: *Let R be a S-ring II, if R has no two sided S-ideals II then we say R is a Smarandache simple ring II (S-simple ring II).*

**Example 3.5.16**: Z is not a S-simple ring II.

**Example 3.5.17**: $Z_6 = \{0, 1, 2, \dots, 5\}$ is a S-simple ring II.

**Example 3.5.18**: Let $Z_{12} = \{0, 1, 2, \dots, 10, 11\}$ be the ring of integers modulo 12. $Z_{12}$ is a S-ring II which is not a S-simple ring II.

**DEFINITION 3.5.15**: *Let R be a S-ring I. I an S-ideal of R. R / I = {a + I /a ∈ R} is a Smarandache quotient ring I (S- quotient ring I) of R related to I.*

**DEFINITION 3.5.16**: *Let R be a S-ring. I a S-pseudo ideal of R; R/I = {a + I / a ∈ R} is a Smarandache pseudo quotient ring (S-pseudo quotient ring) of R related to I.*

**DEFINITION 3.5.17**: *Let R be a S-ring II, I be a S-ideal II. R / I = {a + I / a ∈ R} is defined as the Smarandache quotient ring II (S-quotient ring II) of R.*

<u>**PROBLEMS:**</u>

1. Does $Z_{15}$ have a S-subring?
2. Find S-ideals of $Z_{21}$.
3. Can every ideal of $Z_{28}$ be S-ideal? Substantiate your answer.
4. Prove $Z_{16}$ cannot have S-ideals.
5. Find S-subrings II of $Z_{120}$.
6. Can S-subring I be S-subring II?
7. Give an example of a ring R in which S-subring I and S-subring II are coincident.
8. Let $Z_{12}$ be a S-ring I find a suitable ideal I, so that $Z_{12}$/I is a S-quotient ring I.
9. Is $Z_{11}$ is a S-simple ring? Justify.
10. Is $Z_{13}$ a S-ring I?
11. Is $Z_{19}$ a S-ring II?
12. Can $Z_{23}$ be a S-pseudo simple ring?
13. Find all S-ideals I of $Z_{36}$.
14. Can $Z_{36}$, have S-ideal II?
15. Find in $Z_{36}$, S-pseudo ideal II.
16. Find for the ring $Z_{36}$
    i. S-quotient ring I.
    ii. S-quotient ring II.
    iii. S-pseudo quotient ring.

## 3.6 Smarandache modules

In this section we recall the definition of Smarandache R-module as given by Florentin Smarandache and Padilla Raul and proceed on to define Smarandache module II and Smarandache pseudo module. We illustrate them by examples and give some interesting results about them.

**DEFINITION [73, 60]**: *The Smarandache R-module (S-R module) is defined to be an R-module (A, +, ×) such that a proper subset of A is a S-algebra (with respect to the same induced operations and another '×' operation internal on A)* *where R is a commutative unitary Smarandache ring and S its proper subset which is a field.*

***Example 3.6.1***: Let R[x] be the polynomial ring in the variable x with coefficients from the real field R. Q[x] is a S-R module for it is a S-algebra.

***Example 3.6.2***: Let $R = Q \times Q \times Q$ be the ring. $S = Q \times \{1\} \times \{1\} \subset R$ is a field. A $= Q \times Q \times \{1\}$ is a S-R module over S.

But one may once again recall the definition of a module: "Let A be a ring. An A-module or a left A-module is an additive abelian group M having A as a left operator i.e., $a(x + y) = ax + ay$ for $a \in A$ and $x, y \in M$. Similarly right A-module. If M is simultaneously left and right A-module then we say M is a A-module."

Keeping this in view we can speak of S-modules I first, and then proceed onto define S-module II and S-pseudo modules. Now we have in case of S-modules the following situations:

1. A S-module relative to a subfield B may fail to be a S-module over some other subfield C.
2. Further we may have S-modules to be S-modules over every subfield.

The study of these concepts is innovative and interesting.

**DEFINITION 3.6.1**: *Let R be a S-ring I. A non-empty set B which is an additive abelian group is said to be a Smarandache right (left) module I (S-right(left) module I) relative to a S-subring I, A if $D \subset A$ where D is a field then $DB \subset B$ and $BD \subset B$ i.e. bd (and db) are in B with $b(d + c) = bd + dc$ for all $d, c \in D$ and $b \in B$ $((d + c)b = db + cb)$. If B is simultaneously a S-right module I and S-left module I over the same relative S-subring I then we say B is a Smarandache module I (S-module I).*

***Example 3.6.3***: Let $A = (M_{n \times n}, +, x)$ be the set of $n \times n$ matrices with entries from Q. Now consider R, the set of reals which is a S-ring. Now A is a S-module over the subfield Q. Clearly A is not a S-module over R. Further if we take $B = \{M_{n \times n}, \times, +\}$ the set of all $n \times n$ matrices with entries form Z, then we see B is not a S-module over any subfield of R. Motivated by this example and to overcome this problem we define S-module II.

**DEFINITION 3.6.2**: *Let R be a S-ring II. We say a non-empty set B which is an additive abelian group is said to be a Smarandache right (left) module II (S-right (left) module II) relative to a S-subring II, A if $D \subset A$ where D is a division ring or an integral domain, then $DB \subset B$ and $BD \subset B$; i.e., bd(and db) are in B. with* $b(d + c) = bd + bc \ \forall \ d, \ c \in D$ and $b \in B$ $((d + c) \ b = db + cb)$. If B is simultaneously a S-right module II and S-left module II over the same relative S-subring II then we say B is a Smarandache module II (S-module II).

**Example 3.6.4**: Let Z be a S-ring II, $M = M_{2 \times 2} = \{(a_{ij}) \ / \ a_{ij} \in 2Z\}$. M is a S-module II related to the S-ring II. $A = 2Z$. Clearly M is also a S-module II over the S-subring II, $A_1 = 4Z$ or $A_2 = 8Z$, but M is also S-module II over any $A_p = pZ$. Thus $M_{2 \times 2}$ is a S-module II over any S-subring II of Z

**Example 3.6.5**: Let Z[x] be the S-ring II, $M = Z[x]$, the polynomial ring with only polynomials of even degree. Then M is a S-module II over the S-subring Z but M is not a S-module II over the S-subring, $Y = \{$all polynomial of odd degree over Z$\}$, if we take; $A = \{$all odd degree polynomial with coefficient from 2Z$\}$ as the integral domain. Thus we see in case of S-module we see every S-ideal II is a S-module II.

**DEFINITION 3.6.3**: Let $(A, +, .)$ be a S-ring. B be a proper subset of A $(B \subset A)$ which is a field. A set M is said to be a Smarandache pseudo right (left) module (S-pseudo right(left) module)of A related to B if

1. $(M, +)$ is an additive abelian group
2. For $b \in B$ and $m \in M$ $m.b \in M$ $(b.m \in M)$
3. $(m_1 + m_2)b = m_1b + m_2b$, $(b.(m_1+m_2)=bm_1+bm_2)$ for $m_1, m_2 \in M$ and $b \in B$. If M is simultaneously a S-pseudo right module and S-pseudo left module, we say M is a Smarandache pseudo module (S-pseudo module) related to B.

Here also we wish to state if $M_1$ is a S-pseudo module related to B, $M_1$ need not be S-pseudo module related to some other subfield $B_1$ of A. Thus we see we can have different S-pseudo modules associated with different subfields in a ring.

**Example 3.6.6**: Let $Z_{24} = \{0, 1, \ldots, 23\}$ be the ring of integers modulo 24. I = $\{0, 2, 4, 6, \ldots, 22\}$ is an S-pseudo ideal II as well as, S-pseudo module of $Z_{24}$. For $\{0, 8, 16\}$ is a subfield of characteristic 3. $16^2 \equiv 16 \pmod{24}$, $16 \times 8 \equiv 8 \pmod{24}$. $8 \times 8 \equiv 16 \pmod{24}$. $Z_{24}[x]$ is a S-pseudo module related to the field P = $\{0, 8, 16\} \subset Z_{24}$

**Example 3.6.7**: Let $Z_2S_4$ be the group ring of the symmetric group of degree 4 over the field $Z_2$. $M = \{0, \Sigma g, g \in S_4\}$ ($\Sigma g$ denotes the sum of all elements from $S_4$). M is a S-module II over $Z_2$. M is a S-M-module II over $Z_2A_4$ Clearly M is also a S-ideal II and S-pseudo ideal of $Z_2S_4$.

It is left as an exercise for the reader to find in $Z_2S_4$ :
1. S-right module II.

2. S module II.
3. S-pseudo module II for different fields in $Z_2 S_4$.

## PROBLEMS:

1. For the S-ring $Z_{24}$. Find
   i. S-modules I,
   ii. S-modules II and
   iii. S-pseudo modules.
2. Find for the ring $Z[x]$ (The polynomial ring with coefficient from Z), the S-module II. Can $Z[x]$ have S-module I? Justify your answer.
3. Let $M_{n \times n} = \{(a_{ij}) \ / \ a_{ij} \in Z\}$ be the collection of all $n \times n$ matrices with entries from Z. Can $M_{n \times n}$ have S-pseudo modules? Substantiate your answer.
4. Let $M_{n \times n} = \{(a_{ij}) \ / \ a_{ij} \in Q\}$ be the collection of all $n \times n$ matrices with entries from Q. Can $M_{n \times n}$ have
   i. S-module I?
   ii. S-module II?
   iii. S-pseudo module?.
   Can the same abelian group A be such that it is simultaneously S-module I, S-module II and S-pseudo module?
5. Can the ring in problem 4 have S-right module I over a subfield A which are not S-left module I over the subfield A?
6. Let $ZS_3$ be the group ring of the symmetric group $S_3$ over the ring of integers Z. Can $ZS_3$ have S-module I? Find in $ZS_3$, S-right module II and S-left pseudo module.
7. Let $ZS(4)$ be the semigroup ring of the symmetric semigroup $S(4)$. Find a S-left module II in $ZS(4)$ which is not a S-right module II over the same S-subring II.
8. Does there exist an example of a ring in which no S-ideal I is a S-module I?
9. Does there exist a S-ring II in which every S-ideal II is a S-module II?
10. Give a S-pseudo module for the ring $R = Q \times Q$.
11. Let $R = Q \times Q \times Q \times Q \times Z$ be the ring. Find
    i. S-pseudo module.
    ii. S-module I.
    iii. S-module II of R.
12. For the ring QG where G is the Dihedral group, $G = D_{2n} = \{a, b \ / \ a^2 = b^n = 1;$ bab $= a\}$, Find
    i. S-right module I.
    ii. S-right module II.
    iii. S-right pseudo module.
    iv. S-module II.
    v. An S-ideal II which is a S-module II.

## 3.7 Rings satisfying S-A.C.C and S-D.C.C

In this section we define the concepts of Smarandache A.C.C and Smarandache D.C.C and obtain some interesting results about them. The chapter ends with several problems for the reader to solve.

We know the ring

$$A = \begin{pmatrix} Q & 0 \\ Q & Z \end{pmatrix}$$

is the best known example of a ring that is Noetherian on the right but not Noetherian on the left. The reader is entrusted to find such examples in case of S-Noetherian rings. For very recent work on Artinian modules over a matrix ring refer [64].

**DEFINITION 3.7.1**: *Let R be a ring, we say the ring R satisfies the Smarandache ascending chain condition (S-A.C.C for brevity) if for every ascending chain of S-ideals $I_j$ of R; that is $I_1 \subset I_2 \subset I_3 \subset \ldots$ is stationary in the sense that for some integer $p \geq 1$, $I_r = I_{r+1} = \ldots$. Similarly R is said to have the Smarandache descending chain condition (S-D.C.C for brevity) if every descending chain $N_1 \supset N_2 \supset \ldots \supset N_k \supset \ldots$ of S-ideals $N_j$ of R is stationary. Similarly one can define Smarandache-A.C.C and Smarandache D.C.C for S-right ideals and S-left ideals of a ring.*

**DEFINITION 3.7.2**: *A ring R is said to be Smarandache left Noetherian (or just Smarandache Noetherian) (S-Notherian) if the S-A.C.C on S-left ideals (or on S-ideals) is satisfied.*

**DEFINITION 3.7.3**: *A ring R is said to be Smarandache left Artinian (or just Smarandache Artinian) (S-Artinian) if for the S-left ideals (or S-ideals) of R satisfies the S-D.C.C condition.*

**_Remark_**: It is interesting to note that the matrix ring $A = M_{n \times n}$ over a division ring K is Noetherian as well as Artirian but we do not know whether $M_{n \times n}$ is S-Noetherian or S-Artinian.

**_Example 3.7.1_**: Let $Z_6 = \{0, 1, 2, \ldots, 5\}$ be the ring of integers modulo 6. $Z_6$ is a S-ring but has no S-ideals.

**_Example 3.7.2_**: Let $Z_{12}$ be the ring of integers modulo 12. The ideals of $Z_{12}$ are $\{0\}$, $I_1 = \{0, 2, 4, 6, 8, 10\}$, $I_2 = \{0, 3, 6, 9\}$, $I_3 = \{0, 6\}$, $I_4 = \{0, 4, 8\}$. $I_2$ is not an S-ideal, $I_1$ is an S-ideal for $A = \{0, 4, 8\}$ is a field in $I_1$, so we have $(0) \subset I_1 \subset Z_{12}$ is the S-A.C.C condition on the ring. $I_2$ is not even an S-ideal of $Z_{12}$.

***Example 3.7.3***: $Z_2G$ be the group ring where $G = \langle g / g^{12} = 1 \rangle$. The ideals of $Z_2G$ are $I_0 = \{0, (1 + g + \ldots + g^{11})\}$ which is not an S-ideal, $I_1 =$ Augmentation ideal of $Z_2G$; $I_2$ is a S-ideal for $\{0, g^8 + g^4\}$ is a field of characteristic two. We have $(0) \subset I_2 \subset Z_2G$ so $Z_2G$ satisfies S-A.C.C condition.

<u>PROBLEMS:</u>

1. Find S-ideals of $Z_{60}$.
2. Does the group ring $Z_2S_5$ have S-ideals?
3. Prove all augmentation ideals in $Z_2G$ are S-ideals ($G$ a finite group).
4. Can $Z_3S(4)$ have S-ideals?
5. Give an example of a group ring, which satisfies S-A.C.C.
6. Give an example of a group ring, which is not S-Artinian.
7. Give an example of a semigroup ring, which is S-Noetherian.
8. Find an example of a group ring, which is S-Noetherian.
9. Illustrate by an example a semigroup ring that can be S-Artinian.
10. Find a semigroup ring, which is not S-Noetherian.
11. Is the semigroup ring $Z_{20} S(4)$ S-Noetherian? Justify.
12. Can $Z_{12}S_3$ be S-Artinian? Prove your claim.

## 3.8 Some Special Types of Rings

The main motivation of this section is the introduction of the class of Smarandache semigroup rings, Smarandache group rings and give conditions for group rings and semigroup rings to be Smarandache rings. If RG happens to be group ring which is a S-ring it may still fail to be Smarandache group ring. Likewise a semigroup ring KS may be a S-ring but it may fail to be Smarandache semigroup ring for the semigroup S may not be S-semigroup. Further the concrete class of rings are reals R, rationals Q, integers Z, modulo integers $Z_n$, ring of matrices and polynomial rings but when we get to class of group rings and semigroup rings over those rings with standard well known groups and semigroups we get a very wide class of nice rings with varying properties.

Finally we get only from these ring a class of non-commutative rings apart from the ring of matrices. That is why we have taken special care not only to introduce group rings and semigroup rings in chapter I but also define Smarandache notions of these in this section. This section also discusses about matrix rings.

**Theorem 3.8.1**: *Let R be a field and G any group. The group ring RG is a S-ring.*

*Proof*: Since R ⊂ RG and R is a field; RG is a S-ring.

**Example 3.8.1**: Let $Z_2S_3$ be the group ring. Clearly $Z_2S_3$ is a S-ring.

All group rings are not in general S-rings by an example.

**Example 3.8.2**: Let $Z_4G$ be the group ring where G = ⟨g / g² = 1⟩; clearly the group ring $Z_4G$ is not a S-ring.

**Theorem 3.8.2**: *Let K be a field and S any semigroup with identity; KS the semigroup ring is a S-ring.*

*Proof*: Since K is a field and KS is a ring such that K ⊂ KS, is a S-ring.

All semigroup rings are not in general S-rings. The reader is requested to prove this.

**Definition [73, 60]**: *Let S be any semigroup. We say S is a Smarandache semigroup (S-semigroup) if S has a proper subset A which is a group under the operations of S.*

 We define Smarandache semigroup rings as follows.

**Definition 3.8.1**: *Let S be a semigroup, which is a S-semigroup and K, any field the semigroup ring KS is called a Smarandache semigroup ring (S-semigroup ring). So we see when a semigroup ring contains a group ring as a proper subset we call KS the Smarandache semigroup ring. It is to be noted that when we say KS is a Smarandache semigroup ring we do not demand KS to be a S-ring.*

**Example 3.8.3**: Let S = {0, 1, a, b} be a semigroup given by the following table:

| * | 0 | 1 | a | b |
|---|---|---|---|---|
| 0 | 0 | 0 | 0 | 0 |
| 1 | 0 | 1 | a | b |
| a | 0 | a | 0 | a |
| b | 0 | b | a | 1 |

Hence S is a S-semigroup. For {1, b} is a group in S.

Consider $Z_4S$ the semigroup ring clearly; $Z_4S$ is a S-semigroup ring which is not a S-ring.

**THEOREM 3.8.3**: *All S-semigroup rings in general are not S-rings.*

*Proof*: By an example the semigroup ring $Z_4S$ is not a S-ring but it is a S-semigroup ring.

Group rings are not S-semigroup rings for by the very definition of S-semigroup we take only a semigroup.

**DEFINITION 3.8.2**: *Let G be a group and K a S-ring. The group ring KG is called the Smarandache group ring (S-group ring).*

Note: K is only a S-ring.

Now the group ring KG when K is a field is always a S-ring. We see for the ring Z and G any group. ZG is not a S-group ring for Z is not a S-ring I.

**THEOREM 3.8.4**: *Let K be any commutative ring with 1 or any field. S(n) the symmetric semigroup. KS(n) is a S-semigroup ring.*

*Proof*: We know S(n) for any integer n, is a S-semigroup as $S_n$ is the symmetric group of degree n is a proper subset which is a group. Hence the claim.

**THEOREM 3.8.5**: *$Z_nG$ is a S-group ring for any group G where $Z_n$ is a ring such that there exist a $m \in Z_n$ with $m^2 \equiv m$ (mod n) and $m + m \equiv 0$ (mod n).*

*Proof*: $Z_nG$ is a S-group ring as $Z_n$ becomes a S-ring when $m \in Z_n$, is such that $m^2 \equiv m$ (mod n) and $m + m \equiv 0$ (mod n) as $A = \{0, m\}$ is a subfield of $Z_n$. Hence the theorem.

**THEOREM 3.8.6**: *Let $M_{n \times n} = \{(a_{ij}) / a_{ij} \in F, F$ a field or a S-ring$\}$ be the ring of $n \times n$ matrices. $M_{n \times n}$ is a S-ring.*

*Proof*: Let $A = \{(a_{ij}) / a_{11} \neq 0$ and all $a_{ij}$ are zero $a_{ij} \in F$ if F is a field or $a_{ij} \in B$ if F is a S-ring where $B \subset F$ is the subfield of F$\} \cup \{(0)\}$. $\{(0)\}$ denotes the $n \times n$ zero matrix. It is easily verified that A is a subfield in $M_{n \times n}$. Hence $M_{n \times n}$ is a S-ring.

It is interesting to note when $M_{n \times n}$ takes its entries from $Z_n$ the ring of integers modulo n when n is a composite number, we may have several more results.

This is a non-commutative S-ring, hence we can study S-right ideals, S-left ideals and concepts purely related to the non-commutative rings.

***Example 3.8.4***: Let $M_{2 \times 2} = \{(a_{ij}) \,/\, a_{ij} \in Z_2 = \{0, 1\}\}$, clearly $M_{2 \times 2}$ is a S-ring.

***Example 3.8.5***: Let $ZS(4)$ be the S-semigroup ring. Can $ZS(4)$ have a proper subset which is a field?

***Example 3.8.6***: Find for the group ring $Z_6 S_3$ a proper subset which is a field, apart from the fields $A_1 = \{0, 3\}$ and $A_2 = \{0, 2, 4\}$.

***Example 3.8.7***: Find any proper subset which is a field in the group ring $Z_2 S(3)$ apart from $Z_2 = \{0, 1\}$.

<u>**PROBLEMS:**</u>

1. Prove $Z_{12} S_5$ is a S-group ring.
2. Show $Z_2 S(7)$ is a S-semigroup ring. (Hint: To prove this show $S(7)$ has a proper subset which is a subgroup).
3. $Z_7 S_5$ is a S-ring. Justify.
4. Find all the proper subsets which are fields in the group ring $Z_3 S_4$.
5. Does $M_{2 \times 2} = \{(a_{ij}) \,/\, a_{ij} \in Z_4\}$ have proper subsets which are fields? Is $M_{2 \times 2}$ a S-ring? Justify your answer.
6. Prove $M_{3 \times 3} = \{(a_{ij}) \,/\, a_{ij} \in Z_6\}$ is a S-ring.
7. Find all proper subsets which are fields in $Z_3 G$ where $G = \langle g \,/\, g^7 = 1 \rangle$.
8. How many proper subsets in $Z_3 S(3)$ are fields?
9. Does there exist a S-semigroup ring which is not a S-ring?
10. Give an example of S-semigroup ring of order 64 which is a S-ring.
11. In the matrix ring $M_{n \times n} = \{(a_{ij}) \,/\, a_{ij} \in Z\}$, can we find a subset $P \subset M_{n \times n}$ such that $P$ is a field?

## 3.9 Special elements in S-rings

In this section we introduce the concepts of Smarandache nilpotent elements, Smarandache semi idempotents, Smarandache pseudo commutative pair, Smarandache-pseudo commutative ring, Smarandache strongly regular rings Smarandache quasi-commutative ring and finally the concept of Smarandache nilpotent elements. Several properties enjoyed by these Smarandache notions are proved and some of them are illustrated by examples and several of them are left as an exercise for the reader.

**DEFINITION 3.9.1**: *Let R be a ring. A nilpotent element $0 \neq x \in R$ is said to be a Smarandache nilpotent element (S-nilpotent element) if $x^n = 0$ and there exists a* $y \in R \setminus \{0, x\}$ such that $x^r y = 0$ or $yx^s = 0$, $r, s > 0$ and $y^m \neq 0$ for any integers $m > 1$.

**Example 3.9.1**: Let $Z_{12} = \{0, 1, 2, 3, \ldots, 11\}$ be the ring of integers modulo 12. Clearly $6^2 \equiv 0 \pmod{12}$, $6 \cdot 8 \equiv 0 \pmod{12}$. But $8^m \not\equiv 0 \pmod{12}$ for any m > 0 as $8^3 \equiv 8 \pmod{12}$. Thus 6 is a S-nilpotent element of $Z_{12}$.

**Example 3.9.2**: Let $Z_8 = \{0, 1, 2, 3, \ldots, 7\}$ be ring of integers modulo 8. $2^3 \equiv 0 \pmod 8$, $4^2 \equiv 0 \pmod 8$ and $6^3 \equiv 0 \pmod 8$. These are nilpotents but none of them are S-nilpotents.

In view of this we have the following theorem:

**THEOREM 3.9.1**: *Let R be a ring. Every S-nilpotent element of R is a nilpotent element of R. But in general every nilpotent element of R need not be S-nilpotent element of R.*

*Proof*: By the very definition of S-nilpotent element we see every S-nilpotent element is a nilpotent element of R. But all nilpotents in general need not be S-nilpotents. By example 3.9.2 we see the theorem is evident.

**DEFINITION [24]**: *An element $\alpha \neq 0$ of a ring R is called semi idempotent if and only if $\alpha$ is not in any proper two sided ideal of R generated by $\alpha^2 - \alpha$; i.e., $\alpha \notin R(\alpha^2 - \alpha) R$ or $R = R(\alpha^2 - \alpha)R$. 0 is also counted among semi idempotents.*

Now we proceed onto define Smarandache -semi idempotents.

**DEFINITION 3.9.2**: *Let R be a ring an element $\alpha \in R \setminus \{0\}$ is said to be a Smarandache- semi idempotent I (S-semi idempotent I), if the ideal generated by $(\alpha^2 - \alpha)$ that is $R(\alpha^2 - \alpha) R$ is a S-ideal I and $\alpha \notin R(\alpha^2 - \alpha)R$ or $R = R(\alpha^2 - \alpha) R$. We say $\alpha$ is a Smarandache semi idempotent II (S-semi idempotent II) if the ideal generated by $\alpha^2 - \alpha$ i.e., $R(\alpha^2 - \alpha)R$ is a S-ideal II and $\alpha \notin R(\alpha^2 - \alpha)R$ or $R = R(\alpha^2 - \alpha)R$.*

**THEOREM 3.9.2**: *Every semi idempotent of a ring R in general need not be a S-semi idempotent of R.*

*Proof*: Let $Z_{24} = \{0, 1, 2, \ldots, 23\}$ be the ring of integers modulo 24. $4 \in Z_{24}$ is a semi idempotent. For $\alpha = 4^2 - 4$ generates an ideal $I = \{0, 12\}$. Clearly I is not a S-ideal so 4 is not S-semi idempotent but 4 is only a semi idmepotent. Thus every semi idempotent need not in general be a S-semi idempotent.

**Theorem 3.9.3**: *Let R be a ring every S-semi idempotent I is a semi-idempotent of R.*

*Proof*: We know if $\alpha \in R \setminus \{0\}$ is a S-semi idempotent. Then $\alpha \notin R(\alpha^2 - \alpha)$ R or R $= R(\alpha^2 - \alpha)R$ where $R(\alpha^2 - \alpha)R$ is an S-ideal I of R. But all S-ideals are ideals. Hence the claim.

**Example 3.9.3**: Let $Z_{24} = \{0, 1, 2, \ldots, 23\}$ be the ring of integers modulo 24. $5 \in Z_{24}$ is a S-semi idempotent I of $Z_{24}$. For consider the ideal generated by $\alpha^2 - \alpha = 5^2 - 5 = 20$. $\langle \alpha^2 - \alpha \rangle = \langle 20 \rangle = \{0, 20, 16, 12, 4, 8\} = I$. Clearly $(0, 8, 16) = J \subset I$ is a field isomorphic to the prime field of characteristic 3. $16^2 \equiv 16 \pmod{24}$ acts as unit. $5 \notin I$ so 5 is a S-semi idempotent I of $Z_{24}$.

**Definition [151]**: *Let R be a non-commutative ring. A pair of distinct elements x, y $\in$ R different from the identity of R which are such that xy = yx is said to be a pseudo commutative pair of R if xay = yax for all a $\in$ R. If in a ring R every commutative pair happens to be a pseudo commutative pair of R then R is said to be a pseudo commutative ring.*

Clearly every commutative ring is trivially pseudo commutative.

**Definition 3.9.3**: *Let R be ring with A, a S-subring of R. A pair of elements x, y $\in$ A which are such that xy = yx is said to be a Smarandache pseudo commutative pair (S-pseudo commutative pair) of R if xay = yax for all a $\in$ A. If in a S-subring A, every commuting pair happens to be a S-pseudo commutative pair of A then A is said to be a Smarandache pseudo commutative ring (S-pseudo commutative ring).*

**Theorem 3.9.4**: *Let R be a ring if R is a S-pseudo commutative ring then R is a S-ring.*

*Proof*: Follows from the fact that if R is a S-pseudo commutative ring then R has a S-subring which immediately by the definition of S-subring makes R is a S-ring.

**Definition [151]**: *Let R be a non-commutative ring. A commuting distinct pair of elements x, y $\in$ R is said to be pseudo commutative with respect to a non-empty subset S of R if xsy = ysx for all s $\in$ S.*

**Definition 3.9.4**: *Let R be a non-commutative ring. A commuting distinct pair of elements x, y in R is said to be Smarandache pseudo commutative pair (S-pseudo commutative pair) with respect to a S-subring B of R, if xsy = ysx for all s $\in$ B.*

**THEOREM 3.9.5**: *Let R be a ring having a commuting pair, which is a S-pseudo commutative then R is a S-ring.*

*Proof*: By the very definition of the S-pseudo commutative pair we see the ring R must be a S-ring.

It is left as an exercise for the reader to show if R is a S-ring having a commuting pair still R need not be S-pseudo commutative.

**THEOREM 3.9.6**: *Let $Z_pS_n$ be a group ring of the group $S_n$ over the prime field $Z_p$. $Z_pS_n$ is S-pseudo commutative ring.*

*Proof*: Now $Z_pS_n$ is a S-ring. $A = Z_pB$ where

$$B = \left\{ \begin{pmatrix} 1 & 2 & 3 & 4 & . & . & . & n \\ 1 & 2 & 3 & 4 & . & . & . & n \end{pmatrix}, \begin{pmatrix} 1 & 2 & 3 & 4 & . & . & . & n \\ 2 & 3 & 1 & 4 & . & . & . & n \end{pmatrix} \begin{pmatrix} 1 & 2 & 3 & 4 & . & . & . & n \\ 3 & 1 & 2 & 4 & . & . & . & n \end{pmatrix} \right\}$$

is a subgroup of $S_n$ is a S-subring of $Z_pS_n$. Now take

$$x = \begin{pmatrix} 1 & 2 & 3 & 4 & . & . & . & n \\ 2 & 3 & 1 & 4 & . & . & . & n \end{pmatrix} \text{ and } y = \begin{pmatrix} 1 & 2 & 3 & 4 & . & . & . & n \\ 3 & 1 & 2 & 4 & . & . & . & n \end{pmatrix}$$

we see $xy = yx$ in $Z_pS_n$ and $xay = yax$ for all $a \in Z_pB = A$. Hence the claim.

**THEOREM 3.9.7**: *A S-pseudo commutative ring need not in general be a pseudo commutative ring.*

*Proof*: The example given in theorem 3.9.6 viz. the group ring $Z_pS_n$ is a S-pseudo commutative ring but it is clearly not a pseudo commutative ring, hence the claim.

**THEOREM 3.9.8**: *Let R be a ring. If Z(R) denotes the center of R and Z(R) is a S-subring R, which is nontrivial then R is a S-pseudo commutative ring.*

*Proof*: By the very definition of S-pseudo commutative ring we see R satisfies the conditions in the theorem 3.9.6; hence R is a S-pseudo commutative ring.

**DEFINITION [48]**: *Let R be a ring. For every x, y $\in$ R if we have $(xy)^n = xy$ for some integer, n = n(xy) > 1 then R is called a strongly regular ring.*

**DEFINITION 3.9.5**: *Let R be a ring. We say R is a Smarandache strongly regular ring (S-strongly regular ring) if R contains a S-subring B such that for every x,y in B we have $(xy)^n = xy$ for some integer $n = n(x,y) > 1$.*

We have the following interesting result.

**THEOREM 3.9.9**: *Let R be a ring which is strongly regular ring then R is a S-strongly regular ring provided R has non-trivial S-subring.*

*Proof*: Obvious by the very definition of strongly regular ring and S-strongly regular ring.

**THEOREM 3.9.10**: *Let $Z_p$ be the prime field and S be an ordered semigroup with identity then the semigroup ring $Z_nS$ is a S-strongly regular ring and not a strongly regular ring.*

*Proof*: $Z_pS$ is the semigroup ring. Now $Z_p$ is a S-subring of $Z_pS$. Clearly $Z_pS$ is a S-strongly regular ring.

Now $Z_pS$ is not a strongly regular ring. For given S is an ordered semigroup with 1. Let $\alpha, \beta \in RS$ with $\alpha = \Sigma a_i s_i$ and $\beta = \Sigma \beta_j h_j$ $1 \leq j \leq m$, $1 \leq i \leq n$, $\alpha_i \neq 0$ and $\beta_j \neq 0$. $s_1, \ldots, s_n$ and $h_1, \ldots, h_m$ are assumed to be distinct and

$$s_1 < s_2 < \ldots < s_n$$
$$h_1 < h_2 < \ldots < h_m.$$

It is given S is an ordered semigroup. So in $(\alpha\beta)^p$ we have $(s_1 h_1)^p$ to be the least element and $(s_n h_m)^p$ to be the largest element. Hence $(\alpha\beta)^p \neq \alpha\beta$. $p > 1$. Thus the semigroup ring is not a strongly regular ring only a S-strongly regular ring.

**DEFINITION [39]**: *Let R be a ring, R is called quasi commutative if $ab = b^\gamma a$ for every pair of elements a, b $\in$ R and $\gamma > 1$.*

**THEOREM [130]**: *Let R be a quasi commutative ring. Then for every pair of elements a, b $\in$ R there exists s $\in$ R such that $a^2 b = bs^2$.*

*Proof*: Given R is a quasi commutative ring so $ab = b^\gamma a$ for every pair of elements a,b in R, $\gamma \geq 1$. Now $ab = b^\gamma a$. $a^2 b = ab^\gamma a = ab(b^{\gamma-1}a) = b^\gamma a.b^{\gamma 1}a. = b (b^{\gamma-1}a)^2 = bs^2$ where s $\in$ R.

**THEOREM [130]**: *Let R be a ring in which we have a pair of elements a, b ∈ R such that there exists an s ∈ R with $a^2b = bs^2$ then we need not have $ab = b^\gamma a$ in R for some γ > 1.*

*Proof*: By an example. Let $Z_2 = \{0, 1\}$ be the prime field of characteristic two and

$$S_3 = \left\{ e = \begin{pmatrix} 1 & 2 & 3 \\ 2 & 1 & 3 \end{pmatrix}, \quad p_1 = \begin{pmatrix} 1 & 2 & 3 \\ 1 & 3 & 2 \end{pmatrix}, \right.$$
$$p_2 = \begin{pmatrix} 1 & 2 & 3 \\ 3 & 2 & 1 \end{pmatrix}, \quad p_3 = \begin{pmatrix} 1 & 2 & 3 \\ 2 & 1 & 3 \end{pmatrix},$$
$$\left. p_4 = \begin{pmatrix} 1 & 2 & 3 \\ 2 & 3 & 1 \end{pmatrix} \text{ and } p_5 = \begin{pmatrix} 1 & 2 & 3 \\ 3 & 1 & 2 \end{pmatrix} \right\}.$$

Let $Z_2S_3$ be the group ring of the group $S_3$ over $Z_2$. $p_1p_2 = p_2{}^\gamma p_1$ is not possible for any γ for if γ is even. $p_1p_2 = p_5 \neq ep_1 = p_1$ if γ is odd then $p_2{}^\gamma = p_2$ so $p_2p_1 = p_4$ and $p_4 \neq p_5$. Hence $Z_2S_3$ is not quasi commutative.

**DEFINITION 3.9.6**: *Let R be a ring. We say R is a Smarandache quasi-commutative ring (S-quasi commutative ring) if for any S-subring, A of R we have $ab = b^\gamma a$ for every a, b ∈ A; γ ≥ 1.*

**THEOREM 3.9.11**: *If R is a S-quasi commutative ring then R is S-ring.*

*Proof*: Obvious by the very definition of S-quasi commutative ring.

**THEOREM 3.9.12**: *Every S-ring in general is not S-quasi commutative.*

*Proof*: The ring $Z_6 = \{0, 1, 2, 3, 4, 5\}$ is a S-ring. This has no S-subring so the very concept of S-quasi commutative cannot be defined.

**THEOREM 3.9.13**: *Let G be a torsion free non-abelian group R be any S-ring which is quasi commutative. The group ring RG is S-quasi commutative.*

*Proof*: Given R is quasi commutative and is a S-ring so R is S-quasi commutative. Now R ⊂ RG, so R is a S-subring which is quasi commutative; hence RG is a S-quasi commutative ring.

**DEFINITION [140]**: *An element x of an associative ring R is called semi nilpotent if $x^n - x$ is a nilpotent element of R. If $x^n - x = 0$ we say x is a trivial semi nilpotent.*

**Theorem [140]**: *If $x$ is a nilpotent element of a ring $R$ then $x$ is a semi nilpotent element of $R$.*

*Proof*: Given $x \in R$ is nilpotent so $x^n = 0$ clearly $x^n - x = -x$ so $(-x)^n = 0$ hence our claim.

**Theorem [140]**: *Let $R$ be a ring an unit in $R$ can also be semi nilpotent.*

*Proof*: Let $Z_2 G$ be the group ring of the group $G = \langle g \, / \, g^2 = 1 \rangle$ over the field $Z_2 = \{0, 1\}$. Clearly $g \in Z_2 G$ is such that $g^2 = 1$, so $g$ is a unit of $Z_2 G$ but $g^2 - g = 1 + g$ is nilpotent as $(1 + g)^2 = 0$; hence the claim.

**Theorem [140]**: *Let $R$ be a ring every idempotent in $R$ is semi nilpotent.*

*Proof*: It is left for the reader to prove.

**Theorem [140]**: *Let $K$ be a field of characteristic 0. $G$, a torsion free abelian group. The group ring $KG$ has no nontrivial semi nilpotents.*

*Proof*: Since $KG$ is a domain $KG$ has no zero divisors; so it cannot have semi nilpotents.

**Definition 3.9.7**: *Let $R$ be a ring. An element $x \in R$ is a Smarandache semi nilpotent (S-semi nilpotent) if $x^n - x$ is S-nilpotent.*

**Example 3.9.4**: Can the ring $Z_{24}$ have S-semi nilpotents?

**PROBLEMS:**

1. Find S-nilpotents in the commutative ring $Z_{15}$.
2. Find S-nilpotents of the group ring $Z_2 S_3$.
3. Can the semigroup ring $Z_3 S(4)$ have S-nilpotent? Justify your claim.
4. Find for $Z_{30} = \{0, 1, 2, \ldots, 29\}$ the ring of integers modulo 30, all S-semi idempotents.
5. Does $Z_3 S_5$ have S-idempotents?
6. Find all S-semi idempotents of the semigroup ring $Z_7 S(3)$.
7. Can the ring $M_{5 \times 5} = \{(a_{ij}) \, / \, a_{ij} \in Z_4\}$ have S-semi idempotents? If so find them.
8. Prove $Z_{11} S_5$ is a S-pseudo commutative ring.
9. Give an example of a S-commutative ring.
10. Is the group ring $Z_7 S_3$ a S-strongly regular ring?
11. Is $Z_{25} S_3$ a S-quasi commutative ring?

12.    Can $Z_{27}$ have S-nilpotents?

## 3.10 Special Properties about Smarandache rings

In this section we introduce special properties about Smarandache rings, which are not found in any book. The sole purpose of this section is to define over seventy new Smarandache notions on rings and these concepts illustrated with examples. The vitality of this section is the recollection of several ring theoretical concepts which are interesting and give a Smarandache-ic equivalent of them. Thus this section will not only attract Smarandache researchers but also ring theorist. Finally it ends with 70 problems for the reader to solve to get involved and through with these concepts.

**DEFINITION ([88])**: *A ring R is said to be reduced, if R has no non-zero nilpotent elements.*

***Example 3.10.1***: Z the ring of integers is a reduced ring.

***Example 3.10.2***: $Z_p[x]$ the polynomial ring with coefficients from $Z_p$, p a prime is a reduced ring.

***Example 3.10.3***: $Z_{12} = \{0, 1, 2, \ldots, 11\}$ the ring of integers modulo 12 is not a reduced ring for $6^2 \equiv 0 \pmod{12}$.

**DEFINITION 3.10.1**: *Let R be a ring. R is said to be a Smarandache reduced ring (S-reduced ring) if R has no S-nilpotent elements.*

***Example 3.10.4***: $Z_4 = \{0, 1, 2, 3\}$ the ring of integers modulo 4 is a S-reduced ring. For it has no S-nilpotent elements.

***Example 3.10.5***: $Z_9 = \{0, 1, 2, \ldots, 8\}$ the ring of integers modulo 9 is a S-reduced ring for it has no S-nilpotents only nilpotents, viz. $3^2 \equiv 0 \pmod 9$ and $6^2 \equiv 0 \pmod 9$.

**THEOREM 3.10.1**: *Let R be a reduced ring then R is a S-reduced ring. If R is a S-reduced ring then R need not be a reduced ring.*

*Proof*: If R is a reduced ring R has no nilpotents so R cannot have S-nilpotents so R is a S-reduced ring.

Conversely if R is a S-reduced ring, R need not be a reduced ring. For the rings $Z_4$ and $Z_9$ are S-reduced rings but clearly $Z_4$ and $Z_9$ are not reduced rings.

**DEFINITION [68]**: *A ring R is a zero square ring if $x^2 = 0$ for all $x \in R$.*

**DEFINITION 3.10.2**: *Let R be a ring. We say R is a Smarandache zero square ring (S-zero square ring) if R has S-subring A having a subring B contained in A which is a zero square ring.*

**Example 3.10.6**: Let $Z_{12} = \{0, 1, 2, 3, \ldots, 11\}$ be the ring of integers modulo 12. I = $\{0, 2, 4, 6, 8, 10\}$ is a subring which is a S-subring as $\{0, 8, 4\}$ is a field. Now P = $\{0, 6\}$ is a zero square subring in I so $Z_{12}$ is a S-zero square ring, but clearly $Z_{12}$ is not a zero square ring.

**THEOREM 3.10.2**: *Every zero square ring is never a S-zero square ring.*

*Proof*: Given R is a zero square ring so $a^2 = 0$ for every $a \in R$. So if R has a S-subring say A, then A must have a subset which is a field, so in A we cannot have $a^2 = 0$ for all $a \in A$. Hence the claim.

**THEOREM 3.10.3**: *Every S-zero square ring is never a zero square ring.*

*Proof*: For if R is a S-zero square ring it has a proper subset which is a field and in a field we cannot have $a^2 = 0$ for all $a \in R$. This substantiated by an example 3.10.6 the ring $Z_{12}$ is a S-zero square ring. But clearly $Z_{12}$ is not a zero square ring. For we have several elements in $Z_{12}$ whose square is not zero.

**Example 3.10.7**: Let $Z_{12}$ be the ring and G any group; $Z_{12}G$ be the group ring of G over $Z_{12}$. $Z_{12}G$ is also a S-zero square ring; in view of this we have the following theorem.

**THEOREM 3.10.4**: *Let R be a S-commutative ring of characteristic 0. If R is a S-zero square ring then in R we have xy = 0 for all x, y ∈ A; A the subring of the S-subring B of R.*

*Proof*: We have R is a S-zero square ring; so $x^2 = 0$ for all $x \in A$, A the subring of the S-subring B of R. We have $A \subset R$. Let x, y $\in$ A. Now $(x + y)^2 = 0$ i.e., 2xy = 0 so xy = 0; hence the claim.

**THEOREM 3.10.5**: *Let be a ring. If R is not a S-commutative ring and if R is a S-zero square ring with A the subring of the S-subring B of R is also non-commutative, then every pair in A is anti-commutative.*

*Proof*: Let R be a S-zero square ring i.e., R has a S-subring B such that $A \subset B$ where A is a subring of B, is a zero square ring.

If A is non-commutative but is a zero square ring so we have $x^2 = 0$ $\forall$ $x \in A$. So $(x + y)^2 = 0$ using $x^2 = y^2 = 0$ we have $xy + yx = 0$. So elements in A are anti-commutative.

**THEOREM 3.10.6**: *Let R be a commutative ring with 1 of characteristic 0. G a commutative group (or S a commutative semigroup with 1). RG (RS) is a S-zero square ring if and only if $A^2 = 0$ where A is a subring of a S-subring B of RG (i.e. $A \subset B \subset RG$).*

*Proof*: If RG (RS) is a S-zero square ring then we have $A \subset B \subset RG$ (RS) where A is a subring of B where B is a S-subring of RG. We have A to be a zero square ring by theorem 3.10.4, x . y = 0 for all x, y $\in$ A. Hence $A^2 = 0$. Conversely if $A^2 = 0$ and A is a subring of the S-subring B of RG we have RG to be a S-zero square ring.

Now we leave it as an exercise to the reader the case when G is a non-commutative group.

**DEFINITION [94]**: *Let R be ring. R is called a inner zero square ring if every proper subring of R is a zero square ring.*

**Example 3.10.8**: $Z_4 = \{0, 1, 2, 3\}$ is a inner zero square ring as $\{0, 2\}$ is the only subring, and it is a zero square ring.

Now we proceed on to define Smarandache inner zero square ring.

**DEFINITION 3.10.3**: *Let R be a ring. If every S-subring A of R has a subring B $\subset$ A such that B is an inner zero square ring then we say R is a Smarandache inner square ring (S-inner square ring).*

**Example 3.10.9**: $Z_{12} = \{0, 1, 2, \ldots, 11\}$ is a S-inner zero square ring. For the S-subring $A_1$ of $Z_{12}$, $A_1 = \{0, 2, 4, 6, 8, 10\}$, has B = $\{0, 6\}$ to be an inner zero square ring. Clearly $Z_{12}$ is not an inner zero square ring but is a S-inner square ring.

In view of this we have the following.

**THEOREM 3.10.7**: *Let R be a inner zero square ring then, R in general need not be a S-inner zero square ring. Further if R be a S-inner zero square ring. R is not an inner zero square ring.*

*Proof*: By the above example, now even if R is a inner zero square ring we may not have R to be a S-inner zero square ring for if R is to have S-subring A $\subset$ R then A should contain a field as a proper subset. So if R is a S-inner zero square ring R is never a inner zero square ring.

We define Smarandache weak inner zero square ring.

**DEFINITION 3.10.4**: *Let R be a ring. We say R is a Smarandache weak inner zero square ring (S-weak inner zero square ring) if R has atleast a S-subring $A \subset R$ such that a subring B of A is a zero square ring.*

**THEOREM 3.10.8**: *Let R be a S-inner zero square ring. Then R is a S-weak inner zero square ring.*

*Proof*: Let R be a S-inner zero square ring then obviously by the very definition, R is a S-weak inner zero square ring.

**THEOREM 3.10.9**: *Let G be any group and R a S-inner zero square ring. The group ring RG is a S-weak inner zero square ring.*

*Proof*: Since $R \subset RG$; we have R to be S-inner zero square ring so RG is a S-weak inner zero square ring.

***Example 3.10.10***: Let G be a torsion free abelian group and R a S-inner zero square ring. The group ring RG is only a S-weak inner zero square ring.

The concept of S-inner zero square ring is important as we see a same S-subring; has two subsets one a field one a zero square ring. Except for Smarandache notions, this is an impossibility in the same substructure.

In case of semigroup ring we have the following theorem for which we need to define a new Smarandache notion about semigroups.

**DEFINITION 3.10.5**: *Let S be a multiplicative semigroup with 0, we say S is a Smarandache null semigroup (S-null semigroup) if we have a proper subsemigroup $P \subset S$ such that in P we have*

>    1. *$p^2 = 0$ for every $p \in P$ and*
>    2. *$p_i p_j = p_j p_i = 0$ for every $p_i, p_j \in P$.*

*We say S is a Smarandache strong null semigroup (S-strong null semigroup) if every subsemigroup P of S satisfies 1 and 2.*

***Example 3.10.11***: Let $Z_4 = \{0, 1, 2, 3\}$ be the semigroup under multiplication modulo 4. $Z_4$ is a S-null semigroup; for $\{0, 2\} = P$ is such that $2^2 \equiv 0 \pmod{4}$.

***Example 3.10.12***: Let $Z_6 = \{0, 1, 2, 3, 4, 5\}$ be the semigroup under multiplication modulo 6. $Z_6$ is not a S-null semigroup.

***Example 3.10.13***: Let $Z_8 = \{0, 1, 2, \ldots, 7\}$ be the semigroup under multiplication modulo 8. P = $\{0, \overline{4}\}$ is such that $\overline{4}^2 \equiv 0 \pmod 8$ so $Z_8$ is a S-null semigroup.

**DEFINITION 3.10.6**: *Let R be a ring. R is said to be a Smarandache null ring (S-null ring) if R has a S-subring A and A has a subring P such that in P we have*

       *1. $p^2 = 0$ for all $p \in P$.*
       *2. $p_i p_j = p_j p_i = 0$ for all $p_i p_j \in P$.*

Thus S-null ring localizes the null ring property.

**THEOREM 3.10.10**: *Let R be a commutative ring of characteristic zero. R is a S-zero square ring if and only if R is S-null ring.*

*Proof*: Left for the reader to prove.

It is also important to state, if R is a non-commutative ring, the above result may not in general be true.

***Example 3.10.14***: Let $Z_{24} = \{0, 1, 2, \ldots, 23\}$ be the ring of integers modulo 24. Clearly A = $\{0, 2, 4, 6, 8, \ldots, 20, 22\}$ is a S-subring of $Z_{24}$. $Z_{24}$ is a S-null ring as B = $\{0, 12\}$ is a null ring in A.

**THEOREM 3.10.11**: *The semigroup ring RS is a S-zero square ring if and only if 1 or 2 or 3 is true.*

       *1. R is a S-null ring and S any semigroup.*
       *2. R is any S-zero square commutative ring of characteristic zero and S any commutative semigroup.*
       *3. R any ring and S a S-null semigroup.*

The proof is left as an exercise to the reader as the proof requires vitally only the definitions and a logical use of them.

A recent paper [27] which studies strictly wild algebras with radical square zero may give more innovative ideas when applied to Smarandache concepts.

**DEFINITION [29]**: *A ring R is said to be a p-ring if $x^p = x$ and px = 0 for every x $\in$ R.*

We define here Smarandache p-rings as follows.

**DEFINITION 3.10.7**: *Let R be a ring. R is said to be a Smarandache p-ring (S-p-ring) if R is a S-ring and R has a subring P such that $x^p = x$ and $px = 0$ for every x $\in$ P.*

**THEOREM 3.10.12**: *Let G be the cyclic group of order p-1 and $Z_p$ be the ring of integers modulo p, p a prime. The group ring $Z_pG$ is a S-p-ring.*

*Proof*: $Z_pG$ is obviously a S-ring for we have $Z_p \subset Z_pG$ and $Z_p$ to be such that $x^p = x$ and $px = 0$, so the group ring is a S-p-ring.

**THEOREM 3.10.13**: *Let R be a S-p-ring; R need not be a p-ring.*

*Proof*: Let G be any group and $Z_{12} = \{0, 1, 2, \ldots, 11\}$ be the ring. The group ring $Z_{12}G$ is a S-ring; consider A = $\{0, 4, 8\} \subset \{0, 2, 4, \ldots, 10\} \subset Z_{12}G$. A is a p-ring for p = 3. So $Z_{12}G$ is a S-p-ring which is not a p-ring.

In view of this we have the following.

**THEOREM 3.10.14**: *Let G be any torsion free group. Let R be a S-ring if R is a S-p-ring then the group ring RG is a S-p-ring.*

*Proof*: Since given R is a S-ring which is a S-p-ring we see A $\subset$ R is a subring such that $x^p = x$ and $px = 0$ for all x $\in$ A. Now consider A $\subset$ R $\subset$ RG, so RG is a S-p-ring. Thus we see the group ring RG is not a p-ring but it is a S-p-ring.

**THEOREM 3.10.15**: *Let R be a S-ring which is a S-p-ring. P any semigroup, the semigroup ring RP is a S-p-ring if and only if P has identity.*

*Proof*: Given R is S-p-ring let A $\subset$ R be a S-subring of R. We see in A, $x^p = x$ and $px = 0$ for all x $\in$ A. Now if 1 $\in$ P then, A $\subset$ R. 1 $\subset$ RP; so RP is a S-p-ring. If 1 $\notin$ P then even if R is a S-p-ring. RP in general is not a S-p-ring so RP is a S-p-ring if and only if 1 $\in$ P for any semigroup P.

**DEFINITION [95]**: *A ring R is called an E-ring if $x^{2n} = x$ and $2x = 0$ for every x in R and n a positive integer. The minimal such n is called the degree of the E-ring. It is interesting to note that an E-ring of degree 1 is a Boolean ring.*

Now we proceed on to define Smarandache E-ring.

**DEFINITION 3.10.8**: *Let R be a ring. P a subring of A and A a S-subring of R. if for all x $\in$ P, $x^{2n} = x$ and $2x = 0$ then we say R is a Smarandache E-ring (S-E-ring).*

**Theorem 3.10.16**: *Let R be a S-E-ring then R is a S-ring.*

*Proof*: Follows by the very definition of S-E-ring.

**Theorem 3.10.17**: *Let R be E-ring. If R has S-subring then R is a S-E-ring.*

*Proof*: Obvious by the very definition of S-E-rings and E-rings.

**Example 3.10.15**: Consider the group ring $Z_2S_3$ of the group $S_3$ over the ring $Z_2$. $P = \{0, p_1 + p_2 + p_3, 1 + p_4 + p_5, 1 + p_1 + p_2 + p_3 + p_4 + p_5\}$ be a S-subring of $Z_2S_3$. P is a subring, which is E-ring so $Z_2S_3$ is a S-E-ring. But $Z_2S_3$ is not an E-ring as $(1 + p_1)^2 = 0$ in $Z_2S_3$.

**Theorem 3.10.18**: *Let R be a S-E-ring, then R in general is not an E-ring.*

*Proof*: The above example 3.10.12 is a S-E-ring which is clearly not a E-ring as $(1 + p_1)^2 = 0$ and $1 + p_1 \in Z_2S_3$.

**Definition [46]**: *Let R be ring. R is said to be a pre J-ring if $a^n b = ab^n$ for any pair a, b $\in$ R and n a positive integer.*

To localize this property we now define Smarandache pre J-ring as follows:

**Definition 3.10.9**: *Let R be a ring. P a subring of a S-subring A of R. We say R is a Smarandache pre J-ring (S-pre J-ring) if for every pair a, b $\in$ P we have $a^n b = ab^n$ for some positive integer n.*

**Example 3.10.16**: Let $Z_{12} = \{0, 1, 2, \ldots, 11\}$ be the ring of integers modulo 12. S = $\{0, 2, 4, 6, 8, 10\}$ is a S-subring. But S is a pre J-ring. So $Z_{12}$ is a S-pre J-ring.

**Example 3.10.17**: Let $Z_{12}G$ be the group ring of the group G = $S_3$ over $Z_{12}$. $Z_{12}G$ is a S-pre J-ring.

**Theorem 3.10.19**: *Let R be S-pre J-ring and G any group. The group ring RG is a S-pre J-ring.*

*Proof*: Since R is a S-pre J-ring, we have S $\subset$ R such that S is a S-subring which has a subring to be a pre J-ring. So $S.1 \subset R.1 \subset RG$. Hence for any group G, RG is a S-pre J-ring. It is important to note RG in general is not a pre J-ring.

**Theorem 3.10.20**: *Let R be a S-pre J-ring and P any semigroup with identity. The semigroup ring RP is a S-pre J-ring.*

*Proof*: Obvious from the fact if S ⊂ R and S a S-subring which has a subring P to be a pre J-ring, then S ⊂ R ⊂ RP. So RP is a S-pre J-ring.

**DEFINITION 3.10.10**: *Let R be a ring. R is said to be a Smarandache semi prime ring (S-semi prime ring) if and only if R has no non-zero S-ideal I(II) with square zero.*

We have the following nice theorem about S-ring I.

**THEOREM 3.10.21**: *Let R be a ring, if R has a S-ideal I then R is a S-semi prime ring.*

*Proof*: Follows from the fact that if R has a S-ideal I say A then A has a subset which is a field so $A^2 = (0)$ is impossible. Hence the claim.

In view of this we have the following.

**THEOREM 3.10.22**: *All non-simple S-ring I are S-semi prime.*

*Proof*: Follows from the definition. Left for the reader to prove.

**DEFINITION [42]**: *A commutative ring with 1 is called a Marot ring if each regular ideal of R is generated by a regular element of R. (The author means by a regular element a non-zero divisor and by a regular ideal the elements of the ideal must be non-zero divisors).*

For more about Marot rings please refer [42]. We define Smarandache Marot rings as follows.

**DEFINITION 3.10.11**: *Let R be a ring. If every S-ideal I or (S-ideal II) of R is generated by a regular element and these ideals are regular then we call R Smarandache Marot ring (S-Marot ring).*

**Example 3.10.18**: Z the ring of integers is a S-Marot ring.

**Example 3.10.19**: Let $Z_{10} = \{0, 1, 2, \ldots, 9\}$. $Z_{10}$ is a S-Marot ring. For the only S-ideals of $Z_{10}$ is $\{0, 2, 4, 6, 8\}$ which is regular.

For more about semigroup rings, which are Marot rings please refer [96].

**DEFINITION [102]**: *Let R be a ring. S a proper subring of R. Let I ≠ {0} be a proper subset of S. I is called a subsemi ideal of R, related to the subring S if and only if I is a proper ideal of S and not an ideal of R.*

**DEFINITION [102]**: *The ring R which contains a subsemi ideal is called a subsemi ideal ring.*

**Example 3.10.20**: Let $Z_2 = \{0, 1\}$ be the prime field of characteristic 2. $G = \langle g \, / \, g^4 = 1 \rangle$ The group ring $Z_2G$ is a semi ideal ring. Let $H = \langle 1, g^2 \rangle$ now $I = \langle 0, 1 + g^2 \rangle$ is an ideal of the group ring $Z_2H$ and is not an ideal of $Z_2G$.

**THEOREM [102]**: *Let G be any finite group, having a proper sub group H. Then the group ring KG is a subsemi ideal ring.*

*Proof*: Let $H = \{1, h_1, \ldots, h_n\}$ be the subgroup of G. Then $I = \{0, r \, (1 + h_1 + \ldots + h_n) \, / \, r \in K\}$ is an ideal of KH where KH is a subring of KG and I need not be an ideal of KG. Hence KG is a subsemi ideal ring.

**THEOREM [102]**: *Let G be an infinite group having atleast one element $g \neq e$ of finite order. Then for any field K; the group ring KG is a subsemi ideal ring.*

*Proof*: Given $g \in G$, $g \neq e$ and $g^m = 1$. Let $H = \langle g \rangle$ the cyclic group generated by H. KH is a subring of KG and $I = \langle 0, k(1 + g + \ldots + g^{m-1})/ k \in K \rangle$ is an ideal of KH and not an ideal of KG. Hence the claim.

Now we define Smarandache subsemi ideal and Smarandache subsemi ideal rings.

**DEFINITION 3.10.12**: *Let R be a ring. Let A be a S-subring I or II of R. Let $I \subset A$ be an S-ideal I or II of the S-subring A. Then I is called the Smarandache sub semi ideal I or II (S-sub semi ideal) (I must not be an ideal of R).*

**DEFINITION 3.10.13**: *Let R be a ring if R has a S-subsemi ideal I (II) then we say the ring R is a Smarandache subsemi ideal ring (S-subsemi ideal ring).*

Using the theorem of [102] and the definition of S-subsemi ideal and S-subsemi ideal ring the reader is requested to construct examples and related theorems.

Next we define yet another new notion called Smarandache pre-Boolean ring. A ring R is pre-Boolean ring if xy (x + y) = 0 for every x and y in R.

**DEFINITION 3.10.14**: *Let R be a ring. R is said to be a Smarandache pre Boolean (S-pre Boolean) ring if R is a S-ring and has a subring $A \subset R$ where for all $x, y \in A$ we have xy (x + y) = 0.*

**THEOREM 3.10.23**: *Let R be a S-ring; if $A \subset R$ is a S-subring then A does not satisfy xy (x + y) = 0 for all $x, y \in A$.*

*Proof*: Since if A itself is a S-subring A contains a field so xy (x + y) = 0 may not be possible for all x, y ∈ A unless x + y = 0 for all x, y ∈ A. Hence the claim.

That is why in the definition of S-pre Boolean ring, we demand R to be a S-ring and A a subring not necessarily a S-subring.

**THEOREM 3.10.24**: *Let R be a pre-Boolean ring then R is never a S-pre-Boolean ring*

*Proof*: R is a pre Boolean ring then we have xy (x + y) = 0 for all x, y ∈ R. So R cannot contain a proper subset, A which is a field. For xy (x + y) = 0 forces x + y = 0 or xy is a zero divisor.

**DEFINITION [5]**: *A ring R is called filial if the relation ideal in R is transitive, that is if a subring J is an ideal in a subring I, and I is an ideal in R, then J is an ideal of R.*

We define Smarandache filial ring as follows.

**DEFINITION 3.10.15**: *Let R be a ring. We say R is a Smarandache filial ring (S-filial ring) if the relation S-ideal in R is transitive, that is if a S-subring, J is an S-ideal in a S-subring I and I is a S-ideal of R, then J is an S-ideal of R.*

***Example 3.10.21***: Let R = $Z_2 \times Z_2 \times Z_2$ be a ring. J = $\langle (0, 0, 0)(0, 0, 1) \rangle$ is an ideal in I = $\langle (0, 0, 0), (0, 0, 1), (0, 1, 0), (0, 1, 1) \rangle$ but I is an ideal of R and we see J is an ideal of R. Hence R is a filial ring.

**DEFINITION [100]**: *A ring R is called an n-ideal ring if for every set of n-distinct ideals $I_1, I_2, \ldots, I_n$ of R and for every set of n-distinct elements $x_1, x_2, \ldots, x_n \in R \setminus (I_1 \cup I_2 \cup \ldots \cup I_n)$ we have $\langle x_1 \cup I_1 \cup I_2 \cup \ldots \cup I_n \rangle = \langle x_2 \cup I_1 \cup I_2 \cup \ldots \cup I_n \rangle = \ldots = \langle x_n \cup I_1 \cup I_2 \cup \ldots \cup I_n \rangle$, where $\langle \rangle$ denotes the ideal generated by $x_i \cup I_1 \cup I_2 \cup \ldots \cup I_n$, $1 \leq i \leq n$. (By an ideal we mean only a two sided ideal).*

Now we proceed on to define Smarandache n-ideal rings.

**DEFINITION 3.10.16**: *Let R be a ring. We say R is a Smarandache n-ideal ring (S-n-ideal ring) if for every set of n-distinct S-ideal I (II), $I_1, \ldots, I_n$ of R and for every distinct set of n elements $x_1, x_2, \ldots, x_n \in R \setminus (I_1 \cup I_2 \cup \ldots \cup I_n)$ we have $\langle x_1 \cup I_1 \cup I_2 \cup \ldots \cup I_n \rangle = \langle x_2 \cup I_1 \cup I_2 \cup \ldots \cup I_n \rangle = \ldots = \langle x_n \cup I_1 \cup \ldots \cup I_n \rangle$ denotes the S-ideal generated by $x_1 \cup I_1 \cup I_2 \cup \ldots \cup I_n$; $1 \leq i \leq n$.*

***Example 3.10.22***: Let $Z_{12}$ = {0, 1, 2, …, 11} be the ring of integers modulo 12. $Z_{12}$ is a 3-ideal ring and a 4-ideal ring. $Z_{12}$ is not a S-n-ideal ring for $Z_{12}$ has only one S-ideal.

***Example 3.10.23***: Let $Z_{15}$ = {0, 1, 2, …, 14}. The ideals of $Z_{15}$ are {0, 5, 10} and {0, 3, 6, 9, 12}; clearly $Z_{15}$ is not a S-2 ideal ring.

**DEFINITION [21]**: *A non-empty set S of a ring R is called a generalized left semi ideal of R if S is closed under addition and $x^2 s$ is in S for any $s \in S$ and $x \in R$. Similarly one can define generalized right semi ideal and the generalized semi-ideal when it is both a generalized left and right semi ideal.*

Now we define the concept of Smarandache right (left) generalized semi-ideal and Smarandache generalized semi-ideal.

**DEFINITION 3.10.17**: *Let R be a S-ring I. A generalized right (left) semi ideal I of the S-ring R are called Smarandache generalized right (left) semi ideals (S-generalized right (left) ideals). The Smarandache generalized semi-ideal (S-generalized semi-ideal) is one which is both a S-generalized left and right semi-ideal.*

*If R is a S-ring and R has generalized semi-ideal I, then I is called the Smarandache generalized semi-ideal (S-generalized semi-ideal).*

**THEOREM 3.10.25**: *Let R be a ring; if R is a S-ring having generalized semi-ideals then R has S-generalized semi-ideals.*

*Proof*: By the very definition the result follows.

**THEOREM 3.10.26**: *All rings which are generalized semi-ideal rings need not in general be S-generalized semi-ideal rings.*

*Proof*: By example. $Z_4$ = {0, 1, 2, 3} is not a S-ring. Here $Z_4$ has a generalized semi-ideal but $Z_4$ is not a S-generalized semi-ideal ring.

***Example 3.10.24***: Let $Z_2$ = {0, 1} and G = $\langle g / g^4 = 1 \rangle$ be the group. The group ring $Z_2 G$ is a S-ring. I = {0, 1 + $g^2$} is not a S-ideal but is a S-generalized semi-ideal of $Z_2 G$.

**THEOREM 3.10.27**: *Let K be a non-prime real field of characteristic zero. K has S-generalized semi-ideals.*

*Proof*: K is a non-prime field so K has subfield. Hence K is a S-ring. Take S = $K^+ \cup$ {0} only positive elements. S is closed with respect to addition. For any s ∈ S and x ∈ K; $x^2s$ ∈ S. Hence the claim.

If K is a complex field the result may not be true.

**DEFINITION [25]**: *A ring A is s-weakly regular if for each a ∈ A, a ∈ aAa²A.*

**Example 3.10.25**: Let G = $\langle g \,/\, g^2 = 1 \rangle$ and $Z_2$ = {0, 1}. The group ring $Z_2G$ is not s-weakly regular ring.

**DEFINITION 3.10.18**: *Let R be a ring. A be a S-subring of R. We say R is Smarandache s-weakly regular (S-s-weakly regular) ring if for each a ∈ A. a ∈ aAa²A.*

**DEFINITION [77]**: *Let R be a ring. A right ideal I of R is said to be quasi reflexive if whenever A and B are two right ideals of R with AB ⊂ I then BA ⊂ I.*

A ring R is said to be right quasi reflexive if (0) is a right quasi reflexive ideal of R. Similarly one defines the concept of left quasi reflexive ring. Semi prime rings are left and right quasi reflexive.

One knows the group ring KG is left and right quasi reflexive where K is a field of characteristic 0. The result follows from the fact the group rings KG is semi prime. For more about these results please refer [61, 62]. We just recall: a ring R is semi prime if and only if R contains no non-zero ideal with square zero. We define Smarandache semi prime rings.

**DEFINITION [77]**: *A ring R to be strongly sub commutative if every right ideal of it is right quasi reflexive. (A right ideal I of a ring R is called right quasi reflexive if whenever A and B are right ideals of R with AB ⊂ I then BA ⊂ I).*

We define Smarandache strongly sub commutative rings.

**DEFINITION 3.10.19**: *Let R be a ring. R is said to be a Smarandache strongly sub commutative (S-strongly sub commutative) if every S-right ideal I(II) of it is right quasi reflexive.*

**DEFINITION [7]**: *A commutative ring R is said to be a Chinese ring if given elements a, b ∈ R and ideals I, J ⊂ R such that a ≡ b (I + J) there exists c ∈ R such that c ≡ a (I) and c ≡ b (J), c ≡ a (I) implies $\langle I, a \rangle \equiv \langle I, c \rangle$ i.e., generated by I + a and I + c, for more about Chinese rings refer Aubert.*

**DEFINITION 3.10.20**: *Let R be a ring. R is said to be a Smarandache Chinese ring (S-chinese rings) I(II) if given elements a ,b ∈ R and S-ideal I(II) in R such that ⟨I ∪ J ∪ a⟩ = ⟨I ∪ J ∪ b⟩ there exist an element c ∈ R such that ⟨I ∪ a⟩ = ⟨I ∪ c⟩ and ⟨J ∪ b⟩ = ⟨J ∪ c⟩.*

**Example 3.10.26**: Let $Z_2 = \{0, 1\}$ be the prime field of characteristic two. S = {a, b, 0 / $a^2$ = a, $b^2$ = b, ab = ba = 0}. Clearly $Z_2S$, the semigroup ring is a S-Chinese ring I.

The author has defined group rings, which is a direct sum of subrings.

**DEFINITION [105]**: *A group ring RG is s-decomposable if RG = $S_1$ + … + $S_r$ where $S_i$'s are subrings of RG with $S_i \cap S_j$ = R and every element in RG has a unique representation as a sum.*

**Example 3.10.27**: Let $Z_2S_3$ be the group ring of the group $S_3$ over the field $Z_2 = \{0,1\}$. Let

$$H_1 = \left\{ p_0 = \begin{pmatrix} 1 & 2 & 3 \\ 1 & 2 & 3 \end{pmatrix} , \ p_1 \right\},$$

$$H_2 = \{p_0, p_2\}, H_3 = \{p_0, p_3\} \text{ and } H_4 = \{p_0, p_5, p_4\}$$

be the subgroups of $S_3$. Then $Z_2S_3 = Z_2H_1 + Z_2H_2 + Z_2H_3 + Z_2H_4$ as a direct sum of subrings. $Z_2H_i \cap Z_2H_j = Z_2$, i ≠ j for 1 ≤ i, j ≤ 4.

Now the author defines strongly s-decompsable group rings.

**DEFINITION 3.10.21**: *Let RG the group ring of the group G over the ring R. RG is strongly s-decomposable if RG = $S_1$ + … + $S_r$ with $S_i \cap S_j$ = {0, 1} or {0} and every element of RG has a unique representation as a sum of elements from $S_1$, $S_2$, …, $S_r$.*

Now we define Smarandache s-decompossible and Smarandache strongly s-decomposable ring as follows:

**DEFINITION 3.10.22**: *Let RG be the group ring of the group G over the ring R. We say RG is Smarandache s-decomposable (S-s-decomposable) if RG = $S_1$ + … + $S_r$ where $S_i$ are subrings such that atleast one of the $S_i$ is a S-subring of RG with $S_i \cap S_j$ = R and every element of RG has a unique representation as a sum of elements from $S_1$, $S_2$, …, $S_r$.*

**THEOREM 3.10.28**: *Let RG be the group ring such that R is any field and if RG is s-decomposable then RG is S-s-decomposable.*

*Proof*: Since R is a field we see every subring $S_i$ contains R as a subset so every subring $S_i$ is a S-subring of RG so if RG is s-decomposable then it is a S-s-decomposable.

**DEFINITION 3.10.23**: *The group ring RG is Smarandache strongly s-decomposable (S-strongly s-decomposable) if RG is strongly s-decomposable in which at least one of the $S_i$'s is a S-subring of RG.*

**THEOREM 3.10.29**: *If RG is S-strongly s-decomposable then RG is strongly s-decomposable.*

*Proof*: Obvious by the very definitions of S-strongly s-decomposable and strongly s-decomposable.

**DEFINITION [105]**: *The group ring RG is weakly s-decomposable if we can find subrings $S_1, …, S_r$ of RG with $RG = S_1 + … + S_r$ such that $S_i \cap S_j = G$; $i \neq j$.*

**DEFINITION 3.10.24**: *The group ring RG is Smarandache weakly s-decomposable (S-weakly decomposable) if we can find subrings $S_1, S_2, …, S_r$ of RG of which atleast one of the $S_i$'s is a S-subring of RG with $RG = S_1 + S_2 + … + S_r$ and $S_i \cap S_j = G$; if $i \neq j$.*

The reader is requested to develop relations between weakly s-decomposable and S-weakly s-decomposable group rings. They are also advised to formulate definitions in case of semigroup rings and study them.

**DEFINITION [108]**: *Let R be a ring not necessarily commutative Let L denote the collection of all right ideals of R.*

*If $(A + B) (A + C) (A + D) = A + BC (A + D) + BD (A + C) + DC (A + B)^*$ for all right ideals A, B, C, D in L; where A+B denotes the right ideal generated by $A \cup B$ and AB denotes $A \cap B$. Then we call R a strong right s-ring.*

*\* - This identity is known as the supermodular identity and such lattices are known as supermodular lattices [36].*

The motivation for doing so is from [23], for he called the lattice of ideals, which is distributive to be the strong right D-domain.

We define Smarandache strong right s-ring.

**DEFINITION 3.10.25**: *Let R be a ring not necessarily commutative. Let L denote the collection all right ideals of R if (A + B) (A + C) (A + D) = A + BC (A + D) + BD (A + C) + DC (A + B) for all right ideals A, B, C, D ∈ L of which one of the ideals must be a S-right ideal I(II), then we say R is a Smarandache strong right s-ring (S-strong right s-ring).*

***Example 3.10.28***: Let $Z_2$ = {0, 1} be the ring of integers modulo 2 and S = {1, a, b / $a^2$ = a, $b^2$ = b, ab = a, ba = b, 1.a = 1.a = a and 1.b = b.1 = b} be the multiplicative semigroup. The semigroup ring $Z_2S$ is a S-strong right s-ring, for the ideals of $Z_2S$ are $A_1$ = {0, a}, $A_2$ = {0, b}, $A_3$ = {0, a + b}, $A_4$ = {1, 1 + a, 1 + b, a + b} and $A_5$ = {0, a, b, a + b} of these 5 right ideals any four of them will be S-ideal I, so $Z_2S$ is a S-strong right s-ring.

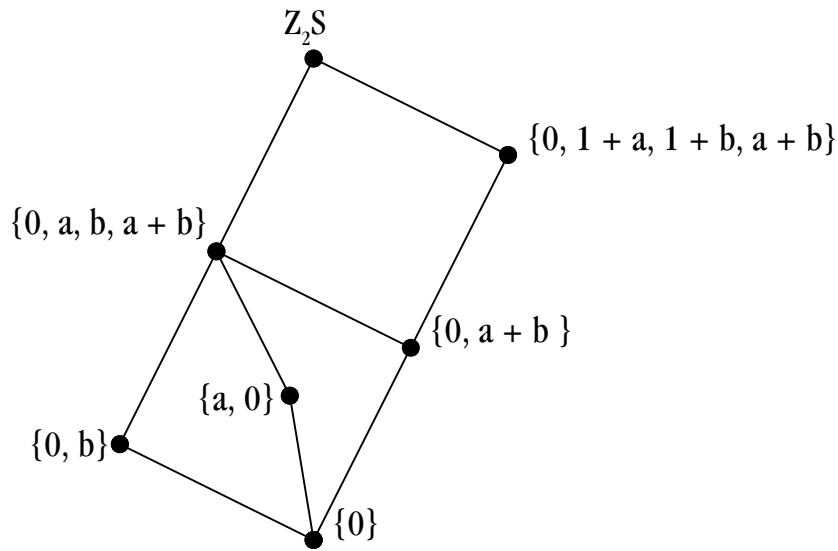

Figure 3.10.1

**DEFINITION 3.10.26**: *Let R be a ring if all the right ideals form a distributive lattice and if in this collection of right ideals we have atleast a right ideal to be a S-ideal I(II) then we call the ring R as the Smarandache strong right D-domain (S-strong right D-domain).*

**THEOREM 3.10.30**: *Every S-strong right-D-ring is a S-strong right s-ring.*

*Proof*: By the very definition; the collection of ideals have S-ideals I(II) so if we put D = C in the supermodular identity we get the result.

**THEOREM 3.10.31**: *If the set of right ideals in a S-strong right s-ring is not a S-strong right D-ring it does not imply the set of two sided ideals of this ring is not a D-ring.*

*Proof*: By an example. Clearly the ideals do not form a S-D-ring. Now consider the set of two sided ideals of $Z_2S$ given in example 3.10.28.

$$B_3 = \{0, 1 + a, 1 + b, a + b\}$$
$$B_2 = \{0, a, b, a + b\},$$
$$B_1 = \{0, a + b\}$$

form a distributive lattice so the ring is a S-strong-D-ring.

**THEOREM 3.10.32**: *If the set of right ideals of a ring R is not a S-strong r-ring still it does not imply the set of two sided ideals is not a S-strong-r-s-ring.*

*Proof*: By an example consider $Z_2 = \{0, 1\}$ and $S = \{a, b, c, 1 / a^2 = a, b^2 = b, c^2 = c, ab = a, ba = b, ca = c, cb = c, ac = a, bc = b, 1.a = a.1 = a, 1.c = c.1 = c, 1.b = b.1 = b\}$ be the multiplicative semigroup. $Z_2S$ be the semigroup ring of S over $Z_2$

Take $A = \{0, 1 + a + b + c, b + c, a + c, a + b, 1 + a, 1 + c, 1 + b\}$, $B = \{0, a\}$, $C = \{0, b\}$, $D = \{0, a + b + c\}$ be the right ideals of $Z_2S$. Clearly (A + B) (A + C) (A + D) = $Z_2S$, A + BC (A + D) + CD (A + B) + DB (A + C) = A. Since $Z_2S \neq AZ_2S$ is not a S-strong-r-s-ring.

Consider the two-sided ideals of $Z_2S$.

$A = \{0, a + b + c\}$, $B = \{0, a + b, a + c, b + c\}$, $C = \{1 + a + b + c, a + b, a + c, b + c, 1 + a, 1 + b, 1 + c, 0\}$ and $D = \{a, b, c, a + b, a + c, b + c, a + b + c, 0\}$. Clearly the set A, B, C, D, $\{0\}$, $Z_2S$ form a supermodular lattice given by the following diagram.

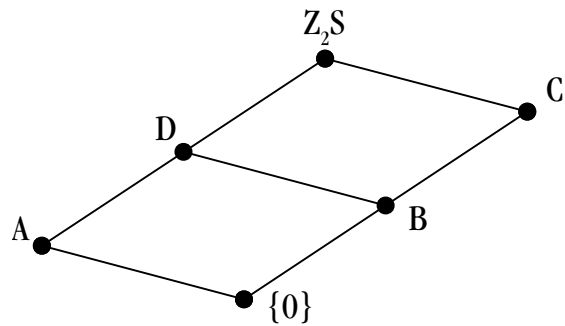

Figure 3.10.2

N.Jacobson calls a ring R to be a J-ring if $x^n = x$ for every $x \in R$, n an integer, n > 1, we motivated by this define Smarandache J-ring as follows:

**DEFINITION 3.10.27**: *Let R be a ring. R is a said to be a Smarandache J-ring (S-J-ring) if R has a S-subring A such that for all a $\in$ A we have $a^n$ = a, n > 1 (n an integer).*

In view this we have the following.

**THEOREM 3.10.33**: *Let R be a J-ring. If R has S-subring then R is a S-J-ring.*

*Proof*: Obvious by the very definition of J-ring and S-J-ring. It is well know that all J-rings are commutative but we see a S-J-ring need not be commutative.

***Example 3.10.29***: Let $Z_2S_3$ be the group ring of the group $S_3$ over $Z_2$. $Z_2S_3$ is a S-J-ring. For $Z_2S_3$ contains a S-subring A = {0, $p_1 + p_2 + p_3$, 1 + $p_4 + p_5$, 1 + $p_1 + p_2 + p_3$ + $p_4 + p_5$}. It is easily verified A is J-ring so $Z_2S_3$ is a S-J-ring. Hence the claim.

**THEOREM 3.10.34**: *Let $Z_n$ be a S-ring. S is a semigroup such that $s_i s_j$ = 0 if i $\neq$ j $s_i s_i$ = $s_i$. Then the semigroup ring $Z_n S$ is a S-J-ring.*

*Proof*: It is left for the reader to verify.

The author has defined for any ring R the strong ideal property as follows.

**DEFINITION [122]**: *Let R be a ring. If every distinct pair of ideals of R generate R, then the set of ideals of R is said to satisfy the strong ideal property.*

**DEFINITION [122]**: *Let R be a ring; if every distinct pair of subrings of R generate R then the set of subrings of R is said to satisfy the strong subring property.*

**DEFINITION [122]**: *Let R be a ring. If {$I_m$} denote the collection of all ideals and {$S_n$} the collection of all subrings and if $\langle S_j I_j / I_j \in \{I_m\}$ and $S_j \in \{S_n\}\rangle$ generate R for every pair ($S_j I_j$) $\in$ {$S_n$} $\times$ {$I_m$} then we say the subrings and ideals of R satisfy the strong subring ideal property.*

***Example 3.10.30***: Let G = $\langle g / g^2 = 1\rangle$ and $Z_2$ = {0, 1} be the ring of integers modulo 2. The group ring $Z_2G$ satisfies strong subring property but does not satisfy strong ideal property. But $Z_2G$ satisfies strong subring ideal property as $S_1$ = {0, 1} and $I_1$ = {0, 1 + g} be the subring and ideal of $Z_2S$. We see $S_1 \cup I_1$ generates $Z_2G$.

***Example 3.10.31***: Let G = $\langle g/ g^3 = 1\rangle$ be the cyclic group of order three and $Z_2$ = {0, 1} be the ring of integers modulo 2. The group ring $Z_2G$ = {1, 0, g, $g^2$, 1 + g, 1 + $g^2$, g + $g^2$, 1 + g + $g^2$}. The subrings of $Z_2G$ are $S_1$ = {0, 1}, $S_2$ = {0, 1 + g + $g^2$}, $S_3$ =

$\{0, g + g^2\}$, $S_4 = \{0, 1 + g, 1 + g^2, g + g^2\}$ and $S_5 = \{0, 1, g + g^2, g + g^2 + 1\}$. The ideals of $Z_2G$ are $I_1 = \{0, 1 + g + g^2\}$, $I_2 = \{0, 1 + g, 1 + g^2, g + g^2\}$; thus clearly $Z_2G$ satisfies strong ideal property. It is easily verified that $Z_2G$ does not satisfy strong subring property further $Z_2G$ does not satisfy the strong subring ideal property.

**THEOREM [122]**: *Let R be a ring. R does not satisfy strong subring property even if a pair of subrings $S_i$, $S_j \in \{S_n\}$ is such that $S_i \subset S_j$ or $S_j \subset S_i$.*

*Proof*: Since $S_i$, $S_j \in \{S_n\}$ if $S_i \subset S_j$ or $S_j \subset S_i$ then $\langle S_i, S_j \rangle = S_j$ if $S_i \subset S_j$ and $\langle S_i, S_j \rangle = S_i$ if $S_j \subset S_i$.

**THEOREM [122]**: *Let R be a ring. R is not a strong ideal ring if there exists a pair of ideals $I_1$, $I_2$ such that $I_1 \subset I_2$ or $I_2 \subset I_1$.*

*Proof*: As in case of subrings.

Now we define S-strong ideal rings, S-strong subring rings and S-strong subring ideal rings.

**DEFINITION 3.10.28**: *Let R be a ring. Let $\{S_i\}$ denote the collection of all S-subrings of R. We say R is a Smarandache strong subring ring (S-strong subring ring) if every pair of S-subrings of R generate R.*

**DEFINITION 3.10.29**: *Let R be a ring. Let $\{I_i\}$ denote the collection of all S-ideals of R. We say R is a Smarandache strong ideal ring (S-strong ideal ring) if every pair of S-ideals of R generate R.*

**DEFINITION 3.10.30**: *Let R be a ring. $\{S_i\}$ and $\{I_j\}$ denote the collection of all S-subrings and S-ideals of R. if every pair $\{S_i, I_j\}$ generate R then we say R is a Smarandache strong subring ideal (S-strong subring ideal) ring. If we do not have S-subrings and S-ideals in a ring then we do not have the concept of S-strong subring ideal or S-strong ideal ring or S-strong subring ring.*

**Example 3.10.32**: Let $Z_6 = \{0, 1, 2, 3, 4, 5\}$ be the ring of integers modulo 6. $S_1 = \{0, 3\}$ and $S_2 = \{0, 2, 4\}$ are subrings as well as ideals. Clearly $Z_6$ is not a S-strong ideal ring and not a S-strong subring ring and not a S-strong subring ideal ring, since this ring has no proper S-subrings or S-ideals.

**Example 3.10.33**: Let $Z_2S_3$ be the group ring of the group $S_3$ over $Z_2$. Clearly $Z_2S_3$ is not a strong subring. For take the two distinct subrings, $S_1 = \{0, 1 + p_1\}$ and $S_2 = \{0, 1 + p_2\}$. $\langle S_1, S_2 \rangle = \{0, 1 + p_1, 1 + p_2, p_1 + p_2, p_4 + p_5, p_2 + 1 + p_1 + p_5, 1 + p_1 + p_2 +$ $p_4$, $p_1 + p_4 + p_5 + p_2$, ...} $\neq Z_2S_3$ as in $\langle S_1, S_2 \rangle$ i.e., the ring generated by $S_1$ and $S_2$ we cannot find elements whose support is odd i.e., single term, sum of three terms or sum of five terms. Hence $Z_2S_3$ is not a strong subring ring. Also $Z_2S_3$ is not a strong ideal ring. For take $I_1 = \{0, 1 + p_1 + p_2 + \ldots + p_5\}$ and $I_2 = \{0, 1 + p_4 + p_5, p_1 + p_2 + p_3, 1 + p_1 + p_2 + p_3 + p_4 + p_5\}$, $\langle I_1 \cup I_2 \rangle$ does not generate $Z_2S_3$. Hence the claim.

Is $Z_2S_3$ a S-strong ideal ring?

$Z_2S_3 = \langle Z_2H_1 \cup Z_2H_3 \rangle$ where $H_1 = \langle 1, p_1 \rangle$ and $H_2 = \langle 1, p_4, p_5 \rangle$ are subgroups of $Z_2S_3$. Clearly both $Z_2H_1$ and $Z_2H_3$ are S-subrings of $Z_2S_3$. The S-subrings of $Z_2S_3$ are A = $\{0, 1 + p_1 + p_2 + p_3 + p_4 + p_5, p_4 + p_5 + 1, p_1 + p_2 + p_3\}$ for $\{0, 1 + p_1 + p_2 + p_3 + p_4 + p_5\}$ acts as the proper subset which is a subfield of A $\subset Z_2S_3$. Also B = $\{0, 1, p_1, 1 + p_1\}$ is a S-subring of $Z_2S_3$. Similarly we have S-subrings $B_1 = \{0, 1, p_2, 1 + p_2\}$, $B_3 = \{0, 1, p_3, 1 + p_3\}$, it is once again easily verified $\langle B_1 \cup B_3 \rangle = Z_2S_3$.

Now the natural question is, will every pair of S-subrings generate $Z_2S_3$. To this end we propose some open problems in chapter 5 and define a weaker Smarandache concept.

**DEFINITION 3.10.31**: *Let R be a ring; we say R is a Smarandache weak ideal (S-weak ideal) ring if there exists a pair of distinct S-ideals $I_1$, $I_2$ in R which generate R i.e., R = $\langle I_1 \cup I_2 \rangle$.*

**DEFINITION 3.10.32**: *Let R be a ring; we say R is a Smarandache weak subring (S-weak subring) ring if there exists a distinct pair of S-subrings $S_1$, $S_2$ in R which generate R. i.e., R = $\langle S_1 \cup S_2 \rangle$.*

**DEFINITION 3.10.33**: *Let R be a ring we say R is a Smarandache weak subring ideal (S-weak subring ideal) ring if there exist an S-ideal I and S-subring A (which is not an S-ideal) such that I $\cup$A generate R i.e., R=$\langle I \cup A \rangle$.*

The following three results can be easily verified; hence the proof is left as an exercise to the reader.

**THEOREM 3.10.35**: *Let R be a ring.*

1. *Every S-strong ideal ring is a S-weak ideal ring.*
2. *Every S-strong subring ring is a S-weak subring ring*
3. *Every S-strong subring ideal ring is a S-weak subring ideal ring.*

Further we leave it as an exercise to the reader to obtain examples to show S-weak structures in general are not S-strong structures.

***Example 3.10.34***: Let Q be the field of rationals. $S_3$ be the group of degree 3. $QS_3$ is a S-weak subring ring. For $A_1 = QH_1$ and $A_2 = QH_2$ are S-subrings and $QS_3 = \langle QH_1 \cup QH_2 \rangle$.

**DEFINITION [123]**: *Let R be a ring we say R is a weakly Boolean ring if $x^{n(\alpha)} = x$ for all $x \in R$ and, for some natural number $n(\alpha) > 1$.*

***Example 3.10.35***: $Z_p = \{0, 1, \ldots, p{-}1\}$ be the prime field of characteristic p. Clearly $Z_p$ is a weakly Boolean ring.

Now we proceed on to define Smarandache weakly Boolean ring.

**DEFINITION 3.10.34**: *Let R be a ring we say R is a Smarandache weakly Boolean ring (S-weakly Boolean ring) if we have a S-subring A of R such that A is a weakly Boolean ring.*

***Example 3.10.36***: Let $Z_{15} = \{0, 1, 2, \ldots, 14\}$ where $G = \langle g / g^2 = 1 \rangle$. Clearly the group ring $Z_{15}G$ is not a weakly Boolean ring. But $Z_{15}G$ is a S-weakly Boolean ring. For take $B = \{0, 5, 10\}$; $BG$ is a S-subring of $Z_{15}G$ which is a S-weakly Boolean ring but $Z_{15}G$ is not a weakly Boolean ring.

**DEFINITION [67]**: *R is a weakly regular ring if for each right (left) ideal I of R; we have $I^2 = I$.*

***Example 3.10.37***: Let $Z_2 = \{0, 1\}$ be the prime field of characteristic two and $G = \langle g / g^3 = 1 \rangle$ be the cyclic group of order 3. The group ring $Z_2G$ is weakly regular. For $I_1 = \{0, 1 + g + g^2\}$, $I_2 = \{0, 1 + g, 1 + g^2, g + g^2\}$ are such that $I_1^2 = I_1$ and $I_2^2 = I_2$. Hence the claim.

Now we see all ideals I in every ring need not satisfy $I^2 = I$.

***Example 3.10.38***: $Z_2 = \{0, 1\}$ be the prime field of characteristic 2 and $G = \langle g / g^2 = 1 \rangle$, the group ring; $Z_2G = \{0, 1, g, 1 + g\}$. The ideal $I = \{0, 1 + g\}$ is such that $I^2 = \{0\}$. So $Z_2G$ is not a weakly regular ring.

**DEFINITION 3.10.35**: *Let R be a ring, we say R is a Smarandache weakly regular ring (S-weakly regular ring) if each S- right (left) ideal I of R satisfies $I^2 = I$.*

**Definition 3.10.36**: *Let R be a ring. If R has atleast one S-ideal I such that $I^2$ = I then we say R is a Smarandache weakly weak regular ring (S-weakly weak regular ring).*

It is easily verified that :

**Theorem 3.10.36**: *Let R be a S-weakly regular ring then R is a S-weakly weak regular ring.*

**Definition [1]**: *A ring R of characteristic p is said to be a pre-p-ring if $x^p y = xy^p$ for every x, y ∈ R.*

**Definition 3.10.37**: *Let R be a ring of characteristic p; R is said to be a Smarandache pre-p-ring (S-pre-p ring) if R has a nontrivial S-subring A such that a subring B of A is a pre-p-ring i.e., in B we have $x^p y = xy^p$ for all x, y ∈ B.*

***Example 3.10.39***: Let $Z_2$ = {0, 1} be the prime field of characteristic 2 and G = ⟨g /$g^6$ = 1⟩. The group ring $Z_2$G is a S-pre-p-ring for take A = {0, 1, $2g^3$, 1 + $g^3$, $2g^3$, 2 + $g^3$, $2g^3$ + 1, $2g^3$ + 2}. So $Z_2$G is a S-pre-p-ring.

**Definition [22]**: *Let R be a commutative ring, an ideal I of a ring R is said to be the multiplication ideal if for each ideal J ⊂ I we have J = IC for some ideal C.*

**Definition 3.10.38**: *Let R be a S-commutative ring. An S-ideal I of R is said to be the Smarandache multiplication ideal (S-multiplication ideal) if for each S-ideal J ⊂ I there is J = IC for some S-ideal C in R.*

**Definition [134]**: *A two sided ideal I of a non-commutative ring R is called a right multiplication ideal, if for each right ideal J ⊆ I there is J = IC for some right ideal C in R.*

**Definition [134]**: *Let R be a ring. A right ideal I of a non-commutative ring is called a right multiplication right ideal if for each right ideal J ⊂ I there is J = IC for some right ideal C in R.*

**Definition [134]**: *Let R be a ring. If every proper two sided ideal of R is a right multiplication ideal of R then we call R as a right multiplication ideal ring.*

**Definition 3.10.39**: *Let R be a ring. If every proper two sided S-ideals of R is a right multiplication ideal of R. We call the ring as a Smarandache right multiplication ideal ring (S-right multiplication ideal ring)*

**DEFINITION [3]**: *Let R be a partially ordered ring without non-zero nilpotents. R is said to be a f-ring if and only if for any a ∈ R, there exists $a_1$, $a_2$ in R with $a_1$ > 0, $a_2$ > 0, a = $a_1$ − $a_2$ and $a_1 a_2 = a_2 a_1 = 0$.*

We define Smarandache f-rings as follows.

**DEFINITION 3.10.40**: *Let R be a ring. Let A be a S-subring of R. R is said to be a Smarandache f-ring (S-f-ring) if and only if A is a partially ordered ring without non-zero nilpotents and for any a ∈ A we have $a_1$, $a_2$ in R with $a_1 \geq 0$, $a_2 \geq 0$, a = $a_1$ − $a_2$ and $a_1 a_2 = a_2 a_1 = 0$.*

**DEFINITION [20]**: *Let R be a ring. R is said to be a chain ring if the set of ideals of R is totally ordered by inclusion.*

For more properties about chain rings please refer [20].

**THEOREM [121]**: *Let $Z_2$ = {0, 1} be the prime field of characteristic 2. G = ⟨g / $g^p$ = 1⟩ and p an odd prime. The group ring $Z_2 G$ is not a chain ring.*

*Proof*: Consider I = {0, 1 + g + … + $g^{p-1}$} and J = the augmentation ideal of $Z_2 G$. Clearly I and J are ideals such that they are not comparable so $Z_2 G$ is not a chain ring.

**THEOREM [121]**: *Let G = ⟨g / $g^{p+1}$ = 1⟩ be a cyclic group of order p + 1, $Z_p$ be the prime field of characteristic p. The group ring $Z_p G$ is not a chain ring.*

*Proof*: Consider the ideals I = {0, n(1 + g + $g^2$ + … + $g^p$)}, 1 ≤ n ≤ p-1 and J the augmentation ideal. Clearly I and J are not comparable; so $Z_p G$ is not a chain ring.

**THEOREM [103]**: *Let $Z_p$ be the prime field and G be a finite group of order n, such that (n, p) = 1. Then the group ring $Z_p G$ is not a chain ring.*

*Proof*: Take J = {0, t(1 + g + … + $g_{n-1}$)}, 1 ≤ t ≤ p −1 and J the augmentation ideal. I and J are incomparable so $Z_p G$ is not a chain ring.

**DEFINITION [103]**: *Let R be a ring. R is said to be a strictly right chain ring only when the right ideals of R is ordered by inclusion.*

**DEFINITION 3.10.41**: *Let R be a ring. If the set of S-ideals of R is totally ordered by inclusion, then we say R is a Smarandache chain ring (S-chain ring).*

**DEFINITION 3.10.42**: *Let R be a ring. If the set of all S-right ideal of R is totally ordered by inclusion then we say R is a Smarandache right chain ring (S-right chain ring).*

**DEFINITION 3.10.43**: *Let R be a ring, A, a S-subring of R. If the set of S-ideals of the S-subring A of R is totally ordered by inclusion then the ring is said to be Smarandache weakly chain ring (S-weakly chain ring).*

**THEOREM 3.10.37**: *Let R be a ring, if R is a S-weakly chain ring then R need not be a S-chain ring.*

*Proof*: By an example the above result can be proved.

**DEFINITION [127]**: *Let R be a ring, $0 \neq I$ be an ideal of R. If for any nontrivial ideal X and Y of R $X \neq Y$ we have $\langle X \cap I, Y \cap I \rangle = \langle X, Y \rangle \cap I$ then I is called the obedient ideal of R.*

***Example 3.10.40***: Let $Z_{12} = \{0, 1, 2, \ldots, 11\}$, $I = \{0, 6\}$ is an obedient ideal of $Z_{12}$ for if we take $X = \{0, 4, 8\}$ and $Y = \{3, 6, 9, 0\}$ two ideals of $Z_{12}$. We see $\langle X \cap I, Y \cap I \rangle = \langle 0, 0, 6 \rangle \{0, 6\} = I$, $\langle X, Y \rangle \cap I = Z_{12} \cap I = \{0, 6\}$. Hence the claim.

**DEFINITION [127]**: *If every ideal I of a ring R is an obedient ideal of R, then we say R is an ideally obedient ring.*

**DEFINITION 3.10.44**: *Let R be a ring. Let I be an S-ideal of R. We say I is a Smarandache obedient ideal (S-obedient ideal) of R if we have two ideals X, Y in R, $X \neq Y$ such that $\langle X \cap I, I \cap Y \rangle = \langle X_I, Y \rangle \cap I$.*

<u>Note</u>: We do not demand that X and Y to be S-ideals of R. It is sufficient if they are distinct ideals but we demand I to be a S-ideal of R.

**DEFINITION 3.10.45**: *Let R be a ring, if every S-ideal I of R is a S-obedient ideal of R then we say R is a Smarandache ideally obedient ring. (S-ideally obedient ring)*

**DEFINITION [49]**: *Let R be an associative ring in which for every x, y in R there exists a positive integer $n = n(x, y) > 1$ such that either $(xy - yx)^n = xy - yx$ or $(xy + yx)^n = xy + yx$.*

*In honour of Lin we call these rings as Lin rings. For more about these structures please refer [49].*

**THEOREM [129]**: *Let F be a field and G any finite non-abelian group. If the group ring FG is a Lin ring then*

1. *It has zero divisors*
2. *FG is a Lin ring having elements of finite order.*

*Proof*: Given FG is a ring which is a Lin ring. Hence $(xy - yx)^n = xy - yx$ or $(xy + yx)^n = xy + yx$. Now $(xy - yx)^n = xy - yx$ implies $(xy - yx) [(xy - yx)^{n-1} - 1] = 0$. Since $xy \neq yx$; we have if FG is a Lin ring, FG has zero divisors provided $(xy - yx)^{n-1} \neq 1$. Clearly if the ring FG has no zero divisors then FG has elements of finite order i.e., $(xy - yx)^{n-1} = 1$ for atleast some pairs $x, y \in$ FG.

**THEOREM [129]**: *If FG is the group ring of a non-commutative group G over the field F. FG has zero divisors or elements of finite order then FG is not in general a Lin Ring.*

*Proof*: By an example. Consider $Z_2 = \{0, 1\}$ be the prime field of characteristic two and $S_3$ be the symmetric group of degree three.

To show $Z_2S_3$ is not a Lin ring, it is sufficient to prove that there exists atleast a pair of elements $x, y$ in $Z_2S_3$ such that $(xy + yx)^n \neq xy + yx$. Consider $p_2, p_4 \in Z_2S_3$. Clearly $p_2p_4 + p_4p_2 = p_1 + p_3$ and $(p_1 + p_3)^2 = p_4 + p_5$. Since $(p_4 + p_5)^2 = p_4 + p_5$. We have $(p_1 + p_3)^n$; for no integer $n>1$ can be equal to $p_1 + p_3$. Hence $Z_2S_3$ is not a Lin ring.

In view of this we can prove.

**THEOREM [129]**: *Let $Z_2 = \{0, 1\}$ and $S_n$ be the permutation group of degree n. The group ring $Z_2S_n$ is not a Lin ring.*

The proof of this theorem is left for the reader as an exercise.

Now we proceed onto define Smarandache Lin rings.

**DEFINITION 3.10.46**: *Let R be a ring. R is said to be a Smarandache Lin ring (S-Lin ring) if R contains a S-subring B such that B is a Lin-ring. We do not demand every element of R to satisfy the Lin identity viz*

$$(xy - yx)^n = xy - yx \ or$$
$$(xy + yx)^n = xy + yx.$$

*but only if elements of B satisfy the Lin identity then it makes R a S-Lin ring.*

**Theorem 3.10.38**: *Let R be a Lin ring having a S-subring then R is a S-Lin ring.*

*Proof*: Follows from the very definition of S-Lin ring and Lin ring.

**Theorem 3.10.39**: *Let R be S-Lin ring then R is a S-ring.*

*Proof*: By the very definition of S-Lin ring it should have a S-subring so R will be a S-ring.

**Definition [160]**: *A ring R with 1 is said to satisfy the right super ore condition if for any r, s ∈ R there is some r' ∈ R such that rs = sr'.*

For more about these concepts please refer [160, 126].

**Theorem [126]**: *Let F be any field or a ring and G = $S_n$, n ≥ 3 be the symmetric group of degree 3. $FS_n$ does not satisfy super ore condition.*

*Proof*: To show $FS_n$ does not satisfy super ore condition it is enough if we show for some x, y ∈ $FS_n$ we have xy = yα is not true for any α ∈ $FS_n$.

Take

$$x = 1 + \begin{pmatrix} 1 & 2 & 3 & 4 & . & . & . & n \\ 3 & 2 & 1 & 4 & . & . & . & n \end{pmatrix} \text{ and } y = 1 + \begin{pmatrix} 1 & 2 & 3 & 4 & . & . & . & n \\ 1 & 3 & 2 & 4 & . & . & . & n \end{pmatrix}.$$

Clearly yx = xτ for no τ ∈ $FS_n$; hence the claim.

**Definition 3.10.47**: *Let R be a ring. We say R satisfies Smarandache super ore condition (S-super ore condition) if R has a S-subring A and for every pair x, y ∈ A we have r ∈ R such that xy = yr.*

**Theorem 3.10.40**: *Let R be a ring which satisfies super ore condition. If R has a S-subring then R satisfies the S-super ore condition.*

*Proof*: By the very definition of S-subring and S-super ore condition the results are easily proved.

**Definition [138]**: *Let R be a ring. R is said to be ideally strong, if every subring of R not containing identity is an ideal of R.*

**Example 3.10.41**: Let $Z_2G$ be the group ring of the group $G = \langle g \mid g^2 = 1 \rangle$ over $Z_2$. $Z_2G$ is an ideally strong ring.

**THEOREM [138]**: *Let $Z_2 = \{0, 1\}$ and $G = \langle g \mid g^{2n} = 1 \rangle$. The group ring $Z_2G$ is not an ideally strong ring.*

*Proof*: The set $S = \{0, 1 + g^n\}$ is a subring of $Z_2G$ but is not an ideal of $Z_2G$. Hence the claim.

**DEFINITION 3.10.48**: *Let R be a ring. We say R is a Smarandache ideally strong (S-ideally strong) ring if every S-subring of R is a S-ideal of R.*

**DEFINITION [137]**: *Let R be a ring; $\{I_j\}$ be the collection of all ideals of R. R is said to be a $I^*$-ring if every pair of ideals $I_1$, $I_2 \in \{I_j\}$ in R and for every $a \in R \setminus (I_1 \cup I_2)$ we have $\langle a \cup I_1 \rangle = \langle a \cup I_2 \rangle$ where $\langle \rangle$ denotes the ideal generated by a and $I_j$, $j = 1, 2$.*

**Example 3.10.42**: $Z_{12} = \{0, 1, 2, \ldots, 11\}$ ring of integers modulo 12, is not a $I^*$ ring for $I_1 = \{0, 6\}$ and $I_2 = \{0, 4, 8\}$ are two ideals of $Z_{12}$ and we have $\langle 3 \cup I_1 \rangle = \{0, 3, 6, 9\}$ where as $\langle 3 \cup I_2 \rangle = Z_{12}$.

**DEFINITION 3.10.49**: *Let R be a ring $\{A_i\}$ be the collection of all S-ideals of R. If for every pair of ideals $A_1$, $A_2 \in \{A_i\}$ we have for every $x \in R \setminus \{A_1 \cup A_2\}$, $\langle A_i \cup x \rangle = \langle A_2 \cup x \rangle$ and they generate S-ideals of R, then we say R is a Smarandache $I^*$-ring ( S-$I^*$ ring).*

**DEFINITION 3.10.50**: *Let R be a ring $\{A_j\}$ be the collection of all S-ideals of R if for $A_1$, $A_2 \in \{A_j\}$, we have some $x \in R \setminus \{A_1 \cup A_2\}$ such that $\langle A_1 \cup x \rangle = \langle A_2 \cup x \rangle$ are S-ideals of R then we say R is a Smarandache weakly $I^*$-ring (S-weakly $I^*$-ring).*

**THEOREM 3.10.41**: *All S-$I^*$ rings are S-weakly $I^*$-rings.*

*Proof*: Follows by the very definitions of S-$I^*$ ring and S-weakly $I^*$-ring.

Here we demand the ideals generated by $\langle A_i \cup x \rangle$ to be S-ideals. Further we have the following.

**THEOREM 3.10.42**: *If R is a S-$I^*$ ring or a S-weakly $I^*$-ring then R is a S-ring.*

*Proof*: Follows from the fact that for R to be a S-I$^*$ ring or S-weakly I$^*$ ring; R must contain nontrivial S-ideals which in turn will imply R is a S-ring.

**DEFINITION [136]**: *Let P and V be any two non-isomorphic finite rings, if there exists non-maximal ideals I of P and J of V such that P/I is isomorphic to V/J. The finite rings P and Q are called Q-rings.*

**Example 3.10.43**: Let $Z_4 = \{0, 1, 2, 3\}$ be the ring of integers modulo 4. $Z_8 = \{0, 1, 2, \ldots, 7\}$ be the ring of integers modulo 8. $Z_8 / J \cong Z_4 / I$ where $J = \{0, 4\}$ and $I = \langle 0 \rangle$. So $Z_4$ and $Z_8$ are Q rings. Similarly we can prove $Z_6$ and $Z_{12}$ are also Q-rings.

**THEOREM [136]**: *Let $Z_n = \{0, 1, 2, \ldots, n - 1\}$ , be the ring of integers modulo n, n not a prime then $Z_n$ is always a Q-ring.*

*Proof*: Left for the reader as an exercise to prove.

**DEFINITION [136]**: *Suppose R is a ring such that all of its ideals are maximal and if we have R /$\langle 0 \rangle$ is isomorphic to some ring then we call R a weakly Q-ring.*

**DEFINITION 3.10.51**: *Let R be a ring. A a S-ideal of R. R/A is defined as the Smarandache quotient ring (S-quotient ring) related to the S-ideal A.*

**DEFINITION 3.10.52**: *Let R and S be two rings; if we have S-ideals A and B of R and S respectively, such that the Smarandache quotient ring R/A is S-isomorphic with the Smarandache quotient ring S/B then we say the ring R is a Smarandache Q-ring (S-Q-ring).*

We assume the S-ideals A and B need not be S-maximal ideals of R and S respectively.

**DEFINITION [13]**: *A ring R is called a F-ring if there is a finite set X of non-zero elements in R such that $aR \cap X \neq \phi$ for any non-zero a in R. If in addition X is contained in the center of R; R is called a FZ-ring.*

**Example 3.10.44**: Let $Z_2 = \{0, 1\}$ be the field and $G = \langle g / g^2 = 1 \rangle$. The group ring $Z_2G$ is a F-ring, $X = \{1 + g\} \subset Z_2G$. Clearly $a.Z_2G \cap X \neq \phi$ for any non-zero a in $Z_2G$.

**THEOREM [135]**: *Let $Z_2 = (0, 1)$ be the field of characteristic 2 and $S_n$ be the symmetric group of degree n. The group ring $Z_2S_n$ is a F-ring.*

*Proof*: Take X = $\left\{ \alpha = \sum_{i=1}^{2m} s_i / s_i \in S_n \right\}$ with |supp $\alpha$| = 2m and 1 < 2m ≤ n!. For any

a $\in$ $Z_2S_n \setminus \{0\}$; $aZ_2S_n \cap X \neq \phi$; hence $Z_2S_n$ is a F-ring.

**DEFINITION 3.10.53**: *Let R be any ring; A a S-subring of R. We say R is a Smarandache F-ring (S-F-ring) if we have a subset X in R and a non-zero b $\in$ R such that bA $\cap$ X $\neq \phi$. It is pertinent to mention here that we need not take X as a subset of A but nothing is lost even if we take X to be a subset of A. Similarly b can be in A or in R.*

**DEFINITION [145]**: *Let R be a ring we say an element a $\in$ R \ {0, 2} is an SS element if $a^2$ = a + a.*

**DEFINITION [145]**: *Let R be a ring if R has atleast one SS-element other than 0 and 2 then we say R is a SS-ring.*

Now we proceed onto define Smarandache SS-elements.

**DEFINITION 3.10.54**: *Let R be a ring an element x $\in$ R is said to be a Smarandache SS element (SSS-element) of R if there exists y $\in$ R \ {x} with x.y = x + y.*

**Example 3.10.45**: Let $Z_{10}$ = {0, 1, 2, …, 9} be the ring of integers modulo 10. 4.8 ≡ 4 + 8 (mod 10). So 4 is a SSS element of $Z_{10}$.

**DEFINITION 3.10.55**: *Let R be a ring if R has atleast one nontrivial Smarandache SS-element we call R a SSS ring.*

**Example 3.10.46**: $Z_{14}$ is a SSS-ring for 4, is a SSS-element as 4 + 6 ≡ 4.6 ≡ 10 (mod 14).

**Example 3.10.47**: Let $Z_{15}$ = {0, 1, 2, …, 14} be the ring of integer modulo 15, 3 is a SSS element for 3.9 ≡ 3 + 9 ≡ 12 (mod 15).

**Example 3.10.48**: $Z_9$ be the ring of integers modulo 9. 3.6 ≡ 3 + 6 ≡ 0 (mod 9), 5.8 ≡ 5 + 8 ≡ 4 (mod 9). This ring has two SSS-elements.

**Example 3.10.49**: $Z_8$ has SSS-element.

Here we solve the problem proposed by [16]: "He asks whether there exists a commutative ring R with the property satisfying the following condition: $R^2$ = R and a

+ a = 0 = $a^2$ for all a $\in$ R. Now if R is a commutative ring in which $x^2$ = 0 for all x $\in$ R then we have xy = 0 for all x, y $\in$ R or characteristic R = 2. For $(x + y)^2$ = 0 but $x^2 + y^2 + 2xy = 0$ so xy = 0 for all x, y $\in$ R or characteristic of R is two. So we prove that such rings do not exist. For the rings R, which satisfy the given condition, cannot contain 1, the identity. Secondly $R^2$ = R is impossible as a + a = $a^2$ = 0 so $R^2$ = {0}. Hence the claim; thus the answer to the question in [16] is answered in the negative.

**DEFINITION [35]**: *Let R be a ring. A subring S $\neq$ {0} of R is said to be a trivial subring of R of $S^2$ = (0).*

***Example 3.10.50***: Let $Z_2$G be the group ring of the group G = $\langle g / g^{2n} = 1 \rangle$ over $Z_2$. This has trivial subrings, for S = {0, 1 + g + … + $g^{2n-1}$} is such that $S^2$ = (0) and $S_1$ = (1 + $g^n$, 0) is such that $S_1^2$ = (0).

**DEFINITION 3.10.56**: *Let R be a ring, we say R has a Smarandache trivial subring (S-trivial subring) if R has a S-subring A such that A has a subring B (B $\subset$ A) with $B^2$ = (0).*

We have to make this form of definition as S-subring A has a subfield in it, so $A^2$ = {0} is impossible, so to overcome this we have to define a subring B $\subset$ A such that $B^2$ = (0).

We leave it for the reader to obtain some interesting results in this direction.

**DEFINITION [65]**: *Let R be a ring, we say R is a $\gamma_n$-ring, n > 1, n an integer if $x^n - x$ is an idempotent for all x $\in$ R.*

***Example 3.10.51***: Let $Z_2$ = {0, 1} be the field and G = $\langle g / g^3 = 1 \rangle$, the group ring $Z_2$G is a $\gamma_n$-ring. For more about $\gamma_n$-ring refer [65].

**DEFINITION 3.10.57**: *Let R be a ring. If $x^n - x$ is a S-idempotent for some integer n > 1, for all x $\in$ R then we say R is a Smarandache $\gamma_n$-ring (S-$\gamma_n$ ring).*

It is left for the reader to obtain some nice examples of S-$\gamma_n$-rings.

**THEOREM 3.10.43**: *Let R be a field of characteristic 0 and G be a torsion free abelian group. The group ring KG is not a S-$\gamma_n$-ring.*

*Proof*: KG is a domain, hence the claim.

The concept of demi modules is a generalization of modules. For every module is a demi module and not conversely. To define demi module we define demi subring.

**DEFINITION [148]**: *Let R be a commutative ring with unit. A non-empty subset V of R is said to be a demi subring of R if V is closed with respect to '+' and '.' of R.*

**Example 3.10.52**: Let Z be the ring of integers, $Z^+$ is a demi subring ($Z^+$ set of positive integers).

**Example 3.10.53**: Let $Z_p = \{0, 1, 2, \ldots, p-1\}$, p a prime, $Z_p$ does not have a demi subring.

**DEFINITION [148]**: *Let R be a commutative ring with unity; P is said to be a demi module over R if P is a semigroup under '+' and '.' and there exists a nontrivial demi subring V of R such that for every $v \in V$ and $p \in P$; vp and pv $\in P$; $0 \in P$ and $v(p_1 + p_2) = vp_1 + vp_2$, $v_1(v_2p) = (v_1v_2)p$. for $p_1, p_2, p \in P$; $v_1, v_2, v \in V$.*

**DEFINITION [148]**: *Let R be a commutative ring and let P be a demi module relative to the demi subring V. Then a non-empty subset T of P is a subdemi module, if T is a demi module for the same demi subring.*

**DEFINITION 3.10.58**: *Let R be a commutative ring with 1. A subset S of R is said to be a Smarandache demi subring (S-demi subring) of R if*

1. *(S, +) is a Smarandache semigroup.*
2. *(S, .) is a Smarandache semigroup.*

**THEOREM 3.10.44**: *All S-demi subrings of the ring are demi subrings.*

*Proof*: Left as an exercise for the reader.

**DEFINITION 3.10.59**: *Let R be a commutative ring with 1. P is said to be Smarandache demi module (S-demi module) over R if*

1. *P is a S-semigroup under + and '.'.*
2. *There exists a nontrivial S-demi subring V of R such that for every $p \in P$ and $v \in V$, vp and pv $\in P$*
3. *$v(p_1 + p_2) = vp_1 + vp_2$.*
4. *$v(v_1p) = (vv_1) p$ for all $p_1, p_2, p \in P$ and $v, v_1 \in V$.*

**THEOREM 3.10.45**: *Let R be a ring; P be a S-demi module relative to V, then P is a demi module relative to V, V a subring of R.*

*Proof*: Left as an exercise for the reader to prove.

**DEFINITION 3.10.60**: *Let R be a ring. P a S-demi module relative the S-demi subring V of R. A non-empty subset T of P is said to be a Smarandache subdemi module (S-subdemi module) if T is a Smarandache demi module for the same S-demi subring.*

The reader is requested to obtain nice results in this direction.

**DEFINITION [63]**: *Let R be a ring, R is said to be locally unitary if for each x ∈ R there exists an idempotent e ∈ R for which ex = xe = x.*

Here we define semiunitary ring using semi idempotents.

**DEFINITION [152]**: *Let R be a ring. R is said to be locally semiunitary if for each x ∈ R there exists a semi idempotent s ∈ R such that xs = sx = x.*

**THEOREM [152]**: *Let R be a locally unitary ring then R is a locally semiunitary ring.*

*Proof*: By the very definition of locally unitary ring and locally semiunitary ring the result follows.

**THEOREM [152]**: *A locally semiunitary ring in general is not locally unitary ring.*

*Proof*: By an example; $Z_2 = \{0, 1\}$ be the field of characteristic two and $G = \langle g / g^2 = 1 \rangle$. The group ring $Z_2G = \{0, 1, g, 1 + g\}$ is locally semiunitary for $(1 + g) g = g(1 + g) = 1 + g$ where g is not an idempotent as $g^2 = 1$. Thus $Z_2G$ is not locally unitary.

**DEFINITION 3.10.61**: *Let R be a ring if for every element x ∈ R there exists a S-idempotent e in R such that xe = ex = x. Then we call the ring R a Smarandache locally unitary ring (S-locally unitary ring).*

**DEFINITION 3.10.62**: *Let R be a ring. If for every x ∈ R there exists a S-semi idempotent s of R such that xs = sx = x then we call the ring R a Smarandache locally semiunitary ring (S-locally semi unitary ring).*

**DEFINITION [151]**: *Let (R, +, .) be a ring. A non-empty subset S of R is called a closed net of R if S is a closed set of R under the operation '.' and is generated by a single element. That is S is a semigroup under multiplication '.'.*

**DEFINITION [151]**: *Let (R, +, .) be a ring. If R is contained in a finite union of closed nets of R then we say the ring R has a closed net.*

**DEFINITION [151]**: *Let R be a ring we say R is a CN-ring if R = $\cup S_i$ where $S_i$'s are closed nets such that $S_i \cap S_j = \phi$ or {1} or {0} if $i \neq j$ and $1 \in R$ and $S_i \cap S_j = S_i$ if i = j and each $S_i$ is a nontrivial closed net of R.*

**Example 3.10.54**: Let $Z_8$ = {0, 1, 2, …, 7} be the ring of integers modulo 8. Clearly R is not a CN-ring for take $S_1$ = {4, 6}, $S_2$ = {1, 3}, $S_3$ = {5, 1}, $S_4$ = {1, 7} and $S_5$ = {0, 2, 4}. Easily verified.

**DEFINITION [151]**: *Let R be a ring. If R $\subset \cup S_i$ with $S_i \cap S_j \neq \phi$ or {1}. Then we say R is a weakly CN-ring.*

We have rings which may not even be a weakly CN-ring.

**Example 3.10.55**: Let $Z_9$ = {0, 1, 2, …, 8} be the ring of integers modulo 9. It is easily verified $Z_9$ is a CN-ring. For $S_1$ = {0, 3}, $S_2$ = {0, 6} and $S_3$ = {2, 4, 8, 7, 5, 1} closed nets of $Z_9$.

**THEOREM [151]**: *Every CN-ring is a weakly CN-ring. But a weakly CN-ring in general is not a CN-ring.*

*Proof*: Left as an exercise to the reader.

**THEOREM [151]**: *$Z_p$ be the prime field of characteristic p. $Z_p$ is not a CN-ring and not even a weakly CN-ring.*

*Proof*: Left for the reader to prove.

Now we proceed onto define Smarandache CN-rings and Smarandache weakly CN-rings to this end we define Smarandache closed net in rings.

**DEFINITION 3.10.63**: *Let R be a ring, we say a subset S of R is said to be Smarandache closed net (S-closed net) if*

1. *S is a semigroup.*
2. *S is a S-semigroup.*

From this we easily see that all S-closed nets are closed nets but every closed net in general is not a S-closed net.

**DEFINITION 3.10.64**: *Let R be a ring. If R is contained in the finite union of S-closed nets of R then we say the ring R has a Smarandache closed net (S-closed net).*

**DEFINITION 3.10.65**: *Let R be a ring, we say R is a Smarandache CN-ring (S-CN-ring) if R = $\cup S_i$ where $S_i$'s are S-closed nets such that $S_i \cap S_j = A$, $i \neq j$, $A \neq S_i$ or $A \neq S_j$ where A is a subgroup of $S_i$.*

**DEFINITION 3.10.66**: *Let R a ring, if $R \subset \cup S_i$ where $S_i$'s are Smarandache closed nets then we say R is a Smarandache weakly CN-ring (S-weakly CN ring).*

**THEOREM 3.10.46**: *Let R be a S-CN-ring. Then R is a weakly CN-ring.*

*Proof*: By the very definition the result is straight forward.

**THEOREM 3.10.47**: *Every S-weakly CN-ring is a weakly CN-ring and not conversely.*

*Proof*: It can be proved by simple working, hence left for the reader as an exercise.

**DEFINITION [153]**: *Let R be a ring. A subset M of R with $|M| \leq 2$, $|M + M| \leq 2$ and $|M^2| \leq 2$ is called the tight subset of R.*

**DEFINITION [153]**: *Let R be a ring. R is said to be a tight ring if we can find a subset M of R which is a tight subset of R.*

***Example 3.10.56***: Let G = $\langle g \,/\, g^2 = 1 \rangle$ and $Z_2 = \{0, 1\}$ the ring of integers modulo 2. The group ring $Z_2 G$ is a tight ring for it has the tight subset M = $\{0, 1 + g\}$ such that $|M| \leq 2$, $|M + M| \leq 2$ and $|M^2| \leq 2$.

**DEFINITION [153]**: *Let R be a ring R; is said to be a strong tight ring if every subset M with $|M| \leq 2$ of R is a tight subset of R.*

***Example 3.10.57***: The group ring R = $Z_2 G$ where G = $\langle g \,/\, g^2 = 1 \rangle$ is a strong tight ring.

**THEOREM [153]**: *Every strong tight ring is a tight ring but every tight ring need not be a strong tight ring.*

*Proof*: By the very definition of tight ring and strong tight ring. To prove the converse we see $Z_8 = \{0, 1, 2, \ldots, 7\}$ is a tight ring which is not a strong tight ring.

**THEOREM [153]**: *No ring of characteristic 0 is a tight ring.*

*Proof*: Left as an exercise for the reader to prove.

**DEFINITION [153]**: *Let R be a ring; R is said to be r-tight ring i.e., $T_r$-ring, $r \geq 2$ if R contains a subset M with $|M| \leq r$ implies $|M + M| \leq r$ and $|M^2| \leq r$. Clearly when r = 2 we get the tight ring.*

*Every $T_r$-ring is a $T_i$-ring for all $i \leq r$. Thus we have T-ring $\subset T_3$-ring $\subset \ldots \subset T_r$-ring.*

**THEOREM [153]**: *Every $T_2$ ring is a $T_3$ ring but all $T_3$ rings need not be $T_2$-rings.*

*Proof*: By the very definition every $T_2$ ring is a $T_3$ ring, to prove the converse we give an example. Consider $Z_9 = \{0, 1, 2, \ldots, 8\}$, is a $T_3$-ring but it is not a $T_2$-ring.

**THEOREM [153]**: *The ring of integers is not a $T_i$-ring for any finite $i \geq 2$.*

*Proof*: Left for the reader as an exercise to prove.

**DEFINITION 3.10.67**: *Let R be a ring we say a non-empty subset M of R is said to be a Smarandache tight set (S-tight set) if*

    1. *M contains a subset S where S is a semigroup under '.' of R.*
    2. *If $|M| \leq 2$, $|M + M| \leq 2$, $|M^2| \leq 2$.*

**DEFINITION 3.10.68**: *Let R be a ring. R is said to be a Smarandache tight ring (S-tight ring) if we can find a subset M of R which is a S-tight set of R.*

**THEOREM 3.10.48**: *Every S-tight set of R is a tight set of R.*

*Proof*: Obvious by the very definitions.

**THEOREM 3.10.49**: *Every S-tight ring R, is a tight ring.*

*Proof*: Left as an exercise to the reader to prove.

**DEFINITION 3.10.69**: *Let R be a ring, R is said to be Smarandache strong tight ring (S-strong tight ring) if every subset M of R is a S-tight set of R.*

**DEFINITION 3.10.70**: *Let R be a ring, R is said to be Smarandache r-tight ring or a S-$T_r$ ring $r \geq 2$ if contains a subset M with $|M| \leq r$, $|M + M| \leq r$, $|M^2| < r$ and M is a S-tight set of R.*

Now we proceed on to define the concept of finite quaternion rings and skew fields as this concept would help in defining in chapter IV Smarandache rings of level II and Smarandache mixed direct product of rings.

**DEFINITION [149]**: *Let $Z_n$ be the ring of integers modulo n. Let $P = \{p_0 + p_1 i + p_2 j + p_3 k / p_0, p_1, p_2, p_3 \in Z_n, n$ finite, $n > 2\}$ Define '+' and '.' on P as follows*

$X = p_0 + p_1 i + p_2 j + p_3 k$ and
$Y = q_0 + q_1 i + q_2 j + q_3 k$ be in P
$X + Y = (p_0 + q_0) + (p_1 + q_1)i + (p_2 + q_2)j + (p_3 + q_3)k.$
$X.Y = (p_0 + p_1 i + p_2 j + p_3 k)(q_0 + q_1 i + q_2 j + q_3 k) = [p_0 q_0 + (n - 1) p_1 q_1 + (n - 1) p_2 q_2 + (n - 1) p_3 q_3] + (p_0 q_1 + p_1 q_0 (n - 1) + p_3 q_2 + p_2 q_3)i + (p_0 q_2 + p_2 q_0 + (n-1)p_1 q_3 + p_3 q_1)j + (p_0 q_3 + p_3 q_0 + p_1 q_3 + (n-1)p_3 q_1)k$

*where $i^2 = j^2 = k^2 = (n - 1) = ijk$ and $ij = (n - 1)ji = k$ where $ji = (n - 1) kj = i$ and $ki = (n - 1) ik = j.$*

*Clearly in P, $0 = 0 + 0i + 0J + 0k$ is the identity with respect to '+' and $1 = 1 + 0i + 0j + 0k$ is the identity with respect to '.'.*

*Now P is a ring called the ring of real quaternion of characteristic n, n a finite prime. If n is a composite number then we have P to be a ring with divisors of zero.*

**THEOREM [149]**: *Let $P = \{\alpha_0 + \alpha_1 i + \alpha_2 j + \alpha_3 k / \alpha_0, \alpha_1, \alpha_2, \alpha_3 \in Z_p = \{0, 1, 2, \ldots, p - 1\}$ be defined as above. Then P is a prime skew field.*

*Proof*: Please refer [149] for proof.

**THEOREM [149]**: *Let $P = \{\alpha_0 + \alpha_1 i + \alpha_2 j + \alpha_3 k / \alpha_0, \alpha_1, \alpha_3, \alpha_2 \in Z_n = \{0, 1, 2, \ldots, n - 1\}\}$ ring of integers modulo n, n a composite number with $i^2 = j^2 = k^2 = n - 1 = ijk$. $ij = (n - 1)ji = k$ and so on. Then P is a ring with divisors of zero.*

*Proof*: Left as an exercise for the reader to prove.

Let P be a ring defined as above. Let $G = \{\alpha_0 + \alpha_1 i + \alpha_2 j + \alpha_3 k / \alpha_0, \alpha_1, \alpha_2, \alpha_3 = m$ where m is a zero divisor in $Z_n$, $\alpha_0, \alpha_1, \alpha_2, \alpha_3 \in Z_n\}$ denote the set of nilpotent elements of P or zero divisors of P.

Let $V = \{\alpha_0 + \alpha_1 i + \alpha_2 j + \alpha_3 k / \alpha_0 + \alpha_1 + \alpha_2 + \alpha_3 = t$ where t is a unit in $Z_n$; $\alpha_0, \alpha_1, \alpha_2, \alpha_3 \in Z_n\}$ denote the set of all units of P. Then it has been verified if $n = p^r$ where p is a prime then $P = G \cup V$. In view of this we propose a few problems in Chapter 5, using the notation P, G and V.

## PROBLEMS:

1. Is $Z_{25}$ the ring of integers modulo 25 a reduced ring?

2. Find all the S-nilpotent elements of the group ring $Z_4 S_3$.

3. Test whether the semigroup ring $Z_6 S(4)$ is a reduced ring.

4. Find all S-nilpotent elements in $QS(5)$, the semigroup ring of the semigroup $S(5)$ over the field of rationals Q.

5. Is $Z_{14}$ a S-zero square ring?

6. Can a S-zero square ring be of order 19?

7. Is the group ring $Z_{12} S_4$ a S-zero square ring?

8. Give an example of a S-null semigroup.

9. Is $Z_{32}$ a S-null ring?

10. Prove $Z_{36}$ is a S-null ring.

11. Prove the group ring $Z_n G$, for any group G is a S-p ring.

12. Show the semigroup ring $Z_3 S(5)$ is a S-p-ring.

13. Give an example of a semigroup ring, which is a S-E ring.

14. Is the group ring $Z_2 S_4$ a S-E-ring?

15. Is $Z_8$ a S-pre J-ring?

16. Is the semigroup ring $Z_{12} S(5)$ a S-pre J-ring?

17. Show $Z_{24}$ is a S-inner zero square ring.

18. Give an example of a S-inner square ring which is not an inner zero square ring.

19. Is $Z_{24}$ a S-weak inner zero square ring?

20. Give an example of S-semi prime ring.

21. Is $Z_{14}$ a S-Marot ring?

22. Prove or disprove: $Z_{15}$ is a S-Marot ring.

23. Is $Z_2 G$, where $G = \langle g \, / \, g^8 = 1 \rangle$ is a S-subsemi ideal ring?

24. Give an example of a S-pre-Boolean ring.

25. Is $R = Z_3 \times Z_3 \times Z_3$ a S-filial ring?

26. Show $Z_{21}$ is a S-ideal ring.

27. Give an example of a S-generalized ideal ring.

28. Is the group ring $Z_2 S_3$ a S-s-weakly regular ring?

29. Illustrate by an example the S-strongly sub commutative ring.

30. Give an example of a S-Chinese ring of finite order.

31. Give an example of group ring which is S-strongly s-decomposable.

32. Illustrate the definition of S-strong right D-domain by an example.

33. Show the group ring $Z_2 S_4$ is a S-J-ring.

34. Can $Z_2 S_3$ be a S-strong subring ring?

35. Is $Z_7 S_2$ a S-strong ideal ring?

36. Prove $Z_7$ is a S-weakly Boolean ring.

37. Give an example of a S-pre p-ring.

38. Illustrate by an example S-multiplication right ideal.

39. Is $Z_{16}$ a S-f-ring?

40. Give a pair of S-obedient ideals in $Z_{26}$.

41. Illustrate by an example the S-obedient ideal ring.

42. Show $Z_2S_4$ is not a S-Lin ring.

43. Can $M_{2\times2} = \{(a_{ij}) / a_{ij} \in Z_4\}$ be a S-Lin ring?

44. Illustrate the definition of S-super ore condition in a ring by an example.

45. Show $Z_2G$, where $G = \langle g / g^8 = 1 \rangle$ is not a S-ideally strong ring.

46. Can $Z_{24}$ be a S-$I^*$-ring?

47. Prove the group ring $Z_2S_6$ is a S-F-ring.

48. Can $Z_{22}$ have SSS-elements?

49. Show by an example, a ring, which is not a S-$\gamma_n$-ring.

50. Give an example of a S-$\gamma_n$-ring.

51. Give a ring R which has a S-demi subring.

52. Is the group ring $Z_7S_4$, S-locally semiunitary?

53. Can the semigroup ring $Z_6S(3)$ be S-locally unitary?

54. Does there exist a locally unitary ring, which is not S-locally unitary?

55. Does there exist a locally semiunitary ring, which is not a S-locally semiunitary ring?

56. Give an example of a weakly CN-ring, which is not a CN-ring.

57. Prove $Z_{11}$ is not even a weakly CN-ring.

58. Is $Z_{24}$ a CN-ring? Justify your answer.

59. Find whether $Z_2S_3$ is a S-CN-ring.

60. Give an example of a CN-ring, which is not a S-CN-ring.

61. Give an example of a weakly CN-ring, which is not a S-weakly CN-ring.

62. Prove in a E-ring every prime ideal is maximal.

63. Is the semigroup ring $Z_2S(5)$ an E-ring?

64. Is the group ring $Z_7S_3$ a E-ring?

65. Give an example of a 5-tight ring which is not a 2-tight ring

66. Show by an example that a tight ring in general is not a S-tight ring.

67. Show by an example that a tight set in general is not a S-tight set of R.

68. Give an example of a semi idempotent in a ring R that is not S-semi idempotent.

69. Give an example of a ring in which every set is a S-tight set.

70. Give an example of a division ring of order 81.

# SOME NEW NOTIONS ON SMARANDACHE RINGS

In this chapter we introduce several new notions and concepts in ring theory to Smarandache rings. This chapter is organized into five sections. In section one we introduce the concept of Smarandache mixed direct product of rings which alone helps us in the building of Smarandache rings of level II, which is dealt in section two. It is noteworthy and important here that several concepts enjoyed by the ring introduced by researchers in different nations have not been consolidated or taken notice of by many books on ring theory. So in this book section three, four and five are especially devoted to the recalling of these concepts and also simultaneously defining the Smarandache analogue of them. The concepts which are taken from different researchers in ring theory are listed in the references/bibliography. Thus at this juncture it is pertinent to mention this book will become an attraction simultaneously to both the ring theorist and Smarandache algebraist.

Section three of this chapter separately gives the introduction and study of elements which enjoy new special properties in a ring like magnifying elements, shrinking elements, semi idempotents and so on. The localization or the Smarandache analaogue is carefully brought out at every stage.

Section four is devoted to the study of new or special properties enjoyed by the substructure of a S-ring like subsets, semigroup (with respect to '+' or '.') subgroups, subrings, ideals etc.

Many new concepts are defined and Smarandache analogue of them are obtained. The importance of these Smarandache analogue or Smarandache notions is that even when a ring fails to enjoy certain property fully, it can enjoy the same property sectionally. So except for these Smarandache notions such local study or a sectional study in a ring would be impossible.

The final section of the chapter entitled miscellaneous properties of Smarandache ring introduces and studies several concepts; the prominent among them are hyperrings, lattice substructures of S-ideals, S-rings etc.

## 4.1 Smarandache Mixed Direct Product Rings

In this section we define what are called Smarandache mixed direct product rings. Only this concept of rings paves way for the introduction of S-rings of level II. For a ring to be a S-ring of level I where we have not mentioned level I we see the ring should contain a subset which is a field. This condition makes the ring of integers Z, the field of rationals to be useless in the analysis of S-rings. So we have used S-ring level II to overcome this problem, which is done by introducing the concept of Smarandache mixed direct products. We just recall the definition of S-ring II, S-subring II and S-ideal II.

**DEFINITION 4.1.1**: *Let $R = R_1 \times R_2$ where $R_1$ is a ring and $R_2$ is an integral domain or a division ring. Clearly this product is called the Smarandache mixed direct product (S-mixed direct product) of two rings which is easily verified to be a ring.*

It is to be noted if both $R_1$ and $R_2$ are just rings then we don't call the direct product as a Smarandache mixed direct product. We extend it to any number of rings and integral domains or division rings.

**DEFINITION 4.1.2**: *Let $R = R_1 \times R_2 \times \ldots \times R_n$ is called the Smarandache mixed direct product of n-rings (S-mixed direct product of n-rings) if and only if at least one or some of the $R_i$'s is an integral domain or a division ring.*

***Example 4.1.1***: Let $R = Z \times Q$. R is a S-mixed direct product of rings. Clearly R is a S-ring as $\{1\} \times Q$ is a field contained in R.

***Example 4.1.2***: Let $R = Z \times R_2$ where $R_2 = \{0, 2\}$ modulo 4. Clearly R is a S-mixed direct product of rings, but R is not a S-ring.

***Example 4.1.3***: Let $R = ZS_3 \times Z$. Clearly R is a S-mixed direct product of rings but is not a S-ring.

**DEFINITION 4.1.3**: *Let R be a ring. R is called a Smarandache ring of level II or in short S-ring II if R contains an integral domain or a division ring.*

***Example 4.1.4***: $ZS_3$ the group ring is a S-ring II, for $Z \subseteq ZS_3$ is an integral domain.

***Example 4.1.5***: Z is a S-ring II and never a S-ring I.

***Example 4.1.6***: $R = Z \times Z \times Z$ is a S-ring II and never a S-ring I.

**THEOREM 4.1.1**: *Every S-ring I is a S-ring II and not conversely.*

*Proof*: By the very definition of S-ring I and S-ring II we see every S-ring I is a S-ring II as all fields are trivially integral domains or division rings.

To prove the converse we give the following example. The ring of integers Z is a S-ring II but is never a Smarandache ring I. Hence the theorem.

**THEOREM 4.1.2**: *Z[x], the polynomial ring is a S-ring II.*

*Proof*: Easily verified.

**THEOREM 4.1.3**: *The class of rings $Z_p$, p a prime are not S-ring I or S-ring II.*

*Proof*: Obvious by the very definition, as $Z_p$ has no non-trivial subfields.

**THEOREM 4.1.4**: *Q the field of rationals is not a S-ring but is a S-ring II.*

*Proof*: Q has no proper subsets which is a field so is not a S-ring, but $Z \subset Q$ is an integral domain so Q is a S-ring II.

**COROLLARY :** *Q[x] is a S-ring I and S-ring II.*

*Proof*: Obvious as $Q \subset Q[x]$.

**DEFINITION 4.1.4**: *Let R be a ring. We say a subset A of R is said to have a Smarandache subring of level II (S-subring II) if*

1. *A is a subring of R.*
2. *A has a proper subset P where P is an integral domain or a division ring under the operations of R.*

***Example 4.1.7***: Let Z be the ring. Z a S-subring of level II for take A = 2Z, A is a subring and P = 8Z is an integral domain contained in A. Hence the claim.

**THEOREM 4.1.5**: *Let R be a ring which has a S-subring II then R is a S-ring II.*

*Proof*: Since R contains a S-subring II, say A and has a proper subset P which is an integral domain or a division ring, we see $P \subset A \subset R$ so $P \subset R$. Hence R is a S-ring II.

**THEOREM 4.1.6**: *Let R be a S-ring II then R need not in general have a S-subring II.*

*Proof*: Consider $Z_6 = \{0, 1, 2, \ldots, 5\}$. Clearly A = $\{0, 2, 4\}$ is a field so $Z_6$ is a S-ring II but $Z_6$ has no proper subring which contains an integral domain or a field or a division ring. Thus $Z_6$ doesn't contain a S-subring II but $Z_6$ is a S-ring II.

***Example 4.1.8***: Let R = $Z_6 \times Z$ clearly R has S-subring II.

***Example 4.1.9***: Let R = Z × Z × $Z_8$, R has S-subring II and no S-subring I. In fact R is not a S-ring only a S-ring II.

**DEFINITION 4.1.5**: *Let R be a ring. A non-empty subset I of R is said to be a Smarandache ideal of level II (S-ideal II) if*

1. *I is a S-subring II.*
2. *ri and ir ∈ I for all r ∈ R and i ∈ I.*

*The notion of Smarandache right/ left ideal of level II can be defined as in the case of right/ left ideals.*

***Example 4.1.10***: Let Z be the ring of integers, Z has S-ideal II.

**THEOREM 4.1.7**: *Let R be a ring if R has a S-ideal II, then R is a S-ring II.*

*Proof*: By the very definition of S-ring II and S-ideal II, the result follows:

**THEOREM 4.1.8**: *Let R be a S-ring II, R need not have S-ideal II.*

*Proof*: By an example. Consider the ring $Z_6 = \{0, 1, 2, 3, 4, 5\}$. Clearly $Z_6$ is a S-ring II but $Z_6$ has no S-ideal II.

It is left as an exercise for the reader to prove the following theorem:

**THEOREM 4.1.9**: *Every S-ideal II is a S-subring II and not conversely.*

It is to be noted that as in the case of S-rings the notion of S-idempotents, S-units and S-zero divisors are defined we do not see any distinction of them in case of S-rings I or S-ring II.

<u>**PROBLEMS:**</u>

1. Show S = $Z_6$× Q is a S-ring II. Find a S-subring.
2. Is R = $Z_8$× $Z_8$ a S-ring II? Substantiate your claim.
3. Give an example of a S-ring II which is not a S-ring I (Examples should be other than the ones discussed in this section).
4. Can $M_{2\times2} = \{(a_{ij}) / a_{ij} \in Z_3\}$, the ring of 3 × 3 matrices with entries from $Z_3$ be a S-ring II?
5. Does $M_{2\times2}$ given in problem 4 have
   a. S-ideal II?
   b. S-subring II?
6. Give an example of a S-subring II which is not an S-ideal II.

7. Can the group ring $Z_8 S_3$ be a S-ring II? Justify your answer.
8. Prove $Z_7 S(3)$ is a S-ring II. Find an S-ideal II of $Z_7 S(3)$.
9. Can $Z_7 S(3)$ have S-ideals I?
10. Find an S-subring II which is not an S- ideal II of $R = Z_6 \times Z_7 \times Z_2$.

## 4.2 Smarandache rings of level II

In the previous section we just defined the concept of S-ring II. Here we discuss some important and interesting properties about them and we illustrate them by examples. We request the reader to find and introduce and study all the properties existing in S-ring I to the case of S-ring II. Though we had introduced S-commutative rings in Chapter III we recall it in this section.

**DEFINITION 4.2.1**: *Let R be a S-ring II. We say R is a Smarandache commutative ring II (S-commutative ring II) if R has a proper subset A where A is a S-subring II and A is a commutative ring.*

***Example 4.2.1***: Let $R = ZS_3$. The ring R is non-commutative. Clearly R is a S-ring II. But R is a S-commutative ring II.

**THEOREM 4.2.1**: *If R is a commutative ring and has a S-subring II then R is a S-commutative ring II.*

*Proof*: Left for the reader to verify as it is an easy consequence of the definition.

**THEOREM 4.2.2**: *Let R be a ring. If R is a S-commutative ring II then R in general need not be a commutative ring.*

*Proof*: The ring $Z_2 S_3 = R$ given in example 4.2.1 is a non-commutative ring but R is clearly a S-commutative ring.

S-mixed direct product of rings will help us to get several examples of such rings.

**DEFINITION 4.2.2**: *Let R be a ring. If every S-subring II of R happens to be a commutative subring then we say R to be a S-strongly commutative ring II. (S-strongly commutative ring II)*

**THEOREM 4.2.3**: *Every S-strongly commutative ring II is a S-commutative ring II and not conversely.*

*Proof*: Follows from the very definitions of these concepts. We prove the converse by an example. Let $R = ZS_4$. Clearly the group ring $ZS_4$ is not a commutative ring but it is a S-commutative ring II. Further $ZS_4$ is not a S-strongly commutative ring II as $ZA_4$ is a S-subring II but $ZA_4$ is not commutative. Hence the claim.

Thus we have the following relational chain. That is all commutative rings with S-subrings II are both S-commutative rings II and S-strongly commutative rings II.

The proof of the following theorem is left as an exercise to the reader.

**THEOREM 4.2.4**: *Let R be a ring. If R has S-ideal II then it need not imply R has a S-ideal I.*

While defining the concept of A.C.C. and D.C.C. to S-ring of level II, i.e. in case of Smarandache A.C.C. (S.A.C.C) on rings of level II we consider only chain of S-ideals II. Similarly for Smarandache D.C.C. (S.D.C.C) on rings of level II we take only S-ideals II. So it is easily seen even if a ring satisfies A.C.C. or D.C.C. on ideals it need not have any relevance for S.A.C.C. or S.D.C.C. of level II on S-ideals II.

**DEFINITION 4.2.3**: *Let R be a ring. A be a S-ideal II of R. We say A is a Smarandache maximal ideal II (S-maximal ideal II) of R if $A \subseteq S \subseteq R$ where S is another S-ideal II of R then either S = A or S = R.*

**DEFINITION 4.2.4**: *Let R be a ring. A be a S-ideal II of R. We say A is a S-minimal ideal II (S-minimal ideal II) of R if for any S-ideal II B of R if $B \subseteq A \subseteq R$ implies B = A or B is empty.*

**DEFINITION 4.2.5**: *Let R be a ring. Let A be a S-ideal II of R; we say A is a Smarandache principal ideal II (S-principal ideal II) of R if A is itself a principal ideal of R.*

**DEFINITION 4.2.6**: *Let R be a ring. A be a S-ideal II of R. A is said to be a S-prime ideal if A is a prime ideal of R.*

**DEFINITION 4.2.7**: *Let R and $R_1$ be two Smarandache rings II. We say a map $\phi$: $R \rightarrow R_1$ is a Smarandache ring homomorphism II (S-ring homomorphism II) if $\phi$ restricted to the integral domain or division rings A and $A_1$ of R and $R_1$ respectively is a integral domain homomorphism or division ring homomorphisms i.e. $\phi (a + b) = \phi (a) + \phi (b)$, and $\phi(ab) = \phi(a) \phi(b)$ for all a, $b \in A$. $\phi$ may or may not be even defined on other elements of R. $\phi$ the Smarandache ring homomorphism II, is a Smarandache ring isomorphism II if $\phi$ : $A \rightarrow A_1$ is an isomorphism from A to $A_1$.*

From this we see certainly the kernel of any homomorphism will be an ideal. "Is it a S-ideal I or S-ideal II?" is an open problem for the reader to solve.

**DEFINITION 4.2.8**: *Let R be a ring. I be a S-ideal II of R. The Smarandache quotient ring II (S-quotient ring II) is defined as R / I (R / I defined in a similar way as that of quotient rings).*

Once again will R / I given in definition 4.2.8 be a S-ring II is an open problem. It may be or it may not be. If I is a maximal ideal and R / I is not a prime ring certainly R / I is a S-ring I as well as S-ring II.

The polynomial rings P[x] will be S-rings II provided P is a S-ring I or S-ring II. The question, when are matrix rings $M_{n \times n}$ S-ring I or S-ring II is yet another interesting study.

<u>**PROBLEMS:**</u>

1. Can $M_{4 \times 4} = \{(a_{ij}) / a_{ij} \in Z\}$ be a S-commutative ring II?
2. Can $M_{4 \times 4}$ in problem 1 have S-subring II?
3. Find conditions on the group ring $Z_p G$ to have
       i. S-subrings II.
       ii. S-ideals II.
   Give at least examples of group rings which have S-subrings II and S-ideals II.
4. Give an example of a S-subring II which is not a S-ideal II.
5. Does there exist a ring in which all S-subrings II are S-ideals II?
6. Does there exist a ring in which all subrings are
       i. S-subrings I?
       ii. S-subrings II?
7. Is $M_{2 \times 2} = \{(a_{ij}) / a_{ij} \in Z_4\}$ a S-ring II?
8. Does the ring $M_{2 \times 2}$ given in problem 7 have
       i. S-subrings II?
       ii. S-ideals II?
9. Can $M_{3 \times 3} = \{(a_{ij}) / a_{ij} \in Z_7\}$ be a S-ring II? Does it have S-subring II which are not S-ideals II?
10. Find all S-subrings II and S-ideal II for the mixed direct product $R = Z_7 \times Z_8 \times Z_{12}$.

## 4.3 Some New Smarandache elements and their properties

This section is completely devoted to the study of properties of elements in rings and their Smarandache analogue. Several properties introduced on elements of a ring are not found in any ring theory texts but only in research papers published by journals. So some of these concepts, which may be really good, do not find adequate importance among researchers in ring theory. So this section uses these and gives a Smarandache analogue so that not only a Smarandache algebraist but any researcher/student in ring theory can find it useful and interesting. We define in rings new notions like super idempotent, shrinkable element, dispotent element, super-related elements, magnifying element, friendly and non-friendly shrinkable and magnifying element, n-like ring and the Smarandache analogue of them.

**DEFINITION 4.3.1**: *Let $R$ be a ring. An element $0 \neq \alpha \in R$ is called a super idempotent of $R$, if $\alpha^2 - \alpha$ is an idempotent of $R$.*

**THEOREM 4.3.1**: *If a ring $R$ has nontrivial super idempotents then it has nontrivial idempotents.*

*Proof*: By the very definition the result follows.

***Example 4.3.1***: Let $Z_2 = \{0, 1\}$ and $G = \langle g / g^3 = 1 \rangle$. $Z_2G$ be the group ring, $(1 + g^2) \in Z_2G$ is a super idempotent for $(1 + g^2)^2 - (1 + g^2) = 1 + g + 1 + g^2 = (g + g^2)$. Now $(g + g^2)^2 = g^2 + g$, hence $1 + g^2$ is a super idempotent which is not an idempotent.

**THEOREM 4.3.2**: *Let $R$ be a ring. Every super idempotent in general need not be an idempotent of $R$.*

*Proof*: By an example; in example 4.3.1, $1 + g^2$ is a super idempotent which is not an idempotent as $(1 + g^2)^2 = 1 + g$.

**THEOREM 4.3.3**: *Let $R$ be a ring. An element $0 \neq \alpha \in R$ is a nontrivial super idempotent if and only if either $\alpha(\alpha^3 - 2\alpha^2 + 1) = 0$ or $\alpha^3 - 2\alpha^2 + 1 = 0$.*

*Proof*: Given $0 \neq \alpha \in R$ is a nontrivial super idempotent of R. So $[(\alpha^2 - \alpha)]^2 = \alpha^2 - \alpha$ i.e. $\alpha^4 - 2\alpha^3 + \alpha^2 - \alpha^2 + \alpha = 0$ i.e. $\alpha^4 - 2\alpha^3 + \alpha = 0$ i.e. $\alpha(\alpha^3 - 2\alpha^2 + 1) = 0$ or $\alpha^3 - 2\alpha^2 + 1 = 0$. Hence the claim.

Conversely if for some $\alpha \neq 0$ in R we have $\alpha(\alpha^3 - 2\alpha^2 + 1) = 0$ or $\alpha^3 - 2\alpha^2 + 1 = 0$. We get $\alpha^4 - 2\alpha^3 + \alpha = 0$ add to this $\alpha^2$ on both sides so that $\alpha^4 - 2\alpha^3 + \alpha + \alpha^2 = \alpha^2$ i.e. $\alpha^4 - 2\alpha^3 + \alpha^2 = \alpha^2 - \alpha$ i.e. $(\alpha^2 - \alpha)^2 = \alpha^2 - \alpha$. Hence the claim.

**THEOREM 4.3.4**: *Let $R$ be a ring; $\alpha \neq 0$ in R be a super idempotent then either $\alpha$ or $\alpha - 1$ is a zero divisor or $\alpha(\alpha - 1) = 1$ is a unit in R.*

*Proof*: From the above theorem we have $\alpha(\alpha^3 - 2\alpha^2 + 1) = 0$ so $\alpha$ is a zero divisor. If $\alpha$ is not a zero divisor then we have $\alpha^3 - 2\alpha^2 + 1 = 0$ that is $\alpha^3 - \alpha^2 + 1 - \alpha^2 = 0$.

$$\alpha^2(\alpha - 1) - (\alpha^2 - 1) = 0$$
$$\alpha^2(\alpha - 1) - (\alpha - 1)(\alpha + 1) = 0$$

i.e. $(\alpha - 1)[\alpha^2 - (\alpha + 1)] = 0$, so either $\alpha - 1$ is a zero divisor or $\alpha^2 - (\alpha + 1) = 0$. If $\alpha^2 - \alpha - 1 = 0$ then we have $\alpha(\alpha - 1) = 1$. Hence the claim.

**THEOREM 4.3.5**: *Let G be torsion free abelian group and R any field. The group ring KG has no nontrivial super idempotents.*

*Proof*: We know by theorem 4.3.4, if KG has super idempotents then it has nontrivial zero divisors or units, but KG is a domain. Hence KG has no super idempotents.

**DEFINITION 4.3.2**: *Let R be a ring $0 \neq \alpha \in R$ is called a Smarandache super idempotent (S-super idempotent) if $\alpha^2 - \alpha$ is a S-idempotent of R.*

Thus we see superidempotents guarantees the existence of zero divisors or units. Obtain analogous results for S-super idempotents as S-idempotents are introduced and studied in chapter 3.

**Example 4.3.2**: Let $Z_{12}$ be the ring of integers modulo 12. $5 \in Z_{12}$ is such that $5^2 - 5 = 25 + 20 = 45 = 9$. Now $(5^2 - 5)^2 \equiv 5^2 - 5 \pmod{12}$, here 5 is unit of $Z_{12}$. 5 is a super idempotent of $Z_{12}$ which is also a unit. All units of $Z_{12}$ are not super idempotents for 7 is a unit but 7 is not a super idempotent of $Z_{12}$.

Now we proceed onto define a new relation in rings called superrelated elements of a ring and also we define Smarandache superrelated elements of a ring. Such relations bring in a interrelation between elements in a ring.

**DEFINITION 4.3.3**: *Let R be a ring. An element $a \in R$ is said to be weakly superrelated if there exists atleast three distinct elements b, c, d in R such that $(a + b)(a + c)(a + d) = a + bc(a + d) + cd(a + b) + bd(a + c)$.*

**Example 4.3.3**: Let $Z_3 = \{0, 1, 2\}$ and $G = \langle g / g^2 = 1 \rangle$. Consider the group ring $Z_3G$. $2 + 2g \in Z_3G$ is weakly superrelated element of $Z_3G$. For take b = 1, c = g + 1 and d = 1 + 2g. We see $2 + 2g$ satisfies the condition for it to be superrelated.

**DEFINITION 4.3.4**: *Let R be a ring. An element $x \in R$ is said to be a superrelated element of R if $(x + a)(x + b)(x + c) = x + bc(x + a) + ab(x + c) + ac(x + b)$ for all a, b, c $\in$ R.*

***Example 4.3.4***: Let $Z_2 = (0, 1)$ be the prime field of characteristic two, G any group. The element 0 of $Z_2G$ is a superrelated element. For $(0 + a)(0 + b)(0 + c) =$ abc, $0 + abc + abc + abc = abc$ as characteristic of $Z_2G$ is two.

**THEOREM 4.3.6**: *Let R be a ring of characteristic two then 0 is a superrelated element of R.*

*Proof*: Left for the reader to prove.

**THEOREM 4.3.7**: *Every superrelated element of R is a weakly superrelated element of R but every weakly superrelated element of R in general need not be a superrelated element of R.*

*Proof*: Follows from the very definition of superrelated element and weakly superrelated element of R.

The reader is requested to prove the converse by giving examples.

**DEFINITION 4.3.5**: *Let R be a ring. R is said to be weakly superrelated ring if every element of R is a weakly superrelated element of R.*

**DEFINITION 4.3.6**: *Let R be a ring. R is said to be a superrelated ring if every element of R is a superrelated element of R.*

**THEOREM 4.3.8**: *Every superrelated ring is a weakly superrelated ring.*

*Proof*: Obvious.

**THEOREM 4.3.9**: *Let ZG be the group ring. ZG is not a weakly superrelated ring.*

*Proof*: $0 \in$ ZG is such that 0 cannot be weakly superrelated as the identity becomes $\alpha\beta\gamma = 3\alpha\beta\gamma$.

Now we proceed to define the Smarandache analogue.

**DEFINITION 4.3.7**: *Let R be a ring an element x in R is said to be a Smarandache weakly superrelated (S-weakly superrelated) in R if there exists $\alpha$, $\beta$, $\gamma \in A$ such that $(x + \alpha)(x + \beta)(x + \gamma) = x + \alpha\beta(x + \gamma) + \alpha\gamma(x + \beta) + \beta\gamma(x + \alpha)$ where A is a S-subring of R. Note if R has no S-subring but R is a S-ring then we say x in R is Smarandache weakly superrelated in R.*

**DEFINITION 4.3.8**: *Let R be a ring. An element x in R is said to be a Smarandache superrelated in R if for all $\alpha$, $\beta$, $\gamma \in A$, where A is a S-subring such* *that $(x + \alpha)(x + \beta)\ (x + \gamma) = x + \alpha\beta(x + \gamma)\ \alpha\gamma(x + \beta) + \beta\gamma(x + \alpha)$. If R has no S-subring but R is a S-ring then we say $x \in R$ is a Smarandache superrelated in R.*

**THEOREM 4.3.10**: *If R is a superrelated ring and if R is a S-ring then R is a S-superrelated ring.*

*Proof*: Easily proved by using definitions and properties of superrelated elements.

**THEOREM 4.3.11**: *If R is a superrelated ring and if R has a S-subring then R is S-superrelated ring.*

*Proof*: Straight forward.

**DEFINITION [28]**: *A ring R is said to be bisimple if it has more than one element and satisfies the following conditions:*

1. *For any $a \in R$ we have $a \in aR \cap Ra$.*
2. *For any non-zero $a, b \in R$ there is some $c \in R$ such that $aR = cR$ and $Rc = Rb$.*

For more about bisimple rings please refer [28].

***Example 4.3.5***: Let $G = \langle g\ /\ g^2 = 1 \rangle$ and $Z_2 = \{0, 1\}$ be the prime field of characteristic two. $Z_2G = \{0, 1, g, 1 + g\}$ is the group ring. For any g, $1 + g \in Z_2G$ there is no c in $Z_2G$ such that $Z_2Gc = Z_2G(1 + g)$, $cZ_2G = gZ_2g$. Thus $Z_2G$ is not bisimple, but for $1, g \in Z_2G$ we have no c in $Z_2G$ such that

$$c \cdot Z_2G = g \cdot Z_2G \text{ and}$$
$$Z_2G \cdot 1 = Z_2G \cdot c$$

**THEOREM 4.3.12**: *Let $G = \langle g/g^n = 1 \rangle$ and $Z_2 = \{0, 1\}$. The group ring $Z_2G$ is not bisimple.*

*Proof*: Left for the reader to prove.

**DEFINITION 4.3.9**: *Let R be a ring, we call R semi bisimple if for any $a, b \in R$ we have $c \in R$ such that $aR = cR$ and $Rb = Rc$.*

**THEOREM 4.3.13**: *Let R be a zero square ring. R is semi bisimple.*

*Proof*: Obvious as in R, $ab = 0$ for all $a, b \in R$.

**DEFINITION 4.3.10**: *Let R be a ring. R is said to be weakly bisimple if for every $a \in R$, $a \in aR \cap Ra$ and for every pair of elements $a, b \in R$ $aR \subset cR$ and $Rb \subset Rc$.*

**THEOREM 4.3.14**: *A weakly bisimple ring need not in general be bisimple.*

*Proof*: Consider the ring $Z_6 = \{0, 1, 2, 3, 4, 5\}$. Clearly for every $a \in Z_6$ as $1 \in Z_6$, $\{(1, 2); c = 5\}$, …, $\{(2, 4); c = 5\}$; it is easily verified $Z_6$ is weakly bisimple.

**_Remark_**: A ring without 1 can be weakly bisimple.

**_Example 4.3.6_**: Let $P = \{0, 2, 4, 6\}$ modulo 8. Clearly P has no unit it is easily verified P is weakly bisimple and not bisimple.

**DEFINITION 4.3.11**: *Let R be a ring, we say R is Smarandache bisimple (S-bisimple) if it has more than one element and satisfies the following conditions:*

1. *For any $a \in A$ we have $a \in aA \cap Aa$ where A is a S-subring.*
2. *For any non-zero $a, b \in A$ (A a S-subring) there is some $c \in A$ such that $aA = cA$ and $Ac = Ab$.*

**DEFINITION 4.3.12**: *Let R be a ring, we call R a Smarandache semi bisimple (S-semi bisimple) if for any $a, b \in A$ where A is a S-subring, we have $c \in A$ such that $aA = cA$ and $Ab = Ac$.*

**DEFINITION 4.3.13**: *Let R be a ring not necessarily commutative. R is said to be a Smarandache weakly bisimple (S-weakly bisimple) if for every $a \in A$, A a S-subring of R we have $a \in aA \cap Aa$ and for every pair of elements $a, b \in A$, $aA \subset cA$ and $Ab \subset Ac$ for some $c \in A$.*

All properties parallel to bisimple rings can also be studied and obtained with modification for S-bisimple rings, S-weakly bisimple rings and S-semi bisimple rings.

Now we proceed onto define trisimple ring and S-trisimple rings.

**DEFINITION 4.3.14**: *Let R be a ring. R is said to be trisimple if R has more than one element and satisfies the following conditions:*

1. *For any $a \in R$, $a \in aR \cap Ra \cap aRa$.*
2. *For any non-zero $a, b \in R$ there is some $c \in R$ such that $aR = cR$ and $Rc = Rb$.*

**THEOREM 4.3.15**: *Let R be a commutative ring with 1. If a ∈ R is such that $a^2$ = 0 then R is not trisimple.*

*Proof*: Left for the reader to verify.

**THEOREM 4.3.16**: *Let $Z_{12}$ = {0, 1} and $S_n$ the permutation group of degree n. The group ring $Z_2 S_n$ is not trisimple.*

*Proof*: Take $1 + p \in Z_2 S_n$ where

$$p = \begin{pmatrix} 1 & 2 & 3 & 4 & . & . & . & n \\ 2 & 1 & 3 & 4 & . & . & . & n \end{pmatrix} \in S_n$$

$p^2 = 1$ and $(1 + p)^2 = 0$ so $Z_2 S_n$ is not trisimple.

Now we define semi trisimple rings.

**DEFINITION 4.3.15**: *Let R be a ring. R is said to be a semi trisimple if for any a ∈ R; a ∈ Ra ∩ aR ∩ aRa the second condition need not be true.*

**DEFINITION 4.3.16**: *Let R be a ring. If for any a ∈ A where A is a S-subring, we have a ∈ aA ∩ Aa ∩ aAa and for any non-zero a, b ∈ A there is some c ∈ R such that aA = cA and Ac = Ab then we call R a Smarandache trisimple ring (S-trisimple ring).*

**DEFINITION 4.3.17**: *Let R be a ring. R is said to be Smarandache semi trisimple (S-semi trisimple) if R has a S-subring A such that for any a ∈ A we have a ∈ aA ∩ Aa ∩ aAa.*

**THEOREM 4.3.17**: *The group ring $Z_2 S_3$ is not S-trisimple ring.*

*Proof*: Follows easily.

**Example 4.3.7**: The ring of integers Z is not S-semi trisimple.

Now we proceed on to define n-like rings.

**DEFINITION [34]**: *A ring R is called a generalized n-like ring if R satisfies: $(xy)^n - xy^n - x^n y + xy = 0$; for all x, y ∈ R. If characteristic of R = n, R is called a n-like ring.*

**Example 4.3.8**: The group ring $Z_2 G$ where G = ⟨g/ $g^2$ = 1⟩ is a 2-like ring.

**Example 4.3.9**: The group ring $Z_2 S_3$ is not a n-like ring.

**THEOREM 4.3.18**: *The group ring KG of a torsion free abelian group G over any field K is not a n-like ring.*

*Proof*: Take g, h ∈ G ⊂ KG clearly if $(gh)^n - g^n h - gh^n + gh = 0$ implies gh $[(gh)^{n-1} - g^{n-1} - h^{n-1} + 1] = 0$ since gh $\neq$ 0 we see $g^{n-1} h^{n-1} - g^{n-1} - h^{n-1} + 1 = 0$; i.e. $g^{n-1} [h^{n-1} - 1] - [h^{n-1} - 1] = 0$ that is $(h^{n-1} - 1)(g^{n-1} - 1) = 0$ which is impossible as G is torsion free. Hence the claim.

**THEOREM 4.3.19**: *Let $Z_2 = \{0\}$ and $G = \langle g \,/\, g^n = 1 \rangle$. The group ring $Z_2 G$ is a (n −1)-like ring.*

*Proof*: Left as an exercise for the reader to verify.

**DEFINITION 4.3.18**: *Let R be a ring we say R is a Smarandache n-like ring (S-n-like ring) if R has a proper S-subring. A of R such that $(xy)^n - xy^n - x^n y + xy = 0$ for all x, y ∈ A.*

**THEOREM 4.3.20**: *R is a n-like ring with a S-subring $A \subsetneq R$ $(A \neq R)$ then R is a S-n-like ring.*

*Proof*: Follows from the very definition.

Construct an example of a S-n- like ring, which is not a n -like ring.

Now we proceed on to define triple identity in rings which is analogous to the identity; $x^n + y^n = z^n$ for x, y, z integers; the famous last theorem of Fermat.

**DEFINITION 4.3.19**: *Let R be a ring. If there exists a triple $\{\upsilon, \nu, \omega\} \in R \setminus \{0\}$ such that $\upsilon, \nu$ and $\omega$ are distinct elements of $R \setminus \{0\}$ which satisfy the identity $\nu^n + \omega^n = \upsilon^n$ $(n > 1)$, we call R a triple identity ring or TI- ring.*

**Example 4.3.10**: Let $Z_6 = \{0, 1, 2, \ldots, 5\}$ be the ring of integers modulo 6. It is easily verified $Z_6$ is a TI- ring.

**Example 4.3.11**: The ring $Z_7$ is a TI -ring for $2^4 + 5^4 = 4^4$, $3^4 = 2^4 + 5^4$ and so on.

**Example 4.3.12**: The group ring $Z_2 S_3$ is a TI-ring. For $p_3^2 + (p_4 + p_5)^2 = (1 + p_4 + p_5)^2$ and $p_1^2 + (p_4 + p_5)^2 = (1 + p_4 + p_5)^2$.

**DEFINITION 4.3.20**: *Let R be a ring we say R is a Smarandache TI-ring (S-TI-ring) if R has a S-subring, A and in A we have 3 distinct elements, $v$, $v$, $\omega$ such that $v^n + \omega^n = v^n$.*

The reader is requested to prove the following theorems:

**THEOREM 4.3.21**: *If R is a S-TI-ring then R is a TI-ring.*

**THEOREM 4.3.22**: $Z_2 S_3$ *is a S-TI-ring.*

Now we proceed onto define the concept of power joined ring and their Smarandache analogue.

Now in case of integers $(a, b) = 1$, $a^n = b^m$ is impossible but we have such identities to be true in ring so we proceed on to define such related elements to be power joined elements.

**DEFINITION 4.3.21**: *Let R be a ring. If for every $a \in R$ there exists atleast one b $\in R$, $(b \neq a)$ such that $a^m = b^n$ for some positive integers m and n then we say R is a power joined ring.*

***Example 4.3.13***: Let $Z_5 = \{0, 1, 2, 3, 4\}$ be the ring of integers modulo 5. $Z_5$ is a power joined ring.

***Example 4.3.14***: Let $Z_2 = \{0, 1\}$ and $G = \langle g/g^3 = 1 \rangle$. The group ring $Z_2 G$ is not power joined as $1 + g + g^2 \in Z_2 G$ cannot be represented as a power of some other element as $(1 + g + g^2)^2 = 1 + g + g^2$. But this does not imply no idempotents can be expressible as power joined elements.

***Example 4.3.15***: Let $Z_3 G$ be the group ring where $G = \langle g / g^2 = 1 \rangle$. Clearly $(2 + 2g)^2 = 2 + 2g$, but we have $2 + 2g = (1 + g)^2$. Hence the claim.

**DEFINITION 4.3.22**: *Let R be a ring if for every $x \in R$ there exists some $y \in R$ such that $x^m = y^n$ for some integers m and n and if the integers m and n happen to be the same for all $x$, $y \in R$ we say R is a (m, n)–power joined ring.*

**THEOREM 4.3.23**: *A Boolean ring can never be a power joined ring.*

*Proof*: Left as an exercise.

**DEFINITION 4.3.23**: *Let R be a ring in which we have for every x ∈ R there exists y ∈ R such that $x^m = y^m$ (x ≠ y) and m ≥ 2. Then we say R is a uniformly power joined ring.*

**THEOREM 4.3.24**: *Every prime field K = $Z_p$ of characteristic p, p ≠ 2 is a power joined ring.*

*Proof*: Obvious by the definition of K = $Z_p$.

Now we proceed onto define Smarandache analogue.

**DEFINITION 4.3.24**: *Let R be a ring. If for every a ∈ A ⊂ R where A is S-subring there exists b ∈ A such that $a^m = b^n$ for some positive integers m and n then we say R is a Smarandache power joined ring (S-power joined ring).*

**DEFINITION 4.3.25**: *Let R be a ring. If for every x ∈ A ⊂ R there exist some y ∈ A ⊂ R such that $x^m = y^n$ where A is a S-subring of R then we say the ring R is a Smarandache-(m, n)-power joined ring (S-(m, n)-power joined ring).*

**DEFINITION 4.3.26**: *Let R be a ring if for every x ∈ A ⊂ R where A is a S-subring there exists some y ∈ A ⊂ R such that $x^m = y^m$ (x ≠ y) and m ≥ 2. Then we say R is a Smarandache uniformly power joined ring (S-uniformly power joined ring).*

**THEOREM 4.3.25**: *If R is a ring which is power joined, then R is a S-power joined ring only if for x ∈ A ⊂ R, y also belongs to A. If y ∉ A then R cannot be a S-power joined ring.*

*Proof*: We seek the proof to be supplied by the reader.

Thus we see the Smarandache notions in this case has made R a locally power joined ring.

Now we discuss about the types of commutativity right commutativity and quasi commutativity and obtain the Smarandache analogue.

**DEFINITION 4.3.27**: *Let R be a ring. If for every pair of elements a, b in R we have ab = $(ba)^r$, r ≥ 1 then we say R is conditionally commutative. If r = 1 for all a, b ∈ R, then R is obviously commutative. Thus every commutative ring R is conditionally commutative.*

*If in a ring we have a pair of elements a, b ∈ R such that ab = $(ba)^r$; r ≥ 1 we say the pair is conditionally commutative.*

We can define a group G to be conditionally commutative if $xy = (yx)^n$, $n \geq 1$ for every $x, y \in G$.

By our Smarandache notions we will localize the property.

**DEFINITION 4.3.28**: *Let R be a ring. We say R is a Smarandache conditionally commutative ring (S-conditionally commutative ring) if for every $x, y \in A$ where A is a S-subring of R, we have $xy = (yx)^n$ for $n \geq 1$.*

**THEOREM 4.3.26**: *If R is a conditionally commutative ring having a S-subring then R is a S-conditionally commutative ring.*

*Proof*: Follows from the very definitions.

However it is left for the reader to construct an example of a S-conditionally commutative ring which is not a conditionally commutative ring. Yet another interesting result is semi right commutativity of rings which leads to give conditions for the existence of zero divisors.

**DEFINITION 4.3.29**: *Let R be a ring. R is said to be a strongly semi right commutative ring if for every triple of elements x, y, z we have atleast one of the following three equalities to be true.*

> *i. $xy = zyx$ (or $yx = zxy$).*
> *ii. $yz = xzy$ (or $zy = xyz$).*
> *iii. $zx = yxz$ (or $xz = yzx$).*

*We can define a strongly semi right commutative triple x, y, z $\in$ R if atleast one of the following three equalities is true.*

> *i. $xy = yxz$ (or $yx = xyz$) or*
> *ii. $yz = zyx$ (or $zy = yzx$) or*
> *iii. $zx = xzy$ (or $xz = zxy$).*

*Similarly we can define strongly semi left commutative ring in a similar way.*

**THEOREM 4.3.27**: *No commutative ring without divisors of zero is strongly semi right commutative.*

*Proof*: Obvious from the very definition, for if $x \neq 1$ or $0$, $y \neq 0$ or $1$, $z \neq 0$ or $1$ where R is commutative. If R is strongly semi right commutative then we have $xy = zyx$ so that $xy = zyx = z(xy)$ as $xy = yx$ we have $(1 - z) xy = 0$ $x \neq 0$, $y \neq 0$ and $z \neq 0$. So $(1 - z) xy = 0$ must be a zero divisor.

**DEFINITION 4.3.30**: *Let R be a ring. R is said to be a Smarandache strongly semi right commutative ring (S-strongly semi right commutative ring) if for every triple of elements x, y, z in A, A a S-subring of R we have atleast one of the following three equalities to be true.*

> *i.    xy  = zyx (or yx = zxy) or*
> *ii.   yz = xzy (or zy  =  xyz) or*
> *iii.  zx = yxz (or xz  =  yzx).*

*Similarly we define Smarandache strongly semi left commutative ring (S-strongly semi left commutative triple) if for every triple x, y, z ∈ A; A a S-subring of R if at least one of the following three equalities is true.*

> *i.    xy  = yxz (or yx  =  xyz) or*
> *ii.   yz = zyx (or zy  =  yzx) or*
> *iii.  zx = xzy (or xz  =  zxy).*

*Finally we define a Smarandache strongly semi right (left) commutative triple (S-strongly semi right (left) commutative triple) only when the triple x, y, z satisfies the above conditions the elements must be only from a proper S-subring A of R.*

Now we proceed on to define right commutativity in rings.

**DEFINITION 4.3.31**: *Let R be a ring. R is said to be strongly right commutative if a(xy) = a(yx) for all a, x, y ∈R.*

*Similarly we define a ring R to be strongly left commutative if (xy)a = (yx)a. for all a, x, y ∈R.*

**THEOREM 4.3.28**: *Every strongly right or left non-commutative ring has nontrivial divisors of zero.*

*Proof*: From the definitions we have for a, x, y ∈ R; xya = (yx)a in both case this implies (xy − yx) a = 0 or a (xy − yx) = 0 as xy ≠ yx and a ≠ 0 we have non-trivial zero divisors in R.

**THEOREM 4.3.29**: *A group ring KG of a group G over any field K can never be a strongly right (left) commutative ring.*

*Proof*: If a, x, y ∈ G ⊂ KG, then axy = ayx forces xy = yx.

Thus it is important to note that this property can only be defined for rings and never for groups. Now we proceed onto define Smarandache analogue.

**DEFINITION 4.3.32**: *Let R be a ring we say R is a Smarandache strongly right commutative (S-strongly right commutative) ring if a(xy) = a(yx) for all a, x, y ∈ A where A is a proper S-subring of R. Similarly we define Smarandache strongly left commutative elements.*

The goodness about the Smarandache structures is that we saw no group rings can be strongly right (left) commutative, but we see the following theorem:

**THEOREM 4.3.30**: *Let R be a S-ring and G any group. RG is a S-strongly right (left) commutative ring provided the following holds good:*

> *i.    R is a right (left) commutative ring.*
> *ii.   R is a S-strongly right (left) commutative ring.*

*Proof*: Clearly by the very definition of Smarandache strongly right (left) commutative ring, we get the theorem to be true under the given conditions.

We see in rings, we can have centre but not notions analogues to commutator in groups; here we proceed onto define a new concept called quasi semi commutator and the quasi semi commutative element.

**DEFINITION 4.3.33**: *Let R be a ring. An element x ∈ R is said to be quasi semi commutative if there exists y ∈ R (y ≠ 0) such that (xy − yx) commutes with every element of x. Trivially if y = 0 then we have xy − yx = 0 which commutes with every element of R.*

**DEFINITION 4.3.34**: *Let R be a ring. For a quasi-semi commutative element x of R we define the quasi semi-commutator to be the set of all p ∈ R such that xp − px commutes with every element of R and denote it by Q(x) i.e. Q(x) = {p ∈ R / xp − px commutes with every element of R}. Clearly Q(x) ≠ φ for 0, 1 ∈ Q(x), if R is a ring with 1.*

**DEFINITION 4.3.35**: *Let R be a ring. R is said to be a quasi semi commutative ring if every element in R is quasi semi commutative. Every commutative ring is obviously quasi semi commutative.*

**DEFINITION 4.3.36**: *Let R be a ring. The quasi semi center of R denoted by Q(x) = {x ∈ R / xp − px is quasi semi commutative}; clearly Q(R) ≠ φ.*

**THEOREM 4.3.31**: *Let R be a non-commutative ring. Z(R) denote the center of R. Then we have Z(R) ⊂ Q(R).*

*Proof*: Clearly Z(R) = {x ∈ R / xy = yx for all y ∈ R} Now Q(R) = {x ∈ R / xy − yx commutes with every element of R}. So Z(R) ⊂ Q(R) as xy − yx = 0 for all x ∈ Z(R).

**DEFINITION 4.3.37**: *Let R be a ring. An element x ∈ A ⊆ R where A is a S-subring of R is said to be a Smarandache quasi semi commutative (S-quasi semi commutative) if there exists y ∈ A (y ≠ 0) such that xy − yx commutes with every element of A.*

**DEFINITION 4.3.38**: *Let R be a ring, x ∈ A ⊂ R (A a S-subring of R) is a Smarandache quasi semi commutative element (S-quasi semi commutative element) of R if x is quasi semi commutative for some y ∈ A. The Smarandache semi commutator (S-semi commutator) of x denoted by SQ(x) = {p ∈ A / xp − px commutes with every element of A}. R is said to be a Smarandache quasi semi commutative ring (S-quasi semi-commutative ring) if for every element in A (A ⊂R) (A a S-subring) is a S-quasi semi commutative.*

**DEFINITION 4.3.39**: *Let R be a ring. The S-quasi semi center (S-quasi semi center) of R denoted by SQ(R) = {x ∈ A / xp − px is S-quasi semi commutative}.*

The reader is requested to derive interesting results about these concepts.

The concept of magnifying and shrinking elements in a ring is an interesting feature. However the notion of magnifying elements was introduced to semigroups by researchers. We introduce them to rings.

**DEFINITION 4.3.40**: *Let R be a ring. v is called left magnifying element of R (v need not be in R) if for some proper subset M of R we have vM = R.*

*Similarly we define right magnifying element of R. In case of commutative rings the notion of right and left magnifying elements coincide. Even if R is a non-commutative ring we may have vM = Mv = R that is v may serve as a magnifying element.*

*If v is in R we say v is a friendly magnifying element of R; if v ∉ R still vM = Mv = R for some proper subset M in R then we say v is a non-friendly magnifying element of R. The concept of friendly and non-friendly magnifying elements plays a vital role only when we define the Smarandache notions of them.*

***Example 4.3.16***: Let Z be the ring of integers Let P = {0, ±2, ±4, …}. Clearly P is a proper subset of Z. Take v = 1/2 clearly v ∉ Z but v . P = P . v = Z so v is a non-friendly magnifying element of Z.

Now the nontrivial question is why should one study only magnifying elements so we introduce the concept of shrinking elements of a ring.

**DEFINITION 4.3.41**: *Let R be a ring. An element x of R is called a shrinking element of R if xR = P where P is a proper subset of R. If x ∈ R we say x is a friendly shrinking element; otherwise we say x is a non-friendly shrinking element of R. The concept of shrinking element is, in a way, just the opposite of magnifying elements.*

Here also the concept of right shrinking, left shrinking and shrinking can be defined as in the case of magnifying elements.

**DEFINITION 4.3.42**: *If in a ring, if every element other than unity shrinks R, we call the ring R as a shrinkable ring (i.e. xR ≠ R for x ∈ R).*

**THEOREM 4.3.32**: *A field has no shrinkable elements other than {0}.*

*Proof*: Obvious by the very definition.

**THEOREM 4.3.33**: *Let KG be the group ring of the group G over the field K. G a finite group. The group ring has shrinkable elements.*

*Proof*: Take $\alpha = (1 + g_1 + \ldots + g_n)$ where $\{1, g_1, \ldots, g_n\} = G$ then $\alpha KG \neq KG$; hence KG has shrinkable elements.

Now we localize this property.

**DEFINITION 4.3.43**: *Let R be a ring. A ⊂ R be a proper S-subring of R. An element v is called Smarandache left magnifying (S-left magnifying) element of R if vM = A for some proper subset M of A, we say v is Smarandache right magnifying (S-right magnifying) if $M_1 v = A$ for some proper subset $M_1$ of A. v is said to be Smarandache magnifying (S-magnifying) if vM = Mv = A for some M a proper subset of A. If v ∈ A, then v is said to be a Smarandache friendly magnifying (S-friendly magnifying) element. If v ∉ A we call v a Smarandache non-friendly magnifying (S-non-friendly magnifying) element even if v ∈ R \ A we still call v a S-non-friendly magnifying element.*

**DEFINITION 4.3.44**: *Let R be a ring. An element x is called a Smarandache left shrinking element (S-left shrinking element) of R if for some S-subring A of R we have proper subset M of R such that xA = M (M ≠ A) or M ≠ R. We define similarly Smarandache right shrinking (S-right shrinking) and Smarandache shrinking (S-shrinking) if xA = Ax = M.*

*If $x \in A$ we call $x$ a Smarandache friendly shrinking element (S-friendly shrinking element); if $x \notin A$ we call $x$ a Smarandache non-friendly shrinking element (S-non-friendly shrinking element).*

Obtain analogues and interesting results about S-shrinking and magnifying elements of a ring R.

Finally we conclude this section by just defining some new concepts viz. semiunit, dispotent elements of a ring and a dispotent ring.

**DEFINITION 4.3.45**: *Let R be a commutative ring with unit 1. An element x of R is said to be a semiunit of R if there exists $y \in R$ such that $(x + 1)(y + 1) = 1$.*

This method can make even zero divisors and idempotents into semiunits hence the study of them is important or to be more specific it makes nilpotent elements into semiunits. For example consider the ring $Z_{12}$.

**Example 4.3.17**: Let $Z_{12}$ be the ring of integers modulo 12. 6 is a zero divisor but 6 is also a semiunit of $Z_{12}$ for $(6 + 1)(6 + 1) \equiv 1 \pmod{12}$.

**THEOREM 4.3.34**: *Let R be a ring. An element x is a semiunit if and only if there exists $y \in R$ with $x + y + xy = 0$, $y \neq 0$.*

*Proof*: $(x + 1)(y + 1) = 1$ forces $xy + x + y = 0$. Now if $xy + x + y = 0$ then we have $x + y + xy + 1 = 1$ forcing $(x + 1)(y + 1) = 1$. Hence the claim.

In case of rings, which are non-commutative, we can also define right semiunit and left semiunit and obtain similar characterizations about them.

**Example 4.3.18**: Let $Z_6 = \{0, 1, 2, 3, 4, 5\}$, $4 \in Z_6$ is an idempotent of R. But 4 is also a semiunit as $(4 + 1)(4 + 1) \equiv 1 \pmod 6$.

Thus we see nilpotents, zero divisors and idempotents can be semiunits of R.

**THEOREM 4.3.35**: *Let K be a field of characteristic 0. Every element is a semiunit.*

*Proof*: Let $x \in K$. Consider

$$(x + 1)\left[\frac{-x}{x + 1} + 1\right] = (x + 1)\frac{1}{x + 1} = 1.$$

Hence the claim.

**DEFINITION 4.3.46**: *Let R be a ring; $x \in R$ is said to be a Smarandache semiunit (S-semiunit) if there exists $y \in R$ such that $x + 1$ and $y + 1$ are S-unit of R.*

The reader is advised to develop interesting results as the notion of S-units are dealt in an entire section in chapter 3 of this book.

[72] had defined the concept of dispotent semigroups. Here we define dispotent rings and their Smarandache analogue.

**DEFINITION [72]**: *A semigroup S is a dispotent semigroup if and only if it has exactly two idempotents.*

**DEFINITION 4.3.47**: *Let R be a ring. R is said to be a dispotent ring if R has exactly two nontrivial idempotents.*

**Example 4.3.19**: Let $Z_2G$ be the group ring where $G = \langle g \,/\, g^2 = 1 \rangle$. This group ring has only two idempotents viz., $1 + g + g^2$ and $g + g^2$.

**Example 4.3.20**: $Z_{18}$ is a dispotent ring.

**Example 4.3.21**: $Z_3S_n$ is not a dispotent ring.

**DEFINITION 4.3.48**: *Let R be a ring if R has a proper S-subring A of R such that the S-subring A has only two S-idempotents then we call R a Smarandache dispotent ring (S-dispotent ring).*

The study of S-idempotents has been carried out in a sole section in chapter 3 of this book. The reader is requested to study and get some interesting results.

**DEFINITION 4.3.49**: *Let R be a S-ring. If every S-subring A of R has exactly two S-idempotents then we say R is a Smarandache strong dispotent ring (S-strong dispotent ring).*

Can we obtain any relation between S-dispotent rings and S-strong dispotent rings.

<u>**PROBLEMS:**</u>

1.   Does the ring $Z_{24}$ have super idempotents?
2.   Find whether the group ring $Z_3A_4$ has super idempotents?
3.   Can the ring $Z_{26}$ have S-super idempotents?
4.   Can the semigroup ring $Z_2S(3)$ have S-super idempotents?

5. Prove a S-superrelated ring in general need not be a superrelated ring.

6. Give an example of a S-weakly superrelated ring which is not a weakly superrelated ring.

7. Give an example of a bisimple ring.

8. Can $Z_n$ be a weakly bisimple ring?

9. Give an example of S-bisimple ring.

10. Can $Z_n$ for any suitable n be S-weakly bisimple? Justify.

11. Give an example of
    i. trisimple ring.
    ii. S-trisimple ring.

12. Find a S-semi trisimple ring which is not a S-trisimple ring.

13. Is $Z_5S(3)$ a
    i. S-trisimple?
    ii. S-semi trisimple?
    iii. Trisimple?
    iv. Semi trisimple?

14. Prove $Z_2G$ where $G = \langle g / g^6 = 1 \rangle$ is a 7-like ring.

15. Give an example of a semigroup ring which is a n-like ring.

16. Can ring of matrices with entries from $Z_2$ be a n-like ring for any suitable n?

17. Prove $Z_2S_4$, the group ring, is a TI-ring.

18. Prove the semigroup ring; $Z_2S(3)$ is a
    i. TI-ring.
    ii. Smarandache TI- ring.

19. Give an example of a power joined ring which is not a S-power joined ring.

20. Is $Z_9$ a S-power joined ring? Justify.

21. Can we say $Z_{15}$ is a (m, n) power joined ring or S-(m, n) power joined ring?

22. Is $Z_2S_4$ a S-conditionally commutative ring? Justify.

23. Can $Z_2S(3)$ be a conditionally commutative ring? Prove your answer.

24. Give an example of a S-conditionally commutative ring which is not a conditionally commutative ring.

25. Prove $Z_2S_3$ has atleast a S-semi commutative triple?

26. Can $Z_2S_3$ be a strongly semi commutative ring?

27. Does $Z_3S(4)$ have a
    i. Strongly semi commutative triple?
    ii. S-strongly semi commutative triple?

28. Give an example of a strongly right (or left) commutative ring.

29. Give an example of a S-strongly right (or left) commutative ring which is not a strongly right (or left) commutative ring.

30. Is $Z_2S_3$ a quasi semi commutative?

31. Can $Z_3S(4)$ be S-quasi commutative? Can the group ring $Z_3S_5$ have
    i. S-shrinking elements?
    ii. magnifying elements?

32. Can the group ring QG have

     i. semiunits?

     ii. S-semiunits? (Q field of rationals)?

33. Can $Z_6S(4)$ have

     i. semiunits?

     ii. S-semiunits?

  Find them if they exist.

34. Is $Z_{2n}$, n a prime be a dispotent ring?

35. Can $Z_{22}$ be a S-dispotent ring? Justify or substantiate your claim.

## 4.4 New Smarandache substructure and their properties

Here we introduce the notions of quasi ordering, semi nilpotent, normal elements in a ring, normal ring, G-ring, S-J ring, n-c-s rings, co-rings, iso-rings and their Smarandache analogues leading to several interesting localized properties on the substructures.

**DEFINITION 4.4.1**: *A sum quasi ordering in a ring R is a subset T of R satisfying the condition T + T ⊂ T.*

**DEFINITION 4.4.2**: *A product quasi ordering in a ring R is a subset U of R satisfying the condition U . U ⊂ U.*

**DEFINITION 4.4.3**: *A quasi ordering in a ring R is a subset I of R which is both a sum quasi ordering and a product quasi ordering.*

***Example 4.4.1***: Let $Z_2G$ be the group ring of the group G = $\langle g \, / \, g^2 = 1 \rangle$ over $Z_2$. I = {0, g} is a sum quasi ordering set which is clearly not a product quasi ordering set. J = {0, 1 + g} is both a sum and a product quasi ordering set.

Now we proceed on to define Smarandache analogue.

**DEFINITION 4.4.4**: *Let R be a ring, we say the set T is a Smarandache sum quasi ordering (S-sum quasi ordering) on R.*

  *a. If T has a proper subset P, (P ⊂ T) and P is a semigroup under addition.*

  *b. P + P ⊂ T.*

**DEFINITION 4.4.5**: *Let R be a ring we say a subset U of R is a Smarandache product quasi ordering (S-product quasi ordering) on R if*

a. U contains a proper subset X such that X is a semigroup under multiplication.

b. X . X ⊂ U.

**DEFINITION 4.4.6**: *Let R be a ring. A subset Y of R is said to be a Smarandache quasi ordering (S-quasi ordering) on R if Y is simultaneous a S-sum quasi ordering and a Smarandache product quasi ordering.*

<u>Note</u>: We can have for Y the set, a proper subset $P \subset Y$, P an additive semigroup and $Z \subset Y$ where Z is a multiplicative semigroup and P in general need not be the same as Z.

***Example 4.4.2***: Let $Z_2S_3$ be the group ring of the group $S_3$ over $Z_2$. Take I = {0, $p_1 + p_2 + p_3$, $1 + p_4 + p_5$, $p_5 + p_4 + p_3 + p_2 + p_1 + 1$}, P = {0, $1 + p_4 + p_5$} is a semigroup under addition, P is also a semigroup under multiplication. Clearly $Z_2S_3$ has a S-quasi ordering in it.

***Example 4.4.3***: Let $Z_4S_n$ be the group ring. Take A = {$\Sigma\alpha_is_i / \alpha_i \in$ {0, 2}}. The set A is both S-quasi sum ordering as well S-quasi product ordering. Thus $Z_4S_n$ has a S-quasi ordering on it.

**THEOREM 4.4.1**: *Let $Z_2S_n$ be the group ring. Then $Z_2S_n$ is a S-sum quasi ordering as well as S-product quasi ordering.*

*Proof*: It is left for the reader to verify.

Now we proceed on to define a new concept called Smarandache semi nilpotent ideals.

**DEFINITION [24]**: *An ideal I of R is semi nilpotent if each ring generated by a finite set of elements belonging to the ideal I is nilpotent. An ideal, which is not nilpotent, is called semi regular.*

*Nilpotent ideals are nil.*

**DEFINITION 4.4.7**: *Let R be a ring. An S-ideal I of R is Smarandache semi nilpotent (S-semi nilpotent) if each ring generated by a finite set of elements belonging to the S-ideal which forms a subring A, contained in I is nilpotent.*

**THEOREM 4.4.2**: *Let K be any field and G a torsion free abelian group. KG has no non- zero S-semi nilpotent ideals.*

*Proof*: KG has no zero divisors, hence no nilpotents as semi nilpotent ideals are nil. So KG has no non-zero semi nilpotent ideals.

**DEFINITION 4.4.8**: *Let R be a ring. If R is an ideal which is not Smarandache semi nilpotent then we call the non-S-semi nilpotent ideal to be Smarandache semi regular (S-semi regular).*

**DEFINITION 4.4.9**: *Let R be a ring. M a proper subring of R. I is called a sub semi ideal of R related to M if and only if I is a proper ideal of M. A ring containing a sub semi ideal is called a sub semi ideal ring.*

An analogue to this is defined for Smarandache rings.

**DEFINITION 4.4.10**: *Let R be a ring. M be a S-subring of R. I is called the Smarandache subsemi ideal (S-subsemi ideal) of the ring R related to the S-subring M if and only if I is a proper S-ideal of M and not an S-ideal of R.*

***Example 4.4.4***: Let $Z_2 = \{0, 1\}$ and $S = \{g, h, k, 1 / g^5 = g, k^2 = k, 1.g = g.1 = g, h^3 = h, gh = g = hg, hk = kh = k gk = kg = k\}$ be a semigroup. $Z_2S$ is the semigroup ring. Take $M = \langle g, h, 1 \rangle \subset S$, $Z_2M$ is a S-subring and $Z_2 I$ where $I = \langle g, 1 \rangle$ is an S-ideal of $Z_2M$.

Now we proceed onto define normal elements in a ring, normal ring and obtain a Smarandache analogue.

**DEFINITION 4.4.11**: *Let R be a ring an element $\alpha \in R \setminus \{0, 1\}$ is called a normal element of R if $\alpha R = R\alpha$.*

**DEFINITION 4.4.12**: *Let R be a ring, if $\alpha R = R\alpha$ for every $\alpha \in R$, we say R is a normal ring.*

Now we just recall the definition of normal semigroups [69].

**DEFINITION [69]**: *Let S be a semigroup, if for every $\alpha \in S$ we have $\alpha S = S\alpha$ then S is called a normal semigroup.*

Using this definition we define Smarandache normal semigroup as follows.

**DEFINITION 4.4.13**: *Let S be a S-semigroup with A a proper subset of S which is a group. If $\alpha A = A\alpha$ for all $\alpha \in S$ then S is a Smarandache normal semigroup (S-normal semigroup).*

**DEFINITION 4.4.14**: *Let R be a ring, X be a S-subring of R. We say R is a Smarandache normal ring (S-normal ring) if αX = Xα for all α ∈ R.*

**DEFINITION 4.4.15**: *Let R be a ring; R is said to be a Smarandache strongly normal ring (S-strongly normal ring) if every S-subring X of R is such that αX = Xα for all α ∈ R.*

**THEOREM 4.4.3**: *Let K be a field and S a normal semigroup then KS the semigroup ring is a normal ring.*

*Proof*: Given αS = Sα for all α ∈ S. Hence αKS = KSα for every α ∈ KS thus KS is a normal ring.

**THEOREM 4.4.4**: *Let R be a ring. Z(R) be the nontrivial center of R and if Z(R) is a S-subring then R is a S-normal ring .*

*Proof*: By simple techniques we can obtain the result.

The author has defined the concept of a G-ring.

**DEFINITION 4.4.16**: *Let R be a ring if for every additive subgroup S of R we have rS = Sr = S for every (r ≠ 0) then we call R a G-ring.*

**DEFINITION 4.4.17**: *Let R be a ring. If for every additive subgroup S of R we have rS = Sr for every r ∈ R (r ≠ 0) then we call R a weakly G-ring.*

***Example 4.4.5***: Let $Z_4$ = {0, 1, 2, 3} and S = {0, 2}; now Sr = rS = S thus $Z_4$ is a G-ring.

***Example 4.4.6***: Let $Z_2$ = {0, 1} and G = ⟨g/$g^2$ = 1⟩. The group ring $Z_2G$ is a weakly G- ring for {0, 1, g, g + 1}, {0, 1}, {0, g} and {0, g + 1} are subgroups of $Z_2G$ under addition. Clearly only S = {0, 1 + g} is such that rS = Sr = S for every r ≠ 0. {0, 1} and {0, g} are such that Sr = rS. Thus $Z_2G$ is a weakly G-ring.

**THEOREM 4.4.5**: Let R be a ring. Every G-ring is a weakly G-ring but a weakly G-ring is not a G-ring.

*Proof*: By definition and example 4.4.6 given above.

Now we proceed on to define Smarandache analogue.

**Definition 4.4.18**: *Let R be a ring. If for every S-semigroup, P under addition we have rP = Pr = P for every r ∈ R(r ≠ 0) then we call R a Smarandache G-ring (S-G-ring).*

**Definition 4.4.19**: *Let R be a ring. If for every additive S-semigroup P of R and for every r ∈ R we have rP = Pr then we call R a Smarandache weakly G- ring (S-weakly G-ring).*

**Theorem 4.4.6**: *Every S-G-ring is a S-weakly G-ring and not conversely.*

*Proof*: Left for the reader to verify.

**Example 4.4.7**: Let $M_{3\times3}$ = { $(a_{ij})$ / $a_i \in Z_4$} be the ring of 3×3 matrices with entries from the ring of integers modulo 4. Take

$$P_{3\times3} = \left\{ \begin{pmatrix} a & b & 0 \\ 0 & 0 & 0 \\ 0 & 0 & 0 \end{pmatrix} / a,b \in Z_4 \right\}$$

$P_{3\times3}$ is a S-semigroup under '+' clearly $M_{3\times3}$ is not a Smarandache G-ring. Take

$$P'_{3\times3} = \left\{ \begin{pmatrix} a & 0 & 0 \\ 0 & b & 0 \\ 0 & 0 & 0 \end{pmatrix} / a,b \in Z_4 \right\}$$

$P'_{3\times3}$ is S-Semigroup ring. $M_{3\times3}$ is not a G-ring.

We define a new property in ring called special identity ring or in short SI -ring.

**Definition 4.4.20**: *Let R be a ring. Let S denote the collection of all proper subrings of R. If $(S_1 + S_2)(S_2 + S_3) = S_1(S_2 + S_3) + S_3(S_1 + S_2) + S_2$ for all $S_1$, $S_2$, $S_3 \in S$ we say R is a SI-ring.*

The Smarandache analogue would be

**Definition 4.4.21**: *Let R be a ring. S denote the set of all proper S-subrings of R. If $(S_1 + S_2)(S_2 + S_3) = S_1(S_2 + S_3) + S_2 + S_3(S_1 + S_2)$ for all $S_1$, $S_2$, $S_3 \in S$ then we say R is a Smarandache SI-ring (S-SI-ring).*

**Example 4.4.8**: $Z_{12}$ = {0, 1, 2, …, 11}, the ring of integers modulo 12 is not a SI-ring.

***Example 4.4.9***: Let $Z_2G$ be the group ring, where $Z_2 = \{0, 1\}$ and $G = \langle g \,/\, g^4 = 1 \rangle$. $Z_2G$ is a SI- ring.

***Example 4.4.10***: Let $Z_2 = \{0, 1\}$ and $S = \{1, a, b \,/\, a^2 = a, \ b^2 = b, ab = a, \ ba = b, \ 1.a = a.1 = a, \ 1.b = b.1 = b\}$. It can be checked $Z_2S$, the semigroup ring is a SI-ring.

***Example 4.4.11***: $Z_{12} = \{0, 1, 2, \ldots, 11\}$ is trivially a S-SI-ring as this ring has only one S-subring viz $S = \{0, 2, 4, 6, 8, 10\}$.

***Example 4.4.12***: $Z_2 = \{0, 1\}$ and $G = \langle g \,/\, g^4 = 1 \rangle$, the group ring $Z_2G$ is not a S-SI-ring as it has no S-subrings.

Now we introduce the concept of n-closed additive subgroups in a ring.

**DEFINITION 4.4.22**: *Let R be a ring, if every nonempty additive subgroup A of R is an n-closed additive subgroup of R i.e., $A^n \subset A$ (n > 1) then we say R is a n-closed additive subgroup ring (n-c- s ring).*

***Example 4.4.13***: Let $Z_4 = \{0, 1, 2, 3\}$ be the ring of integers modulo 4. $S = \{0, 2\}$ is an additive subgroup such that $S^2 \subset S$ so $Z_4$ is a n-c-s ring.

***Example 4.4.14***: Let $Z_2 = \{0, 1\}$ and

$$S_3 = \{\begin{pmatrix} 1 & 2 & 3 \\ 1 & 2 & 3 \end{pmatrix} = 1, p_1, p_2, p_3, p_4, p_5\}$$

be the symmetric group of degree 3. $Z_2S_3$ is the group ring which is not a n-c-s ring for $A = \{0, p_1 + p_5\}$ is a group but $A^n \not\subset A$, n > 1.

In view of this we have the following theorem:

**THEOREM 4.4.7**: *Let $Z_2 = \{0, 1\}$ and $S_n$ be the symmetric group of degree n. The group ring $Z_2S_n$ is not a n-c-s ring.*

*Proof*: By taking $S = \{0, p + q\}$ where

$$p = \begin{pmatrix} 1 & 2 & 3 & 4 & . & . & . & n \\ 2 & 3 & 1 & 4 & . & . & . & n \end{pmatrix}$$

and

$$q = \begin{pmatrix} 1 & 2 & 3 & 4 & . & . & . & n \\ 3 & 1 & 2 & 4 & . & . & . & n \end{pmatrix},$$

we have S is an additive subgroup. But $S^r \subset S$ for $r > 1$.

**THEOREM 4.4.8**: $Z_2 S_m$ contains n-c-s subgroups for n = 2, 4, …, m.

*Proof*: Let S = {0, $\alpha$} where $\alpha = 1 + s_i$ where $s_i$ permutes only an even number of elements, $(1 + s_i)^1 = 0$. So $S^i \subseteq S$. If $s_i$ permutes odd number of elements then $(1+s_i)^{1+1} = 1 + s_i$, thus if S = {0, $(1 + s_i)$} we have $S_i^{1+1} \subset S_i$. Hence the claim.

**THEOREM 4.4.9**: *Let $Z_2$ = {0, 1} and G be a torsion free group. No subgroup of the form {0, g / g $\in$ G} is a n-c-s subgroup of $Z_2 G$.*

*Proof*: Follows from the fact G is a torsion free group.

Now we hint at the Smarandache analogue of these definitions.

**DEFINITION 4.4.23**: *Let R be a ring. If for every additive Smarandache semigroup A of R we have $A^n \subset A$ (n > 1) then we say R is a Smarandache n-closed additive subgroup ring. (S-n-closed additive subgroup ring).*

**Example 4.4.15**: Let $Z_2 G$ be the group ring where G = $\langle g / g^6 = 1 \rangle$. Take A = {0, $g^4$, $g^2 + g^4$, $1 + g^4$, $1 + g^2 + g^4$, $1 + g^2$, 1, $g^2$}; yet clearly A is a S-semigroup, $A^n \subset A$ (n > 1). Hence A is S-additive subgroup of $Z_2 G$. It is easily verified $Z_2 G$ is not a S-n-closed additive subgroup ring.

We introduce yet another new concept called co-rings.

**DEFINITION 4.4.24**: *Let R be a ring with identity 1. Two subrings A and B of same order in R is said to be conjugate if there exists some x $\in$ R such that A = $xBx^{-1}$.*

**DEFINITION 4.4.25**: *Let R be a ring with 1. R is said to be a conjugate ring (co-ring) if every distinct pair of subrings of same order are conjugate.*

**DEFINITION 4.4.26**: *Let R be a ring, R is said to be a weak co-ring if there is atleast one pair of distinct subrings of same order which are conjugate to each other.*

**THEOREM 4.4.10**: *Every co-ring is a weak co-ring.*

*Proof*: Obvious by the very definition.

**DEFINITION 4.4.27**: *Let R be a ring. A ring in which every pair of distinct subrings of same order are isomorphic is called an iso ring.*

**DEFINITION 4.4.28**: *Let R be a ring. R is called a weak-iso-ring if there exists atleast a pair of distinct subrings of same order which are isomorphic.*

**THEOREM 4.4.11**: *Every iso-ring is a weak iso ring.*

*Proof*: Obvious by the very definition.

Recall from [2].

**DEFINITION [2]**: *An arbitrary group G is called a B-group if any two subgroups of same order are conjugate and G is a iso group if any two subgroups of same order are isomorphic.*

***Example 4.4.16***: Let $G = \langle g / g^2 = 1 \rangle$ and $Z_3 = \{0, 1, 2\}$ be the prime field of characteristic 3. $Z_3G$ be the group ring. It is easily verified $Z_3G$ is an iso-ring but is not a co-ring.

In view of this example we propose open problems in chapter 5 of this book.

**THEOREM 4.4.12**: *Let $Z_2 = \{0, 1\}$ and $S_n$ be the symmetric group of degree n. The group ring $Z_2S_n$ is a weak co-ring and a weak-iso-ring.*

*Proof*: To prove this we have to find two subrings in $Z_2S_3$ which are isomorphic and two subrings which are conjugate. To this end consider the subrings $A = \{0, 1 + p_2\}$ and $B = \{0, 1 + p_1\}$ where

$$p_1 = \begin{pmatrix} 1 & 2 & 3 & 4 & . & . & . & n \\ 1 & 3 & 2 & 4 & . & . & . & n \end{pmatrix} \text{ and } p_2 = \begin{pmatrix} 1 & 2 & 3 & 4 & . & . & . & n \\ 3 & 2 & 1 & 4 & . & . & . & n \end{pmatrix}.$$

Clearly A and B are isomorphic as subrings. Take $X = \{1 + q_1, 0\}$ and $Y = \{0, 1 + q_2\}$ where

$$q_1 = \begin{pmatrix} 1 & 2 & 3 & 4 & 5 & 6 & . & . & . & n \\ 2 & 1 & 4 & 3 & 5 & 6 & . & . & . & n \end{pmatrix}$$

and

$$q_2 = \begin{pmatrix} 1 & 2 & 3 & 4 & 5 & 6 & . & . & . & n \\ 3 & 4 & 1 & 2 & 5 & 6 & . & . & . & n \end{pmatrix}.$$

$$Y = \begin{pmatrix} 1 & 2 & 3 & 4 & 5 & . & . & . & n \\ 1 & 3 & 2 & 4 & 5 & . & . & . & n \end{pmatrix} X \begin{pmatrix} 1 & 2 & 3 & 4 & 5 & . & . & . & n \\ 1 & 3 & 2 & 4 & 5 & . & . & . & n \end{pmatrix}$$

Thus X and Y are conjugate subrings; hence $Z_2S_n$ is a weak co-ring.

**THEOREM 4.4.13:** *Let F be any ring or a field and G a B-group then the group ring FG is a weak co-ring and a weak iso-ring.*

*Proof*: Left as an exercise for the reader to prove.

Now we proceed onto define Smarandache co-ring and Smarandache iso-ring.

**DEFINITION 4.4.29**: *Let R be a ring if two S-subrings of R of same order are conjugate then we say R is a Smarandache weak co-ring (S-weak co-ring).*

**DEFINITION 4.4.30**: *Let R be a ring if every pair of S-subrings of same order are conjugate then we say R is a Smarandache co-ring (S-co-ring).*

**DEFINITION 4.4.31**: *Let R be a ring, if every pair of S-subrings of same order are isomorphic then we say R is a Smarandache iso ring (S-iso-ring).*

**DEFINITION 4.4.32**: *Let R be a ring if R has a pair of S-subrings of same order that are isomorphic then we say R is a Smarandache weak iso ring (S-weak iso-ring).*

The following two theorems are left for the reader to prove as an exercise.

**THEOREM 4.4.14**: *Every S-co-ring is a S-weak co-ring.*

**THEOREM 4.4.15**: *Every S-iso-ring is a S-weak iso-ring.*

We propose simple problems at the end of this section as well as difficult problems in the last chapter for the reader to solve.

**DEFINITION [24]**: *A ring is e-primitive if every nonzero ideal in R contains a nonzero idempotent element.*

***Example 4.4.17***: Let $Z_2 = \{0, 1\}$ and $S = \{0, 1, a, b / a^2 = 0, b^2 = 1, ab = ba = a\}$. The semigroup ring $Z_2S$ is not e-primitive.

**Example 4.4.18**: The ring $Z_{12} = \{0, 1, 2, \ldots, 11\}$ is not e-primitive. For ideals of $Z_{12}$ are $I_1 = \{0, 4, 8\}$, $I_2 = \{0, 2, 4, 6, 8, 10\}$, $I_3 = \{0, 3, 6, 9\}$ and $I_4 = \{0, 6\}$. All the ideals $I_1$, $I_2$ and $I_3$ are e-primitive where as $I_4 = \{0, 6\}$ is not e-primitive.

**DEFINITION 4.4.33**: *Let R be ring. If R has atleast one ideal which has an idempotent in it then we say R is weakly e-primitive.*

**Example 4.4.19**: The ring $Z_{12}$ is weakly e-primitive.

**DEFINITION 4.4.34**: *Let R be a ring. If every nonzero S-ideal in R contains a nonzero S-idempotent then R is Smarandache e-primitive (S-e-primitive).*

**DEFINITION 4.4.35**: *Let R be a ring. If R has atleast one nonzero S-ideal in R, which contains a nonzero S-idempotent then we say R is a Smarandache weakly e-primitive (S-weakly e-primitive).*

**THEOREM 4.4.16**: *Every S-e-primitive ideal is S-weakly e-primitive.*

*Proof*: Left as an exercise for the reader to prove.

**THEOREM 4.4.17**: *Let R be a field of characteristic 0 and G be a torsion free abelian group. The group ring KG is never*

       *1. Weakly e-primitive.*
       *2. S-weakly e-primitive.*

*Proof*: Obvious from the fact that KG has no nontrivial idempotents; as KG is a domain hence has no S-idempotents.

**DEFINITION [161]**: *Let R be a ring. If $M \neq 0$ is an additive subgroup of a ring R with zero divisors then M is an SV-group, in case $x \cdot y = 0$ for all $x, y \in M$ and $L_R(M) \cap R_R(M) \subset M$ ($L_R$ and $R_R$ are the left and right annihilators).*

We weaken this concept and define weakly SV-group as follows.

**DEFINITION 4.4.36**: *Let R be a ring. If $M \neq 0$ is an additive subgroup of a ring R with zero divisors then we call M a weakly SV-group.*

**Example 4.4.20** : $Z_6 = \{0, 1, 2, 3, 4, 5\}$ is not even a weakly SV-group.

**Example 4.4.21**: $Z_{12} = \{0, 1, 2, \ldots, 11\}$ is a weakly SV group. For $M = \{0, 2, 4, 6, 8\}$ is such that

$$2.6 \equiv 0 \pmod{12}$$
$$4.6 \equiv 0 \pmod{12}$$
$$8.6 \equiv 0 \pmod{12}.$$

Now we proceed on to define Smarandache parallel.

**DEFINITION 4.4.37**: *Let R be a ring. If M ≠ 0 be a S-semigroup under addition of the ring R. M is called a Smarandache SV-group (S-SV group) if in case x . y = 0 for all x, y ∈ M and $L_R(M) \cap R_R(M) \subseteq M$.*

**DEFINITION 4.4.38**: *Let R be a ring. If M ≠ 0 be a S-semigroup under addition with S-zero divisors then we call M a Smarandache weakly SV-group (S-weakly SV-group).*

**THEOREM 4.4.18**: *Let R be a SV group then R is a S-weakly SV group only when every divisor in M is a S-zero divisor.*

*Proof*: Left for the reader to verify using definitions.

**THEOREM 4.4.19**: $Z_p G$ *where $Z_p$ is a prime field of characteristic p and G = $\langle g/g^{2p} = 1 \rangle$ be the group ring. Then $Z_p G$ is a weakly SV-group and a S-weakly SV-group only when all zero divisors in the subgroup are S-zero divisor.*

*Proof*: Take H = $\{1, g^2, g^4, \ldots, g^{2p-2}\}$. Now $Z_p H$ is a weakly SV-group. $Z_p G$ is a S-weakly SV-group only if all zero divisor in ZH are S-zero divisor.

Now we define radix for rings which once again uses additive subgroup.

**DEFINITION [19]**: *An additive subgroup S of a commutative ring R is called a Radix provided $tx^3$ and $(t^3 - t)x^2 + t^2 x$ are in S for every x in S and for every t in R.*

**Example 4.4.22**: Let $Z_2 = (0, 1)$ be the field of integers. G = $\langle g/g^2 = 1 \rangle$, $\{0, 1 + g\}$ is not a radix of $Z_2 G$. $\{0, 1\}$ is not a radix of $Z_2 G$.

**DEFINITION 4.4.39**: *Let R be a non-commutative ring. Let S be an additive subgroup of R. S is said to be a left radix of R if $tx^3$, $(t^3 - t) x^2 + t^2 x$ are in S for every t ∈ R and every x ∈ S. Similarly we define right radix of R if $x^3 t$, $x^2(t^2 - t) + xt^2$ are in S for every x ∈ S and t ∈ R. S is called a radix, if S is both a left and right radix of R.*

**Example 4.4.23**: Every right ideal I of a ring R is a right radix of R.

**THEOREM 4.4.20**: *Let $Z_2 = \{0, 1\}$ be the field of integers modulo 2 and G be a cyclic group of even order then $Z_2 G$ has a radix which is not an ideal of $Z_2 G$.*

*Proof*: Let $G = \langle g \mid g^{2n} = 1 \rangle$ and $Z_2 = \{0, 1\}$ Take $S = \{0, 1 + g^2 + g^4 + \ldots + g^{2n-1}\}$. Clearly S is a radix which is not an ideal of $Z_2 G$.

**THEOREM 4.4.21**: *Let R be a ring then x be an element which annihilates every element of R. Then $S = \langle \{0, x\} \rangle$ is a radix of R.*

*Proof*: Obvious.

Now we proceed in to define Smarandache radix for a ring.

**DEFINITION 4.4.40**: *Let R be a commutative ring. An additive S-semigroup S of R is said to be a Smarandache radix (S-radix) of R is $x^3 t$, $(t^2 - t)x^2 + xt^2$ are in S for every $x \in S$ and $t \in R$. If R is a non-commutative ring then for any S-semigroup S of R we say R has Smarandache left radix (S-left radix) if $tx^3$, $(t^2 - t)x^2 + t^2 x$ are in S for every $x \in S$ and $t \in R$. Similarly we define Smarandache right radix (S-right radix) of R. If S is simultaneously a S-left radix and S-right radix of a non-commutative ring then we say R has a S-radix.*

**THEOREM 4.4.22**: *Let R be a ring and S a radix of R. S in general is not a S-radix.*

*Proof*: By an example. Consider $Z_2 G$ where $Z_2 = \{0, 1\}$ and $G = \langle g \mid g^8 = 1 \rangle$. The group ring has $H = \{0, 1 + g^2 + g^4 + g^6\}$ to be radix of $Z_2 G$ but H is not a S-radix as H is not a S-semigroup.

**THEOREM 4.4.23**: *Let R be a ring, if H is a S-radix of R then H is a radix of R.*

*Proof*: Clear from the very definitions of the radix and S-radix.

[4] has defined rings which has $\gamma$-semigroups and obtained some interesting results about them.

**DEFINITION [4]**: *A multiplicative semigroup M of a ring R is a $\gamma$-semigroup if for each $a \in M$, the additive subgroup of R generated by a is contained in M.*

**Example 4.4.24**: Let $Z_2 = \{0, 1\}$ and $G = \langle g \mid g^2 = 1 \rangle$. The group ring $Z_2 G$ has $\gamma$-semigroup. For $M = \{0, 1 + g\}$ is a $\gamma$-semigroup of $Z_2 G$.

**THEOREM 4.4.24**: *Let G be a group having an element of finite order and $Z_2 = \{0, 1\}$. The group ring $Z_2G$ has $\gamma$-semigroup.*

*Proof*: Let $g \in G$ such that $g^n = 1$. Then $M = \{0, 1 + g + g^2 + \ldots + g^{n-1}\}$ is a $\gamma$-semigroup.

**THEOREM 4.4.25**: *Let S(n) be the semigroup and $Z_2 = \{0, 1\}$. The semigroup ring $Z_2S(n)$ has $\gamma$-semigroup.*

*Proof*: We know $S_n \subset S(n)$, by the above theorem (4.4.24). $S(n)$ has elements of finite order; hence $Z_2S(n)$ has $\gamma$-semigroup.

**THEOREM 4.4.26**: *Every group ring KG is a $\gamma$-semigroup.*

*Proof*: Given KG is the group ring. Two cases arise

    1. G has elements of finite order.
    2. G has no elements of finite order.

If G has elements of finite order then we see by theorem 4.4.25; KG has $\gamma$-semigroup. If G has no elements of finite order take $g \in G$ and let H be the infinite cyclic group generated by g. Then KH is a $\gamma$-semigroup of KG.

**DEFINITION 4.4.41**: *Let R be a ring. A multiplicative S-semigroup M of R is a Smarandache $\gamma$-semigroup (S-$\gamma$-semigroup) if for each $a \in M$ the additive subgroup of R generated by a is contained in M.*

**THEOREM 4.4.27**: *Let R be a ring; if M is a S-$\gamma$-semigroup of R, then M is a $\gamma$-semigroup of R.*

*Proof*: Easily derived from the definitions.

**THEOREM 4.4.28**: *Let R be a ring; if M is a $\gamma$-semigroup of R then M need not in general be a S-$\gamma$ semigroup of R.*

*Proof*: By an example. Let $Z_2G$ be the group ring with $G = \{g/g^2 = 1\}$; $M = \{0, 1 + g\}$ is a $\gamma$-semigroup of $Z_2G$, but M is not a S-$\gamma$-semigroup of $Z_2G$.

Now we proceed on to define yet another new concept called $\delta$-semigroups.

**DEFINITION 4.4.42**: *Let R be a ring. A non-empty multiplicative semigroup M not containing 1 is called a δ-semigroup if for every a ∈ M; the ideal of R generated by a is such that aR ⊂ M and Ra ⊂ M.*

**Example 4.4.25**: Let $Z_4$ = {0, 1, 2, 3} be the ring of integers modulo 4. Clearly M = {0, 2} is a δ-semigroup.

**DEFINITION 4.4.43**: *Let R be a ring. R is called a δ-semigroup ring (δ-s-ring) if for every multiplicative semigroup M of R not containing 1, M is a δ-semigroup.*

**Example 4.4.26**: Let $Z_4$ = {0, 1, 2, 3}, $Z_4$ is a δ-s-ring.

**Example 4.4.27**: $Z_2$G where G = ⟨g / $g^2$ = 1⟩ is a δ-s-ring.

**THEOREM 4.4.29**: *A field can never have δ-semigroups.*

*Proof*: Obvious by the very structure of fields.

**Example 4.4.28**: $Z_{12}$ = {0, 1, 2, …, 11}, $Z_{12}$ is not a δ-s-ring.

**THEOREM 4.4.30**: *Let G be a finite group and K any field. The group ring KG has δ-semigroup.*

*Proof*: Take M = {0, cΣ$g_i$ / c ∈ K} where cΣ$g_i$ = c(1 + $g_1$ + $g_2$ + … + $g_n$) such that G = {1, $g_1$, …, $g_n$} ; clearly M is a δ-semigroup.

**DEFINITION 4.4.44**: *Let R be a ring. A non-empty multiplicative S-semigroup M not containing 1 of R (unit of R) is called a Smarandache δ-semigroup (S-δ-semigroup) if for every a ∈ M, the ideal generated by a is contained in M.*

**Example 4.4.29**: $Z_{12}$ = {0, 1, 2, 3, …, 11} be the ring of integers modulo 12. M = {0, 2, 4, 6, 8, 10} is a S-δ-semigroup for A = {4, 8} is a group with 4 as identity. M is a S-δ-semigroup.

**THEOREM 4.4.31**: *Let R be a ring. If M is a S-δ-semigroup then M is a δ-semigroup.*

*Proof*: Obvious by the very definition.

**THEOREM 4.4.32**: *Let R be a ring, if M is a δ-semigroup; M in general is not a S-δ-semigroup.*

*Proof*: By an example, $Z_4 = \{0, 1, 2, 3\}$ is in fact a δ-s-ring but no δ semigroup of $Z_4$ is a S-δ semigroup of $Z_4$.

**DEFINITION 4.4.45**: *Let R be a ring. R is called a Smarandache δ-semigroup ring (S-δ-s-ring) if every S-δ-semigroup under multiplication of R is a S-δ-semigroup.*

Now we recall the concept of SG-rings which also makes use of multiplicative semigroup.

**DEFINITION 4.4.46**: *Let R be a ring. R is said to be a SG-ring if R = $\cup S_i$ where $S_i$'s are multiplicative semigroups and $S_i \cap S_j = \phi$, if $i \neq j$.*

***Example 4.4.30***: Let $Z_2G$ be the group ring where $G = \langle g \, / \, g^2 = 1 \rangle$; $Z_2G = S_1 \cup S_2$ where $S_1 = \{0, 1 + g\}$ and $S_2 = \{1, g\}$, so $Z_2G$ is a SG-ring.

**DEFINITION 4.4.47**: *Let R be a ring. R is said to be a weakly SG-ring if R = $\cup S_i$ and $S_i \cap S_j \neq \phi$ even if $i \neq j$.*

***Example 4.4.31***: Let $Z_8 = \{0, 1, 2, \ldots, 7\}$ be the ring of integers modulo 8. $Z_8$ is a weakly SG-ring. For $Z_8 = \{0, 2, 4\} \cup \{1, 3\} \cup \{1, 5\} \cup \{1, 7\} \cup \{0, 2, 4, 6\} = \{1, 3\} \cup \{1, 5\} \cup \{1, 7\} \cup \{0, 2, 4, 6\}$. Hence the claim.

**THEOREM 4.4.33**: *Every SG-ring is a weakly SG-ring but not conversely.*

*Proof*: Obvious by the very definition. $Z_8$ is not a SG-ring but only a weakly SG-ring.

***Example 4.4.32***: Let $Z_7 = \{0, 1, 2, \ldots, 6\}$ be the ring of integers modulo 7. $Z_7$ is not a weakly SG-ring.

**DEFINITION 4.4.48**: *Let R be a ring. R is said to be a Smarandache SG-ring (S-SG-ring) if R = $\cup S_i$ where $S_i$ are multiplicative S-semigroups such that $S_i \cap S_j = \phi$ if $i \neq j$.*

**DEFINITION 4.4.49**: *Let R be a ring. R is said to be a Smarandache weakly SG-ring (S-weakly SG-ring) if R = $\cup S_i$ where $S_i$'s are S-multiplicative semigroups and $S_i \cap S_j \neq \phi$ even if $i \neq j$.*

**THEOREM 4.4.34**: *Let R be a SG-ring. Then R need not in general be a S-SG-ring.*

*Proof*: By an example; clearly the group ring $Z_2G$ is a SG-ring but is not a S-SG-ring.

**THEOREM 4.4.35**: *If R is a S-SG-ring then R is a SG-ring*.

*Proof*: By the very definition of these concepts the result follows.

Another interesting property about multiplicative semigroups is:

**DEFINITION 4.4.50**: *Let R be a ring. Let M be a multiplicative semigroup; we say R has a 0-semigroup if $S^2 = \{0\}$ and idempotent semigroup if $S^2 = S$.*

*A ring which has every multiplicative semigroup to be either a 0-semigroup or an idempotent semigroup is called a ZI-ring. If R has atleast a 0-semigroup and an idempotent semigroup then we call R a weak ZI-ring.*

**THEOREM 4.4.36**: *Every ZI-ring is a weak ZI-ring*.

The proof is left as an exercise to the reader.

**DEFINITION 4.4.51**: *Let R be a ring if in R every multiplicative semigroup M which is a S-semigroup is such that $M^2 = M$ or if M has a sub semigroup N such that $N^2 = 0$, then we call R a Smarandache ZI-ring (S-ZI-ring). If R has atleast a S-semigroup which is such that $M^2 = M$ and has a S-semigroup which has a subsemigroup N such that $N^2 = 0$ then we say R is a Smarandache weakly ZI-ring (S-weakly ZI-ring).*

**THEOREM 4.4.37**: *If R is a Boolean ring. Then R has multiplicative semigroups M such that $M^2 = M$.*

*Proof*: Obvious by the very definition of Boolean ring R; $a^2 = a$ for all $a \in R$.

**DEFINITION 4.4.52**: *Let R be a ring, if R has a multiplicative semigroup M, which has nontrivial idempotents or nontrivial nilpotents or both nontrivial idempotents and nilpotents or $S^2 = S$ or $S^2 = \{0\}$ then we call the ring R pseudo ZI-ring. We say M is a Smarandache pseudo ZI-ring (S-pseudo ZI-ring) if M is a multiplicative Smarandache semigroup which has nilpotents or S-idempotents, or $S^2 = S$ or $S^2 = \{0\}$ or has both nilpotents and S-idempotents.*

**THEOREM 4.4.38**: *If R is a ZI-ring then R is a pseudo ZI ring.*

*Proof*: Left for the reader to verify.

***Example 4.4.33***: Let $Z_{12} = \{0, 1, 2, \ldots, 11\}$. Now $M_1 = \{0, 4, 8\}$, $M_2 = \{0, 2, 4, 6, 8, 10\}$, $M_3 = \{0, 6\}$ and $M_4 = \{0, 3, 6, 9\}$ are semigroups under multiplication where $M_3^2 = \{0\}$; $M_1$, $M_2$ and $M_4$ are S-semigroups. We see $Z_{12}$ is a S-pseudo ZI-ring as all the S-semigroups $M_1$, $M_2$ and $M_4$ have both nilpotents and idempotents, for $6^2 \equiv 0$ (mod 12), $9^2 \equiv 9$ (mod 12) is such tht $3^2 \equiv 9$ (mod 12) and $3.9 \equiv 3$ (mod 12) and $4^2 \equiv 4$ (mod 12), $8^2 \equiv 4$ (mod 12) and $4.8 \equiv 8$ (mod 12).

Now we proceed to define square sets in rings.

**DEFINITION 4.4.53**: *Let R be ring. A non-empty subset A of R, $|A| > 1$ is said to be a square set of R if for every $a \in A$ there exists atleast one $b \in R$ ($b \neq a$) such that $a = b^2$.*

***Example 4.4.34***: Let $Z_4 = \{0, 1, 2, 3\}$ be the ring of integers modulo 4. $A = \{0, 1\}$ is a square set of $Z_4$ as $0 \equiv 2^2 (\text{mod } 4)$ and $1 \equiv 3^2 (\text{mod } 4)$, $0^2 \equiv 0$ and $1^2 = 1$ are called trivial forms.

***Example 4.4.35***: Let $Z_p$ be the prime field of characteristic p. $A = \{1\}$ is not a square set though we have $(p - 1)^2 \equiv 1$ (mod p).

**THEOREM 4.4.39**: *Let Z be the ring of integers the square set A of Z is non-empty.*

*Proof*: $A = \{n^2 / n \in Z\} \neq \phi$. For $\{9, 4, 25, 36\} = A$. $3 \in Z$ is such that $3^2 = 9$, $3 \in Z$ is such that $2^2 = 4$, $5 \in Z$ is such that $5^2 = 25$, $6 \in Z$ is such that $6^2 = 36$ and so on.

**THEOREM 4.4.40**: *Let K be a field of characteristic 0 and G a torsion free abelian group. The group ring KG has a square set which is non-empty.*

*Proof*: Left for the reader to prove using the fact KG is a domain.

**DEFINITION 4.4.54**: *Let R be a ring. R is said to have a Smarandache square set (S-square set) A if $|A| > 1$ and A is an additive S-semigroup and $a \in A$ is such that there exist $r \in R$ with $a = r^2$.*

***Example 4.4.36***: Let $Z_{16} = \{0, 1, 2, \ldots, 14, 15\}$ be the ring of integers modulo 16. $\{1, 4\} = A$, is a square set for $9^2 \equiv 1 (\text{mod } 16)$ and $14^2 \equiv 4$ (mod 16). It is easily verified $7^2 \equiv 1$ (mod 16) and $10^2 \equiv 4$ (mod 16). Thus a single element can have several representations. But we see $Z_{16}$ has no S-square set.

**THEOREM 4.4.41**: *Let R be a ring, if A is a S-square set of R then A is a square set of R and not conversely.*

*Proof*: Given A ⊂ R is a S-square set for A is a semigroup under addition and A is a square set.

The square set in general need not be a S-square set for in $Z_{16}$ given in example 4.4.36; A = {1, 4} is a square set but is not a S-square set.

**DEFINITION [66]**: *If R is a ring and $0 \neq r \in R$ then a non-empty subset X of R is said to be a (right) insulator for r in R if the right annihilator $r_R$ = {rx / x ∈ X} = 0. If every non-zero element of R has a finite insulator the author calls the ring R to be (right) strongly prime i.e., A non-zero ring R is said to be (right) strongly prime if every non-zero element of R has finite insulator.*

**THEOREM 4.4.42**: *Let G be a torsion free abelian group and K any field of characteristic zero. No element in KG have insulators.*

*Proof*: Follows from the fact KG is a domain.

**DEFINITION 4.4.55**: *Let R be a ring. We say $0 \neq r \in R$ is called a Smarandache insulator (S-insulator) if for r there exists a non-empty subset X of R where X is a S-semigroup under addition and the right annihilator $r_R$ = ({rx / x ∈ X}) = 0. A non-zero ring R is said to be Smarandache strongly prime (S-strongly prime) if every non-zero element of R has a finite S-insulator.*

Obtain interesting results about them.

**DEFINITION 4.4.56**: *Let R be a commutative ring and P an additive subgroup of R.P is called the n-capacitor group of R if $x^n P \subset P$ for every $x \in R$ and $n \geq 1$ and n a positive integer.*

**Example 4.4.37**: Let $Z_4$ = {0, 1, 2, 3} be the ring of integers modulo 4. P = {0, 2} is a n-capacitor group of $Z_4$.

**Example 4.4.38**: Let R be a ring. Every ideal I of R is called the n- capacitor group of R.

**Example 4.4.39**: Let $Z_2 G$ be the group ring of the group G = ⟨g / $g^3$ = 1⟩. K = {0, 1 + g} and I = {0, 1 + $g^2$} are 3k-capacitor group of the group ring, k > 1.

**THEOREM 4.4.43**: *Let $Z_2 G$ be the group ring where G = ⟨g / $g^p$ = 1⟩. The group ring $Z_2 G$ has pk-capacitor group for k = 1, 2, 3, … .*

*Proof*: Let I = {0, 1 + g + … + $g^{p-1}$} is pk capacitor group for k = 1, 2, 3, … .

**THEOREM 4.4.44**: *Let F be a field of characteristic two and G = ⟨g/g$^n$ = 1⟩ be the cyclic group of degree n. The group ring KG has n-capacitor groups which are not ideals, {(n, 2) = 1 and n not a prime}.*

*Proof*: Left for the reader to verify.

**THEOREM 4.4.45**: *Let G be a torsion free abelian group and K a field of characteristic zero. ZG the group ring has no n capacitor groups other than the ideals.*

*Proof*: We know all ideals are n-capacitor groups. But KG has no n-capacitor group. For if P is a subgroup; for all x ∈ KG we have x$^n$P ⊄ P.

Now we define relative Smarandache notions.

**DEFINITION 4.4.57**: *Let R be a commutative group and P an additive S-semigroup of R. P is called a Smarandache n-capacitor group (S-n-capacitor group) of R if x$^n$p ⊆ P for every x ∈ R and n ≥ 1 and n a positive integer.*

**THEOREM 4.4.46**: *Let R be a commutative ring. If R has S-n-capacitor group then R has n-capacitor group.*

*Proof*: Follows by the definition.

Just for the sake of completeness we give the definition of semiring as it used to define semiorder in rings.

**DEFINITION [50]**: *Let S be a non-empty set on which is defined two binary operation '+' and '.' such that the following are true.*

1. *(S, +) is a monoid with 0 as the additive idenitity.*
2. *(S, .) is a semigroup under '.'.*
3. *s(a + b) = sa + sb and (a + b)s = as + bs for all a, b, s ∈ S*

*then (S, +, .) is a semiring. A proper subset A of S is a subsemiring if (A, +, .) is itself a semiring.*

**DEFINITION 4.4.58**: *Let R be a ring with identity. A non-empty subset S of R is a semi order in R if S satisfies the following conditions:*

1. *(S, +) is a semigroup with identity with respect to the operation + of R.*
2. *(S, .) is a semigroup '.' the operation of R that is (S, +, .) is a semiring in R which we call as the subsemiring of R.*

3. *For every nonzero divisor of S the inverse whenever exists is in R.*
4. *Any x = R is of the form x = sy or x = zs$^I$ for some s, s$^I$ ∈ S and some y, z in R.*

*If such a nontrivial S in R exists we call R a semi order ring or a so-ring.*

***Example 4.4.40***: Let Z be the ring of integers Z = $Z^+$ ∪ $Z^-$ ∪ {0}. S = $Z^+$ ∪ {0} is a semiring which is a semi order in Z. So Z is a so-ring.

**Theorem 4.4.47**: *Every field F of characteristic zero is a so-ring.*

*Proof*: Obvious by the definition.

**Definition 4.4.59**: *Let R be a ring. A nonempty subset S of R is a Smarandache semi order (S-semi order) in R if S satisfies for the following conditions.*

1. *(S, +) is a S-semigroup with identity with respect to '+' the operation in R*
2. *(S, .) is a S-semigroup '.' the operation of R and (S, +, .) is a S-semiring in R.*
3. *For every nonzero divisor of S the inverse whenever exists is in R.*
4. *Any x ∈ R is of the form x = sy or x = zs$^I$ for some s, s$^I$ ∈ S and some y, z ∈ R.*

For more about S-semigroups and S-semirings please refer [154, 156, 157].

**Theorem 4.4.48**: *Every S-semi order is a semi order and not conversely.*

*Proof*: Follows easily by the definitions; hence left for the reader as an exercise.

**Definition 4.4.60**: *Let R be a ring. A subset I of R which is closed under '+' of R is called the square ring ideal of R if $r^2i$ ∈ I and $ir^2$ ∈ I for all i ∈ I and ∀ r ∈ R.*

***Example 4.4.41***: Let $Z_2G$ be the group ring where G = ⟨g / $g^2$ = 1⟩; I = {0, g} is a square ideal ring of $Z_2G$.

**Definition 4.4.61**: *Let R be a ring. A subset I of R which is closed under '+' is called a n-ring ideal of R if $r^ni$ ∈ I and $ir^n$ ∈ I for all i ∈ I and r ∈ R (n > 1, n a positive integer).*

It is important to note that even if n = 1 we see I is not an ideal of R.

**THEOREM 4.4.49**: *Let $Z_2 = \{0, 1\}$ and $G = \langle g \,/\, g^{2n} = 1 \rangle$. The group ring $Z_2G$ is a 2n-ring ideal of R.*

*Proof*: Take any I = {0, g}; such that g ∈ G. Clearly $s^{2n}i \in I$ and $is^{2n} \in I$ for all i ∈ I and s ∈ $Z_2$G.

**THEOREM 4.4.50**: *Let G be a torsion free abelian group and K any field. No finite subset of KG even closed under addition is a n-ring ideal for any n.*

*Proof*: Easily evident from the fact that G is torsion free abelian.

Now we proceed on to define the Smarandache analogue.

**DEFINITION 4.4.62**: *Let R be a ring. A subset I of R which is a S-semigroup under '+' of R is called the Smarandache square ring ideal (S-square ring ideal) of R if $ir^2$ and $r^2i \in I$ for all i ∈ I and r ∈ R.*

**DEFINITION 4.4.63**: *Let R be a ring. A subset I of R which is a S-semigroup under '+' is called the Smarandache n-ring ideal (S-n-ring) of R if for all i ∈ I and r ∈ R we have $r^ni$ and $ir^n \in I$.*

Obtain examples and some interesting results. The following theorem can be easily proved.

**THEOREM 4.4.51**: *Let R be a ring. If I is a S-square ideal of R then I is a square ideal of R.*

**DEFINITION 4.4.64**: *Let R be a ring. An ideal I of R is said to be quasi nilpotent if I does not contain any semi idempotent elements.*

***Example 4.4.42***: Let $Z_2$G be the group ring where G = $\langle g\,/\, g^2 = 1 \rangle$; clearly the ideal I = {0, 1 + g} has no semi idempotents so I is quasi nilpotent.

**DEFINITION 4.4.65**: *Let R be a ring. A S-ideal I of R is said to be Smarandache quasi nilpotent (S-quasi nilpotent) if I does not contain any S-semi idempotents.*

***Example 4.4.43***: $Z_6$G be the group ring of the group G = $\langle g/g^2 = 1 \rangle$ over the ring of integers modulo 6; $Z_6$ = {0, 1, 2, …, 5}; clearly $Z_6$ has S-ideals HG where H = {0, 3} which is also a S-quasi nilpotent of $Z_6$G.

**DEFINITION [24]**: *Gray defined a radical ideal as follows: A subset P of a ring R is a radical if*

1. *P is an ideal.*
2. *P is a nil ideal.*
3. *R/P has no nonzero nilpotent right ideals.*

*The sum of all ideals in R satisfying 1) and 2) is the upper radical of R and is denoted by $\cup(R)$. The intersection of all those ideals in R satisfying 1) and 3) is the lower radical of R. L(R).*

For more about radical ideals please refer [24].

***Example 4.4.44***: The group ring $Z_2G$ where $G = \langle g \, / \, g^2 = 1 \rangle$ is a radical ideal. It is important to note that even the subideal of a radical ideal in general need not be a radical ideal.

***Example 4.4.45***: Let $Z_2G$ be the group ring where $G = \langle g \, / \, g^3 = 1 \rangle$. $Z_2G$ has no proper radical ideals as $Z_2G$ has no nilpotent elements so it cannot have nil ideals hence the claim.

**DEFINITION 4.4.66**: *Let R be a ring. The Smarandache radical ideal (S-radical ideal) P of R is defined as follows:*

1. *P is a S-ideal of R.*
2. *$S \subset P$ where S is a subideal of P is a nil ideal.*
3. *R / P has no non-zero nilpotent right ideals.*

*The sum of all S-ideal of R satisfying 1) and 2) is called the Smarandache upper radical (S-upper radical) of R and is denoted by $S(\cup(R))$. The intersection of all those S-ideals in R satisfying 1) and 3) is the Smarandache lower radical (S-lower radical) of R denoted by $S(L(R))$.*

The reader is requested to study radical ideals and S-radical ideals of a ring R.

**THEOREM 4.4.52**: *The group ring KG of the torsion free abelian group G over any field K has no*

1. *radical ideals.*
2. *S-radical ideals.*

*Proof*: Follows from the fact that KG has no nilpotent or zero divisors as KG is a domain.

**DEFINITION 4.4.67**: *Let R be a ring. A be an additively closed subset of R. For $a, b \in R$, $a, b \notin A$. We say a is right A related to b if $a \in Ab$; a is said to be left A related to b if $a \in bA$. If a is both right and left A related to b then we say a is A* *related to b. Obviously in case of commutative rings the notion of right and left related coincides.*

**Example 4.4.46**: Let $Z_2 S_3$ be the group ring. Take A = $\{0, p_5\}$; $p_1 \in Ap_2 = \{0, p_1\}$ So $p_1$ is right A related to $p_2$ but $p_1$ is not left A related to $p_2$.

**DEFINITION 4.4.68**: *Let R be a ring. A be a semigroup with respect to '+'. For a, b $\in$ R (a, b $\in$ A) we say a is both way related to b (or a and b related with respect to A or relative to A) if a $\in$ Ab and b $\in$ Aa.*

**THEOREM 4.4.53**: *Every prime field of characteristic p is relation free.*

*Proof*: Left for the reader to prove.

**THEOREM 4.4.54**: *Let R the field of reals. No pair in R can be related.*

*Proof*: Left for the reader to verify.

**DEFINITION 4.4.69**: *Let R be a ring. A be a S-semigroup under '+'. For a, b $\in$ R and a, b $\notin$ A; we say a is Smarandache right related (S-right related) to b if a $\in$ Ab. a is said to be Smarandache left related (S-left related) if a $\in$ bA; if a is both Smarandache right and left related to b then we say a is Smarandache A related to b (S-A related to b).*

**DEFINITION 4.4.70**: *Let R be a ring. A be a S-semigroup under addition. For a, b $\in$ R, a, b $\notin$ A we say a is both way related to b (or a and b are related with respect to A) with respect to A if a $\in$ Ab and b $\in$ Aa. The pair (a, b) is called also as a Smarandache related pair (S-related pair).*

Obtain interesting results about S-related pairs. On similar lines when we replace the semigroup under addition by a subring we define a relation called subring relation on R. When a ring has this relation the ring has nontrivial divisors of zero (and) or units.

**DEFINITION 4.4.71**: *Let R be a ring. A pair of elements x, y $\in$ R is said to have a subring right link relation if there exists a subring M of R in R \ {x, y} i.e., M $\subseteq$ R \ {x, y} such that x $\in$ My and y $\in$ Mx. Similarly subring left link relation if x $\in$ yM and y $\in$ xM. If it has both a left and a right link relation for the same subring M then we have x and y have a subring link relation and is denoted by xMy.*

**Example 4.4.47**: $Z_4$ = $\{0, 1, 2, 3\}$ be the ring of integers modulo 4. No pair of elements in $Z_4$ has a subring link relation.

**Theorem 4.4.55**: *Let R be a ring. M a subring such that x, y ∈ R \ M are subring link related. Then R has nontrivial divisors of zero or a unit.*

*Proof*: Let x, y ∈ R with x and y ∉ M, where M is a subring such that x ∈ yM and y ∈ xM that is x = yt and y = xu for some u, t ∈ M so that x = xut leading to x(1 − ut) = 0. The two possibilities are either x(1 − ut) = 0 is a zero divisor or ut = 1, then R has a unit.

**Example 4.4.48**: Let $Z_6$ = {0, 1, 2, …, 5} be the ring of integers modulo 6. M = {0, 2, 4} is a subring of $Z_6$. $Z_6$ \ {0, 4, 2} has no pair which is linked.

**Example 4.4.49**: Let $Z_2S_3$ be the group ring. The element $p_4$ and $p_5$ cannot be subring related through any subring.

**Definition 4.4.72**: *Let R be a ring. We say a pair x, y in R has a weakly subring link with a subring P in R \ {x, y} if either y ∈ Px or x ∈ Py, 'or' in the strictly mutually exclusive sense and we have subring Q, Q ≠ P such that y ∈ Qx (or x ∈ Qy).*

**Definition 4.4.73**: *Let R be a ring. We say a pair x, y ∈ R is said to be one way weakly subring link related if we have a subring P ⊂ R \ {x, y} such that x ∈ Py and for no subring S ⊂ R \ {x, y} we have y ∈ Sx.*

**Definition 4.4.74**: *Let R be a ring a pair x, y ∈ R is said to have a Smarandache subring right link relation (S-subring right link relation) if there exists a S-subring P in R \ {x, y} such that x ∈ Px and y ∈ Py. Similarly Smarandache subring left link relation (S-subring left link relation) if x ∈ yP and y ∈ xP. If it has both a Smarandache left and right link relation for the same S-subring P then we say x and y have a Smarandache subring link (S-subring link).*

*We say x, y ∈ R is Smarandache weak subring link (S-weak subring link) with a S-subring P in R \ {x, y} if either x ∈ Py or y ∈ Px ('or' in strictly mutually exclusive sense) we have a S-subring Q ≠ P such that y ∈ Qx (or x ∈ Qy). We say a pair x, y ∈ R is said to be Smarandache one way weakly subring link related (S-one way weakly subring link related) if we have a S-subring P ⊂ R \ {x, y} such that x ∈ Py and for no subring Q ⊆ R \ {x, y} we have y ∈ Qx.*

Thus we see that subring link relation between a pair of elements in a ring leads to either zero divisor or units leading to the following:

**Theorem 4.4.56**: *Let KG be the group ring of a torsion free abelian group G and K any field, the group ring KG has no pair which is subring linked.*

*Proof*: Obvious from the fact a pair is subring link related forces zero divisors or units and as KG has no zero divisors or units. KG cannot have a pair which is subring linked.

The same result holds good in case of S-subring related pairs. Now the subrings of a ring are studied but no inter relation between them are studied. Here we define a concept called essential subrings and we feel the study of Smarandache essential subrings will throw more light on the S-subrings of a ring. With this view we just define the concept of essential subrings.

**Definition 4.4.75**: *Let R be a ring. A subring A of R is said to be an essential subring of R, if intersection of A with every other subrings of R is zero. By subring we mean only proper subrings.*

**Definition 4.4.76**: *Let R be a ring if every subring of R is an essential subring of R then we call R an essential ring.*

**Definition 4.4.77**: *Let R be a ring. A be a S-subring of R. A is said to be a Smarandache essential subring (S-essential subring) of R if the intersection of every other S-subring is zero. By S-subring we mean only proper S-subrings.*

**Definition 4.4.78**: *Let R be a ring. If every S-subring of R is S-essential S-subring then we call R a Smarandache essential ring (S-essential ring).*

**Definition 4.4.79**: *Let R be a ring. If for a pair of subrings P and Q of R there exists a subring T of R (T ≠ R) such that the subrings generated by PT and TQ are equal i.e. $\langle PT \rangle = \langle TQ \rangle$, then we say the pair of subrings are stabilized subrings and T is called the stabilizer subring of P and Q.*

**Definition 4.4.80**: *Let R be a ring. A pair of subrings A and B of R is said to be a stable pair if there exists a subring C of R (C ≠ R) such that $C \cap A = C \cap B$ and $\langle C \cup A \rangle = \langle C \cup B \rangle$ where $\langle \rangle$ denote generated by $C \cup A$ and $C \cup B$. C is called the stability subring for the stable pair of subrings.*

**Theorem 4.4.57**: *Let R be a ring. If the subring A, B of R is a stable pair then A, B is a stabilized pair and not conversely.*

*Proof*: Follows by the very definition of these two concepts. To prove the converse is not true, consider the ring $Z_{12} = \{0, 1, 2, \ldots, 11\}$ be the ring of integers modulo 12. $S_1 = \{0, 6\}$, $S_2 = \{0, 6, 3, 9\}$, $S_3 = \{0, 4, 8\}$ and $S_4 = \{0, 2, 4, 6, 8, 10\}$. The subrings $S_3$ and $S_4$ is a stabilized pair but it is not a stable pair. Hence the claim.

**DEFINITION 4.4.81**: *Let R be a ring. If every pair of subrings of R is a stable pair then we say R is a stable ring.*

**DEFINITION 4.4.82**: *Let R be a ring. If for a pair of S-subrings P and Q of R there exists a S-subring T of R (T ≠ R) such that the S-subrings generated by PT and TQ are equal i.e. ⟨PT⟩ = ⟨TQ⟩ then we say the pair P and Q is a Smarandache stabilized pair (S-stabilized pair) and T is called the Smarandache stabilizer (S-stabilizer) of P and Q.*

**DEFINITION 4.4.83**: *Let R be a ring. A pair of S-subrings A, B of R is said to be a Smarandache stable pair (S-stable pair) if there exists a S-subring C of R (C ≠ R) such that C ∪ A = C ∪ B and ⟨C ∪ A⟩ = ⟨C ∪ B⟩ where ⟨⟩ means the subring generated by C ∪ A and C ∪ B; C is called the Smarandache stability S-subring (S-stability S-subring) for the Smarandache stable pair (S-stable pair) A and B.*

**DEFINITION 4.4.84**: *Let R be a ring if every pair of S-subrings of R is S-stable pair then we say R is a Smarandache stable ring (S-stable ring).*

It is left for the reader to prove that the following theorem:

**THEOREM 4.4.58**: *Every S-stable ring is a S-stabilized ring.*

**PROBLEMS:**

1. Does $Z_8 S(5)$ have a S-quasi ordering?
2. Can $Z_5 S_3$ have a S-product quasi ordering?
3. Find a set A in $Z_7 S_5$ which has S-sum quasi ordering but A is not product quasi ordering.
4. Find the S-semi regular ideal of $Z_5 S_3$.
5. Can the semigroup ring $Z_8 S(5)$ have
   i.   S-semi nilpotent ideals?
   ii.  Semi nilpotent ideals?
   iii. S-Semi regular ideals?
   iv.  Semi regular ideals?
6. Find whether the group ring $Z_7 D_{26}$ have
   i.  S-semi nilpotent ideals.
   ii. S-semi regular ideals where $D_{26} = \{a, b/a^2 = b^6 = 1, \ bab = a\}$.
7. Is the semi group ring $Z_3 S(5)$ a
   i.  S-Sub semi ideal ring?
   ii. Sub semi ideal ring?
8. Can the group ring $Z_7 S_4$ be a
   i. S-sub semi ideal ring?

ii.   Sub semi ideal ring?

9.   Can the semigroup ring $Z_6S(4)$ be a
   i.   normal ring?
   ii.   S-normal ring?

10.   Is the ring $Z_{22}$ a
   i.   normal ring?
   ii.   S-normal ring?

11.   Is $Z_{19}$ a
   i.   normal ring?
   ii.   S-normal ring?

12.   Give an example of S-weakly G-ring.

13.   Is $Z_8S_3$ a S-G-ring?

14.   Is a Boolean ring a weakly G-ring?

15.   Is $Z_6S_3$ a n-closed subgroup ring?

16.   Can $Z_3G$ where $G = \langle g / g^3 = 1 \rangle$ be a n-closed subgroup ring?

17.   Give an example of a S-n-closed subgroup ring.

18.   Give an example of a S-n-closed subgroup ring which is not a n-closed subgroup ring.

19.   Give an example of a ring which is a co-ring.

20.   Is $Z_{20}$ a weak co-ring?

21.   Can $Z_2S_3$ be a S-weak iso-ring?

22.   Give an example of a S-weak iso-ring which is not a S-iso-ring.

23.   Give an example of S-co-ring.

24.   Is the semigroup ring $Z_4S(3)$ a
   i.   S-co-ring?
   ii.   S-iso-ring?
   iii.   S-weak iso-ring?

25.   Is the ring $Z_{27}$ S-e-primitive?

26.   Can $Z_{27}$ be atleast e-primitive?

27.   Give an example of a weakly e-primitive ring which is not e-primitive.

28.   Can the group ring $Z_3S_4$ be
   i.   e-primitive?
   ii.   S-e-primitive?
   iii.   S-weakly primitive?
   iv.   Weakly e-primitive?

29.   Give an example of a ring which is a S-weakly SV-group.

30.   Is $Z_{24}$ a S-SV-group?

31.   Can $Z_{27}$ be a SV-group?

32.   Give an example of a ring which has S-radix.

33.   Can $Z_{28}$ have a radix?

34.   Give an example of ring which has a radix which is not a S-radix.

35.   Prove $Z_2S_4$ has a γ-semigroup.

36.   Prove $Z_3S(3)$ has a γ-semigroup.

37. Give an example of a ring which has no $\gamma$-semigroup.

38. Can $M_{3\times3} = \{(a_{ij}) \,/\, a_{ij} \in Z_{12}\}$ the ring of $3 \times 3$ matrices have $\gamma$-semigroup?

39. Give an example of a ring which has $\delta$-semigroups.

40. Does $Z_{24}$ have
     i. $\delta$-semigroups?
     ii. S-$\delta$ semigroups?

41. Is $Z_{24}$ a SG-ring?

42. Give an example of a S-SG ring.

43. Can $Z_{15}$ be a weakly SG-ring?

44. Show $Z_{20}$ is a pseudo ZI ring?

45. Is the ring $Z_{20}$ a S-pseudo ZI-ring?

46. Will $(ZG)^2 = ZG$ when $G$ is torsion free abelian?

47. Prove $(KG)^2 = KG$, if $K$ is a field of characteristic zero and $G$ is torsion free abelian.

48. Give an example of a ring R which has S-square set.

49. Give an example of a ring R in which every square set is a S-square set.

50. Find an example of a ring which has nontrivial insulators.

51. Give a nontrivial example of a ring which has S-insulators.

52. Can the group ring $Z_5(S_4)$ be a
     i. S-semi order ring?
     ii. Semi order ring?

53. Can the group ring $Z_3G$ where $G = \langle g \,/\, g^9 = 1 \rangle$ have
     i. S-n-capacitor group?
     ii. n-capacitor group?

54. Give an example of a semi order ring which is not a S-semi order ring.

55. Give an example of a ring which has n-ideals but no S-n-ideals.

56. Can $Z_{24}$ be a square ideal ring?

57. Is $Z_{24}$ a S-square ideal ring?

58. Does the group ring $Z_2G$ where $G = \langle g \,/\, g^{2n} = 1 \rangle$ have a
     i. quasi nilpotent ideals?
     ii. S-quasi nilpotents ideal?

59. Does the semigroup ring $Z_3(S(4))$ have
     i. radical ideal?
     ii. upper radical ideal?
     iii. S-radical ideal?
     iv. lower S-radical ideal?

60. Give an example of a ring which has only radical ideals and does not contain S-radical ideals.

61. Can the group ring $Z_2S_n$ have
     i. Smarandache related pairs?
     ii. Related pair? (related to any semigroup under '+').

62. Can $Z_{28}$ have a pair which is

    i. subring linked?

    ii. S-subring linked?

63. Does $Z_2S_n$ contain at least one essential subring?

64. Is $Z_2G$, where $G = \langle g / g^{27} = 1 \rangle$, a S-essential ring?

65. Can $Z_2S(n)$ have at least one S-essential subring?

66. Does there exist a stabilized pair of subrings in $Z_{24}$?

67. Can $Z_{12}$ be a S-stable ring?

68. Give an example of a stable ring which is not a S-stable ring.

69. Is the semigroup ring $Z_2S(n)$ a stable ring?

70. Can the semigroup ring $Z_2S(n)$ be a S-stable ring?

## 4.5 Miscellaneous properties about Smarandache Rings

In this section we introduce several important properties to Smarandache rings like hyperrings, Hamiltonian rings, $J_k$-group rings, structure of fixed support in group rings, quasi distributivity, the lattice of S-ideals and S-subrings of a ring and show the lattice of S-ideals of a ring in general is not a modular lattice.

In this section several new very recently introduced concepts like, integrally closed semigroup (or rings), system of local units, $\pi$-regular rings, generalized stable ring are given together with their Smarandache analogue.

[18] introduces the concept of hypergroups using modulo integers that is for, groups under addition. Here we introduce hyperrings and Smarandache hyperrings I and Smarandache hyperrings II.

**DEFINITION 4.5.1**: *Let $Z_n$ be the ring of integers modulo n. The hyperring $(Z_n, q)$ $(q \leq n)$ obtained from $Z_n$ by defining $x + y = (x + y, x + y + q)$ and $x . y = (x.y, x.y.q)$ denoted by $(Z_n, +)$ and $(Z_n, q, .)$ respectively is a subset of $Z_{n \times n}$. We say the hyperring $(Z_n, q)$ has a ring structure only when $(Z_n, q, +)$ $[(Z_n, q, .)]$ which is a subset of $Z_{n \times n}$ is a ring under component wise '+' and ' .' modulo n. $(Z_n, q, +)$ may partition $Z_n \times Z_n$ or $(Z_n, q, .)$ may partition $Z_n \times Z_n$ for varying $q \in Z_n$.*

**_Example 4.5.1_**: Let $Z_4 = \{0, 1, 2, 3\}$ be a ring of integers modulo 4. The hyperrings for all $q \in Z_4$ under '+' are

$$\{Z_4, 3, +\} = \{(0, 3), (1, 0), (2, 1), (3, 2)\}$$
$$\{Z_4, 2, +\} = \{(0, 2), (1, 3), (2, 0), (3, 1)\}$$
$$\{Z_4, 1, +\} = \{(0, 1), (1, 2), (2, 3), (3, 0)\}$$

$$\{Z_4, 0, +\} = \{(0, 0), (1, 1), (2, 2), (3, 3)\}$$

$\{Z_4, q, '+'\}$ partitions $Z_4 \times Z_4$. $\{Z_4, +, 0\}$ is a subring all others are not even closed under '+'.

The hyperrings for all $q \in Z_4$ under '.' are

$$(Z_4, 3, .) = \{(0, 0), (1, 3), (2, 2), (3, 1)\}$$
$$(Z_4, 2, .) = \{(0, 0), (1, 2), (2, 0), (3, 2)\}$$
$$(Z_4, 1, .) = \{(0, 0), (1, 1), (2, 2), (3, 3)\}$$
$$(Z_4, 0, .) = \{(0, 0), (1, 0), (2, 0), (3, 0)\}$$

Thus we see $Z_4 \times Z_4$ is not properly partitioned by $(Z_4, q, '.')$ defined by x.y = $\{(x.y, x.y.r)/ r \in Z_4\}$ and $\{Z_4, 1, .\}$ and $\{Z_4, 0, .\}$ are the only subrings of $Z_4 \times Z_4$.

**THEOREM 4.5.1**: *Let $Z_n$ be a ring of integers modulo n.*

1. *$(Z_n, 1, .), (Z, 0, .)$ and $\{Z_n, 0, +\}$ are the only subrings of $Z_n \times Z_n$.*
2. *$Z_n \times Z_n$ is never partitioned by the operation '.'.*
3. *$(Z_n, 1, .) = \{Z_n, 0, +\}.$*

*Proof*: $(Z_n, 1, .) = \{(x.y, x.y.1) / x, y \in Z_n\}$. It is easily verified $(Z_n, 1, .)$ is a subring for $(Z_n, 1, .) = \{(0, 0), (1, 1), (2, 2), \ldots, (n-1, n-1)\}$. $(Z_n, 0, .) = \{(x.y, x.y.0) / x, y \in Z_n\}$ is a subring of $Z_n \times Z_n$. It is easily verified as $(Z_n, 0, .) = \{(0, 0), (1, 0), (2, 0), \ldots, (n-1, 0).\}$. Now $\{Z_n, +, 0\} = \{(0, 0), (1, 1), (2, 2), \ldots, (n-1, n-1)\} = (Z_n, 1, .)$ is a subring of $Z_n \times Z_n$.

**DEFINITION 4.5.2**: *Let $Z_n$ be a ring with A to be S-subring of $Z_n$. Define the Smarandache hyperring I (S-hyperring I) to be a subring of $A \times A$ given by for any $q \in A$. $(A, q, +) = \{ (a_1 + a_2, a_1 + a_2 + q) / a_1, a_2 \in A\}$ and $(A, q, .) = \{(a_1 . a_2, a_1 . a_2 . q)/ a_2, a_1 \in A\}$. Similarly we define Smarandache hyperring II for any S-subring II of B.*

[40] defines a ring R to be a generalized Hamiltonian ring if every non-zero subring of R includes a non-zero ideal of R.

***Example 4.5.2***: Let Z be the ring of integers. Z is a generalized Hamiltonian ring.

**THEOREM 4.5.2**: *Let KG be the group ring of the group G over any field K. The group ring KG is not a generalized Hamiltonian ring.*

*Proof*: For K ⊂ KG is a subring which cannot include any ideal. Hence the claim.

**Theorem 4.5.3**: *Suppose the group ring RG is Hamiltonian then we see R is Hamiltonian.*

*Proof*: Follows from the fact that R ⊂ RG is a subring of RG so R should be a ring in which every non-zero subring includes a non-zero ideal. Hence the claim.

**Definition 4.5.3**: *Let R be a ring. We say R is a Smarandache Hamiltonian ring (S-Hamiltonian ring) if every S-subring includes a non-zero S-ideal.*

**Definition 4.5.4**: *Let R be a ring we say R is a Smarandache Hamiltonian ring II (S-Hamiltonian ring II) if every S-subring II includes a non-zero S-ideal II.*

**Theorem 4.5.4**: *Every S-Hamiltonian ring I is a S-Hamiltonian ring II and not conversely.*

*Proof*: By the very definition of S-subrings I and S-subrings II and S-ideals I and S-ideals II. The result is true. For the converse consider the ring Z. Clearly Z is a Smarandache Hamiltonian II and is not Smarandache Hamiltonian I.

We just recall in a group ring KG or in a semigroup ring KS we define for any $\alpha \in$ KG (or KS) |supp $\alpha$| = {$g_i / \alpha_i \neq 0$} where $\alpha = \sum \alpha_i g_i$ and |supp $\alpha$| denotes the number of elements in $\alpha$ which has non-zero coefficients. It is a subset of G(or S). For more about these ideas please refer [61, 62].

**Definition 4.5.5**: *Let RG be a group ring of a group G over the ring R. Let N = {$\alpha \in$ RG / |supp $\alpha$| = n}, n a fixed positive integer. If 0 is adjoined with N and if N $\cup$ {0} becomes a subring of RG we call N a fixed support subring of the group ring RG or n-subring of RG. (The same holds good for semigroup rings).*

***Example 4.5.3***: Let $Z_2$G be the group ring of the group G = $S_3$ over $Z_2$. N = {$p_4 + p_5$, $1 + p_4$, $1 + p_5$ / |supp $\alpha$| = 2}; N $\cup$ {0} is a 2-subring of $Z_2$G. Similarly M = {$1 + p_4 + p_5$ / |supp $\alpha$| = 3} adjoined with 0 gives a 3-subring of $Z_2$G.

**Definition 4.5.6**: *Let RG be the group ring of the group G over the ring R, if A = {$\alpha$ / |supp $\alpha$| = n}, n a fixed number is a semigroup under multiplication then N is called the fixed support semigroup of the group ring RG or the n-subsemigroup of RG.*

***Example 4.5.4***: Let $Z_2G$ be the group ring of the group $G = \langle g / g^3 = 1 \rangle$. Clearly P = $\{1 + g, 1 + g^2, g + g^2\}$ and $P_1 = \{1 + g + g^2\}$ are 2-subsemigroup and 3-subsemigroup of $Z_2G$ respectively.

**DEFINITION 4.5.7**: *Let RG be the group ring of the group G over the ring R. Let H = {α / |supp α| = m}, m a fixed integer. If H is a subgroup under multiplication after adjoining the identity 1 then we call H the fixed support subgroup of the group ring RG or m-subgroup of RG.*

**THEOREM 4.5.5**: *Every group ring RG has a 1-fixed support subgroup.*

*Proof*: Take $G = \{g / g \in G, |g| = 1\}$. Clearly G is a 1-fixed support subgroup.

**THEOREM 4.5.6**: *Let KS be the semigroup ring. We have S = {α / |supp α| = 1 and α ∈ S} is a 1-fixed subsemigroup.*

*Proof*: Obvious by the very definition of KS.

By using Smarandache notions we can combine the concept of fixed support of subgroup and fixed support of subsemigroup.

**DEFINITION 4.5.8**: *Let RG be the group ring. Let N = {α ∈ RG/ |supp α| = n}, n a fixed positive integer. If 0 adjoined in N becomes a S-subring of RG we call N a Smarandache fixed support subring (S-fixed support subring) of the group ring RG or S-n-subring of RG.*

**DEFINITION 4.5.9**: *Let RG be the group ring of the group G over the ring R if S = {α / |supp α| = n}, n a fixed number is a S-semigroup under multiplication then we call the set S to be the Smarandache fixed support of the group ring (S-fixed support of the group ring).*

Obtain interesting results about these concepts and study them.

**DEFINITION [44]**: *A ring R is said to be semi-connected if the center of R contains a finite number of idempotents.*

**THEOREM 4.5.7**: *A group ring KG of a finite group G over a field K is semi-connected.*

*Proof*: Obvious from the fact that K is a field and KG has idempotents as G is finite. Hence KG is semi-connected.

**DEFINITION 4.5.10**: *Let R be a ring. We say R is Smarandache semi-connected (S-semi-connected) if the center of R contains a finite number of S-idempotents.*

Once again as the main motive of this book is for researchers to develop Smarandache concepts we leave it for the reader to study this concept and get some nice results and examples of rings which are Smarandache semi-connected.

**THEOREM 4.5.8**: *Let R be a S-semi-connected ring then R is semi-connected.*

*Proof*: Follows from the very definitions of these notions.

[70] had defined the concept of $J_k$-ring and has studied them.

**DEFINITION [70]**: *Let R be a ring, k a fixed positive integer. We say R is a $J_k$-ring if for each $x_1, x_2, \ldots, x_k$ of R there exists a $n = n(x_1, x_2, \ldots, x_k) > 1$ such that $(x_1 x_2 \ldots x_k)^n = x_1 \ldots x_k$.*

For more about $J_k$ rings please refer [70].

**THEOREM 4.5.9**: *The group ring $Z_2 S_n$ is not a $J_k$-ring (k > 1).*

*Proof*: Left for the reader to verify. For $1 + p_1 \in Z_2 S_n$ where

$$p_1 = \begin{pmatrix} 1 & 2 & 3 & 4 & \ldots & n \\ 2 & 1 & 3 & 4 & \ldots & n \end{pmatrix}.$$

Clearly $(1 + p_1)^2 = 0$. $[(1 + p_1)(1 + p_1) \ldots (1 + p_1)]^n = 0 \neq 1 + p_1$. Hence the claim.

**THEOREM 4.5.10**: *Let G be a torsion free abelian group and K any field. The group ring KG is not a $J_k$-ring.*

*Proof*: Obvious from the fact for, if we take $g_1, g_2, \ldots, g_n \in G \subset KG$ then for no k > 1 we have $(g_1 \ldots g_n)^K = g_1 \ldots g_n$ as G is torsion free abelian.

**DEFINITION 4.5.11**: *Let R be a ring. We say R is a Smarandache $J_k$-ring (S-$J_k$ ring) if R contains a S-subring A, ($A \neq R$ but $A \subset R$) such that A is a $J_k$-ring.*

**THEOREM 4.5.11**: *Let R be a ring if R is a $J_k$ ring and has a non-trivial S-subring. Then R is a S-$J_k$-ring.*

*Proof*: Obvious by the very definitions.

Now we proceed onto find the nature of the lattice of the substructure of a ring. We know the set of all two-sided ideals of a ring form a modular lattice. Now we are interested in studying the following:

1.  Let R be a finite ring. M denote the collection of all S-subrings of R including {0} and R. What is the lattice structure of M?

2.  If we replace S-subring of R by S-subrings II of R, what is the structure of the lattice? Will it be distributive? modular? or non-modular?

3.  Let R be a ring M = {Set of all S-ideals of R}. What is the lattice structure of M?

We assume {0} and the ring R are trivially S-ideals, S-subrings, S-ideal II and S-subring II which act as the least and the greatest element of the lattice respectively.

***Example 4.5.5***: Let R = $Z_7 \times Z_9$ be a ring, the S-subrings of R are {{0}, R, $Z_7 \times \{0\}$, $Z_7 \times \{0, 3, 6\}$}. The lattice diagram is a 4 element chain lattice which is distributive and hence modular.

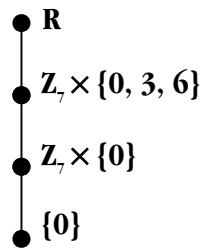

**R**

$Z_7 \times \{0, 3, 6\}$

$Z_7 \times \{0\}$

{0}

Figure 4.5.1

Clearly these S-subrings are also S-ideals of R.

***Example 4.5.6***: Let R = $Z_3 \times Z_{12} \times Z_7$ be the S-mixed direct product of rings. The S-subrings of R, are

$$
\begin{aligned}
A_1 &= Z_3 \times \{0\} \times Z_7 \\
A_2 &= Z_3 \times Z_{12} \times \{0\} \\
A_3 &= Z_3 \times \{0, 6\} \times \{0\} \\
A_4 &= Z_3 \times \{0, 4, 8\} \times \{0\} \\
A_5 &= Z_3 \times \{0, 3, 6, 9\} \times \{0\} \\
A_6 &= Z_3 \times \{0, 2, 4, \ldots, 10\} \times \{0\} \\
A_7 &= Z_3 \times \{0, 6\} \times Z_7 \\
A_8 &= Z_3 \times \{0, 3, 6, 9\} \times Z_7 \\
A_9 &= Z_3 \times \{0, 4, 8\} \times Z_7
\end{aligned}
$$

$$A_{10} = Z_3 \times \{0, 2, \ldots, 10\} \times Z_7$$
$$A_{11} = \{0\} \times Z_{12} \times \{0\}$$
$$A_{12} = \{0\} \times \{0, 2, 4, \ldots, 10\} \times \{0\}$$
$$A_{13} = \{0\} \times Z_{12} \times Z_7$$
$$A_{14} = \{0\} \times \{0, 6\} \times Z_7$$
$$A_{15} = \{0\} \times \{0, 4, 8\} \times Z_7$$
$$A_{16} = \{0\} \times \{0, 3, 9, 6\} \times Z_7$$
$$A_{17} = \{0\} \times \{0, 2, \ldots, 10\} \times Z_7$$
$$A_{18} = R \text{ and}$$
$$A_0 = \{0\}$$

Thus we get $S = \{A_0, A_1, \ldots, A_{18}\}$ the collection of S-subrings of R. This is easily verified to be also S-ideals of R. The lattice representation of them is as follows:

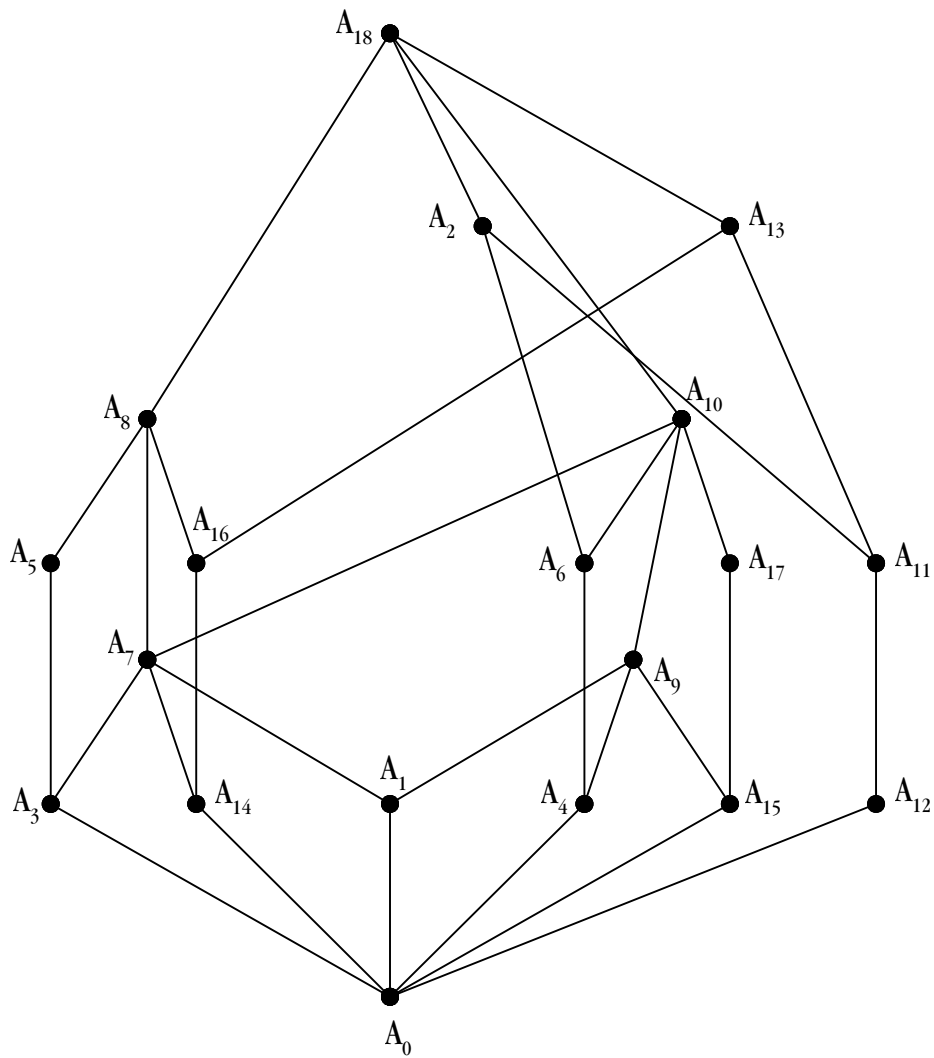

Figure 4.5.2

The set N = {{0}, A₃, A₇, A₁₅, A₁₀} forms a pentagon lattice which is non-modular.

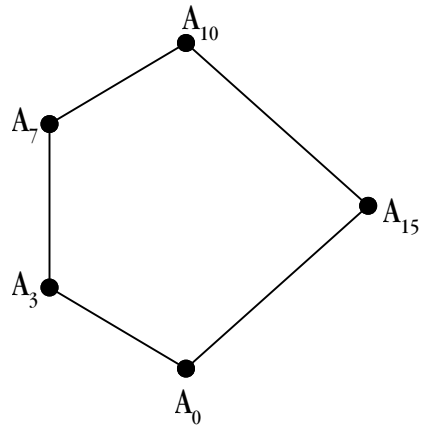

Figure 4.5.3

Hence N = {{0}, A₇, A₃, A₁₅, A₁₀} and forms a sublattice which is a pentagon lattice.

Thus in case of S-rings the set of S-ideals in general will not form a modular lattice which is a marked difference between ideals of a ring and S-ideals of a ring.

**DEFINITION [57]**: *Let L be a lattice. L is said to be a quasi distributive lattice if for all x, y, z, u in L, satisfies the following:*

*$(x \cup y) \cap (z \cup u) = \{x \cap (z \cup u)\} \cup \{y \cap (z \cup u)\} \cup \{z \cap (x \cup y)\} \cup \{u \cap (x \cup y)\}$ and $(x \cap y) \cup (x \cap u) = \{x \cup (z \cap u)\} \cap \{y \cup (z \cap u)\} \cap \{u \cup (x \cap y)\} \cap \{z \cup (x \cup y)\}$.*

***Example 4.5.7***: Let $Z_{12}$ = {0, 1, …, 11} be the ring. For S = H₀ = {0}, H₁ = {0, 6}, H₂ = {0, 3, 6, 9}, H₃ = {0, 4, 8}, H₄ = {0, 2, 4, 8, 6, 10} and H₅ = $Z_{12}$.

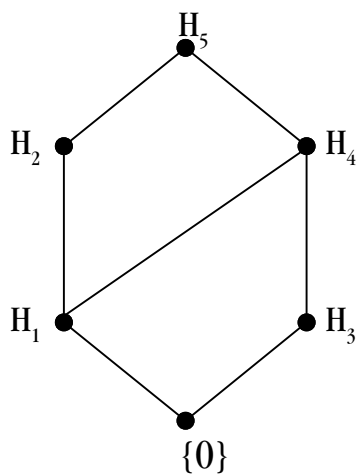

Figure 4.5.4

The reader is advised to verify whether the ideals form a quasi distributive lattice.

**Example 4.5.8**: The lattice of ideals of $Z_{16}$ is a quasi distributive lattice in fact a chain lattice. Left for the reader to draw the lattice diagram.

**Example 4.5.9**: Given the lattice diagram of the ring given by the S-mixed direct product of rings; R= $Z_3 \times Z_{12}$. The S-subrings of R are $A_0 = \{0\} \times \{0\}$, $A_1 = \{0\} \times \{0, 2, 4, 6, 8, 10\}$, $A_2 = \{0\} \times Z_{12}$, $A_3 = Z_3 \times \{0, 6\}$, $A_4 = Z_3 \times \{0, 4, 8\}$, $A_5 = Z_3 \times \{0, 3, 6, 9\}$, $A_6 = Z_3 \times \{0, 2, 6, 4, 8\ 10\}$, $A_7 = Z_3 \times Z_{12} = R$.

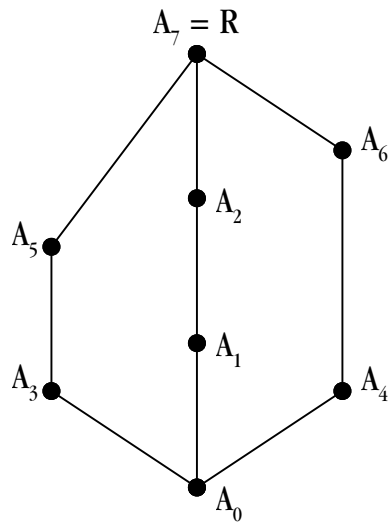

Figure 4.5.5

Test whether this lattice satisfies the quasi distributive identity

**DEFINITION [26]**: *An element $a \in R$, R a ring is called clean if it can be expressed as the sum of an idempotent and a unit in R.*

*A ring is called a clean ring if every element of R is clean.*

It has been shown by [26] that "If e is an idempotent in a ring R such that both eRe and $(1 - e)R(1 - e)$ are clean rings then R is a clean ring.

**DEFINITION 4.5.12**: *Let R be a ring. An element $a \in A$ where A is a S-subring of R is said to be a Smarandache clean (S-clean) element of R if it can be expressed as the sum of an idempotent and a unit in R. A ring R is called a Smarandache clean ring (S-clean ring) if every element of R is S-clean.*

*Thus for a ring R to be a S-clean ring it is sufficient if R has a S-subring A which is a S-clean ring. We do not demand the whole ring to be clean but it is localized.*

**Theorem 4.5.12**: *Let R be a ring. If R has a S-subring and R is a clean ring then R is a S-clean ring.*

*Proof*: Obvious by the very definitions of clean ring and S-clean ring.

**Example 4.5.10**: Let $Z_6 = \{0, 1, 2, \ldots, 5\}$. In this ring every element other than 0 and 1 are clean. But $Z_6$ is not a S-clean ring.

Further all clean rings R need not be S-clean for that ring R may not have a S-subring. Now we proceed onto define a concept viz. Smarandache strongly clean rings.

**Definition 4.5.13**: *Let R be a ring. We say $a \in R$ is a Smarandache strongly clean (S-s-clean) element if a can be written as a sum of a S-idempotent and a S-unit in R. If every $a \in R$ is S-s-clean then we call R a Smarandache strongly clean ring (S-strongly clean ring).*

Here it is important to note that R need not be a S-ring.

Further we leave the following theorems for the reader to prove.

**Theorem 4.5.13**: *If R is S-s-clean ring then it is a clean ring.*

**Theorem 4.5.14**: *Every S-s-clean ring need not be a S-clean ring.*

**Definition [53]**: Let S be a semigroup. S is integrally closed if $n\alpha \in S$ for some integer $n \in N$ implies $\alpha \in S$.

[53] has studied the integral closure of semigroup rings. We now proceed onto define Smarandache integrally closed rings.

**Definition 4.5.14**: *Let S be a S-semigroup. We say S is a Smarandache integrally closed semigroup (S-integrally closed semigroup) if S has a S-subsemigroup A which is such that whenever $n\alpha \in A$, n some integer and $\alpha \in A$.*

**Definition 4.5.15**: *Let R be a ring we say R is an integrally closed ring if R has a subset M, such that M is a multiplicatively closed semigroup and M is an integrally closed semigroup.*

**Definition 4.5.16**: *Let R be a ring. We say the ring R is Smarandache integrally closed (S-integrally closed) ring if R contains a S-semigroup M, $M \subset R$ under multiplication and S is an S-integrally closed semigroup.*

**DEFINITION [58]**: *A subset E of a semigroup S is called a system of local units if and only if the following conditions are satisfied:*

1. *E consists of commuting idempotents.*
2. *For any x ∈ S there exists e ∈ E such that xe = ex = x.*

We proceed onto define the concept of Smarandache local units of a semigroup, local units of a ring and Smarandache local units of a ring.

**DEFINITION 4.5.17**: *Let S be a S-semigroup. A subset M of S is called a Smarandache system of local units (S-system of local units) if and only if the following conditions are satisfied:*

1. *M consists of commuting S-idempotents*
2. *For any x ∈ S there exists e ∈ M such that ex = xe = x.*

**DEFINITION 4.5.18**: *Let R be a ring. A subset P of R is said to be a system of local units if and only if the following conditions are satisfied:*

1. *P consists of commuting idempotents.*
2. *For any r ∈ R there exists p ∈ P such that px = xp = x.*

We define Smarandache system of local units.

**DEFINITION 4.5.19**: *Let R be a S-ring. A subset M of R is called a Smarandache system of local units (S-system of local units) if and only if the following conditions are satisfied:*

1. *M consists of commuting S-idempotents.*
2. *For any s ∈ R there exists e ∈ M such se = es = s.*

The reader is advised to develop Morita equivalence on semigroups with systems of local units.

A ring (or semigroup) R is said to be π-regular if some power of every element is von Neumann regular. If a power of every element in R belongs to a subgroup of R, R is said to be uniformly π regular [8]. This paper [8] is a piece of nice research and the greatness of the paper lies in its extensive bibliography and the reader is advised to develop the Smarandache concepts about regularity.

**DEFINITION [162]**: *Let S be an additive subgroup of a finite ring R and suppose that either S is a subring or S is semisimple. Each of the following conditions are equivalent to S being a subideal of R.*

1. *S is a subideal of the ring generated by S and r for all r ∈ R*
2. *S ∩ T is a subideal of T for all 2-generator subrings T of R.*
3. *For all s ∈ S and r ∈ R there is a positive integer n such that both $(sr)^n$ and $r(sr)^n$ lie in S.*

**DEFINITION [162]**: *A subring S of R is said to be a subideal if there is a finite chain.*

*$S = R_m \subseteq R_{m-1} \subseteq \ldots \subseteq = R$ such that $R_i$ is an ideal of $R_{i-1}$ for i = 1, 2, …, m.*

Now we proceed onto define the Smarandache analogue.

**DEFINITION 4.5.20**: *Let A be an additive S-semigroup of a finite ring R and suppose that either $P \subset A$ (P a subgroup of A) is a S-subring or P is semisimple. Each of the following conditions are equivalent to P being a Smarandache subideal (S-subideal) of R.*

1. *P is a S-subideal of the ring generated by P and r for all r ∈ R.*
2. *P ∩ T is a S-subideal of T for all 2 generator S-subrings T of R.*
3. *For all s ∈ P and r ∈ R there is a positive integer n such that $(sr)^n$ and $r(sr)^n$ lie in P.*

**DEFINITION [163]**: *A ring R is called weakly periodic if for every x in R can be written x = a + b where a is nilpotent and b potent in the sense that $b^{n(b)}$ for some integer n(b) > 1.*

**DEFINITION 4.5.21**: *A S-subring A of R is said to be a S-subideal if there is a finite chain.*

*$A = R_m \subseteq R_{m-1} \subseteq \ldots \subseteq R$ such that $R_i$ is an S-ideal of $R_{i-1}$ for i = 1, 2, …, m.*

**DEFINITION 4.5.22**: *A S-ring R is called Smarandache weakly periodic (S-weakly periodic) if every x in R can be written in the form x = a + b where a is S-nilpotent and b potent in the sense that $b^{n(b)} = b$ for some integer n(b) > 1.*

The study of stable and stabilizer in rings was introduced earlier. Now we give here the concept of generalized stable ring as given by [12].

**DEFINITION [12]**: *Let R be an associative ring with 1. K(R) be the set {x ∈ R / there exists s, t in R such that sxt = 1}. The author defines R to be a generalized stable ring provided that aR + bR = R with a, b ∈ R implies a + by ∈ K(R) for some y ∈ R.*

**DEFINITION [12]**: *A ring R satisfies n-stable range condition for whenever $a_1$, $a_2$, ..., $a_{n+1} \in R$ with $a_1 R + ... + a_{n+1} R = R$ there exists elements $b_1$, $b_2$, ..., $b_n$ in R such that $(a_1 + a_{n+1} b_1)R + ... + (a_n + a_{n+1} b_n)R = R$.*

**DEFINITION [12]**: *R has stable range 1 if and only if whenever a, b $\in$ R with ab and ba strongly $\pi$-regular the Drazin inverses of ab and ba are conjugate via a unit of R.*

For more about stable range please refer [12].

Now when we replace R by a S-ring we get the corresponding results but for stable range 1 the Smarandache analogue is as follows:

**DEFINITION 4.5.23**: *Let R be a ring. R is said to be a Smarandache stable range I (S-stable range I) if and only if whenever a, b $\in$ R with ab and ba strongly $\pi$-regular the Drazin inverses of ab and ba are conjugate via a S-unit of R.*

Thus we request the reader to develop all these concepts and do research on Smarandache notions on $\pi$-regular elements as it has not been carried out in this book.

**PROBLEMS:**

1. Find a S-hyperring II of the ring $Z_{24}$.
2. Find the hyperring of $Z_{22}$.
3. Can the group ring $Z_{18}$ have S-hyperring II which is not S-hyperring I? Justify your claim.
4. Does the semigroup ring $Z_3 S(4)$ have
   i. Fixed support subring?
   ii. Fixed support subsemigroup?
   iii. S-fixed support subring?
5. Find a fixed support subring of $Z_2 S_4$.
6. Can $Z_5 S_4$ have S-fixed support subsemigroup?
7. Give an example of a ring which is semi-connected but not S-semi-connected.
8. Is the group ring $Z_4 S_7$
   i. semi-connected?
   ii. S-semi-connected? Justify your claim.
9. Give a non-trivial example of a $J_k$-ring.
10. Give an example of a S-$J_k$ ring which is not a $J_k$ ring.
11. Find the lattice of S-ideals and S-subrings for the ring $R = Z_8 \times Z_3 \times Z_{16} \times Z_7$.

12. For the S-mixed direct product ring $R = Z_{10} \times Z_7$ draw the lattice of S-ideals. Does it satisfy quasi-distributive lattice identity?

13. Give an example of a clean ring of order 18.

14. Show by an example a clean ring need not be a S-clean ring.

15. Is $Z_{20} = \{0, 1, 2, \ldots, 19\}$ the semigroup under multiplication
    - i. integrally closed?
    - ii. S-integrally closed?

16. Can the group ring $Z_3 S_5$ be
    - i. integrally closed?
    - ii. S-integrally closed?

17. Give an example of a semigroup which has system of local units.

18. A semigroup which can never have a system of local units (Will $Z^+$ under multiplication be a system of local units).

19. Can the ring $Z_5 S_3$ have system of
    - i. local units?
    - ii. S-local units?

20. Give an example of a weakly periodic ring.

21. Show by an example a S-ring which is not weakly periodic.

22. Is $Z_{40}$ a clean ring?

23. Can $Z_{25}$ be a S-clean ring?

24. Give an example of a S-clean ring which is not clean.

25. Give an example of a clean ring which is not S-clean.

**Chapter five**

# SUGGESTED PROBLEMS

This section is completely devoted to introducing several problems for researchers both in ring theory and Smarandache ring theory. Some of the problems are relatively simple and easily solvable whereas many problems can be treated as serious research problems. This chapter has 203 problems which are engrossing and an innovative researcher would certainly find them interesting.

Except for the classical zero divisor conjecture (problem) for group rings (1940) we have not repeated or included any of the open problems from other texts. Several problems are termed as characteristly by which we mean only to obtain a necessary and sufficient condition for the results to be true. If the student/ researcher has solved all problems at the end of each section in each chapter then certainly the researcher will not only find problems interesting but may solve them. As the reader is advised to have a good background in algebra in general and ring theory in particular. Several of the problems are characterization of S-group rings and S-semigroup rings.

Finally we state that this book gives importance to problems not only related to units, idempotents, zero divisors but those special elements like semi-idempotents, semi-units, super-idempotents etc. which indirectly guarantee the existence of units, zero-divisors, idempotents, etc. Likewise, not only study of ideals or subrings are studied and their related problems discussed but importance is given to substructure like additive/ multiplicative subgroups of ring, additive/ multiplicative semigroups and S-semigroups, S-ideals, S-subrings and so on.

Thus this chapter will be a boon to all researchers in Smarandache algebra and in ring theory.

**Problems:**

1. Obtain a necessary and sufficient condition for a unit to be a S-unit in a ring R (R is not a field).

2. Characterize those rings R in which every unit is a S-unit (R is not a field).

3. Characterize those rings R in which no unit of R is a S-unit of R.

4. Find conditions on the group G so that the group ring KG has S-units (K any field).

5. Find conditions on the semigroup ring FS so that FS has units which are S-units.

6. Characterize those fields of characteristic 0 in which every unit is not a S-unit.

7. Characterize those group rings in which every zero divisor is a S-zero divisor.

8. Characterize those rings which have zero divisiors but no S-zero divisors.

9. Characterize those semigroup rings.

      a.   In which every zero divisor is a S-zero divisor.
      b.   No zero divisor is a S-zero divisor.

10. Characterize those S-integral domains which are not integral domains.

11. Does there exist S-division rings which are not division rings?

12. Let $G$ be a torsion free non-abelian group, $R$ any field of characterize 0; can $KG$ have S-idempotents? (The existence of S-idempotents in $KG$ will force one to settle the problem of zero divisor conjecture, for it would imply the existence of zero divisor in $KG$ which is an open problem from the year 1940). This is an equivalent formulation of the zero divisor conjecture.

13. Characterize those rings in which every idempotent is a S-idempotent.

14. Characterize those rings in which no idempotent is a S-idempotent.

15. Prove in any ring $Z_n$ (Ring of integers modulo n). If a is a S-idempotent and b a S-co-idempotent of a then $a + b \equiv n \pmod{n}$.

16. Characterize those rings $R$ in which every S-idempotent has a unique S-co-idempotent.

17. Can $Z_{p^n} = \{0, 1, 2, \ldots, p^n - 1\}$ the ring of integers modulo $p^n$; p a prime, $n \geq 2$ have S-idempotents? Characterize them.

18. Let $Z_n$ be the ring of integers modulo n. $n = p_1 p_2 p_3$ ($p_1, p_2, p_3$ are 3 distinct primes) of which atleast one is an even prime.

      a.   Can $Z_n$ have only 6 idempotents of which 5 are S-idempotents?
      b.   If $p_1$, $p_2$ and $p_3$ are all odd primes, can we prove $Z_n$ has 6 idempotents all of which are S-idempotents?

19. Find the number of idempotents in $Z_n = \{0, 1, 2, \ldots, n-1\}$ (where $n = p_1^{\alpha_1} p_2^{\alpha_2} \ldots p_m^{\alpha_m}$, $p_i$ are distinct primes; $\alpha_I > 1$) which are S-idempotents.

20. Let $Z_n G$ be the group ring of the group G over the ring $Z_n$.

        a. Characterize those group rings $Z_n G$ (by giving conditions on $Z_n$ and/or on G) so that every idempotent in $Z_n G$ is an S-idempotent of $Z_n G$.
        b. No idempotent in $Z_n G$ is a S-idempotent.

21. Let $Z_n S$ be the semigroup ring of a semigroup S over the ring $Z_n$

        a. Characterize those semigroup so that in $Z_n S$, every idempotent is a S-idempotent.
        b. No idempotent in $Z_n S$ is a S-idempotent.

22. Characterize those rings R in which every ideal is a S-ideal I.

23. Characterize those rings in which no ideal is a S-ideal I.

24. Obtain conditions on rings R so that the concept of S-ideal I and S-subring I coincide.

25. Can we say the concept of S-rings I and S-rings II coincide on finite rings?

26. Can we say all rings $Z_n$ are both S-ring I and S-ring II?

27. Characterize all rings R which are S-ring II and not S-ring I. (Note – Don't take Z or Z[x] or matrices over Z or direct product of Z).

28. Determine those S-rings which have only S-pseudo ideals and no S-ideals I or S-ideals II.

29. Describe mathematically those group rings which are not

        a. S-simple rings I.
        b. S-simple rings II.
        c. S-pseudo simple rings.

30. Characterize those rings R in which an S-module I M, related to a field $F \subset R$ is S-module I for all fields contained in R.

31. Do there exist rings R in which all modules M which is a S-module I is also a S-module II?

32. In a ring R, can the concept of S-module II coincide with S-pseudo module?

33. Does there exist rings for which S-module II can never be defined? If such rings exists, characterize them.

34. Characterize those class of rings which are S-strong right S-rings.

35. Characterize those class of rings which are precisely S-strong right D-rings.

36. Does there exist a two sided ideal of order $p^{n!/2}$ in the group ring $Z_p S_n$?

        a. When p is a prime and p/n!
        b. When p is a prime and (p, n!) = 1.
        c. When p is a composite number and (p, n!) = 1.
        d. When p is a composite number and (p, n!) = d.

(Remark. The group ring $Z_2 S_3$ has no two sided ideals of order 8 but has a right ideal of order 8. Further the group ring $Z_2 S_3$ has two sided ideals of order 2,4,16 and 32 where the order of $Z_2 S_3$ is 64 but has no two sided ideals of order 8). Now we propose the following:

37. Does there exist a two-sided S-ideal I and II of order $p^{n!/2}$ in the group ring $Z_p S_n$?

The conditions mentioned as four cases in problems 36 should be discussed in the case of S-ideals.

38. Study the same problems given in 36 and 37 in case of right ideals, S-right ideals I (II), left ideals and S-left ideals I (II).

39. Characterize those group rings which are S-J-rings. [we see $Z_2 S_3$ is a S-J-rings can we say all rings $Z_p S_n$ will be S-J-ring. Find conditions on p and n, so that the group ring $Z_p S_n$ is a S-J-ring].

40. Obtain conditions on the semigroup S and on the ring R so that the semigroup ring RS is a S-J-ring.

41. Let $Z_n$ be the ring of integers modulo n. $S_m$ be the symmetric group of degree m. Let $Z_n S_m$ be the group ring of the group $S_m$ over $Z_n$. Is $Z_n S_m$ a S-strong subring? Discuss the cases

        a. When n= m = p (p a prime).
        b. When n= m = non-prime.
        c. When (n,m) = 1.
        d. When (n,m) = p, n > m p − prime.
        e. When (n,m) = d, d any non-prime integer.

42. Study problem 41 for S-ideals; i.e., is $Z_n S_m$ a S-strong ideal ring?

43. When will $Z_n S_m$ be a S-strong subring ideal ring?

44. Find conditions on $Z_n$ and $S_m$, so that the group ring $Z_n S_m$ is a S-weakly Boolean ring.

45. Let $Z_n$ be the ring of integers modulo n which is a S-ring. Characterize those rings $Z_n$ which are S-weakly Boolean.

(Hint: $Z_6 = \{0,1,2,3,4,5\}$ is a S-ring which is not a S-weakly Boolean ring.)

The group ring $Z_{15}$ G where G=$\{g/g^2=1\}$, $Z_{15} = \{0, 1, 2, 3, \ldots, 14\}$ is a S-weakly Boolean ring for it has the S-subring BG where B=$\{0,5,10\}$.

46. Characterize those group rings $Z_n S_m$ which are S-right multiplication ideal ring.

47. Characterize those semigroup rings $Z_n S(m)$ (S(m) − symmetric semigroup) which are S-right multiplication ideal ring. Discuss atleast the cases

        a. When m = n,
        b. (m, n) = 1,
        c. n odd prime, m a multiple of n.
        d. (m, n) = d (d, a non-prime).

48. Characterize those semigroup rings $Z_n S(m)$ which are S-weakly Boolean rings?

49. Characterize those group rings $Z_n S_m$ which are S-pseudo commutative.

50. Characterize those group rings $Z_n S_n$ which have a pair which is S-pseudo commutative relative to a S-subring.

(Hint: Should discuss when n is a prime and n is a non-prime m is prime or a power of a prime and m a non-prime and finally (n, m) = n or (n, m) = m, (n, m) = 1 and (n, m) = d, d < n, d < m).

51. Characterize those semigroup rings $Z_n S(m)$ which are S-pseudo commutative.

52. Do we have semigroup rings $Z_n S(m)$ which has a pair which is S-pseudo commutative relative to a S-subring? Classify and characterize those semigroup rings.

53. Let $Z_n$ be a prime field or a ring of integers modulo n. $S_m$ be the symmetric group of degree n. Characterize the group ring $Z_n S_m$ so that it is

        a.  Strictly right chain ring.
        b.  Chain ring.

(Hint: Use conditions on n and m as (n, m)=1, n prime, m composite, (m, n) = 1 both m and n prime n ≠ m, (n, m) = n or m).

54. Characterize those rings which are ideally obedient rings.

55. Characterize those rings which do not have any obedient ideals.

56. Characterize those rings which are S-ideally obedient rings.

57. Characterize those ring which has

        a.  S-obedient ideals.
        b.  which has S-ideals but none of them are S-obedient ideals.

58. Characterize those group rings which are S-strongly clean rings.

59. Characterize those semigroup rings which are

        a.  S-ideally obedient rings.
        b.  Ideally obedient rings.
        c.  S-clean rings.
        d.  Clean rings.

60. Characterize those group rings

        a.  Which are Lin rings.
        b.  Which are S-Lin rings.

61. Characterize those semigroup rings

        a. Which are Lin ring.
        b. Which are not Lin rings.

62. Characterize those rings R which are S-Lin rings (R not a group rings a semigroup ring).

63. Find the class of group rings KG (by giving conditions on the group G and on the field or ring K) so that they are S-Lin rings.

64. Find the class of semigroup rings FS (S a semigroup, F a commutative ring with 1 or a field) which are S-Lin rings.

65. Let G be torsion free non-abelian group and F any field or a commutative ring with 1. Can the group ring FG satisfy super ore conditions?

66. Characterize those rings which does not satisfy super ore condition but satisfies S-super ore condition.

67. Obtain a necessary and sufficient condition for a group ring FG to be an ideally strong ring.

68. Characterize those semigroup rings RS where R is a ring and S a semigroup which are
        a. Ideally strong.
        b. Which are not ideally strong.

69. Characterize those rings which are S-ideally strong rings.

70. Find conditions on the group G and on the ring R so that the group ring RG is both ideally strong and a S-ideally strong ring.

71. Study the same problem of characterization when the group G is replaced by a semigroup.

72. Is $Z_p S_m$ a Q-ring? Discuss all the possible cases when

        a. p a prime.
        b. p/m.
        c. (p, q) = 1.
        d. p = 2 and m a prime.

73. Will the semigroup ring $Z_p S_m$ be a Q-ring

        a. p a prime m non-prime (p, m) = 1.

b. (p, m) = d (p need not be a prime).

c. p = 2 m any integer.

d. (p, q) = 1 p and q primes ?

74. Characterize those rings which are Q-rings.

75. Characterize those group rings which are S-Q-rings.

76. Characterize those semigroup rings which are both Q-rings and S-Q-rings.

77. Can the ring $M_{n \times n} = \{ (a_{ij}) / a_{ij} \in Z_m$, m not a prime$\}$ be

  a. a Q-ring (m = n and m ≠ n)?

  b. a S-Q-ring (m = n, and m ≠ n)?

78. Is the group ring $Z_p S_n$ a F-ring, p an odd prime?

79. Classify those group rings which are S-F-rings.

80. Can the semigroup ring $Z_p S(n)$ be a F-ring, p an odd prime?

81. Classify those semigroup rings which are S-F rings and F-rings.

82. Can a group ring KG where K is a field of characteristic 0 and G a torsion free non-abelian group have semi nilpotent elements. (The solution to this problem is equivalent to the zero divisor conjecture in group rings proposed in 1940; can KG have zero divisors if G is a torsion free non-abelian group?)

83. Let $Z_p$ be the prime field of characteristic p (p > 2) and $G = \langle g / g^q = 1 \rangle$ be a cyclic group of order q.

  a. If (p, q) = 1, can the group ring $Z_p G$ have nontrivial semi nilpotent elements?

  b. If p/q, can the group ring $Z_p G$ have nontrivial semi nilpotent elements?

(We see the group ring $Z_2 G$ has semi nilpotents where $G = \langle g/g^2 = 1 \rangle$. We also observe that if $G = \langle g/g^3 = 1 \rangle$ then the group ring $Z_2 G$ has no nontrivial semi nilpotents).

84. Let $Z_n = \{0, 2, 3, \ldots, n - 1\}$ be the ring of integers modulo n. Find conditions on n so that $Z_n$ is a SSS ring. Find the maximum number of SS-elements in $Z_n$. (Hint: $Z_9 = \{0, 1, 2, \ldots, 8\}$ has 4-SSS elements or 2 pairs of SS elements viz (3,6) and (5,8)).

85. Can we have a ring R in which all elements in R \ {0,1} are SS-elements?

86. Let G be a torsion free non-abelian group, K any field. Can KG have nontrivial S-subrings?

87. Let G be a torsion free non-abelian group and K any field of characteristic zero. Can KG be

        a. Locally semiunitary?
        b. Locally unitary?

88. Characterize those group rings KG which are

        a. S-locally unitary.
        b. S-locally semiunitary.
        c. Locally semiunitary.

89. Characterize those semigroup rings KS which are

        a. S-locally unitary.
        b. S-locally semiunitary.
        c. Locally unitary.

90. Let G be a torsion free non-abelian group and K a field of characteristic 0.

        a. Can KG be S-semiunitary?
        b. Can KG be S-unitary?

91. Let G be a torsion free group and K any field. Can the group ring KG be a

        a. CN ring?
        b. Weakly CN ring?

92. Characterize those group rings KG which are

        a. S-weakly CN ring.
        b. S-CN ring.
        c. Weakly CN-ring.

93. Characterize those group rings which are

        a. Tight rings.
        b. S-Tight rings.

94. Characterize those rings R which are

a. Never a tight ring.
b. Never a S-tight ring.
c. S-tight ring.

95. Let G be a torsion free non-abelian group and K any field of characterize 0. Can the group ring KG be a $\gamma_n$-ring?

96. Let G be a torsion free non-abelian group and K any field of characterize 0. Can the group ring KG be a S-$\gamma_n$-ring?

97. Characterize those semigroup rings KS that are S-$\gamma_n$-rings.

98. What is the condition on the group G and on the field K so that the group ring KG has

a. S-demi subrings?
b. Demi subrings?

(Study the above problem in case of semigroup rings KS).

99. Give any nice characterization theorem for the S-demimodules to exists for the group ring KG.

100. Let P={$\alpha_0 + \alpha_1 i + \alpha_2 j + \alpha_3 k / \alpha_0, \alpha_1, \alpha_2, \alpha_3 \in Z_n$; n a composite number}; Can we have P=G $\cup$ V?

101. Can P have idempotents other than the ones given in $Z_n$?

102. Can we say by using $Z_p$, p a prime we can construct P which is a finite division ring of dimension $p^4$, p an odd prime?

103. Use P given in problem 100 to define some interesting and innovative Smarandache notions on rings.

104. Characterize those S-group rings which are S-Artinian.

105. Characterize those S-semigroup rings which are S-Artinian.

106. Determine the class of group rings which are S-Noetherian.

107. Does there exist a class of semigroup rings which are S-Noetherian?

108. Let the group ring RG be any S-ring. Depending on G and on R is it possible to find the number of proper subsets of RG which are fields. (Hint: Study this in case of $Z_nS_m$ and $QS_m$).

109. Let RS be a S-semigroup ring. R any field and S a S-semigroup.

    a. Does there exist a method by which the number of proper subsets which are fields can be found out?
    b. Characterize those S-semigroup rings which has no proper subset which is a field (For ZS(3) is S-semigroup ring which has no proper set which is a field).

110. Characterize the ideals in the ring $Z_n$, n a composite number so that $Z_n$ has

    a. S-ideals.
    b. S-A.C.C. condition is satisfied.
    c. no S-ideals.

111. Characterize those group rings $Z_nS_n$ which satisfy S-A.C.C. condition.

112. Characterize those semigroup rings $Z_nS(m)$ which satisfy S-A.C.C. II on ideals.

113. Characterize those group rings (semigroup rings) which satisfy both A.C.C. and S.A.C.C. II.

114. Characterize those group rings (semigroup rings) which satisfy both S-D.C.C. II and D.C.C.

115. Does there exist a ring R which cannot be S-quasi ordered?

116. Let K be any field. G a torsion free non-abelian group. Can KG have

    a. non-trivial semi nilpotent ideals?
    b. non-trivial S-semi nilpotent ideals?

117. Can the group ring $Z_pG$ (G-finite group) have non-trivial S-semi nilpotent ideals?

118. Can the semigroup ring $Z_pS(n)$ have non-trivial S-semi nilpotent ideals?

119. Let $Z_pG$ be the group ring and G is a p-group. Does $Z_pG$ have non-trivial S-semi nilpotent ideals not including $\omega(Z_pG)$?

120. Characterize those rings which are

a. S-subsemiideal rings.
        b. subsemiideal rings.

121. Give a necessary and sufficient condition for the group ring to be

        a. S-subsemi ideal rings.
        b. subsemi ideal rings

122. Characterize those semigroup rings which are

        a. S-subsemi ideal rings.
        b. subsemi ideal rings.

123. Obtain conditions on a ring R so that every subsemi ideal ring is a S-subsemi ideal.

124. Can KG the group ring where K is a field and G a torsion free non-abelian group have non-trivial super-idempotents?

125. Obtain conditions on the group G and on the ring R so that the group ring RG has non-trivial super idempotents.

126. Does the group ring $Z_pG$ when G is a cyclic group of order q have non-trivial super-idempotents when:

        a. $(p, q) = 1$.
        b. $p/q$.
        c. $p = q$.
        d. $(p, q) = d$.

127. Can the semigroup ring $Z_pS(n)$ have non-trivial S-super-idempotents?

128. Study problems 124 to 126 in the context of S-super idempotents.

129. Characterize those normal rings which are not S-normal rings.

130. Characterize those normal rings which are S-normal rings.

131. Classify those group rings which are

        a. normal rings.
        b. S-normal rings.

132. Characterize those semigroup rings which are

a. normal rings.
b. S-normal rings.

133. Characterize those group rings which are

a. S-SI rings.
b. SI rings.

134. Study problem 133 in case of semigroup rings.

135. Classify those SI-rings which are S-SI-rings and those S-SI-rings which are SI-rings.

136. Can KG where K is a field and G a torsion free non-abelian group be a

a. n-c-s-ring.
b. S-n-c-s ring.

137. Characterize those group rings $Z_n S_m$ which are

a. n-c-s ring.
b. S-n-c-s ring.

by varying n and m as

a. $(n, m) = 1$.
b. $n/m$.
c. $(n, m) = d$.
d. n prime, m prime.
e. n prime, m a non-prime.

138. Characterize those semigroup rings $Z_n S(m)$ which are

a. n-c-s rings.
b. S-n-c-s rings.

under the conditions mentioned in problem 137.

139. Let G be a cyclic group of order p, p a prime and $Z_g$ be the prime field of characteristic q such that $(p, q) = 1$. Is $Z_p G$ an iso-ring? Can $Z_p G$ be a co-ring? If p/q will $Z_p G$ be an iso-ring and co-ring? If q is not a prime will $Z_p G$ be a co-ring and an iso ring?

140. Obtain conditions for the group ring $Z_pS_n$ to be a

      a.  iso-ring.
      b.  co-ring.
      c.  not a co-ring.
      d.  not an iso-ring.

(by imposing conditions p and n).

141. Characterize those semigroup rings $Z_nS(m)$ which are

      a.  S-co-rings.
      b.  S-weak co-rings.
      c.  S-iso-rings.
      d.  S-weak iso-rings.

142. Find condition on n and m so in the group ring $Z_nS_m$ we have

      a.  S-co-rings.
      b.  S-weak co-rings.
      c.  S-iso-rings.
      d.  S-weak iso-rings.

143. Let K be a field of characteristic zero and G a torsion free non-abelian group. Can KG the group ring be

      a.  S-e-primitive?
      b.  At least S-weakly e-primitive?
      c.  e-primitive?
      d.  At least weakly e-primitive?

144. Characterize those group rings $Z_pG$ which are

      a.  S-weakly e-primitive
      b.  weakly e-primitive.
      c.  e-primitive.
      d.  S-e-primitive.

145. Characterize those semigroup rings $Z_pS(m)$ which are

      a.  e-primitive.
      b.  S-e-primitive.
      c.  S-weakly e-primitive.
      d.  weakly e-primitive.

146. Obtain conditions on the ring R and the group G so that the group ring RG can have

        a. SV-group.
        b. Weakly SV-group.
        c. S-SV-group.
        d. S-Weakly SV-group.

147. Let R be a ring and S a semigroup. Obtain conditions on R and S so that the semigroup ring is a

        a. S-SV group and
        b. not a SV group.

148. Let K be a field and G be a torsion free non-abelian group

        a. Can KG have SV − group?
        b. Can KG have WSV groups?
        c. S-SV-groups?
        d. S-weakly SV-groups?

149. Let K be a field of characteristic 0 and G be a torsion free non-abelian group. Can KG have S-radix? Can KG have radix?

150. Classify those rings R in which no radix is a S-radix.

151. Classify those rings which has no S-radix and no radix.

152. Classify those rings in which every $\delta$-semigroup (under multiplication) is a S-$\delta$-semigroup.

153. Characterize those group rings which are

        a. SG-ring.
        b. weakly SG-ring.
        c. S-SG-ring.
        d. S-weakly SG-ring.

154. Study and characterize those semigroup rings which are SG-rings/ S-SG rings/ S-weakly SG ring/ weakly SG-rings.

155. Classify those rings which are SG-rings and also S-SG rings.

156. Let KG be the group ring of the torsion free non-abelian group G and K the field of characteristic 0. Can KG be a

     a. ZI ring?
     b. S-ZI ring?
     c. Weakly ZI ring?
     d. S-weakly ZI ring?
     e. pseudo ZI ring?

157. Characterize or classify all semigroup rings $Z_n S(P)$ that are

     a. pseudo ZI ring.
     b. S-pseudo ZI ring.
     c. ZI ring.
     d. S-ZI ring.

158. Let K be a field of characteristic 0 and G a torsion free non-abelian group. Can the group ring KG have non-empty square sets?

159. Let K be a prime field of characteristic zero and G a torsion free non-abelian group. Can the group ring KG have insulators?

160. Let $Z_n$ be the ring of integers modulo n and G be an abelian group. When does the group ring $Z_n G$ have n-capacitor groups

     a. If $(|G|, n) = 1$.
     b. If $(n, |G|) \neq 1$.
     c. If n/|G|.

161. Characterize those group rings (semigroup rings) in which all n-capacitor groups are S-n-capacitor groups.

162. If R is a ring without nilpotent elements of order 2. Does it imply R is trisimple?

163. $Z_p G$ be the group ring where $G = \langle g \, / \, g^n = 1 \rangle$. Can $Z_p G$ be trisimple if $(p, n) = 1$, $(p, n) = d$ and $(p, n) = p$?

164. Characterize those group rings and semigroup rings which are

     a. Trisimple.
     b. S-trisimple.
     c. S-semi trisimple.

165. Can Z be a S-semi-order ring?

166. Give a complete characterization of group rings (semigroup rings) which are

      a.  so-ring.
      b.  S-so-ring.

167. Let $Z_p$ be the prime field of characteristic p. $G = \langle g/ g^q = 1 \rangle$. $Z_pG$ be the group ring. For what values of p and q will the group ring $Z_pG$ be a

      a.  Square ideal ring?
      b.  S-square ideal ring?
      c.  S-n-ideal ring?
      d.  n-ideal ring?

168. Let G be a torsion free non-abelian group and K any field. Can the group ring KG be a

      a.  Square ideal ring?
      b.  S-square ideal ring?
      c.  S-n-ideal ring?
      d.  n-ideal ring?

169. Characterize those group rings (semigroup rings or rings) in which all square ideals are S-square ideals and all n-ideals are S-n-ideals.

170. Let G be a torsion free non-abelian group. K any field of characteristic 0 or p. Can the group ring KG be a n-like ring for any n?

171. Let $Z_p$ be a prime field of characteristic p, $p \neq 2$ and $G = \langle g / g^q = 1 \rangle$. Is the group ring $Z_pG$ a n-like ring

      a.  when (p, q) = 1?
      b.  p = q?
      c.  p is a multiple of q or q is multiple of p.

172. Let G be a finite group. K any field of characteristic 0. Can KG the group ring be a

      a.  TI ring?
      b.  S-TI ring?

173. Let G be a torsion free group and K any field characteristic zero or p. Can the group ring KG be a

      a.  TI ring?

b. S-TI-ring?

174. Characterize those rings which are power joined is also S-power joined.

175. Obtain a necessary and sufficient condition for the (m, n) (uniformly) power joined ring to be S-(m, n) (uniformly) power joined.

176. Characterize those group rings and semigroup rings which

        a. has quasi nilpotent ideals.
        b. S-quasi nilpotent ideals

177. Characterize those group rings (semigroup rings) which have

        a. radical ideals.
        b. S-radical ideals.

178. Characterize those group rings (semigroup rings) in which

        a. radical ideal coincides with upper radical.
        b. S-radical ideals which coincides with S-upper radical.

179. Does there exist a method by which we can find whether the ring contains at least a

        a. related pair?
        b. S-related pair?

180. For what group G, the group ring QG (Q the field of rationals) has related pairs.

181. Can we ever find a ring R in which subring link relation happens to be an equivalence relation?

182. Can reals or ring of integers have pairs which are subring linked?

183. Characterize those group rings (semigroup rings) in which at least a pair can be

        a. Subring linked.
        b. S-subring linked.

184. Characterize those rings in which both the notions of

        a. stable and stabilized pair of subrings coincide.
        b. S-stable and S-stabilized pair of subring coincides.

185. Classify those rings which are

        a.   stable rings.
        b.   S-stable rings.

186. Can a torsion free non-abelian group be conditionally commutative?

187. Let G is a conditionally commutative group and R a conditionally commutative ring. Can the group ring RG be a conditionally commutative ring?

188. Let RG be a group ring; obtain a necessary and sufficient condition so that the group ring RG is a generalized Hamiltonian ring.

189. Characterize those semigroup rings RS so that they are generalized Hamiltonian rings.

190. Classify those group rings and semigroup rings which are

        a.   S-Hamiltonian.
        b.   S-Hamiltonian II.

191. Classify those rings in which every S-Hamiltonian II is also S-Hamiltonian I.

192. Characterize those group rings (semigroup rings) which has fixed support subring.

193. Classify those group rings (semigroup rings) which has fixed support semigroup.

194. Classify those group rings (semigroup rings) which have S-fixed support subring.

195. Classify those group rings (semigroup rings) which are

        a.   semi connected.
        b.   S-semi-connected.

196. Classify those rings which are

        a.   $J_k$ ring.
        b.   S-$J_k$ ring.

197. Find conditions on the ring so that every S-$J_k$-ring is a $J_k$ ring and every $J_k$ ring is a S-$J_k$-ring.

198. Find conditions on the ring R to have the S-ideals I and S-ideals II to be modular.

199. Find conditions on the ring so that all S-ideals II are S-subrings.

200. Find those S-rings whose S-subrings forms a quasi distributive lattice.

201. Characterize those S-mixed direct product rings using only modulo rings which has the

       a.  lattice of S-ideals to be a quasi distributive lattice.
       b.  lattice of S-subrings to be a quasi distributive lattice.
       c.  those rings in which all S-ideals form a modular lattice.

202. Does there exist a ring R in which every triple $x, y, z \in R \setminus \{0\}$ satisfies the identity $x^n + y^n = z^n$; $n > 1, 2$.

203. Find Smarandache analogue of classical theorems in ring theory.

# REFERENCES

1.  Abian Alexander and McWorter William, *On the structure of pre p-rings*, Amer. Math. Monthly, Vol. 71, 155-157, (1969).

2.  Adaoula Bensaid and Robert. W. Vander Waal, *Non-solvable finite groups whose subgroups of equal order are conjugate*, Indagationes Math., New series, No.1(4), 397-408, (1990).

3.  Allan Hayes, *A characterization of f-ring without non-zero nilpotents*, J. of London Math. Soc., Vol. 39, 706-707, (1969).

4.  Allevi.E, *Rings satisfying a condition on subsemigroups*, Proc. Royal Irish Acad., Vol. 88, 49-55, (1988).

5.  Andruszkiewicz.R, *On filial rings*, Portug. Math., Vol. 45, 136-149, (1988).

6.  Atiyah.M.F and MacDonald.I.G, *Introduction to Commutative Algebra*, Addison Wesley, (1969).

7.  Aubert.K.E and Beck.L, *Chinese rings*, J. Pure and Appl. Algebra, Vol.24, 221-226, (1982).

8.  Bogdanovic Stojan, Ciric Miroslov and Petkovic Tatjana, *Uniformly $\pi$-regular semigroups: a survey*, Zb. Rad., Vol. 9, 5-82, (2000).

9.  Bovdi Victor and Rosa A.L., *On the order of unitary subgroup of a modular group algebra*, Comm. in Algebra, Vol. 28, 897-1905, (2000).

10. Chen Huanyin, *On generalized stable rings*, Comm. in Algebra, Vol. 28, 1907-1917, (2000).

11. Chen Huanyin, *Exchange rings having stable range one,* Inst. J. Math. Sci., Vol. 25, 763-770, (2001).

12. Chen Huanyin, *Regular rings with finite stable range*, Comm. in Algebra, Vol. 29, 157-166, (2001).

13. Chen Jain Long and Zhao Ying Gan, *A note on F-rings*, J. Math. Res. and Expo. Vol. 9, 317-318, (1989).

14. Connel.I.G, *On the group ring*, Canad. J. Math. Vol. 15, 650-685, (1963).

15.   Corso Alberto and Glaz Sarah, *Guassian ideals and the Dedekind- Merlens lemma*, Lecture notes in Pure and Appl. Math., No. 217, 113-143, Dekker, New York, (2001).

16.   Cresp.J and Sullivan.R.P, *Semigroup in rings*, J. of Aust. Math. Soc., Vol. 20, 172-177, (1975).

17.   De Bruijin.N.G and Erdos.P, *A colour Problem for infinite graphs and a problem in theory of relations*, Nederl Akad. Wetensch Proc., Vol. 54, 371-373, (1953).

18.   De Maria Francesco, *Hypergroups ($Z_q$, q)*, Riv. Math. Pure App., No.5, 49-58, (1989).

19.   Douglas Costa and Gordon Keller, *Normal subgroups of SL (2, A),* BAMS., Vol. 24, 131-135, (1991).

20.   Dubrovin.N.I, *Chain rings*, Rad. Math., Vol. 5, 97-106, (1989).

21.   Dutta. T.K, *On generalized semi-ideal of a ring*, Bull. Calcutta Math. Soc., Vol. 74, 135-141, (1982).

22.   Erdogdu.V, *Regular multiplication rings*, J. Pure and Appl. Algebra, Vol. 31 , 55-59 (1989).

23.   Grater Joachim, *Strong right-D-domains*, Monatsh Math., Vol. 107, 189-205, (1989).

24.   Gray.M, *A Radical Approach to Algebra*, Addison Wesley, (1970).

25.   Gupta.V, *A generalization of strongly regular rings*, Acta. Math. Hungar. , Vol. 43, 57-61, (1984).

26.   Han Juncheol and Nicholson W.K., *Extensions of clean rings*, Comm. in Algebra, Vol. 29, 2589-2595, (2001).

27.   Han Yang, *Strictly wild algebras with radical square zero*, Arch-Math. (Basel), Vol. 76, 95-99, (2001).

28.   Hannah John, *Regular bisimple rings*, Proc. Edin. Math. Soc., Vol. 34, 89-97, (1991).

29.   Herstein.I.N, *Topics in Algebra*, John Wiley and Sons, (1964).

30. Herstein.I.N, *Non-commutative Rings*, M.A.M, (1968).

31. Herstein.I.N, *Topics in Ring theory*, Univ. of Chicago Press, (1969).

32. Higman.D.G, *Modules with a group of operators*, Duke Math. J., Vol. 21, 369-376, (1954).

33. Higman.G, *The units of group rings*, Proc. of the London Math. Soc., Vol. 46, 231-248, (1940).

34. Hirano Yasuyuki and Suenago Takashi, *Generalizations of von Neuman regular rings and n-like rings*, Comment Math. Univ. St. Paul., Vol. 37, 145-149, (1988).

35. Hirano Yasuyuki, *On π-regular rings with no infinite trivial subring*, Math. Scand., Vol. 63, 212-214, (1988).

36. Iqbalunnisa and Vasantha W.B., *Supermodular lattices*, J. Madras Univ., Vol. 44, 58-80, (1981).

37. Istvan Beck, *Coloring of commutative rings,* J. of Algebra, Vol. 116, 208-226, (1988).

38. Jacobson.N, *Theory of rings*, American Mathematical Society, (1943).

39. Jacobson.N, *Structure of ring*, American Mathematical Society, (1956).

40. Jin Xiang Xin, *Nil generalized Hamiltonian rings*, Heilongiang Daxue Ziran Kexue Xuebao, No. 4, 21-23, (1986).

41. Johnson P.L, *The modular group ring of a finite p-group*, Proc. Amer. Math. Soc., Vol. 68, 19-22, (1978).

42. Katsuda.R, *On Marot Rings*, Proc. Japan Acad., Vol. 60, 134-138, (1984).

43. Kim Nam Kyin and Loe Yang, *On right quasi duo-rings which are π-regular,* Bull. Korean Math. Soc., Vol. 37, 217-227, (2000).

44. Kishimoto Kazuo and Nagahara Takas, *On G-extension of a semi-connected ring*, Math. J. of Okayama Univ., Vol. 32, 25-42, (1990).

45. Krempa.J, *On semigroup rings*, Bull. Acad. Poland Sci. Ser. Math. Astron. Phy., Vol. 25, 225-231, (1977).

46.     Lah Jiang, *On the structure of pre J-rings*, Hung-Chong Chow, 65[th] anniversary volume, Math. Res. Centre, Nat. Taiwan Univ., 47-52, (1962).

47.     Lang.S, *Algebra*, Addison Wesley, (1984).

48.     Ligh.S and Utumi.Y, *Direct sum of strongly regular rings and zero rings*, Proc. Japan Acad., Vol. 50, 589-592, (1974).

49.     Lin Jer Shyong, *The structure of a certain class of rings*, Bull. Inst. Math. Acad. Sinica, Vol. 19, 219-227, (1991).

50.     Louis Dale, *Monic and Monic free ideals in a polynomial semiring*, Proc. Amer. Math. Soc., Vol.56, 45-50, (1976).

51.     Louis Halle Rowen, *Ring theory*, Academic Press, (1991).

52.     Mason.G, *Reflexive ideals*, Comm. in Algebra, Vol. 17, 1709-1724, (1981).

53.     Matsuda Ryoki, *Note on integral closure of semigroups*, Tamkang J. Math., Vol. 31, 137-144, (2000).

54.     Mc Coy.N.H, and Montgomery.D, *A representation of generalized Boolean Rings*, Duke Math. J. Vol. 3, 455-459, (1937).

55.     Mitsch.H, *Subdirect products of E-inverse semigroups*, J. Aust. Math. Soc., Vol. 48, 66-78, (1990).

56.     Munn.W.D., *Bisimple rings*, Quat. J. of Math. Oxford, Vol. 32, 181-191, (1981).

57.     Negru. I.S, *Quasi distributivity of* lattices, Math. Issled, No. 66, 113-127, (1982).

58.     Neklyudova V.V. Morita, *Equivalence of semigroups with systems of local units,* Fundam. Prikl. Math., Vol. 5, 539-555, (1999).

59.     Northcott.D.G, *Ideal theory*, Cambridge Univ. Press, (1953).

60.     Padilla Raul, Smarandache Algebraic structures, Bull. of Pure and Appl. Sci., Vol. 17E, 119-121, (1998).

61.     Passman.D.S, *Infinite Group Rings*, Pure and Appl. Math., Marcel Dekker, (1971).

62.  Passman.D.S, *The Algebraic Structure of Group Rings*, Inter-science Wiley, (1977).

63.  Peric Veselin, *Commutativity of rings inherited by the location of Herstein's condition*, Rad. Math., Vol. 3, 65-76, (1987).

64.  Pimenov K.L. and Yakovlev. A. V., *Artinian Modules over a matrix ring, Infinite length modules*, Trends Math. Birkhauser Basel, Bie. 98, 101-105, (2000).

65.  Putcha Mohan.S and Yaqub Adil, *Rings satisfying a certain idempotency condition*, Portugal Math. No.3, 325-328, (1985).

66.  Raftery.J.G, *On some special classes of prime rings*, Quaestiones Math., Vol. 10, 257-263, (1987).

67.  Ramamurthi.V.S, *Weakly regular rings*, Canada Math. Bull., Vol. 166, 317-321, (1973).

68.  Richard.P.Stanley, *Zero square rings,* Pacific. J. of Math., Vol. 30, 8-11, (1969).

69.  Schwarz.S, *A theorem on normal semigroups*, Check. Math. J., Vol. 85, 197, (1960).

70.  Searcold.O.Michael, *A Structure theorem for generalized J rings*, Proc. Royal Irish Acad., Vol. 87, 117-120, (1987).

71.  Shu Hao Sun, *On the least multiplicative nucleus of a ring*, J. of Pure and Appl. Algebra, Vol. 78, 311-318, (1992).

72.  Shung.T.M.S, *On invertible dispotent semigroup*, Bull. Fac. Sci. Engr. Chuo. Univ. Ser. J. Math., Vol. 34, 31-43, (1991).

73.  Smarandache Florentin, *Special Algebraic Structures,* in Collected Papers, Abaddaba, Oradea, Vol. III, 78-81, (2000).

74.  Stephenaion.W, *Modules whose lattice of submodules is distributive*, Proc. London Math. Soc., No. 28, 291-310, (1974).

75.  Stojan B and Miroslav.C, *Tight semigroups*, Publ. Inst, Math., Vol. 50, 71-84, (1991).

76.  Ursul.N.I, *Connectivity in weakly Boolean Topological rings*, IZV Acad. Nauk Moldor. Vol. 79, 17-21, (1989).

77. Van Rooyen.G.W.S, *On subcommutative rings*, Proc. of the Japan Acad., Vol. 63, 268-271, (1987).

78. Vasantha Kandasamy W.B., *On zero divisors in reduced group rings over ordered groups*, Proc. of the Japan Acad., Vol. 60, 333-334, (1984).

79. Vasantha Kandasamy W.B., *On semi idempotents in group rings*, Proc. of the Japan Acad., Vol. 61, 107-108, (1985).

80. Vasantha Kandasamy W.B., *A note on the modular group ring of finite p-group*, Kyungpook Math. J., Vol. 25, 163-166, (1986).

81. Vasantha Kandasamy W.B., *Zero Square group rings*, Bull. Calcutta Math. Soc., Vol. 80, 105-106, (1988).

82. Vasantha Kandasamy W.B., *On group rings which are p-rings*, Ganita, Vol. 40, 1-2, (1989).

83. Vasantha Kandasamy W.B., *Semi idempotents in semi group rings*, J. of Guizhou Inst. of Tech., Vol. 18, 73-74, (1989).

84. Vasantha Kandasamy W.B., *Semigroup rings which are zero square ring*, News Bull. Calcutta Math. Soc., Vol. 12, 8-10, (1989).

85. Vasantha Kandasamy W.B., *A note on the modular group ring of the symmetric group $S_n$*, J. of Nat. and Phy. Sci., Vol. 4, 121-124, (1990).

86. Vasantha Kandasamy W.B., *Idempotents in the group ring of a cyclic group*, Vikram Math. Journal, Vol. X, 59-73, (1990).

87. Vasantha Kandasamy W.B., *Regularly periodic elements of a ring*, J. of Bihar Math. Soc., Vol. 13, 12-17, (1990).

88. Vasantha Kandasamy W.B., *Semi group rings of ordered semigroups which are reduced rings*, J. of Math. Res. and Expo., Vol. 10, 494-493, (1990).

89. Vasantha Kandasamy W.B., *Semigroup rings which are p-rings*, Bull. Calcutta Math. Soc., Vol. 82, 191-192, (1990).

90. Vasantha Kandasamy W.B., *A note on pre J-group rings*, Qatar Univ. Sci. J., Vol. 11, 27-31, (1991).

91. Vasantha Kandasamy W.B., *A note on semigroup rings which are Boolean rings*, Ultra Sci. of Phys. Sci., Vol. 3, 67-68, (1991).

92. Vasantha Kandasamy W.B., *A note on the mod p-envelope of a cyclic group*, The Math. Student, Vol.59, 84-86, (1991).

93. Vasantha Kandasamy W.B., *A note on units and semi idempotents elements in commutative group rings*, Ganita, Vol. 42, 33-34, (1991).

94. Vasantha Kandasamy W.B., *Inner Zero Square ring*, News Bull. Calcutta Math. Soc., Vol. 14, 9-10, (1991).

95. Vasantha Kandasamy W.B., *On E-rings*, J. of Guizhou. Inst. of Tech., Vol. 20, 42-44, (1991).

96. Vasantha Kandasamy W.B., *On semigroup rings which are Marot rings*, Revista Integracion, Vol.9, 59-62, (1991).

97. Vasantha Kandasamy W.B., *Semi idempotents in the group ring of a cyclic group over the field of rationals*, Kyungpook Math. J., Vol. 31, 243-251, (1991).

98. Vasantha Kandasamy W.B., *A note on semi idempotents in group rings*, Ultra Sci. of Phy. Sci., Vol. 4, 77-78, (1992).

99. Vasantha Kandasamy W.B., *Filial semigroups and semigroup rings*, Libertas Mathematica, Vol.12, 35-37, (1992).

100. Vasantha Kandasamy W.B., *n-ideal rings*, J. of Southeast Univ., Vol. 8, 109-111, (1992).

101. Vasantha Kandasamy W.B., *On generalized semi-ideals of a groupring*, J. of Qufu Normal Univ., Vol. 18, 25-27, (1992).

102. Vasantha Kandasamy W.B., *On subsemi ideal rings*, Chinese Quat. J. of Math., Vol. 7, 107-108, (1992).

103. Vasantha Kandasamy W.B., *On the ring $Z_2 S_3$*, The Math. Student, Vol. 61, 246-248, (1992).

104. Vasantha Kandasamy W.B., *Semi group rings that are pre-Boolean rings*, J. of Fuzhou Univ., Vol. 20, 6-8, (1992).

105. Vasantha Kandasamy W.B., *Group rings which are a direct sum of subrings*, Revista Investigacion Operacional, Vol. 14, 85-87, (1993).

106.     Vasantha Kandasamy W.B., *On strongly sub commutative group ring*, Revista Ciencias Matematicas, Vol. 14, 92-94, (1993).

107.     Vasantha Kandasamy W.B., *Semigroup rings which are Chinese ring*, J. of Math. Res. and Expo., Vol.13, 375-376, (1993).

108.     Vasantha Kandasamy W.B., *Strong right S-rings*, J. of Fuzhou Univ., Vol. 21, 6-8, (1993).

109.     Vasantha Kandasamy W.B., *s-weakly regular group rings*, Archivum Mathematicum, Tomus. 29, 39-41, (1993).

110.     Vasantha Kandasamy W.B., *A note on semigroup rings which are pre p-rings*, Kyungpook Math. J., Vol.34, 223-225, (1994).

111.     Vasantha Kandasamy W.B., *A note on the modular semigroup ring of a finite idempotent semigroup*, J. of Nat. and Phy. Sci., Vol. 8, 91-94, (1994).

112.     Vasantha Kandasamy W.B., *Coloring of group rings*, J. Inst. of Math. and Comp. Sci., Vol. 7, 35-37, (1994).

113.     Vasantha Kandasamy W.B., *f-semigroup rings*, The Math. Edu., Vol. XXVIII, 162-164, (1994).

114.     Vasantha Kandasamy W.B., *J-semigroups and J-semigroup rings*, The Math. Edu., Vol. XXVIII, 84-85, (1994).

115.     Vasantha Kandasamy W.B., *Lie ideals of $Z_2 S_3$*, Acta Technica Napocensis, Vol. 37, 113-115, (1994).

116.     Vasantha Kandasamy W.B., *On a new type of group rings and its zero divisor*, Ult. Sci. of Phy. Sci., Vol. 6, 136-137, (1994).

117.     Vasantha Kandasamy W.B., *On a new type of product rings*, Ult. Sci. of Phy. Sci., Vol.6, 270-271, (1994).

118.     Vasantha Kandasamy W.B., *On a problem of the group ring $Z_p S_n$* , Ult. Sci. of Phy. Sci., Vol.6, 147, (1994).

119.     Vasantha Kandasamy W.B., *On pseudo commutative elements in a ring*, Ganita Sandesh, Vol. 8, 19-21, (1994).

120.     Vasantha Kandasamy W.B., *On rings satisfying $A^{\gamma} = b^s = (ab)^t$*, Proc. Pakistan Acad. Sci., Vol. 31, 289-292, (1994).

121.    Vasantha Kandasamy W.B., *On strictly right chain group rings*, Hunan. Annele Math., Vol. 14, 47-49, (1994).

122.    Vasantha Kandasamy W.B., *On strong ideal and subring of a ring*, J. Inst. Math. and Comp. Sci., Vol.7, 197-199, (1994).

123.    Vasantha Kandasamy W.B., *On weakly Boolean group rings*, Libertas Mathematica, Vol. XIV, 111-113, (1994).

124.    Vasantha Kandasamy W.B., *Regularly periodic elements of group ring*, J. of Nat. and Phy. Sci., Vol. 8, 47-50, (1994).

125.    Vasantha Kandasamy W.B., *Weakly Regular group rings*, Acta Ciencia Indica., Vol. XX, 57-58, (1994).

126.    Vasantha Kandasamy W.B., *Group rings which satisfy super ore condition*, Vikram Math. J., Vol. XV, 67-69, (1995).

127.    Vasantha Kandasamy W.B., *Obedient ideals in a finite ring*, J. Inst. Math. and Comp. Sci., Vol. 8, 217-219, (1995).

128.    Vasantha Kandasamy W.B., *On group semi group rings*, Octogon, Vol. 3, 44-46, (1995).

129.    Vasantha Kandasamy W.B., *On Lin group rings*, Zesztyty Naukowe Poli. Rzes., Vol. 129, 23-26, (1995).

130.    Vasantha Kandasamy W.B., *On Quasi-commutative rings*, Caribb. J. Math. Comp. Sci. Vol.5, 22-24, (1995).

131.    Vasantha Kandasamy W.B., *On semigroup rings in which $(xy)^n = xy$*, J. of Bihar Math. Soc., Vol. 16, 47-50, (1995).

132.    Vasantha Kandasamy W.B., *On the mod p-envelope of $S_n$*, The Math. Edu., Vol. XXIX, 171-173, (1995).

133.    Vasantha Kandasamy W.B., *Orthogonal sets in group rings*, J. of Inst. Math. and Comp. Sci., Vol.8, 87-89, (1995).

134.    Vasantha Kandasamy W.B., *Right multiplication ideals in rings*, Opuscula Math., Vol.15, 115-117, (1995).

135.    Vasantha Kandasamy W.B., *A note on group rings which are F-rings*, Acta Ciencia Indica, Vol. XXII, 251-252, (1996).

136. Vasantha Kandasamy W.B., *Finite rings which has isomorphic quotient rings formed by non-maximal ideals*, The Math. Edu., Vol. XXX, 110-112, (1996).

137. Vasantha Kandasamy W.B., *I\*-rings*, Chinese Quat. J. of Math., Vol. 11, 11-12, (1996).

138. Vasantha Kandasamy W.B., *On ideally strong group rings*, The Math. Edu., Vol. XXX, 71-72, (1996).

139. Vasantha Kandasamy W.B., *Gaussian Polynomial rings*, Octogon, Vol.5, 58-59, (1997).

140. Vasantha Kandasamy W.B., *On semi nilpotent elements of a ring*, Punjab Univ. J. of Math. , Vol. XXX, 143-147, (1997).

141. Vasantha Kandasamy W.B., *On tripotent elements of a ring*, J. of Inst. of Math. and Comp. Sci., Vol. 10, 73-74, (1997).

142. Vasantha Kandasamy W.B., *A note on f-group rings without non-zero nilpotents*, Acta Ciencia Indica, Vol. XXIV, 15-17, (1998).

143. Vasantha Kandasamy W.B., *Inner associative rings*, J. of Math. Res. and Expo., Vol. 18, 217-218, (1998).

144. Vasantha Kandasamy W.B., *On a quasi subset theoretic relation in a ring*, Acta Ciencia Indica, Vol. XXIV, 9-10, (1998).

145. Vasantha Kandasamy W.B., *On SS-rings*, The Math. Edu., Vol. XXXII, 68-69, (1998).

146. Vasantha Kandasamy W.B., *Group rings which have trivial subrings*, The Math. Edu., Vol. XXXIII, 180-181, (1999).

147. Vasantha Kandasamy W.B., *On E-semi group rings*, Caribbean J. of Math. and Comp. Sci., Vol. 9, 52-54, (1999).

148. Vasantha Kandasamy W.B., *On demi- modules over rings*, J. of Wuhan Automotive Politechnic Univ., Vol. 22, 123-125, (2000).

149. Vasantha Kandasamy W.B., *On finite quaternion rings and skew fields*, Acta Cienca Indica, Vol. XXIV, 133-135, (2000).

150. Vasantha Kandasamy W.B., *On group rings which are $\gamma_n$ rings*, The Math. Edu., Vol. XXXIV, 61, (2000).

151.    Vasantha Kandasamy W.B., *CN rings*, Octogon, Vol.9, 343-344, (2001).

152.    Vasantha Kandasamy W.B., *On locally semi unitary rings*, Octogon, Vol.9, 260-262, (2001).

153.    Vasantha Kandasamy W.B., *Tight rings and group rings,* Acta Ciencia Indica, Vol. XXVII, 87-88, (2001).

154.    Vasantha Kandasamy W.B., *Smarandache Semigroups*, American Research Press, Rehoboth, NM, (2002).

155.    Vasantha Kandasamy W.B., *On Smarandache pseudo ideals in rings*, (2002).
        http://www.gallup.unm.edu/~smaranandache/pseudoideals.pdf

156.    Vasantha Kandasamy W.B., *Smarandache Semirings and Smarandache Semifields*, Smarandache Notions Journal, American Research Press, Vol. 13, 88-91, (2002).

157.    Vasantha Kandasamy W.B., *Smarandache Semirings, Semifields and Semivector spaces*, American Research Press, Rehoboth, NM, (2002).

158.    Vasantha Kandasamy W.B., *Smarandache Zero divisors*, (2002).
        http://www.gallup.unm.edu/~smarandache/ZeroDivisor.pdf

159.    Vasantha Kandasamy W.B., *Finite zeros and finite zero-divisors*, Varahmihir J. of Math. Sci., Vol. 2, (To appear), 2002.

160.    Victoria Powers, *Higher level orders on non-commutative rings*, J. Pure and Appl. Algebra, Vol. 67, 285-298, (1990).

161.    Vougiouklis Thomas, *On rings with zero divisors strong V-groups*, Comment Math. Univ. Carolin J., Vol. 31, 431-433, (1990).

162.    Wilson John. S., *A note on additive subgroups of finite rings*, J. Algebra, No. 234, 362-366, (2000).

163.    Yakub Adil, *Structure of weakly periodic rings with potent extended commutators*, Int. J. of Math. Sci., Vol. 25, 299-304, (2001).

164.    Zariski.O and Samuel.P, *Commutative Algebra*, Van Nostrand Reinhold, (1958).

165.    Zhang-Chang quan, *Inner Commutative Rings*, Sictiuan Daxue Xuebao, Vol. 26, 95-97, (1989).

# INDEX

## ABOUT THE AUTHOR

Dr. W. B. Vasantha is an Associate Professor in the Department of Mathematics, Indian Institute of Technology Madras, Chennai, where she lives with her husband Dr. K. Kandasamy and daughters Meena and Kama. Her current interests include Smarandache algebraic structures, fuzzy theory, coding/ communication theory. In the past decade she has completed guidance of seven Ph.D. scholars in the different fields of non-associative algebras, algebraic coding theory, transportation theory, fuzzy groups, and applications of fuzzy theory to the problems faced in chemical industries and cement industries. Currently, six Ph.D. scholars are working under her guidance. She has to her credit 241 research papers of which 200 are individually authored. Apart from this she and her students have presented around 262 papers in national and international conferences. She teaches both undergraduate and post-graduate students at IIT and has guided 41 M.Sc. and M.Tech projects. She has worked in collaboration projects with the Indian Space Research Organization and with the Tamil Nadu State AIDS Control Society. She is currently authoring a ten book series on Smarandache Algebraic Structures in collaboration with the American Research Press.

She can be contacted at vasantha@iitm.ac.in
You can visit her on the web at: http://mat.iitm.ac.in/~wbv